\documentclass[a4paper,11pt, oneside]{book}

\usepackage{geometry}
\usepackage{pdfpages} 

\usepackage{todonotes}

\usepackage{changepage} 

\usepackage{mathrsfs}

\date{\today}
\usepackage{authblk}

\usepackage[italian, english]{babel} 
\usepackage[OT2, T1]{fontenc} 
\usepackage[utf8]{inputenc} 

\usepackage{latexsym}
\usepackage{amssymb,amsmath}
\usepackage{amsthm,graphicx}
\usepackage{enumerate}
\usepackage{array}

\usepackage{xcolor} 
\usepackage{color} 

\definecolor{uniba}{RGB}{0, 64, 112}

\usepackage[colorlinks = true,
 linkcolor = uniba,
 urlcolor = uniba,
 citecolor = uniba,
 anchorcolor = uniba
]{hyperref}

\usepackage{bm}

\usepackage[hyperpageref]{backref} 
\usepackage{hyperref}

\usepackage{esint}

\usepackage{centernot}

\usepackage{lipsum}

\usepackage[noadjust]{cite} 

\usepackage{sectsty} 

\chapterfont{\color{uniba}} 
\sectionfont{\color{uniba}}
\subsectionfont{\color{uniba}}
\subsubsectionfont{\color{uniba}}

\usepackage{titlesec}
\usepackage{microtype}
\usepackage{tikz}

\titleformat{\chapter}[display]
 {\normalfont\bfseries\color{uniba}}
 {\filleft%
 \begin{tikzpicture}
 \node[
 outer sep=0pt,
 text width=2.5cm,
 minimum height=3cm,
 fill=uniba,
 font=\color{white}\fontsize{80}{90}\selectfont,
 align=center
 ] (num) {\thechapter};
 \node[
 rotate=90,
 anchor=south,
 font=\color{black}\Large\normalfont
 ] at ([xshift=-5pt]num.west) {\textls[180]{\textsc{\chaptertitlename}}}; 
 \end{tikzpicture}%
 }
 {10pt}
 {\titlerule[2.5pt]\vskip3pt\titlerule\vskip4pt\LARGE\sffamily}

\usepackage[]{frontespizio}

\usepackage{emptypage} 

\usepackage{lmodern}

\usepackage{ifthen}
\usepackage{fancyhdr} 
\newenvironment{abstract}{\cleardoublepage
\vspace*{3\baselineskip}\center
\huge{\textbf{\tu{Abstract}} \normalsize}\endcenter
\quote}{\endquote\clearpage}

\allowdisplaybreaks 

\def\claim#1.{\noindent {\bf #1.}}

\def\flushright#1{{\unskip\nobreak\hfil\penalty50\hskip2em\hbox{}\nobreak\hfil%
#1\parfillskip=0pt\finalhyphendemerits=0\par}}

\def\bull{\vrule height 1.8ex width 1.0ex depth .1ex }
\def\QED{\ifmmode\eqno\hbox{$\bull$}\else\flushright{\hbox{$\bull$}}\fi}

\newcommand{\tnorm}[1]{{\left\vert\kern-0.25ex\left\vert\kern-0.25ex\left\vert #1 
		\right\vert\kern-0.25ex\right\vert\kern-0.25ex\right\vert}}

\usepackage{amsmath}
\usepackage{amsthm}
\usepackage{amscd}
\usepackage{amssymb}
\usepackage{amsfonts}
\usepackage{graphics}
\usepackage{color}
\usepackage{comment}
\usepackage[dvipsnames]{xcolor}

\usepackage[all,cmtip]{xy}
\usepackage[normalem]{ulem}

\newlength{\defbaselineskip }
 \long\def\salta#1{\relax}

 \theoremstyle{plain}
\newtheorem{theorem}{Theorem}[section]
\newtheorem{proposition}[theorem]{Proposition}
\newtheorem{lemma}[theorem]{Lemma}

\newtheorem{corollary}[theorem]{Corollary}

\renewcommand{\a}{{\alpha}}

\theoremstyle{definition}
\newtheorem{definition}[theorem]{Definition}

\newtheorem{remark}[theorem]{Remark}

\newcommand{\ja}{J_{\alpha}}
\newcommand{\la}{\lambda}
\newcommand{\xm}{X_-^m}
\newcommand{\xp}{X_+^m}
\newcommand{\te }{\theta}

\newcommand{\R}{\mathbb{R}}
\newcommand{\N}{\mathbb{N}}
\newcommand{\Sm}{\mathbb{S}}
\newcommand{\bu}{\bar{u}}
\newcommand{\bv}{\bar{v}}
\newcommand{\bz}{\bar{z}}

\newcommand{\na}{\mathbb{N}}
\newcommand{\ze}{\mathbb{Z}}

\newcommand{\re}{\mathbb{R}}

\newcommand{\rnk}[1]{\operatorname{\mathrm{rank}}{\,#1}}

\def\N{\mathbb{N}}

\def\de{\delta}

\def\wr{W^{-1,r'}(\Omega)}
\def\t1p0{T^{1,p}_{0}(\Omega)}

\def\m2{M^{\frac{N(p-1)}{N-1}}(\Omega)}

\def\div{\text{div}}

\def\sobr{W^{1,r}_{0}(\Omega)}
\def\sobq{W^{1,q}_{0}(\Omega)}
\def\sob{W^{1,p}_{0}(\Omega)}

\def\into{\int_{\Omega}}

\def\w-1p'{W^{-1,p'}(\Omega)}
\def\pw-1p'u{L^{p'}(0,1;W^{-1,p'}(\Omega))}

\def\lp'n{(L^{p'}(\Omega))^{N}}

\newcommand{\tu}[1]{\textcolor{uniba}{#1}}

\usepackage{biolinum}

\usepackage{minitoc} 

\usepackage[a-2b]{pdfx} 

\begin{document}

\pagenumbering{Roman}

\thispagestyle{empty}

\newgeometry{bottom=-6in}

\begin{adjustwidth}{-1in}{-1in}

\begin{tikzpicture}[remember picture,overlay]
\node [opacity=1,scale=1] at (current page.center) {\includegraphics{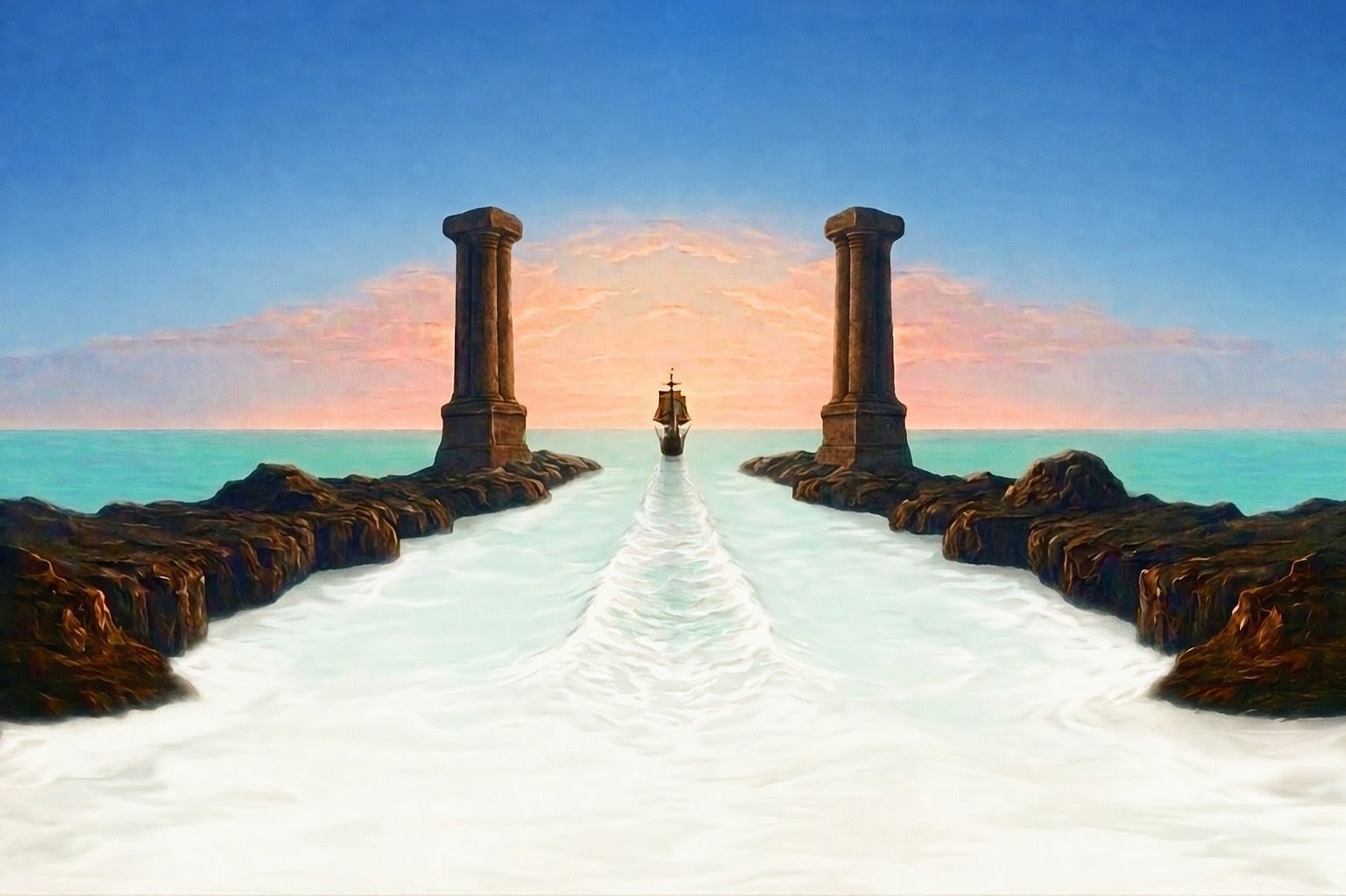}};
\end{tikzpicture}

\biolinum{

\begin{flushleft}

\begin{center}
\vspace{-2.5cm}
\Large{\textbf{Ph.D. Thesis}}
\end{center}

\begin{center}

{\Huge \textbf{Variational and Geometric Analysis for Quasilinear Elliptic Equations and Systems \\}}

\end{center}

\vspace{4cm}

\begin{center}
\huge{
\textbf{Natalino Borgia}
}
\end{center}

\vspace{10cm}

\begin{center}
\textbf{
\Large{Supervised by}}\\
\vspace{0.5cm}
\textbf{
\LARGE{
Prof. Silvia Cingolani}
}
\end{center}

\vskip4cm

\begin{center}
\includegraphics*[width=0.42\textwidth]{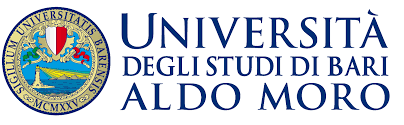}
\end{center}

\end{flushleft} 

}

\end{adjustwidth}

\restoregeometry

\cleardoublepage
\thispagestyle{empty}
\vspace*{\stretch{2}}
{\emph{Many people think that mathematical ideas are static. They think that the ideas originated at some time in the historical past and remain unchanged for all future times. There are good reasons for such a feeling. After all, the formula for the area of a circle was $\pi r^2$ in Euclid's day and at the present time is still $\pi r^2$. But to one who knows mathematics from the inside, the subject has rather the feeling of a living thing. It grows daily by the accretion of new information, it changes daily by regarding itself and the world from new vantage points, it maintains a regulatory balance by consigning to the oblivion of irrelevancy a fraction of its past accomplishments.}}

\vspace{1cm}
\rightline{ \textasciitilde \emph{ Philip J. Davis } \textasciitilde }
\vspace{\stretch{2}}


\newpage

\thispagestyle{empty}

\begin{figure*}
\begin{center}
\includegraphics*[width=0.35\textwidth]{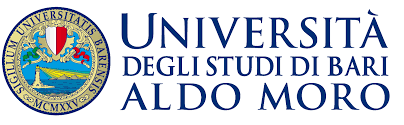}
\end{center}
\vspace{-4em}
\end{figure*}
\begin{center}

\rule{13.5cm}{0.4pt}

\bigskip
\bigskip
\bigskip

\bigskip
\bigskip
\bigskip

{\large Università degli Studi di Bari Aldo Moro
\\}

\medskip

{\large Department of Mathematics
\\}

\bigskip
\medskip

{\large Ph.D. in Computer Science and Mathematics\\
}

\medskip

{
\normalsize
XXXVIII Cycle\\}

\bigskip
\medskip

{
Mathematical Analysis
MATH-03/A\\}

\bigskip
\bigskip
\bigskip

\tu{
{\huge \textbf{Variational and Geometric Analysis for\\
 Quasilinear Elliptic Equations and Systems}}
}

\end{center}

\bigskip
\bigskip
\bigskip
\bigskip
\bigskip

\begin{minipage}{0.5\textwidth}
Ph.D. Candidate:\\
\textbf{Natalino Borgia}
\end{minipage}
\begin{minipage}{0.5\textwidth}
Supervisor:\\
\textbf{Prof. Silvia Cingolani}
\end{minipage}
\vskip0.5cm
\begin{minipage}{0.5\textwidth}
$\,$\\
$\,$
\end{minipage}
\vskip1.7cm
\begin{minipage}{0.5\textwidth}
$\,$\\
$\,$
\end{minipage}
\begin{minipage}{0.5\textwidth}
Ph.D. Coordinator:\\
\textbf{Prof. Giovanna Castellano}\\
\end{minipage}

\bigskip
\bigskip
\bigskip
\bigskip
\bigskip
\bigskip
\bigskip

\begin{center}
\noindent\rule{8cm}{0.4pt}\\

\medskip

Ph.D. Thesis, 2026\\
\noindent\rule{8cm}{0.4pt}
\end{center}

\newpage

\cleardoublepage
\thispagestyle{empty}

\begin{abstract} 

\bigskip
\noindent
In this thesis, we employ variational and geometric analysis to investigate quasilinear elliptic equations and systems arising from nonlinear phenomena with a significant environmental impact.  We focus on quasilinear elliptic  systems driven by various nonlinear operators, such as the $p$-Laplacian, and nonlinear sources that are allowed to exhibit both subcritical and critical growth. We aim to establish the existence of solutions for perturbation of specific  eigenvalue problems, by employing variational and topological methods. 
To establish existence results for autonomous systems of quasilinear PDEs in the spirit of the celebrated paper by Amann and Zehnder \cite{amann_zehnder}, we develop a local Morse theory for functional associated to quasilinear elliptic systems. Our objective is to bridge relationship between the local behaviour of the functional near its critical points and differential notions like the Morse index. To this aim, by refining topological arguments introduced by Cingolani and Degiovanni in Banach product spaces in \cite{CD1}, we establish the finiteness of the critical groups and we derive a Poincaré-Hopf formula in a Banach product space, in presence of both subcritical and critical nonlinear coupling. Furthermore, we address the lack of regularity in the associated Euler functional by employing a penalization method combined with uniform $L^{\infty}$-estimates. We also establish uniform boundedness results for anisotropic quasilinear systems, that are of interest within regularity theory. To show existence results for non-autonomous systems of quasilinear PDEs in the spirit of the celebrated paper \cite{LANDLAZ} of  Landesman and Lazer, we consider the eigenvalue problem for quasilinear elliptic systems introduced by de Thélin in \cite{DT}. Leveraging the intrinsic geometric structure of the eigenvalue problem, we prove the simplicity and isolation of the first eigenvalue $\lambda_1$. Furthermore, we show the existence of a sequence of eigenvalues by employing a suitable deformation lemma proved by Bonnet in \cite{BON}. Subsequently, we analyze a non-autonomous perturbation of the eigenvalue problem. In particular, the established spectral properties enable us to derive new sufficient Landesman-Lazer type conditions within the framework of quasilinear elliptic systems.
While the quasilinear elliptic systems under consideration typically involve the $p$-Laplacian for $1 < p < N$, where $N$ denotes the dimension of the space $\mathbb{R}^N$, we also investigate the geometric role of the $N$-Laplacian in the $N$-dimensional Euclidean Onofri inequality. This inequality was established by Del Pino and Dolbeault \cite{DD3} for smooth functions with compact support. Specifically, after extending this inequality to a suitable weighted Sobolev space via a density argument, we exploit its connection with the Liouville equation on $\mathbb{R}^N$ to prove an equivalence with the sharp logarithmic Moser-Trudinger inequality on the unit ball $B_1 \subset \mathbb{R}^N$,  implicitly contained in the work \cite{CarlesonChang} of Carleson and Chang. For the specific case $N=2$, we further highlight a geometric link to the unit sphere $\mathbb{S}^2$ through the stereographic projection. 
\end{abstract}

\newpage

\cleardoublepage
\thispagestyle{empty}

\begin{abstract} 

\bigskip
\noindent
In questa tesi, utilizziamo l'analisi variazionale e geometrica per studiare equazioni e sistemi quasilineari ellittici che nascono da fenomeni non lineari con un impatto significativo sull'ambiente.  Ci concentriamo su sistemi quasilineari ellittici guidati da diversi operatori non lineari come il $p$-Laplaciano, e sorgenti non lineari che possono presentare crescita sia sottocritica che critica. Il nostro obiettivo è stabilire l'esistenza di soluzioni per una perturbazione di specifici problemi agli autovalori,  utilizzando metodi sia variazionali che topologici. 
Per stabilire risultati di esistenza per sistemi autonomi di PDE quasilineari nello spirito del celebre lavoro di Amann e Zehnder \cite{amann_zehnder}, sviluppiamo una teoria di Morse locale per funzionali associati a sistemi quasilineari ellittici. Il nostro obiettivo è creare un collegamento tra il comportamento locale del funzionale vicino i suoi punti critici e nozioni differenziali come l'indice di Morse. A tal fine, affinando gli argomenti topologici introdotti da Cingolani e Degiovanni in \cite{CD1} per spazi di Banach prodotto, stabiliamo la finitezza dei gruppi critici e deriviamo una formula di  Poincaré-Hopf in uno spazio di Banach prodotto, in presenza di accoppiamenti non lineari sia sottocritici che critici. Inoltre, affrontiamo la mancanza di regolarità nel funzionale associato utilizzando  un metodo di penalizzazione combinato con  stime $L^{\infty}$-uniformi. Stabiliamo inoltre risultati di limitatezza uniforme per sistemi quasilineari anisotropi, che rivestono un certo interesse nell'ambito della teoria della regolarità. Per mostrare risultati di esistenza per sistemi non-autonomi di PDE quasilineari nello spirito del celebre lavoro  \cite{LANDLAZ} di Landesman e Lazer, consideriamo il problema agli autovalori per sistemi quasilineari ellittici introdotto da de Thélin in \cite{DT}. Facendo leva sulla struttura geometrica intrinseca del problema agli autovalori, proviamo la semplicità e l'isolatezza del primo autovalore $\lambda_1$. Inoltre, mostriamo l'esistenza di una successione di autovalori utilizzando un opportuno lemma di deformazione, dimostrato da Bonnet in \cite{BON}. Successivamente, analizziamo una perturbazione non-autonoma del problema agli autovalori. In particolare, le proprietà spettrali ottenute ci permettono di derivare nuove sufficienti condizioni di tipo  Landesman-Lazer all'interno del framework di sistemi quasilineari ellittici.
Mentre i sistemi in considerazione tipicamente coinvolgono il $p$-Laplaciano per $1 < p < N$, dove $N$ denota la dimensione dello spazio $\mathbb{R}^N$, il ruolo geometrico  dell' $N$-Laplaciano lo investighiamo nella disuguaglianza di Onofri euclidea $N$-dimensionale. Questa disuguaglianza è stata ottenuta da  Del Pino e Dolbeault in  \cite{DD3} per funzioni regolari a supporto compatto. Nello specifico, dopo aver esteso la disuguaglianza ad un opportuno spazio di Sobolev pesato tramite un argomento di densità, sfruttiamo la connessione con l'equazione di Liouville in  $\mathbb{R}^N$ per dimostrare un risultato di equivalenza con la ottimale disuguaglianza logaritmica di Moser-Trudinger   sulla palla unitaria $B_1 \subset \mathbb{R}^N$, contenuta implicitamente nel lavoro \cite{CarlesonChang} di Carleson e Chang. Nel caso specifico di $N=2$, sottolineiamo il collegamento con la sfera unitaria due-dimensionale $\mathbb{S}^2$ attraverso la proiezione stereografica. 

\end{abstract}

\newpage

\cleardoublepage
\thispagestyle{empty}

\tu{
{\LARGE \textbf{Notation}}}\\

\begin{itemize}
\item We denote by $N \in \mathbb{N}$ the dimension of the space $\mathbb{R}^N$. In particular, $N \geq 2$ unless otherwise stated. \medskip

\item We denote by $\Omega \subset \mathbb{R}^N$ a bounded domain (open connected set) with smooth boundary $\partial \Omega$. \medskip

\item We denote by $B_1$ the open unit ball of $\mathbb{R}^N$, and we write $V_N$ for its volume and $\omega_{N-1}$ for the measure of its surface. \medskip

\item Given a set $A$, we denote by $|A|$ its Lebesgue measure. \medskip

\item Given a function $f: U \to \mathbb{R}$ with $U\subseteq \mathbb{R}^N$, for any $i=1,...,N$ we write the first partial derivative as $ \displaystyle \frac{ \partial f}{ \partial x_i}$. We denote by $\nabla f(x) $ the gradient of $f$  evaluated at $x$.  \medskip

\item For $w: \Omega \to \mathbb{R}$ and $\Psi: \mathbb{R}^N \to \mathbb{R}$, we denote by $\nabla \Psi (\nabla w(x))$ the gradient of $\Psi$ evaluated at $\nabla w(x)$. \medskip

\item For any $r > 1$, we denote by $L^{r}(\Omega)$  the standard Lebesgue space of measurable functions $f: \Omega \to \mathbb{R}$ such that $\into |f|^r \, dx < + \infty.$ We denote by  $\| \cdot \|_{r}$ the usual norm in $L^r(\Omega)$. \medskip

\item For any $r > 1$, we denote by $W_0^{1,r}(\Omega)$  the standard Sobolev space of functions in $L^r(\Omega)$, vanishing on the boundary, whose first-order weak derivatives also belong to $L^r(\Omega)$.\\
We denote by $\| \cdot \|_{1,r}$  the usual gradient norm in $W^{1,r}_0(\Omega)$.\\
Moreover, we denote by $ \wr$ the dual space of $W_0^{1,r}(\Omega)$.\medskip

\item For any $r > 1$, we denote by $r^*:=Nr/(N-r)$ the critical Sobolev exponent of $r$, and by $r':=r/(r-1)$ the conjugate exponent. \medskip

\item Given $1<p,q<N$, we denote by $X$ the Banach product space $X:=W_0^{1,p}(\Omega) \times W_0^{1,q}(\Omega)$. We equip the space $X$ with  the norm

$$\|z\|= \|u\|_{1,p} + \|v\|_{1,q}, \qquad z=(u,v) \in X.$$

\medskip
\noindent
In particular, $X$ is not uniformly convex under the norm defined above (see \cite{CLARKSON}). \medskip

\item We denote by $(\mathbb{S}^2,g_0)$ the two-dimensional unit sphere $\mathbb{S}^2$ endowed with the standard metric $g_0$, which is represented by the matrix

$$
g_0(x):=\begin{pmatrix}
\lambda(x) & 0\\
0 & \lambda(x)
\end{pmatrix}, \qquad \text{where} \quad \lambda:\mathbb{R}^2 \to \mathbb{R}, \quad \lambda(x):= \frac{4}{(1+|x|^2)^2}.$$ 

\item We denote by $W^{1,2}(\mathbb{S}^2, g_0)$ the Sobolev space consisting of functions that, along with their distributional first derivatives, are square-integrable with respect to the Riemannian volume form $d\nu_{g_0}:= \lambda(x) dx$. \medskip

\item We denote by $\nabla_{g_0} f$ the gradient of $f$ with respect to the metric $g_0$.

\end{itemize}

\newpage

\frontmatter 

\newpage

\tableofcontents 

\chapter*{Introduction
}
\addcontentsline{toc}{chapter}{Introduction}
\markboth{Introduction}{Introduction}

Several nonlinear partial differential equations (PDEs) and systems arise from nonlinear phenomena with significant environmental impacts. In this thesis, we investigate quasilinear elliptic systems driven by different operators, involving isotropic and anisotropic ones, which come from several nonlinear physical processes, such as non-Newtonian mechanics, nonlinear elasticity, glaciology, combustion theory, and population biology (see for instance \cite{diazthelin, glow, manamaw, marcellini1,LWO,PR, RZ}). While isotropic operators  act uniformly in all directions, anisotropic operators exhibit a dependency on the spatial direction, introducing a bias or preference toward specific axes or orientations. The most common nonlinear isotropic operator is the $p$-Laplacian $\Delta_p$, which generalizes the standard Laplacian $\Delta$  and is defined for any $1<p < \infty$ as $\Delta_p u = \text{div}(|\nabla u|^{p-2} \nabla u).$ In particular, $\nabla u$ identifies the direction of steepest ascent, while its magnitude $|\nabla u|$ quantifies the slope's intensity. The term $|\nabla u|^{p-2}$ acts as a weight, or stiffness factor, determined by the gradient itself. The divergence operator, ``$\text{div}$'', characterizes the flux from a point. The case $p=2$ recovers the standard linear Laplacian $\Delta$. For $p>2$, the operator is degenerate, exhibiting increased stiffness in regions with large gradients and becoming `softer' where the gradient is small. Conversely, for $1<p<2$, the operator is singular and tends to blow up as the gradient approaches zero.
From a physical standpoint, the most prominent application lies in the transition from Newtonian to non-Newtonian fluid dynamics \cite{AM}. While the standard Laplacian $\Delta$ describes Newtonian flows, the $p$-Laplacian $\Delta_p$ characterizes non-Newtonian behavior, where viscosity is a function of the applied stress. Specifically, for $1 < p < 2$, the fluid is pseudoplastic (shear-thinning), exhibiting decreased viscosity as the shear rate increases. Conversely, for $p > 2$, the fluid is dilatant (shear-thickening), showing increased resistance in response to greater applied force. Specifically, this thesis focuses on the analytical and geometrical aspects of properties arising from the $p$-Laplacian, which represents the core link between the different problems addressed here.  From the mathematical point of view, the study of quasilinear PDEs presents several specific difficulties and challenging theoretical problems. 

\medskip

A first target of this thesis is to develop a local Morse theory for functional associated to quasilinear elliptic systems in presence of subcritical and critical nonlinear coupling. We aim to establish a relationship between the local behaviour of the functional near its critical points and differential notions like the Morse index. In particular, we aim to derive a  Poincaré-Hopf formula in a Banach product space. We recall that in the Hilbert space setting, the celebrated Poincaré-Hopf formula has been extended to functionals $f$ whose gradient $\nabla f$ is a compact perturbation of the identity. This extension allows for the application of Leray-Schauder degree theory, rather than the Brouwer degree typically used for $C^1$ functions $f: \mathbb{R}^N \to \mathbb{R}$ (see, e.g., \cite{mawhin_willem1989, chang, li_li_liu2005, bartsch_dancer2009, benci1991}). When $Y$ is a Banach space (not necessarily Hilbert), the derivative $f'$ is naturally defined from $Y$ into its dual $Y'$. Consequently, this does not fit the standard framework of Leray-Schauder degree, but rather that of the Browder degree developed in \cite{browder1983} for demicontinuous maps of class $(S)_+$. In our context, the functional setting is defined by the Banach product space $X := W_0^{1,p}(\Omega) \times W_0^{1,q}(\Omega)$, where $1 < p, q < N$ and $\Omega \subset \mathbb{R}^N$ ($N \geq 2$) is a smooth bounded domain. Specifically, we consider the following quasilinear elliptic system:

\begin{equation}\label{INTRODUZIONEpqPoincareHopf}
	\begin{cases}
		\begin{array}{ll}
			-\text{\rm div}  \; [ \nabla \Psi_{a,p} \left( \nabla u  \right)] =F_s(u,v) & x\in\Omega,
			\medskip \\
			-\text{\rm div}  \; [ \nabla \Psi_{b,q} \left( \nabla v  \right)]=F_t(u,v) & x\in \Omega, \medskip \\
			u=v=0  & x\in \partial\Omega,
		\end{array}
	\end{cases} \tag{I.1}
\end{equation}

\medskip
\noindent
where, given $c \geq 0$ and  $r >1$, we denote by $\Psi_{c,r}: \mathbb{R}^N \to \mathbb{R}$ the function defined as

\begin{equation*}
	\displaystyle \Psi_{c,r}(\xi):= \frac{1}{r} \biggl[ \left( c +
	\left| \xi \right|^2\right)^{ \frac{r}{2}} - c^{\frac{r}{2}} \biggr].
\end{equation*}	

\medskip
\noindent
Hence, in \eqref{INTRODUZIONEpqPoincareHopf} we consider $\Psi_{a,p}$ and $\Psi_{b,q}$ for $1<p,q<N$ and $a,b \geq 0$. The nonlinearity $F \in C^1(\mathbb{R}^2)$ is allowed to exhibit both subcritical and critical growth. We show that the Euler $C^1$-functional $J_{a,b}: X \to \mathbb{R}$ associated with \eqref{INTRODUZIONEpqPoincareHopf} has a Fréchet derivative $J_{a,b}': X \to X'$ that is locally a map of class $(S)_+$. By employing topological arguments introduced by Cingolani and Degiovanni in \cite{CD1}, we establish the finiteness of the critical groups of $J_{a,b}$, thereby extending results previously obtained for quasilinear elliptic equations at critical growth to the case of quasilinear systems at critical growth. Furthermore, the use of the Browder degree \cite{browder1983} enables us to establish a Poincaré–Hopf formula in the setting of quasilinear elliptic systems.

\medskip
An other target of this thesis is to show existence results for autonomous systems of quasilinear PDEs in the spirit of the celebrated paper by Amann and Zehnder \cite{amann_zehnder}. The main idea in \cite{amann_zehnder} is that if the "linearized" behavior at the origin is sufficiently different from the behavior at infinity (measured by their respective Morse indices), then a non-trivial solution must exist. To extend this approach to systems within a Banach framework, we employ several innovative techniques. First, we establish a general abstract Linking Theorem, extending the result of Cingolani, Degiovanni and Vannella \cite[Theorem 7.1]{CDV} concerning symmetric cones to a broader setting of $(p,q)$-homogeneous sets. Second, the nonlocal coupling induces a lack of regularity in the functional. This phenomenon can be overcome by using a penalization method and uniform regularity estimates. A Chapter of the thesis is therefore devoted to establish uniform $L^\infty$-estimates results for systems involving anisotropic quasilinear elliptic equations. Consequently, we investigate the following non-autonomous system:

\begin{equation}\label{INTRODUZIONEpqLinfinitoboundedness}
	\begin{cases}
		\begin{array}{ll}
			-\text{\rm div}  \; [ \nabla \Psi_1 \left( \nabla u  \right)] = H_s(\delta,x, u,v) & \hbox{in} \ \Omega,
			\medskip \\
			-\text{\rm div} \; [ \nabla \Psi_2 \left( \nabla v  \right)]= H_t(\delta,x,u,v) & \hbox{in} \ \Omega, \medskip\\
			u=v=0  & \hbox{on} \ \partial\Omega,
		\end{array}
	\end{cases} \tag{I.2}
\end{equation}

\medskip
\noindent
where $\Omega$ and  $N$ are defined as above. The nonlinearity $H$ is allowed to exhibit both subcritical and critical growth and may depend on a parameter $\delta \in I \subseteq \mathbb{R}$. Moreover, $\Psi_1, \Psi_2: \mathbb{R}^N \to \mathbb{R}$ are convex functions satisfying quite general conditions for $1 < p, q < N$, as detailed in Chapter \ref{SezReg}; these conditions allow for the recovery of the previously introduced operators $\Psi_{a,p}$ and $\Psi_{b,q}$. We  establish the $L^{\infty}$-boundedness of all solutions to \eqref{INTRODUZIONEpqLinfinitoboundedness} and  we prove that these solutions are uniformly bounded with respect to $\delta$ within any sufficiently small ball of the functional space setting $X$. To address the singularities arising when $1 < p < 2$ or $1 < q < 2$, we employ a suitable generalization of Stampacchia's Lemma (see \cite{BDO}). Furthermore, we handle the possibly critical growth of $H$ by proving results that hold up to the critical threshold, extending those established in \cite{Vann1}. We stress that in quasilinear elliptic systems of the form \eqref{INTRODUZIONEpqLinfinitoboundedness}, an additional difficulty in proving such results is the possible coupling of the functions $u$ and $v$ in the nonlinearity $H$.

\medskip
Building upon these preliminary results, we pass to detect nontrivial solutions for quasilinear systems where the nonlinearity interacts with the spectrum of an associated eigenvalue problem, in the spirit of the celebrated paper \cite{amann_zehnder} of Amann and Zehnder. We analyze the following autonomous quasilinear system:

\begin{equation}\label{INTRODUZIONEpqNontrivialSolutions}
	\begin{cases}
		\begin{array}{ll}
			-\text{\rm div}  \; [ \nabla \Psi_{a,p} \left( \nabla u  \right)] = G_s(u,v) & x\in\Omega,
			\medskip \\
			-\text{\rm div}  \; [ \nabla \Psi_{a,q} \left( \nabla v  \right)]= G_t(u,v) & x\in \Omega, \medskip	\\
			u=v=0  & x\in \partial\Omega,
		\end{array}
	\end{cases} \tag{I.3}
\end{equation}

\medskip
\noindent
where $\Omega$ is as before, $2\leq p,q<N$, $a\geq 0$,   and $G \in C^{1}(\re^2, \re)$ is a $(p,q)$-asymptotically linear term at infinity, with $ \nabla G(0,0)=(0,0)$, and it is of class $C^2$ in a suitable neighborhood of the origin. Since $(0,0) \in X$ is an a priori trivial solution of \eqref{INTRODUZIONEpqNontrivialSolutions}, we aim to find nontrivial solutions in both the asymptotic resonant and nonresonant cases, relative to the following eigenvalue problem:

\begin{equation}\label{INTROeigen}
	\begin{cases}
		\begin{array}{ll}
			- \Delta_p u  = \lambda |u|^{p-2} u +
			\frac{\lambda}{\beta +1}|u|^{\alpha} |v|^{\beta} v & x\in\Omega,
			 \medskip \\
			- \Delta_q v = \lambda |v|^{q-2} v +
			\frac{\lambda}{\alpha +1}|u|^{\alpha} |v|^{\beta} u, & x\in\Omega, \medskip \\
			u=v=0,  & x\in \partial\Omega,
		\end{array}
	\end{cases} \tag{I.4}
\end{equation}

\medskip
\noindent
where $\Omega,N,p,q$ are as before and $\alpha, \beta>0$ are real numbers satisfying
$(\alpha +1)/ p + (\beta +1)/q =1$.

As mentioned, since the functional is defined on the product space $X = W_0^{1,p}(\Omega) \times W_0^{1,q}(\Omega)$, we exploit recent Morse-theoretic techniques in Banach spaces (see \cite{CV, CV2, CV4, U}); moreover, we derive a Saddle Linking Theorem specifically tailored to $(p,q)$-homogeneous sets. To prove the non-triviality of the solution, we compute the critical groups at the origin using differential notions. However, applying Morse arguments in Banach spaces presents severe difficulties, as the classical Morse Lemma and the generalized Morse Lemma \cite{gromoll_meyer1969-t} are not directly applicable in this context. Moreover, the energy functional $J_a$ associated with \eqref{INTRODUZIONEpqNontrivialSolutions} is not of class $C^2$; hence the classical notion of Morse index cannot be defined. To overcome this lack of regularity, we introduce a suitable quadratic form $Q : X \to \mathbb{R}$, such that $m_0$ and $m_0^*$, defined as the supremum of the dimensions of the subspaces of $X$ where $Q$ is negative definite and negative semidefinite respectively, could serve as the Morse index and the large Morse index. Finally, we relate these differential notions $m_0$ and $m_0^*$ to the behavior of a penalized version of $J_a$, which satisfies $C^2$ regularity only within a ball centered at the origin. It is still open the singular cases and the mixed interaction of singular and degenerate quasilinear PDEs.

\medskip
In the third Chapter we pass to study spectral properties of the following eigenvalue problem for quasilinear elliptic systems, originally introduced by de Thélin in \cite{DT}: 

\begin{equation}\label{INTRODUZIONEEigenDeThélin0}
	\begin{cases}
		\begin{array}{ll}
			- \Delta_p u = \lambda |u|^{\alpha-1} |v|^{\beta+1} u
			&  \text{ in }\Omega,
			\medskip \\
			- \Delta_q v  = \lambda  |u|^{\alpha+1} |v|^{\beta - 1}v 
			&  \text{ in } \Omega, \medskip \\
			u=v=0  & \text{ on }  \partial\Omega,
		\end{array}
	\end{cases} \tag{I.5}
\end{equation}

\medskip
\noindent
where $\Omega$ is a bounded domain of $\mathbb{R}^N$ with smooth boundary, $N \geq 2$, $\lambda \in \R$,  $1 < p,q<N$ and $\alpha, \beta>0$ are real numbers satisfying $(\alpha +1)/ p + (\beta +1)/q =1$.
In \cite[Theorem 1]{DT} de Thélin proved that problem \eqref{INTRODUZIONEEigenDeThélin0} admits a smallest eigenvalue $\lambda_1>0$. In this thesis we show that the first eigenvalue $\lambda_1$ is simple in a suitable sense, a result largely contained in \cite{DT}, and subsequently we prove that $\lambda_1$ is isolated, following the approach in \cite[Theorem 1]{AH} and employing Picone’s identity (see \cite[Lemma 24]{AG}). Furthermore, we provide a variational characterization for a sequence of eigenvalues $\{\lambda_k\}_k$ by applying a suitable deformation lemma for $C^1$-submanifolds proved by Bonnet in \cite{BON}. In addition, we use spectral properties on $\lambda_1$ and $\lambda_2$   to address, in the spirit of the pionering paper \cite{LANDLAZ} of  Landesman and Lazer, the following non-autonomous perturbation of problem \eqref{INTRODUZIONEEigenDeThélin0}:

\begin{equation}\label{INTRODUZIONESyst0}
	\begin{cases}
		\begin{array}{ll}
			-\Delta_p u = \lambda_1 |u|^{\alpha-1} |v|^{\beta+1} u + \frac{1}{\alpha + 1} [F_s(x,u,v) - h_1(x) ]& x\in\Omega,
			\bigskip
			\\
			-\Delta_q v= \lambda_1  |u|^{\alpha+1} |v|^{\beta - 1}v  + \frac{1}{\beta + 1} [ F_t(x,u,v) - h_2(x)] & x\in \Omega,
			\bigskip
			\\
			u=v=0  & x\in \partial\Omega,
		\end{array}
	\end{cases} \tag{I.6}
\end{equation}

\medskip
\noindent
where $\Omega,N,p,q,\alpha$ and $\beta$ are defined as above. The nonlinearity  $F:\Omega \times \mathbb{R}^2 \to \mathbb{R}$ is a $C^1$-Carathéodory function, while $h_1 \in L^{p'}(\Omega)$ and  $h_2 \in L^{q'}(\Omega)$. Inspired from the work \cite{AO} of Arcoya and Orsina, we prove the existence of a weak solution for the perturbed problem \eqref{INTRODUZIONESyst0}. This is achieved under sufficient Landesman-Lazer conditions involving $p, q, F, h_1, h_2$, and the component-wise positive eigenfunction $(\varphi_1, \psi_1)$ associated with $\lambda_1$, normalized such that $\lVert \varphi_1 \rVert_{1,p} + \lVert \psi_1 \rVert_{1,q} = 1$. These conditions, that will be specified in detail in Chapter \ref{CapitolodeThelin}, ensure that the "nonlinear interference due to $F$" and the "external force" $(h_1,h_2)$ balanced out in a way that prevents the solution from sliding off to infinity. In particular, while the Euler functional associated with \eqref{INTRODUZIONESyst0} always satisfies the Palais-Smale (PS) condition, it is coercive only in certain cases. On the other hand, when coercivity is lacking, verifying that the functional exhibits the geometric structure required by the Saddle Point Theorem is not straightforward. Indeed, since the set of eigenfunctions associated with $\lambda_1$ for quasilinear systems is generally not a vector subspace of $X$, the standard splitting argument employed in \cite{AO} cannot be directly applied here. To overcome this difficulty, we exploit the isolation of $\lambda_1$ and the variational characterization of $\lambda_2$ to recover the saddle point geometry and thereby ensure the existence of a solution.

\medskip

In the quasilinear systems discussed previously, we focused on the $p$-Laplacian (or its generalizations) where $p < N$. However, the $N$-Laplacian operator $\Delta_N$ also plays a fundamental role in quasilinear elliptic equations with an underlying geometric structure. For instance, the $N$-Laplacian arises in the following $N$-dimensional Liouville equation

\begin{equation}\label{INTRODUZIONELiouvilleEquation}
\begin{cases}
\displaystyle - \Delta_N U(x) = e^{U(x)} & x \in \mathbb{R}^N, \bigskip\\
\displaystyle \int_{\mathbb{R}^N}   e^{U(x)} \, dx  < + \infty.
\end{cases} \tag{I.7}
\end{equation}

\medskip
\noindent
Denoting with $c_N:=N^N \left( \frac{N}{N-1} \right)^{N-1}$, this problem has an explicit solution given by

$$ \displaystyle U(x)= \ln \left(  \frac{c_N}{ \left( 1+ |x|^{\frac{N}{N-1}} \right)^N  }   \right).$$

\medskip
\noindent
Specifically, defining $\mu_N(x) := \frac{1}{V_N \left( 1+ |x|^{\frac{N}{N-1}} \right)^N}$, the function $v_N(x) := \ln \mu_N(x)$ is connected to Liouville equation \eqref{INTRODUZIONELiouvilleEquation} because

\begin{equation*}
 \displaystyle  -\Delta_N v_N(x) = c_N V_N \, e^{v_N(x)}.
\end{equation*}

\medskip
\noindent
The functions $\mu_N(x)$ and $v_N(x)$ arise in the $N$-dimensional Euclidean Onofri inequality established by Del Pino and Dolbeault in \cite{DD3}. By considering the endpoint of the optimal Gagliardo–Nirenberg interpolation inequalities, first discovered in \cite{DD1} and later extended in \cite{DD2}, they proved the $N$-dimensional Euclidean Onofri inequality for any $u \in C_0^{\infty}(\mathbb{R}^N)$:

\begin{equation}\label{INTRODUZIONEEuclideanOnofriNdimensional}
\displaystyle  \ln \left( \int_{\mathbb{R}^N} e^u \, d\mu_N \right) \leq  \frac{1}{\widetilde{\omega_N}} \int_{\mathbb{R}^N} H_N(u,\mu_N) \, dx + \int_{\mathbb{R}^N} u \, d\mu_N, \tag{I.8}
\end{equation}

\medskip
\noindent
where $ \displaystyle \widetilde{\omega_N}:= N^N \left( \frac{N}{N-1} \right)^{N-1} \omega_{N-1}$ and

\begin{equation*}
 H_N(u,\mu_N) :=| \nabla v_N + \nabla u |^N - |\nabla v_N|^N - N |\nabla v_N |^{N-2} \nabla v_N \cdot \nabla u. 
\end{equation*}

\medskip
\noindent
The function $v_N$ is involved in the definition of the aforementioned operator; specifically, the connection with the Liouville equation \eqref{INTRODUZIONELiouvilleEquation} is established via the term $N |\nabla v_N |^{N-2} \nabla v_N \cdot \nabla u$ and the Gauss-Green formula.
The $N$-dimensional Euclidean Onofri inequality was achieved with different techniques also by Agueh, Boroushaki and Ghoussoub in \cite{ABG} for any 
$u \in W^{1,N}(\mathbb{R}^N)$.  In this thesis, we are interested in determining the natural space in which such inequality remains valid.  We then introduce the  weighted Sobolev space 
\begin{align*}
\displaystyle 
W_{\mu_N}(\mathbb{R}^N):= \{ u \in L^1(\mathbb{R}^N,d\mu_N): \quad  |\nabla u| \in L^N(\mathbb{R}^N,dx),                     \quad |\nabla u|^2 |\nabla v_N|^{N-2} \in L^1(\mathbb{R}^N,dx) \},
\end{align*}

\medskip
\noindent
endowed with the  norm
\begin{equation*}
\lVert u \rVert_{\mu_N}:= \int_{\mathbb{R}^N} |u| \, d\mu_N + \lVert \nabla u \rVert_N + \left( \int_{\mathbb{R}^N}  |\nabla u|^2 |\nabla v_N|^{N-2} \, dx \right)^{\frac{1}{2}}.
\end{equation*}

\medskip
\noindent
We extend inequality \eqref{INTRODUZIONEEuclideanOnofriNdimensional} to this setting by proving that $C_0^{\infty}(\mathbb{R}^N)$ is dense in $W_{\mu_N}(\mathbb{R}^N)$ with respect to the norm $\lVert \cdot \rVert_{\mu_n}$. Furthermore, we  establish an equivalence with the following sharp logarithmic Moser-Trudinger inequality on the unit ball $B_1 \subset \mathbb{R}^N$, which can be derived from the work of Carleson and Chang \cite{CarlesonChang}:

\begin{equation}\label{CCIntroThesis}
\displaystyle \ln \left( \frac{1}{|B_1|} \int_{B_1} e^{u}   dx  \right) < \frac{1}{\widetilde{\omega_N}} \int_{B_1} |\nabla u|^N  dx + \sum_{k=1}^{N-1} \frac{1}{k} \qquad \forall \, u \in  W^{1,N}_0(B_1). \tag{I.9}
\end{equation}

\medskip
\noindent
The case $N=2$ is of particular interest from both mathematical and physical perspectives. In fact inequality \eqref{INTRODUZIONEEuclideanOnofriNdimensional} becomes

\begin{equation*}
\displaystyle  \ln \left(\int_{\mathbb{R}^2} e^u \, d \mu_2 \right) \leq  \frac{1}{16 \pi} \int_{\mathbb{R}^2} | \nabla u |^2 \, dx + \int_{\mathbb{R}^2} u \, d \mu_2
\end{equation*}

\medskip
\noindent
for any $ u \in  W_{\mu_2}(\mathbb{R}^2):= \{ u \in L^1(\mathbb{R}^2,d\mu_2) \, : \, |\nabla u| \in L^2(\mathbb{R}^2,dx) \}$. Furthermore, via stereographic projection, this inequality is equivalent to the Onofri inequality (see \cite{O}) on the two-dimensional unit sphere $\mathbb{S}^2$ equipped with the standard metric $g_0$:

\begin{equation}\label{ONOFRI2DIMENSIONALE}
\ln \left(\frac{1}{4\pi}\int_{\Sm^2} e^u \, dv_{g_0} \right) \leq  \frac{1}{16 \pi} \int_{\Sm^2} | \nabla_{g_0} u |_{g_0}^2 \, dv_{g_0} + \frac{1}{4 \pi} \int_{\Sm^2}u \, dv_{g_0} \qquad \forall \, u \in W^{1,2}(\mathbb{S}^2,g_0). \tag{I.10}
\end{equation}  

\medskip
\noindent
Onofri inequality plays an important role in spectral analysis for the Laplace-Beltrami operator thanks  to Polyakov's formula (see \cite{OPS,OPS2,Poly1,Poly2}) and it was proved in 1982 by Onofri in \cite{O}  using conformal invariance and an earlier result by Aubin  \cite{Aubin}.
Inequality \eqref{ONOFRI2DIMENSIONALE} can be obtained via several methods, including various limiting procedures, mass transportation methods in the radial case, and rigidity methods associated with a related nonlinear flow (we refer to the papers by Dolbeault, Esteban, Jankowiak, and Tarantello \cite{DEJ, DET}, and the references therein for full proofs and further details). The Onofri inequality also has important applications; for instance, in chemotaxis it appears in the study of the Keller-Segel model (see \cite{CALVCORR, GAJZAK}).

\medskip
\noindent
Concerning inequality \eqref{CCIntroThesis}, in dimension $N=2$ it becomes

\begin{equation}\label{INTROLOGccDUEDIM}
\displaystyle \ln \left( \frac{1}{\pi} \int_{B_1} e^{u} \,   dx   \right) < \frac{1}{16 \pi} \int_{B_1} |\nabla u|^2  dx + 1 \qquad \forall \, u \in  W^{1,2}_0(B_1). \tag{I.11}
\end{equation}

\medskip
\noindent
A generalization of previous inequality was provided in \cite{CLMP} (see also \cite{CCL}). Specifically, if $\Omega \subseteq \mathbb{R}^2$ is a bounded domain with smooth boundary, then for any $u \in W_0^{1,2}(\Omega)$ we have 
 
\begin{equation}\label{INTROTHESIOnofriDomain}
\ln \left(\frac{1}{|\Omega|}\int_{\Omega} e^u \,dx \right) \le  \frac{1}{16\pi}\int_{\Omega} |\nabla u|^2 dx  + 1 + {4\pi \sup_{\Omega}\gamma(x) + \ln \frac{\pi}{|\Omega|}}, \tag{I.12}
\end{equation}

\medskip
\noindent
where $\gamma$ is the Robin's function of $\Omega$. This inequality plays a crucial role in the study of the mean field equation 

\begin{equation}\label{INTROmean}
-\Delta u  = \frac{\rho e^u}{\int_{\Omega}e^u \,dx} \quad \text{ in }\Omega, \tag{I.13}
\end{equation}

\medskip
\noindent
which appears in the statistical mechanic description of vortex formations in 2d-models for turbulent flows in fluid dynamics (see e.g. \cite{JM,CLMP}). We also mention \cite{CL1,CL2,Malchiodi,Bart,BarMal} for existence results  for  \eqref{INTROmean} and related problems.\\

The thesis is organized as follows. Chapter \ref{CAPITOLOPOINCARE} is devoted to show the finiteness of the critical groups and a local Poincar\'e-Hopf type formula for functionals associated to quasilinear elliptic systems. In Chapter \ref{SezReg} we present new results concerning $L^{\infty}$-estimates for solutions to non-autonomous quasilinear systems involving operators in divergence form and nonlinearities up to critical growth. Building upon results accomplished in previous Chapters, in Chapter \ref{CAPAMANNZEHNDER} we deal with Amann-Zehnder type results for quasilinear elliptic systems, in order to detect a nontrivial weak solution in both resonant and nonresonant case. Chapter \ref{CapitolodeThelin} is devoted to Landesman-Lazer type conditions for quasilinear elliptic systems, after proving spectral properties on the de Thélin eigenvalue problem associated. Finally, Chapter \ref{CAPITOLOONOFRI} is devoted to extending the $N$-dimensional Euclidean Onofri inequality to the weighted Sobolev space $W_{\mu_N}(\mathbb{R}^N)$. Furthermore, we establish its equivalence with the sharp logarithmic Moser–Trudinger inequality on the unit ball $B_1 \subset \mathbb{R}^N$, providing an insight into the case $N=2$ that highlights the underlying geometric link with the two-dimensional unit sphere $\mathbb{S}^2$.\\

This thesis is based on the papers \cite{ABC,BCM,BCM2,BCV,BCV2,BCV3}.

\newpage

\mainmatter

\chapter{A Poincar\'e-Hopf formula for functional associated to Quasilinear Elliptic Systems}\label{CAPITOLOPOINCARE}

\medskip
\noindent
In this Chapter we present some new results about the finiteness of the critical groups and a local Poincar\'e-Hopf type formula for functionals associated to quasilinear elliptic systems, involving  $p$-Laplacian or $p$-area operators. In particular, we show that the Fréchet derivative of such a functionals is locally of class $(S)_+$ when the nonlinearity exhibits critical growth, so that we extend the results proved in \cite{CD1} by Cingolani and Degiovanni for quasilinear elliptic equations at critical growth.

\medskip
This Chapter is based on the paper \cite{BCV2}.

\bigskip

\section{Browder degree in Banach spaces}

\bigskip
\noindent
In the context of Hilbert spaces, the celebrated Poincar\'e-Hopf formula has been extended for functionals $f$ such that the map $\nabla f$ is a compact perturbation of the identity. 
This allows to use the Leray-Schauder degree theory, instead of the Brouwer degree for $C^1$ functions $f:\re^N \to \re$
(see e.g. \cite[Theorem~8.5]{mawhin_willem1989},
~\cite[Theorem~3.2]{chang},
\cite[Theorem~3.2]{li_li_liu2005}, \cite[Theorem~1.3]{bartsch_dancer2009} and  \cite{benci1991}).

Conversely if $Y$ is a Banach (not Hilbert) space, then
the derivative $f'$ is naturally defined from $Y$ into $Y'$ and this is not the usual setting of Leray-Schauder degree.
In \cite{kim_wang1989},  Kim and Wang  obtained a Poincar\'e-Hopf type result in the special case  in which the Banach space is densely embedded in a Hilbert space and the functional $f$ induces a vector field $\nabla f$ on the Banach space through the inner product of the Hilbert space. They applied extensions of the Leray-Schauder degree for maps defined from the Banach space into itself.

In  \cite{CD1}, Cingolani and Degiovanni established a Poincar\'e-Hopf formula in a general reflexive Banach space $Y$. Their main idea was to 
replace the Leray-Schauder degree by the  Browder degree developed in~\cite{browder1983} for maps that are demicontinuous (namely continuous from the strong topology of the domain to the weak topology od the codomain) and of class~$(S)_+$, that is:

\begin{definition} 
Let $(Y, \lVert \cdot \rVert)$ be a reflexive Banach space, whose dual is denoted by $Y'$. Let $D\subseteq Y$ and let us consider a map $F:D\longrightarrow Y'$. $F$ is said to be \emph{of class~$(S)_+$} if, for every sequence $\{w_k\}_k$ in $D$ weakly convergent to $w$ in $Y$ with
		\[
		\limsup_k \,\langle F(w_k),w_k-w\rangle \leq 0\,,
		\]
we have $\|w_k-w\|\to 0$.
\end{definition}

We notice that if $Y$ is a Hilbert space and $Y'$ is identified with $Y$ in
the standard way, then any compact perturbation of the identity
is of class~$(S)_+$. However, the class~$(S)_+$ allows to consider further
interesting cases, also in the Hilbert setting.
\par
As noticed in \cite{CD1,degiovanni2009},  any compact perturbation of a map of
class~$(S)_+$ is still of class~$(S)_+$. 
Therefore the class includes the perturbations of the
$p$-Laplace operator by terms with subcritical growth. 
We refer to \cite[Theorem~10]{dinca_jebelean_mawhin2001} where $-\Delta_p:W^{1,p}_0(\Omega)\longrightarrow W^{-1,p'}(\Omega)$
is continuous and of class~$(S)_+$. 
In \cite{CD1} the case with critical growth has been covered, which is interesting also for $p=2$.
We also mention the paper \cite{ad}  where Almi and Degiovanni showed how the degree for maps of class $(S)_+$ can be used to
define, by a suitable approximation technique, a degree for quasilinear elliptic
equations with natural growth conditions.

In this Chapter we  prove the finiteness of the critical groups and a local Poincar\'e-Hopf type formula for the functional $J_{a,b}: X \to \mathbb{R}$ defined for any $z=(u,v) \in X$ as
  
\begin{align}\label{functionalPoincareHopf}
	J_{a,b}(z)  =  \frac{1}{p} \into  \left(a+
	|\nabla u|^2\right)^{\frac p 2} \ dx + \frac{1}{q}\into \left(b +
	|\nabla v|^2\right)^{\frac q 2} \ dx  - \into
	F(u,v) \ dx,
\end{align}

\medskip
\noindent
where   $1<p,q <N$ and  $a,b\geq 0$.  Moreover, we assume that $F\in C^{1}(\mathbb{R}^2)$ and that there exists $C>0$ such that for any $(s,t) \in \re^2$ we have

\begin{equation}\label{primaPoincareHopf}
	| F_s(s,t)|\leq C
	\left(|s|^{p^*-1}+|t|^{q^{*}\frac{p^*-1}{p^{*}}}+1\right),
\end{equation}

\begin{equation}\label{secondaPoincareHopf}
	|F_t(s,t)| \leq C
	\left(|s|^{p^*\frac{q^{*}-1}{q^{*}}}+|t|^{q^*-1}+ 1 \right).
\end{equation} 

\medskip
\noindent
Previous assumptions and Nemytskii operator theory give that $J_{a,b}$ is of class $C^1$ on $X$.

\medskip
As in \cite{CD1} for the $p$-Laplace operator (see also \cite{ad,CDV}), we show that also the case with critical growth for the functional $J_{a,b}$, associated to a quasilinear system, gives rise locally to an operator of class $(S)_+$. 
However, in our setting we have to face some additional difficulties in order to show that the Fréchet derivative $J'_{a,b}$ of the Euler functional associated to a quasilinear system is locally an operator  of class $(S)_+$. More precisely, we generalize to a $(p,q)$-area type operator which returns us the $(p,q)$-Laplacian as a special case. About the nonlinearity $F$, in a quasilinear system setting  new difficulties arise in computations due to the coupling of $u$ with $v$, and we are able to manage them also when $F$ is allowed to grow critically. Therefore, we partially answer to question $(i)$ in the final comments of \cite{boccardodefiguerido}, where nonlinearity that exhibits critical growth is not treated because of a loss of compactness.

Precisely, we will prove the following.

\begin{theorem}\label{JabS+}

Let $\Omega \subset \mathbb{R}^N$ be a smooth bounded domain.
Let $1<p,q<N$ and $F:\mathbb{R}^2 \to \mathbb{R}$ a $C^1$-function satisfying \eqref{primaPoincareHopf} and \eqref{secondaPoincareHopf}.
Let $a,b \geq 0$ and  $J_{a,b}$ the $C^1$-functional defined in \eqref{functionalPoincareHopf}.
Then there exists a radius $R:=R(N,p,q,C)$ such that, for every $z_0=(u_0,v_0) \in X$, the map $J'_{a,b}$ is of class $(S)_+$ on
 
$$ \displaystyle \overline{B_R(z_0)}:=\left\lbrace z \in X \; : \; \lVert z - z_0 \rVert \leq R \right\rbrace.$$

\end{theorem}

\medskip
\noindent
Let us now recall that  for any $q \in \mathbb{N}$, the $q$-th critical group of a functional $f$ at $\bar{u} \in \text{Dom} f$ is defined as

\[
C_q(f,\bar{u}) =
H_q\left(\{f\leq f(\bar{u})\},\{f\leq f(\bar{u})\}\setminus\{\bar{u}\}\right)\,,
\]

\medskip
\noindent
where $H_*$ denotes singular homology (see e.g.~\cite{chang1981,chang, mawhin_willem1989,rybakowski1987,spa}).

\medskip
By applying  Theorem \ref{JabS+} and retracing the same ideas contained in \cite[Theorem 1.1]{CD1},   we also obtain the finiteness of the critical groups and a local Poincar\'e-Hopf  formula for the functional $J_{a,b}$, that is:

\begin{theorem}\label{TheorFinGrupCrit}
If $\bar{z}$ is an isolated critical point of the functional $J_{a,b}$ defined in \eqref{functionalPoincareHopf}, then $C_*(J_{a,b},\bar{z})$ is of finite type and we have

\begin{equation*}
\mathrm{deg}(J_{a,b}',V,0) = \sum_{m\geq 0}
(-1)^m\,\rnk{C_m(J_{a,b},\bar{z})}
\end{equation*}

\medskip
\noindent
for every sufficiently small neighborhood $V$ of $\bar{z}$.
\end{theorem}

We mention \cite{CCMV} for the computation of the critical groups of $J_{a,b}$ when $a=b$ and $2 < p,q <N$. Furthermore, in Chapter \ref{CAPAMANNZEHNDER} the computation of the critical groups is needed to obtain Amann-Zender type results for quasilinear systems. It is still an open problem the evaluation of the critical groups for the functional $J_{a,b}$ when $p$ or $q$ are less than $2$, or $a \neq b$. Let us refer the reader to   \cite{AFTPAC,ASSCING,CD,CDS,CDV,CV,CV2} for critical groups estimates for quasilinear elliptic equations (see also \cite{PAO}).\\

The following Chapter is organized as follows: in Section \ref{Sezpareaqarea} we show that $p$-area and $q$-area operators are of class $(S)_+$, even when we consider their simultaneous action in the product space $X$. Section \ref{SezFunzPoincare} is devoted to the proof of Theorem \ref{JabS+}. Finally, Section 
\ref{finitezzagruppicritici} is devoted to the proof of Theorem \ref{TheorFinGrupCrit}.

\bigskip

\section{$(p,q)$-area operators}\label{Sezpareaqarea}

\bigskip
\noindent
We begin to recognize that  the $r$-area operator is of class $(S)_+$ for $r>1$.
More precisely, for any  $c \geq 0$ and $r > 1$, we consider the map $H_{c,r}\!:\!\sobr\to \wr$ defined for any $w \in \sobr$ by

$$	H_{c,r}(w) \, := \, -\div\left(\left(c +|\nabla w|^2\right)^{\frac{r-2}{2}}\nabla w\right),$$

\medskip
\noindent
where the action of the map $H_{c,r}(w): \sobr \to \mathbb{R}$ is given for any $\phi \in \sobr$ by

$$ \displaystyle \langle H_{c,r}(w), \phi \rangle := \into  \left(c +|\nabla w|^2\right)^{\frac{r-2}{2}}\nabla w \cdot \nabla \phi \, dx.     $$

\begin{lemma}\label{s+PoincareHopf}
For any  $c\geq 0$ and $r > 1$, the map $H_{c,r}: \sobr\to \wr$  is of class $(S)_+$.
\end{lemma}

\begin{proof}
To get the thesis, we will apply \cite[Theorem 3.5]{ad}.\\ 
Let  $c\geq 0$ and $r > 1$ and let us consider the map $d_{c,r}:\re^N \to \re^N$ defined as

\begin{equation*}
\displaystyle  d_{c,r}(\xi)=(c+|\xi|^2)^{\frac{r-2}{2}}\xi. 
\end{equation*}

\medskip
\noindent
We know that the map $d_{c,r}$ is strictly monotone for any $c \geq 0$ and for any $r>1$ (see Lemma \ref{smc}  for the proof). Now, let us assume $r \geq 2$. It is immediate to see that for any $\xi \in \mathbb{R}^N$ we have

$$ \displaystyle d_{c,r}(\xi) \cdot \xi \geq |\xi|^r.$$

\medskip
\noindent
On the other hand, we have

$$ \displaystyle \lim_{|\xi| \to + \infty} \frac{\left| d_{c,r}(\xi) \right|}{|\xi|^{r-1}} =  \lim_{|\xi| \to + \infty} \frac{(c+|\xi|^2)^{\frac{r-2}{2}}|\xi|}{|\xi|^{r-1}}  = 1, $$ 

\medskip
\noindent
hence there exists a constant $C$, depending on $c$, such that for any $\xi \in \mathbb{R}^N$ we have

$$ \displaystyle \left| d_{c,r}(\xi) \right| \leq C + C|\xi|^{r-1}.$$

\medskip
\noindent
If $1<r<2$, then it is immediate to see that for any $\xi \in \mathbb{R}^N$ we have

$$ \displaystyle \left| d_{c,r}(\xi) \right| \leq |\xi|^{r-1}.$$

\medskip
\noindent
On the other hand, we have

$$ \displaystyle \lim_{|\xi| \to + \infty} \frac{ d_{c,r}(\xi) \cdot \xi }{|\xi|^{r}} =  \lim_{|\xi| \to + \infty} \frac{(c+|\xi|^2)^{\frac{r-2}{2}} |\xi|^2}{|\xi|^{r}}  = 1, $$ 

\medskip
\noindent
hence there exists a constant $C$, depending on $c$, such that for any $\xi \in \mathbb{R}^N$ we have

$$ \displaystyle d_{c,r}(\xi) \cdot \xi \geq C |\xi|^r + C.$$

\medskip
\noindent
In any case, for any $c \geq 0$ and for any $r>1$ assumptions of \cite[Theorem 3.5]{ad} are satisfied, by which we get the thesis.

\end{proof}

\medskip
\noindent
Now, for any  $a, b \geq 0$ and $p,q > 1$,  we consider the map $H_{a,p,b,q}: X \to X'$ defined for any $z=(u,v) \in X$ as

\begin{align}\label{Hapbq}
\displaystyle
	 H_{a,p,b,q}(z) = \left[ H_{a,p}(u) , H_{b,q}(v) \right], 
\end{align}

\medskip
\noindent
where the action of the map $H_{a,p,b,q}(z): X \to \mathbb{R}$ is given for any $ \bar{z}=(\bar{u},\bar{v}) \in X$ by

\begin{align*}
\displaystyle
 \langle H_{a,p,b,q}(z), \bar{z} \rangle  =  \into  \left(a +|\nabla u|^2\right)^{\frac{p-2}{2}}\nabla u \cdot \nabla \bar{u} \, dx  + \into  \left(b +|\nabla v|^2\right)^{\frac{q-2}{2}}\nabla v \cdot \nabla \bar{v} \, dx.
\end{align*}

\begin{proposition}\label{HapbqS+}
The map $H_{a,p,b,q}: X \to X'$  is of class $(S)_+$.
\end{proposition}

\begin{proof}

Let $\{z_k\}_k=\{(u_k,v_k)\}_k \subset X$ and $z=(u,v) \in X$ such that $z_k \rightharpoonup z$ in $X$ and

\begin{equation}\label{ipotesi1}
 \displaystyle \limsup_{k \to \infty} \, \langle  H_{a,p,b,q}(z_k), z_k - z  \rangle \leq 0.
\end{equation}

\noindent
We need to prove $z_k \to z$ in $X$.\\
Since $v_k \rightharpoonup v$ in $W_0^{1,q}(\Omega)$ and $H_{b,q}(v) \in W_0^{-1,q'}(\Omega)$,  we have

\begin{equation}\label{conto1}
\displaystyle \into  \left(b +|\nabla v|^2\right)^{\frac{q-2}{2}}\nabla v \cdot \nabla \left(v_k - v \right) \, dx  \to 0.
\end{equation}

\noindent
By using the monotonicity of the map $d_{b,q}:\re^N\to \re^N$ defined in \eqref{definizionemappadcr},
we deduce

\begin{align*}
\displaystyle
0 
& \leq \into \left[ d_{b,q}(\nabla v_k)-d_{\beta,q}(\nabla v) \right] \cdot \nabla \left( v_k - v \right) \, dx \\
& =\into  \left(b +|\nabla v_k|^2\right)^{\frac{q-2}{2}}\nabla v_k \cdot \nabla \left(v_k - v \right) \, dx  - \into  \left(b +|\nabla v|^2\right)^{\frac{q-2}{2}}\nabla v \cdot \nabla \left(v_k - v \right) \, dx ,
\end{align*}

\noindent
hence

$$ \displaystyle \into  \left(b +|\nabla v|^2\right)^{\frac{q-2}{2}}\nabla v \cdot \nabla \left(v_k - v \right) \, dx \leq \into  \left(b +|\nabla v_k|^2\right)^{\frac{q-2}{2}}\nabla v_k \cdot \nabla \left(v_k - v \right) \, dx  $$

\medskip
\noindent
and so 

\begin{align}
\displaystyle
  & \into  \left(a+|\nabla u_k|^2\right)^{\frac{p-2}{2}}\nabla u_k \cdot \nabla \left(u_k - u \right) \, dx \nonumber +  \into  \left(b +|\nabla v|^2\right)^{\frac{q-2}{2}}\nabla v \cdot \nabla \left(v_k - v \right) \, dx  \nonumber \\
\leq & \into  \left(a +|\nabla u_k|^2\right)^{\frac{p-2}{2}}\nabla u_k \cdot \nabla \left(u_k - u \right) \, dx  +  \into  \left(b +|\nabla v_k|^2\right)^{\frac{q-2}{2}}\nabla v_k \cdot \nabla \left(v_k - v \right) \, dx  \\
= &  \langle  H_{a,p,b,q}(z_k), z_k - z  \rangle. \nonumber
\end{align}

\noindent
Taking into account last inequality, \eqref{ipotesi1} and \eqref{conto1}, we deduce

$$ \displaystyle \limsup_{k \to \infty} \, \langle  H_{a,p}(u_k), u_k - u  \rangle \leq  \limsup_{k \to \infty} \, \langle  H_{a,p,b,q}(z_k), z_k - z  \rangle \leq 0, $$

\noindent
hence, by Lemma \ref{s+PoincareHopf}, we deduce $u_k \to u$ in $W_0^{1,p}(\Omega)$. Similarly, we deduce $v_k \to v$ in $W_0^{1,q}(\Omega)$.\\
Finally, we have $z_k \to z$ in $X$.

\end{proof}

\bigskip

\section{$(S)_+$-condition}\label{SezFunzPoincare}

\medskip
\noindent
Let us now consider the following quasilinear system

\begin{equation}\label{pqPoincareHopf}
	\begin{cases}
		\begin{array}{ll}
			-\text{\rm div} \left( (a+|\nabla u|^2)^{\frac{p-2}{2}})\nabla u\right) =F_s(u,v) & x\in\Omega,
			\medskip \\
			-\text{\rm div} \left( (b +|\nabla v|^2)^{\frac{q-2}{2}})\nabla v\right)=F_t(u,v) & x\in \Omega, \medskip \\
			u=v=0  & x\in \partial\Omega,
		\end{array}
	\end{cases}
\end{equation}

\medskip
\noindent
where $N \geq 2$, $1<p,q<N$, and $a,b\geq 0$.
We recall that the Euler functional $ J_{a,b}: X \to \mathbb{R}$ associated to \eqref{pqPoincareHopf} is defined for any $z=(u,v) \in X$ as

\begin{align*}
	J_{a,b}(z)  =  \frac{1}{p} \into  \left(a+
	|\nabla u|^2\right)^{\frac p 2} \ dx + \frac{1}{q}\into \left(b +
	|\nabla v|^2\right)^{\frac q 2} \ dx  - \into
	F(u,v) \ dx.
\end{align*}

\medskip
\noindent
Moreover,  $J_{a,b}$ is of class $C^1$ on $X$ and for any $z_0=(u_0,v_0), z=(u, v) \in X$ it results

\begin{align*}
	\langle J'_{a,b}(z_0), z \rangle
& = \displaystyle  \into (a+|\nabla
	u_0|^2)^{\frac {p-2}{2}}\nabla u_0\nabla u \, dx +\into
	(b+|\nabla v_0|^2)^{\frac {q-2}{2}}\nabla v_0 \nabla v \, dx \\
& -\displaystyle\into ( F_s(u_0,v_0) u + F_t(u_0,v_0) v)
	dx.
\end{align*}

\medskip
\noindent
In this Section we show that the map $J_{a,b}':X \to X'$ is locally of class $(S)_+$.

\begin{remark}\label{TroncataRresto}
For every $h \in \mathbb{N},$ $h \geq 1$ we define $\mathcal{T}_h,\mathcal{R}_h:\mathbb{R} \to \mathbb{R}$ as 

\begin{equation}\label{ThRh}
\displaystyle \mathcal{T}_h(s):= \min\{\max\{s,-h\},h\} \qquad \text{and}  \qquad \mathcal{R}_h(s):=s-\mathcal{T}_h(s).
\end{equation}

\noindent
If $w \in W_0^{1,r}(\Omega)$ for some $r>1$,  then $\mathcal{T}_h(w), \, \mathcal{R}_h(w) \in W_0^{1,r}(\Omega)$. Moreover, we have

\begin{equation*}
\qquad \mathcal{T}_h(w(x))=
\begin{cases}
   h       & \text{ a.e. in } \{ x \in \Omega \; : \; w(x)\geq  h \},\\
    w(x)             & \text{ a.e. in } \{ x \in \Omega \; : \; - h < w(x) < h \},\\
 - h        & \text{ a.e. in } \{ x \in \Omega  \; : \; w(x) \leq - h \};
\end{cases}
\end{equation*}

\medskip

\begin{equation*}
 \nabla \mathcal{T}_h(w(x))=
\begin{cases}
 \nabla w(x)        & \text{ a.e. in } \{ x \in \Omega \;  : \; |w(x)| <  h \}, \medskip\\
0                   & \text{ a.e. in } \{ x \in \Omega \;  : \; |w(x)| \geq h \};
\end{cases}
\end{equation*}

\bigskip

\begin{equation*}
\qquad \qquad \mathcal{R}_h(w(x))=
\begin{cases}
w(x) -  h       & \text{ a.e. in } \{ x \in \Omega \; : \; w(x)\geq  h \},\\
0               & \text{ a.e. in } \{ x \in \Omega \; : \; - h < w(x) < h \},\\
w(x) + h        & \text{ a.e. in } \{ x \in \Omega  \; : \; w(x) \leq - h \};
\end{cases}
\end{equation*}

\medskip

\begin{equation*}
 \nabla \mathcal{R}_h(w(x))=
\begin{cases}
 \nabla w(x)        & \text{ a.e. in } \{ x \in \Omega \;  : \; |w(x)| >  h \}, \medskip\\
0                   & \text{ a.e. in } \{ x \in \Omega \;  : \; |w(x)| \leq h \};
\end{cases}
\end{equation*}

\medskip
\noindent
and 

\begin{equation}\label{GradhRh}
\qquad \qquad \nabla w(x)=
\begin{cases}
 \nabla \mathcal{T}_h(w(x))        & \text{ a.e. in } \{ x \in \Omega \; : \; |w(x)|< h \},\\
 \nabla \mathcal{R}_h(w(x))                  & \text{ a.e. in } \{ x \in \Omega \; : \; |w(x)| > h \},\\
 0       & \text{ a.e. in } \{ x \in \Omega  \; : \; |w(x)|=h \}.
\end{cases}
\end{equation}

\end{remark}

\medskip
\noindent
We will use the following result.

\begin{lemma}\label{convergenzasopralivelli}
 Let $\Omega$ be a bounded domain of $\mathbb{R}^N$.
Let $\{w_k\}_k \subset L^1(\Omega)$ and $w \in L^1(\Omega)$ such that $w_k(x) \to w(x)$ a.e. $x \in \Omega$ as $k \to \infty$. Let $h \in \mathbb{N}$, $h \geq 1$ and let us define, almost everywhere, the sets

\begin{equation*}
\displaystyle
\Omega^{w_k}_{h,+}:= \{ x \in \Omega \; : \;  |w_k(x)| \;  \geq \;  h   \} \qquad \text{and} \qquad \Omega^{w}_{h,+}:= \{ x \in \Omega \; : \;  |w(x)| \;  \geq \;  h   \}.
\end{equation*} 

\medskip
\noindent
Then 

$$ \displaystyle \limsup_{k \to \infty} \left| \Omega^{w_k}_{h,+}  \right| \leq \left| \Omega^{w}_{h,+}  \right|.$$ 

\end{lemma}

\begin{proof}
Since $w_k(x) \to w(x)$ a.e. $x \in \Omega$ as $k \to \infty$,  $w_k$ converges to $w$ also in measure, that is

\begin{equation}\label{convmisura}
 \forall \varepsilon >0, \qquad \left| \{ x \in \Omega \, : \, \left| w_k(x) - w(x) \right| \geq \varepsilon    \} \right| \to 0 \qquad \text{ as } k \to \infty.
\end{equation}

\medskip
\noindent
Let $\varepsilon >0$ and let us define
 
\begin{align*}
\displaystyle
D_{k,h} & := \{x \in \Omega \, : \, \left| w_k(x) - w(x) \right|  + \left|w(x) \right| \geq h    \}; \\
A_{k,h}^{\varepsilon} & := \{x \in \Omega \, : \, \left|w(x)\right| \geq h - \varepsilon, \quad\left| w_k(x) - w(x) \right|  + \left|w(x) \right| \geq h    \}; \\
B_{k,h}^{\varepsilon} & := \{x \in \Omega \, : \, \left|w(x)\right| < h - \varepsilon, \quad \left| w_k(x) - w(x) \right|  + \left|w(x) \right| \geq h    \}; \\
\mathcal{A}_{h}^{\varepsilon} & := \{x \in \Omega \, : \, \left|w(x)\right| \geq h - \varepsilon  \}; \\
\mathcal{B}_{k,h}^{\varepsilon} & := \{x \in \Omega \, : \,  \left| w_k(x) - w(x) \right| \geq  \varepsilon    \}. 
\end{align*}

\noindent
Since

$$ \displaystyle \Omega^{w_k}_{h,+} \subseteq D_{k,h} = A_{k,h}^{\varepsilon} \cup B_{k,h}^{\varepsilon} \subseteq \mathcal{A}_{h}^{\varepsilon} \cup \mathcal{B}_{k,h}^{\varepsilon},$$

\noindent
for any $k,h \in \mathbb{N}^*$ and for any $\varepsilon>0$ we have

$$ \displaystyle \left| \Omega^{w_k}_{h,+} \right| \leq \left| \mathcal{A}_{h}^{\varepsilon} \right| + \left| \mathcal{B}_{k,h}^{\varepsilon} \right|.$$

\noindent
By \eqref{convmisura}, we deduce

$$ \displaystyle  \left| \mathcal{B}_{k,h}^{\varepsilon} \right| \to 0 \qquad \text{ as } k \to \infty \qquad \forall \varepsilon >0,$$

\noindent
by which

\begin{equation}\label{c}
 \displaystyle   \limsup_{k \to \infty} \left| \Omega_{k,h}^{+} \right| \leq 
\left| \mathcal{A}_{h}^{\varepsilon} \right|, \qquad \forall \varepsilon >0.
\end{equation}

\medskip
\noindent
Since $\left(\mathcal{A}_{h}^{1/n} \right)_{n \in \mathbb{N}}$ is a decreasing sequence of subsets with $\left| \mathcal{A}_{h}^{1} \right| < \infty$ and 

$$ \displaystyle \bigcap_{n=1}^{\infty} \mathcal{A}_{h}^{1/n} = \Omega^{w}_{h,+},$$

\medskip
\noindent
by continuity of Lebesgue measure we infer

$$ \displaystyle \lim_{n \to \infty} \left| \mathcal{A}_{h}^{1/n} \right| = \left| \Omega^{w}_{h,+} \right|,$$

\medskip
\noindent
by which, together with \eqref{c}, we deduce

$$\displaystyle \limsup_{k \to \infty} \left| \Omega^{w_k}_{h,+} \right| \leq \left| \Omega^{w}_{h,+} \right|.$$

\end{proof}

\medskip
\noindent
We will frequently use also the following result of functional analysis (see \cite{BREZIS}).

\begin{proposition}\label{weakconvLp}
Let $1< r < \infty$ and $(\Omega,\mathcal{M}, \mu)$ a measure space with $\sigma$-finite measure $\mu$.\\
Let $\{f_k\}_k \subset L^r(\Omega)$ a sequence such that:

\begin{itemize}
\item[i)] $\{f_k\}_k$ is bounded in $L^r(\Omega)$; \medskip

\item[ii)] $f_k(x) \to f(x)$ a.e. $x \in \Omega$. 
\end{itemize}

\medskip
\noindent
Then 

$$f_k \rightharpoonup f \text{ in } L^r(\Omega).$$
\end{proposition}

\medskip
\noindent
{\mbox {\it Proof of Theorem~\ref{JabS+}.~}} Let $a,b \geq 0$ and $z_0=(u_0,v_0) \in X$. Let $R>1$ and let us consider a sequence $\{z_k\}_k=\{(u_k,v_k)\}_k \subset \overline{B_R(z_0)}$  weakly convergent to some $z=(u,v)$ in $X$ with

\begin{equation}\label{ipotesi}
\displaystyle \limsup_{k\to \infty} \, \langle J_{a,b}'(z_k), z_k-z\rangle\leq 0.
\end{equation} 

\noindent
By means of the map $H_{a,p,b,q}$ defined in \eqref{Hapbq}, we can rewrite the preceding assumption as follows:

\begin{align}\label{ipotesiHapbq}
\displaystyle 
\limsup_{k \to \infty}
\, \left\lbrace \langle H_{a,p,b,q}(z_k),z_k-z    \rangle 
 - \into \left[ F_s(u_k,v_k) (u_k - u) + F_t(u_k,v_k)(v_k - v) \right]
dx  \right\rbrace \leq 0.
\end{align}

\medskip
\noindent
Let $c \geq0$, $r>1$ and let us consider the functional $f:W_0^{1,r}(\Omega) \to \mathbb{R}$ defined as

\begin{equation}\label{funzionalegenerico}
\displaystyle f(w):= \frac{1}{r} \into \left( c + | \nabla w |^2 \right)^{\frac{r}{2}} \, dx,
\end{equation} 

\noindent
Since $f$ is convex (see Lemma \ref{convexity}), for any $w,w_k \in W_0^{1,r}(\Omega)$ we have

$$ \displaystyle f(w) \geq f(w_k) + \langle f'(w_k),w-w_k\rangle,$$

\medskip
\noindent
by which

\begin{align*}
 \displaystyle \into \left( c + | \nabla w_k |^2\right)^{\frac{r-2}{2}} \nabla w_k \cdot \nabla( w_k - w) \, dx  
\geq \frac{1}{r} \into \left( c + | \nabla w_k|^2\right)^{\frac{r}{2}} \, dx  
 - \frac{1}{r} \into \left( c + | \nabla w |^2 \right)^{\frac{r}{2}} \, dx.
\end{align*}

\noindent
Hence, we get

\begin{align}\label{contoPSlocale}
\displaystyle 
0  \geq \limsup_{k \to \infty} & \; \langle J_{a,b}'(z_k), z_k-z\rangle \\
 \geq \limsup_{k \to \infty} & \left\lbrace  \frac{1}{p} \into \left( a + | \nabla u_k|^2 \right)^{\frac{p}{2}} \, dx  - \frac{1}{p} \into \left( a + | \nabla u |^2 \right)^{\frac{p}{2}} \, dx \right. \nonumber \\
& +  \frac{1}{q} \into \left( b+ | \nabla v_k|^2 \right)^{\frac{q}{2}} \, dx  - \frac{1}{q} \into \left( b + | \nabla v |^2 \right)^{\frac{q}{2}} \, dx \nonumber \\ 
& \left. - \into \left[ F_s(u_k,v_k) (u_k - u) + F_t(u_k,v_k)(v_k - v) \right] \, dx \right\rbrace . \nonumber 
\end{align}

\noindent
We know that $u_k \rightharpoonup u$ in $W_0^{1,p}(\Omega)$ and $v_k \rightharpoonup v$ in $W_0^{1,q}(\Omega)$.  Moreover, we will use the fact that  $\{u_k\}_k$ is bounded in $L^{p^*}(\Omega)$ and
 $\{v_k\}_k$ is bounded in $L^{q^*}(\Omega)$. For every $h \in \mathbb{N}$, $h \geq 1$, let us consider the functions $\mathcal{T}_h,\mathcal{R}_h:\mathbb{R} \to \mathbb{R}$ defined in \eqref{ThRh}. We have:

\begin{itemize}

\item[$i)$] $\mathcal{T}_h(u_k) \rightharpoonup \mathcal{T}_h(u) \text{ in } W_0^{1,p}(\Omega)$, and $\mathcal{T}_h(v_k) \rightharpoonup \mathcal{T}_h(v) \text{ in } W_0^{1,q}(\Omega);$ \medskip

\item[$ii)$] $\mathcal{T}_h(u_k) \rightarrow \mathcal{T}_h(u) \text{ in }  L^{\sigma}(\Omega)$ if $ 1 \leq \sigma < \infty$, and $\mathcal{T}_h(v_k) \rightarrow \mathcal{T}_h(v) \text{ in }  L^{\tau}(\Omega)$ if $ 1 \leq \tau < \infty$. \medskip

\item[$iii)$] $\mathcal{R}_h(u_k) \rightharpoonup \mathcal{R}_h(u) \text{ in } W_0^{1,p}(\Omega)$, and $\mathcal{R}_h(v_k) \rightharpoonup \mathcal{R}_h(v) \text{ in } W_0^{1,q}(\Omega);$ \medskip

\item[$iv)$] $\mathcal{R}_h(u_k) \rightarrow \mathcal{R}_h(u) \text{ in }  L^{\sigma}(\Omega)$ if $ 1 \leq \sigma < p^*$, and $\mathcal{R}_h(v_k) \rightarrow \mathcal{R}_h(v) \text{ in }  L^{\tau}(\Omega)$ if $ 1 \leq \tau < q^*$.
\end{itemize}

\medskip
\noindent
Taking into account Remark \ref{TroncataRresto}, we have

\begin{align*}
\displaystyle
 \into \left( a + | \nabla u|^2 \right)^{\frac{p}{2}} \, dx 
 & = \int_{ |u| > h  }  \left( a+ | \nabla \mathcal{R}_h (u) |^2 \right)^{\frac{p}{2}} \, dx  + \int_{ |u| < h  }  \left( a + | \nabla \mathcal{T}_h (u) |^2 \right)^{\frac{p}{2}} \, dx + a^{\frac{p}{2}}  \int_{ |u| = h  } \, dx \\
 & = \into \left( a + | \nabla \mathcal{R}_h(u)|^2 \right)^{\frac{p}{2}} \, dx + \into \left( \alpha + | \nabla \mathcal{T}_h(u)|^2 \right)^{\frac{p}{2}} \, dx  - a^{\frac{p}{2}} \left| \Omega   \right|. 
\end{align*}

\noindent
Similarly, we have 

\begin{align*}
\displaystyle
 \into \left( a + | \nabla u_k|^2 \right)^{\frac{p}{2}} \, dx = \into \left( a + | \nabla \mathcal{R}_h(u_k)|^2 \right)^{\frac{p}{2}} \, dx + \into \left( a + | \nabla \mathcal{T}_h(u_k)|^2 \right)^{\frac{p}{2}} \, dx  - a^{\frac{p}{2}} \left| \Omega   \right|,
\end{align*}

\begin{align*}
\displaystyle
 \into \left( b+ | \nabla v_k|^2 \right)^{\frac{q}{2}} \, dx = \into \left( b + | \nabla \mathcal{R}_h(v_k)|^2 \right)^{\frac{q}{2}} \, dx + \into \left( b + | \nabla \mathcal{T}_h(v_k)|^2 \right)^{\frac{q}{2}} \, dx  - b^{\frac{q}{2}} \left| \Omega   \right|,
\end{align*}
 
\begin{align*}
\displaystyle
 \into \left( b + | \nabla v|^2 \right)^{\frac{q}{2}} \, dx = \into \left( b + | \nabla \mathcal{R}_h(v)|^2 \right)^{\frac{q}{2}} \, dx + \into \left( b + | \nabla \mathcal{T}_h(v)|^2 \right)^{\frac{q}{2}} \, dx  - b^{\frac{q}{2}} \left| \Omega   \right|.
\end{align*}

\medskip
\noindent
Considering these relations in \eqref{contoPSlocale}, we get

\begin{align}\label{contoPSlocale2}
\displaystyle 
0  \geq  \limsup_{k \to  \infty} & \;  \left\lbrace \frac{1}{p} \into \left( a + | \nabla \mathcal{R}_h(u_k)|^2 \right)^{\frac{p}{2}} \, dx - \frac{1}{p} \into \left( a + | \nabla \mathcal{R}_h(u)|^2 \right)^{\frac{p}{2}} \, dx \right.  \\
& + \frac{1}{q} \into \left( b + | \nabla \mathcal{R}_h(v_k)|^2 \right)^{\frac{q}{2}} \, dx - \frac{1}{q} \into \left( b + | \nabla \mathcal{R}_h(v)|^2 \right)^{\frac{q}{2}} \, dx \nonumber \\
& - \into \left[ F_s(u_k,v_k) (\mathcal{R}_h(u_k) - \mathcal{R}_h(u)) + F_t(u_k,v_k)(\mathcal{R}_h(v_k) - \mathcal{R}_h(v)) \right] \, dx \nonumber \\
& - \into \left[ F_s(u_k,v_k) (\mathcal{T}_h(u_k) - \mathcal{T}_h(u)) + F_t(u_k,v_k)(\mathcal{T}_h(v_k) - \mathcal{T}_h(v)) \right] \, dx \nonumber \\
& + \frac{1}{p} \into \left( a + | \nabla \mathcal{T}_h(u_k)|^2 \right)^{\frac{p}{2}} \, dx - \frac{1}{p} \into \left( a + | \nabla \mathcal{T}_h(u)|^2 \right)^{\frac{p}{2}} \, dx \nonumber \\
& + \left. \frac{1}{q} \into \left( b + | \nabla \mathcal{T}_h(v_k)|^2 \right)^{\frac{q}{2}} \, dx - \frac{1}{q} \into \left( b + | \nabla \mathcal{T}_h(v)|^2 \right)^{\frac{q}{2}} \, dx \right\rbrace  \nonumber
\end{align}

\medskip
\noindent
Considering $ii)$ and the boundedness of $\{u_k\}_k $ in $ L^{p^*}(\Omega)$ and of $\{v_k\}_k$ in $L^{q^*}(\Omega)$, we get

\begin{align*}
\displaystyle
 \into F_s(u_k,v_k) (\mathcal{T}_h(u_k) - \mathcal{T}_h(u)) \, dx \to 0, \qquad \text{and} \qquad  \into F_t(u_k,v_k) (\mathcal{T}_h(v_k) - \mathcal{T}_h(v)) \, dx \to 0. 
\end{align*}

\medskip
\noindent
Moreover, by $i)$ and the weak lower semicontinuity of the functional $f$ defined in \eqref{funzionalegenerico}, we get

\begin{align*}
 \frac{1}{p} \into \left( a + | \nabla \mathcal{T}_h(u) |^2 \right)^{\frac{p}{2}} \, dx \leq \liminf_{k \to  \infty} \left[ \frac{1}{p} \into \left( a + | \nabla \mathcal{T}_h(u_k) |^2 \right)^{\frac{p}{2}} \, dx \right],
\end{align*}

\medskip
\noindent
and

\begin{align*}
 \frac{1}{q} \into \left( b + | \nabla \mathcal{T}_h(v) |^2 \right)^{\frac{q}{2}} \, dx \leq \liminf_{k \to  \infty} \left[ \frac{1}{q} \into \left( b + | \nabla \mathcal{T}_h(v_k) |^2 \right)^{\frac{q}{2}} \, dx \right] . 
\end{align*}

\medskip
\noindent
Considering also that $ \displaystyle \limsup_{k \to \infty} \left( a_k + b_k  \right) \geq \limsup_{k \to \infty} a_k + \liminf_{k \to \infty} b_k $, by \eqref{contoPSlocale2}  we deduce

\begin{align*}
\displaystyle 
0  \geq  \limsup_{k \to \infty} & \left\lbrace \frac{1}{p} \into \left( a + | \nabla \mathcal{R}_h(u_k)|^2 \right)^{\frac{p}{2}} \, dx - \frac{1}{p} \into \left( a + | \nabla \mathcal{R}_h(u)|^2 \right)^{\frac{p}{2}} \, dx \right. \\
& + \frac{1}{q} \into \left( b + | \nabla \mathcal{R}_h(v_k)|^2 \right)^{\frac{q}{2}} \, dx - \frac{1}{q} \into \left( b + | \nabla \mathcal{R}_h(v)|^2 \right)^{\frac{q}{2}} \, dx \\
& \left. - \into \left[ F_s(u_k,v_k) (\mathcal{R}_h(u_k) - \mathcal{R}_h(u)) + F_t(u_k,v_k)(\mathcal{R}_h(v_k) - \mathcal{R}_h(v)) \right] \, dx \right\rbrace . 
\end{align*}

\medskip
\noindent
Let us now define, almost everywhere, the following sets:

\begin{align*}
\displaystyle
\Omega^{u_k}_{h,+}:= \{ x \in \Omega \; : \;  |u_k(x)| \;  \geq  \;  h   \} \qquad \text{and} \qquad \Omega^{u_k}_{h,-}:= \{ x \in \Omega \; : \;  |u_k(x)| \;   <     \;  h   \}. 
\end{align*}

\noindent
Taking into account  Remark \ref{TroncataRresto}, we observe that 

\begin{align*}
\displaystyle
   \into \left( a + | \nabla \mathcal{R}_h(u_k)|^2 \right)^{\frac{p}{2}} \, dx 
& =    a^{\frac{p}{2}} \left| \Omega^{u_k}_{h,-}  \right|   +  \int_{ \Omega^{u_k}_{h,+}} \left( a + | \nabla \mathcal{R}_h(u_k)|^2 \right)^{\frac{p}{2}} \, dx \\
& \geq  a^{\frac{p}{2}} \left| \Omega^{u_k}_{h,-}  \right|   +  \int_{ \Omega^{u_k}_{h,+}} | \nabla \mathcal{R}_h(u_k)|^p \, dx \\
& = a^{\frac{p}{2}} \left| \Omega \right| +  \into | \nabla \mathcal{R}_h(u_k)|^p \, dx - a^{\frac{p}{2}} \left| \Omega^{u_k}_{h,+}  \right|. 
\end{align*}

\noindent
By Lemma \ref{convergenzasopralivelli} we have  

$$ \displaystyle \limsup_{k \to \infty} \left| \Omega^{u_k}_{h,+}  \right| \leq \left| \Omega^{u}_{h,+}  \right|.$$

\noindent
Similarly, we have 

\begin{align*}
\displaystyle
   \into \left( b + | \nabla \mathcal{R}_h(v_k)|^2 \right)^{\frac{q}{2}} \, dx  \geq b^{\frac{q}{2}} \left| \Omega \right| +  \into | \nabla \mathcal{R}_h(v_k)|^q \, dx - b^{\frac{q}{2}} \left| \Omega^{v_k}_{h,+}  \right|
\end{align*}

\noindent
with

$$ \displaystyle \limsup_{k \to \infty} \left| \Omega^{v_k}_{h,+}  \right| \leq \left| \Omega^{v}_{h,+}  \right|.$$

\noindent
Hence, at now we have proved

\begin{align*}
\displaystyle
\limsup_{k \to \infty}  & \left( \frac{1}{p} \into | \nabla \mathcal{R}_h(u_k)|^p \, dx  +  \frac{1}{q} \into | \nabla \mathcal{R}_h(v_k)|^q \, dx  \right. \\
&  \left. - \into F_s(u_k,v_k) \left( \mathcal{R}_h(u_k)-\mathcal{R}_h(u) \right) \, dx - \into F_t(u_k,v_k) \left( \mathcal{R}_h(v_k)-\mathcal{R}_h(v) \right)  \, dx  \right) \\
& \leq \frac{1}{p} \into \left( a+ | \nabla \mathcal{R}_h(u)|^2 \right)^{\frac{p}{2}} \, dx - \frac{1}{p} a^{\frac{p}{2}} \left| \Omega \right| + \frac{1}{p} a^{\frac{p}{2}} \left| \Omega^{u}_{h,+}  \right| \\
& + \frac{1}{q} \into \left( b + | \nabla \mathcal{R}_h(v)|^2 \right)^{\frac{q}{2}} \, dx -  \frac{1}{q} b^{\frac{q}{2}} \left| \Omega \right| + \frac{1}{q} b^{\frac{q}{2}} \left| \Omega^{v}_{h,+}  \right| . 
\end{align*}

\noindent
Now, on the one hand by Chebyshev inequality we have 

$$ \displaystyle \left| \Omega^{u}_{h,+}  \right| \to 0 \quad \text{and} \quad \left| \Omega^{v}_{h,+}  \right| \to 0 \qquad \text{ as } h \to  \infty, $$

\medskip
\noindent
and on the other hand by dominated convergence theorem we have

$$ \displaystyle \into \left( \alpha+ | \nabla \mathcal{R}_h(u)|^2 \right)^{\frac{p}{2}} \, dx \to \alpha^{\frac{p}{2}} \left| \Omega \right|  \quad \text{ and }  \quad \into \left( \beta+ | \nabla \mathcal{R}_h(v)|^2 \right)^{\frac{q}{2}} \, dx \to \beta^{\frac{q}{2}} \left| \Omega \right|  \qquad \text{ as } h \to  \infty. $$

\noindent
Therefore, denoting with $o_h(1)$ a quantity that goes to $0$ as $h \to \infty$, we infer that

\begin{align*}
\displaystyle
\limsup_{k \to \infty} &  \left( \frac{1}{p} \into | \nabla \mathcal{R}_h(u_k)|^p \, dx  +  \frac{1}{q} \into | \nabla \mathcal{R}_h(v_k)|^q \, dx  \right. \\
&  \left. - \into F_s(u_k,v_k) \left( \mathcal{R}_h(u_k)-\mathcal{R}_h(u) \right) \, dx - \into F_t(u_k,v_k) \left( \mathcal{R}_h(v_k)-\mathcal{R}_h(v) \right)  \, dx  \right) \\
& \leq o_h(1).
\end{align*}

\medskip
\noindent
By continuity of $F_s$ and the convergence a.e. $x \in \Omega$ of $u_k(x)$ to $u(x)$ and of $v_k(x)$ to $v(x)$, we deduce that $F_s(u_k(x),v_k(x)) \to F_s(u(x),v(x))$ a.e. $x \in \Omega$. By \eqref{primaPoincareHopf} and the boundedness of $\{u_k\}_k$ in $L^{p^*}(\Omega)$ and of $\{v_k\}_k$ in $L^{q^*}(\Omega)$, we also deduce that the sequence $\{F_s(u_k,v_k)\}_k$ is bounded in $L^{(p^*)'}(\Omega)$. Hence we can apply Proposition \ref{weakconvLp} and say that $F_s(u_k,v_k) \rightharpoonup F_s(u,v)$ in $L^{(p^*)'}(\Omega)$, and with similar argument also $F_t(u_k,v_k) \rightharpoonup F_t(u,v)$ in $L^{(q^*)'}(\Omega)$.\\
Since $\mathcal{R}_h(u) \in L^{p^*}(\Omega)$ and $\mathcal{R}_h(v) \in L^{q^*}(\Omega)$,  we infer that

$$ \displaystyle \into F_s(u_k,v_k)\mathcal{R}_h(u) \, dx \to \into F_s(u,v)\mathcal{R}_h(u) \, dx, \quad \into F_t(u_k,v_k)\mathcal{R}_h(v) \, dx \to \into F_t(u,v)\mathcal{R}_h(v) \, dx,$$

\medskip
\noindent
whence

\begin{align*}
\displaystyle
\limsup_{k \to \infty} 
& \left( \frac{1}{p} \into \left| \nabla \mathcal{R}_h(u_k)  \right|^p \, dx +  \frac{1}{q} \into \left| \nabla \mathcal{R}_h(v_k)  \right|^q \, dx  \right. \\
&  \left. - \into F_s(u_k,v_k) \mathcal{R}_h(u_k)  \, dx - \into F_t(u_k,v_k) \mathcal{R}_h(v_k)   \, dx  \right) \\
& \leq  - \into F_s(u,v) \mathcal{R}_h(u)  \, dx - \into F_t(u,v) \mathcal{R}_h(v)  \, dx + o_h(1). 
\end{align*}

\medskip
\noindent
By  $iv)$ we know that

$$ \displaystyle \into \left| \mathcal{R}_h(u_k)  \right| \, dx \to \into \left|  \mathcal{R}_h(u)  \right| \, dx \quad \text{and} \quad \into \left| \mathcal{R}_h(v_k)  \right| \, dx \to \into \left|  \mathcal{R}_h(v)  \right| \, dx \qquad \text{ as } k \to \infty.$$

\medskip
\noindent
Considering also \eqref{primaPoincareHopf} and \eqref{secondaPoincareHopf}, we get 

\begin{align*}
\displaystyle
\limsup_{k \to \infty} 
& \left( \frac{1}{p} \into \left| \nabla \mathcal{R}_h(u_k)  \right|^p \, dx +  \frac{1}{q} \into \left| \nabla \mathcal{R}_h(v_k)  \right|^q \, dx  \right. \\
& - C \into |u_k|^{p^*-1} \left| \mathcal{R}_h(u_k) \right| \, dx - C \into |v_k|^{q^{*}\frac{p^*-1}{p^{*}}} \left| \mathcal{R}_h(u_k) \right| \, dx  \\
& \left. - C \into |u_k|^{p^*\frac{q^{*}-1}{q^{*}}} \left| \mathcal{R}_h(v_k) \right| \, dx - C \into |v_k|^{q^*-1} \left| \mathcal{R}_h(v_k) \right| \, dx  \right) \\
& \leq  C \into |u|^{p^*-1} \left| \mathcal{R}_h(u) \right| \, dx + C \into |v|^{q^{*}\frac{p^*-1}{p^{*}}} \left| \mathcal{R}_h(u) \right| \, dx  \\
& + C \into |u|^{p^*\frac{q^{*}-1}{q^{*}}} \left| \mathcal{R}_h(v) \right| \, dx + C \into |v|^{q^*-1} \left| \mathcal{R}_h(v) \right| \, dx \\
& + 2C \into \left| \mathcal{R}_h(u) \right| \, dx  + 2C \into \left| \mathcal{R}_h(v) \right| \, dx + o_h(1).
\end{align*}

\medskip
\noindent
By Gagliardo–Nirenberg–Sobolev inequality there exist $S(p,N)>0$ and $S(q,N)>0$ such that

$$ \displaystyle S(p,N) \lVert \bar{u} \rVert_{p^*}^{p} \leq \into | \nabla \bar{u}|^p \, dx  \quad \text{and} \quad S(q,N) \lVert \bar{v} \rVert_{q^*}^{q} \leq \into | \nabla \bar{v}|^q \, dx$$

\medskip
\noindent
for any $\bar{u} \in W_0^{1,p}(\Omega)$ and for any $\bar{v} \in W_0^{1,q}(\Omega)$, by which for any $h$ and $k$ we get

\begin{equation*}
\displaystyle \frac{1}{p} \into | \nabla \mathcal{R}_h(u_k)|^p \, dx \geq \frac{S(p,N)}{p} \lVert \mathcal{R}_h(u_k) \rVert_{p^*}^{p}  \quad \text{and} \quad
\displaystyle \frac{1}{q} \into | \nabla \mathcal{R}_h(v_k)|^q \, dx \geq \frac{S(q,N)}{q} \lVert \mathcal{R}_h(v_k) \rVert_{q^*}^{q}. 
\end{equation*}

\medskip
\noindent
Hence

\begin{align}\label{calcoloPalaisSmale3}
\displaystyle
\limsup_{k \to \infty} 
& \left( \frac{S(p,N)}{p} \lVert \mathcal{R}_h(u_k) \rVert_{p^*}^{p} +  \frac{S(q,N)}{q} \lVert \mathcal{R}_h(v_k) \rVert_{q^*}^{q} \right. \\
& - C \into |u_k|^{p^*-1} \left| \mathcal{R}_h(u_k) \right| \, dx - C \into |v_k|^{q^{*}\frac{p^*-1}{p^{*}}} \left| \mathcal{R}_h(u_k) \right| \, dx  \nonumber \\
& \left. - C \into |u_k|^{p^*\frac{q^{*}-1}{q^{*}}} \left| \mathcal{R}_h(v_k) \right| \, dx - C \into |v_k|^{q^*-1} \left| \mathcal{R}_h(v_k) \right| \, dx  \right) \nonumber \\
& \leq C \into |u|^{p^*-1} \left| \mathcal{R}_h(u) \right| \, dx + C \into |v|^{q^{*}\frac{p^*-1}{p^{*}}} \left| \mathcal{R}_h(u) \right| \, dx  \nonumber \\
& + C \into |u|^{p^*\frac{q^{*}-1}{q^{*}}} \left| \mathcal{R}_h(v) \right| \, dx + C \into |v|^{q^*-1} \left| \mathcal{R}_h(v) \right| \, dx  \nonumber \\
& + 2C \into \left| \mathcal{R}_h(u) \right| \, dx  + 2C \into \left| \mathcal{R}_h(v) \right| \, dx + o_h(1). \nonumber
\end{align}

\noindent
Let us now consider the term

$$ \into |u_k|^{p^*-1} \left| \mathcal{R}_h(u_k) \right| \, dx.$$

\noindent
By H\"older inequality with  conjugated exponents $p^*/p$ and $p^*/(p^* - p)$, we have

\begin{align*}
\displaystyle
\into |u_k|^{p^*-1} \left| \mathcal{R}_h(u_k) \right| \, dx
& = \into |u_k|^{p^*-p} \left| \mathcal{T}_h(u_k) + \mathcal{R}_h(u_k)  \right|^{p-1} \left| \mathcal{R}_h(u_k) \right| \, dx \\
& \leq c(p) \left( \into |u_k|^{p^*-p} \left| \mathcal{T}_h(u_k)  \right|^{p-1} \left| \mathcal{R}_h(u_k) \right| \, dx + 
\into |u_k|^{p^*-p} \left| \mathcal{R}_h(u_k)  \right|^{p}  \, dx \right) \\
& \leq c(p) \left( \into |u_k|^{p^*-p} \left| \mathcal{T}_h(u_k)  \right|^{p-1} \left| \mathcal{R}_h(u_k) \right| \, dx +  \into |u_0|^{p^*-p} \left| \mathcal{R}_h(u_k)  \right|^{p}  \, dx \right) \\
& + c(p) \into |u_k - u_0|^{p^*-p} \left| \mathcal{R}_h(v_k)  \right|^{p}  \, dx \\
& \leq  c(p) \left( \into |u_k|^{p^*-p} \left| \mathcal{T}_h(u_k)  \right|^{p-1} \left| \mathcal{R}_h(u_k) \right| \, dx  + \into |u_0|^{p^*-p} \left| \mathcal{R}_h(u_k)  \right|^{p}  \, dx \right) \\
& + c(p) \lVert u_k - u_0 \rVert_{p^*}^{p^* - p} \lVert \mathcal{R}_h(u_k) \rVert_{p^*}^p.
\end{align*}

\noindent
By Proposition \ref{weakconvLp} we infer that $ \left| \mathcal{R}_h(u_k)  \right|^{p} \rightharpoonup \left| \mathcal{R}_h(u)  \right|^{p}$ in $L^{\frac{p^*}{p}}(\Omega)$. Since $|u_0|^{p^*-p} \in L^{ \frac{p^*}{p^*-p}}(\Omega)$, we have

$$ \displaystyle \into |u_0|^{p^*-p} \left| \mathcal{R}_h(u_k)  \right|^{p}  \, dx \to \into |u_0|^{p^*-p} \left| \mathcal{R}_h(u)  \right|^{p}  \, dx. $$

\medskip
\noindent
On the other hand, by dominated convergence theorem we have 

$$ \displaystyle \into |u_k|^{p^*-p} \left| \mathcal{T}_h(u_k)  \right|^{p-1} \left| \mathcal{R}_h(u_k) \right| \, dx \to \into |u|^{p^*-p} \left| \mathcal{T}_h(u)  \right|^{p-1} \left| \mathcal{R}_h(u) \right| \, dx.$$

\medskip
\noindent
Summing up, denoting with $o_k(1)$ a quantity the goes to zero as $k \to \infty$, we have shown that 

\begin{align*}
\displaystyle
\into |u_k|^{p^*-1} \left| \mathcal{R}_h(u_k) \right| \, dx 
& \leq  c(p) \left( \into |u|^{p^*-p} \left| \mathcal{T}_h(v)  \right|^{q-1} \left| \mathcal{R}_h(u) \right| \, dx  + \into |u_0|^{p^*-p} \left| \mathcal{R}_h(u)  \right|^{p}  \, dx \right) \\
& + c(p) \lVert u_k - u_0 \rVert_{p^*}^{p^* - p} \lVert \mathcal{R}_h(u_k) \rVert_{p^*}^p + o(1).
\end{align*}

\medskip
\noindent
Similarly, we have

\begin{align*}
\displaystyle
\into |v_k|^{q^*-1} \left| \mathcal{R}_h(v_k) \right| \, dx 
& \leq  c(q) \left( \into |v|^{q^*-q} \left| \mathcal{T}_h(v)  \right|^{q-1} \left| \mathcal{R}_h(v) \right| \, dx  + \into |v_0|^{q^*-q} \left| \mathcal{R}_h(v)  \right|^{q}  \, dx \right) \\
& + c(q) \lVert v_k - v_0 \rVert_{q^*}^{q^* - q} \lVert \mathcal{R}_h(v_k) \rVert_{q^*}^q + o(1).
\end{align*}

\medskip
\noindent
Let us now consider the term

$$ \displaystyle \into  |u_k|^{p^*\frac{q^{*}-1}{q^{*}}} \left| \mathcal{R}_h(v_k) \right| \, dx. $$

\noindent
Denoting with

$$ \displaystyle \gamma:= p^*\frac{q^{*}-1}{q^{*}} - \frac{p}{q'}=\frac{p^* \, q^* - p^* -  \frac{p}{q'}q^*}{q^*}, $$

\noindent
we can write 

\begin{align*}
\displaystyle
& \into  |u_k|^{p^*\frac{q^{*}-1}{q^{*}}} \left| \mathcal{R}_h(v_k) \right| \, dx \\
= & \into  |u_k|^{\gamma} \left| \mathcal{T}_h(u_k) + \mathcal{R}_h(u_k) \right|^{\frac{p}{q'}} \left| \mathcal{R}_h(v_k) \right| \, dx \\
\leq & c(p,q) \left( \into  |u_k|^{\gamma} \left| \mathcal{T}_h(u_k) \right|^{\frac{p}{q'}} \left| \mathcal{R}_h(v_k) \right| \, dx  + \into  |u_k|^{\gamma} \left| \mathcal{R}_h(u_k) \right|^{\frac{p}{q'}} \left| \mathcal{R}_h(v_k) \right| \, dx      \right) \\
\leq &   c(p,q) \left( \into  |u_k|^{\gamma} \left| \mathcal{T}_h(u_k) \right|^{\frac{p}{q'}} \left| \mathcal{R}_h(v_k) \right| \, dx  +  \into  |u_0|^{\gamma} \left| \mathcal{R}_h(u_k) \right|^{\frac{p}{q'}} \left| \mathcal{R}_h(v_k) \right| \, dx   \right. \\
+ &   \left. \into  |u_k - u_0|^{\gamma} \left| \mathcal{R}_h(u_k) \right|^{\frac{p}{q'}} \left| \mathcal{R}_h(v_k) \right| \, dx \right).  \\
\end{align*}

\noindent
By applying the extension of H\"older inequality with conjugated exponents $\displaystyle \frac{p^*}{\gamma}$, $\displaystyle \frac{q'}{p}p^*$ and $q^*$, we have

$$ \displaystyle \into  |u_k - u_0|^{\gamma} \left| \mathcal{R}_h(u_k) \right|^{\frac{p}{q'}} \left| \mathcal{R}_h(v_k) \right| \, dx \leq \lVert  u_k - u_0 \rVert_{p^*}^\gamma \lVert \mathcal{R}_h(u_k) \rVert_{p^*}^{\frac{p}{q'}} \lVert \mathcal{R}_h(v_k) \rVert_{q^*}.$$

\medskip
\noindent
Moreover, by Young's inequality with conjugated exponents $q'$ and $q$, we have 

$$ \displaystyle \lVert \mathcal{R}_h(u_k) \rVert_{p^*}^{\frac{p}{q'}} \lVert \mathcal{R}_h(v_k) \rVert_{q^*} \leq \frac{1}{q'} 
\lVert \mathcal{R}_h(u_k) \rVert_{p^*}^{p} + \frac{1}{q} \lVert \mathcal{R}_h(v_k) \rVert_{q^*}^q,$$

\noindent
whence

\begin{align*}
\displaystyle 
c(p,q) \into  |u_k - u_0|^{\gamma} \left| \mathcal{R}_h(u_k) \right|^{\frac{p}{q'}} \left| \mathcal{R}_h(v_k) \right| \, dx 
& \leq  \, c(p,q) \lVert  u_k - u_0 \rVert_{p^*}^\gamma \left( \lVert \mathcal{R}_h(u_k) \rVert_{p^*}^{p} +  \lVert \mathcal{R}_h(v_k) \rVert_{q^*}^q  \right).
\end{align*}

\medskip
\noindent
We obtained

\begin{align*}
\displaystyle
\into  |u_k|^{p^*\frac{q^{*}-1}{q^{*}}} \left| \mathcal{R}_h(v_k) \right| \, dx  
\leq & c(p,q) \left(\into  |u_k|^{\gamma} \left| \mathcal{T}_h(u_k) \right|^{\frac{p}{q'}} \left| \mathcal{R}_h(v_k) \right| \, dx  +  \into  |u_0|^{\gamma} \left| \mathcal{R}_h(u_k) \right|^{\frac{p}{q'}} \left| \mathcal{R}_h(v_k) \right| \, dx   \right) \\
+ &  \lVert  u_k - u_0 \rVert_{p^*}^\gamma c(p,q)\left( \lVert \mathcal{R}_h(u_k) \rVert_{p^*}^{p} +  \lVert \mathcal{R}_h(v_k) \rVert_{q^*}^q  \right) .  \\
\end{align*}

\medskip
\noindent
Observe now that by dominated convergence theorem we have 

$$ \displaystyle \into  |u_k|^{\gamma} \left| \mathcal{T}_h(u_k) \right|^{\frac{p}{q'}} \left| \mathcal{R}_h(v_k) \right| \, dx  \to \into  |u|^{\gamma} \left| \mathcal{T}_h(u) \right|^{\frac{p}{q'}} \left| \mathcal{R}_h(v) \right| \, dx.$$

\medskip
\noindent
By Proposition \ref{weakconvLp}, we deduce that $\left| \mathcal{R}_h(u_k) \right|^{\frac{p}{q'}} \left| \mathcal{R}_h(v_k) \right| \rightharpoonup  \left| \mathcal{R}_h(u) \right|^{\frac{p}{q'}} \left| \mathcal{R}_h(v) \right|$ in $\displaystyle L^{\left( \frac{p^*}{\gamma} \right)'}(\Omega).$\\
Since $|u_0| \in L^{\frac{p^*}{\gamma} }(\Omega)$, we have

$$  \into  |u_0|^{\gamma} \left| \mathcal{R}_h(u_k) \right|^{\frac{p}{q'}} \left| \mathcal{R}_h(v_k) \right| \, dx  \to  \into  |u_0|^{\gamma} \left| \mathcal{R}_h(u) \right|^{\frac{p}{q'}} \left| \mathcal{R}_h(v) \right| \, dx. $$

\noindent
Hence, 

\begin{align*}
\displaystyle
\into  |u_k|^{p^*\frac{q^{*}-1}{q^{*}}} \left| \mathcal{R}_h(v_k) \right| \, dx  
& \leq  c(p,q) \left(\into  |u|^{\gamma} \left| \mathcal{T}_h(u) \right|^{\frac{p}{q'}} \left| \mathcal{R}_h(v) \right| \, dx  +  \into  |u_0|^{\gamma} \left| \mathcal{R}_h(u) \right|^{\frac{p}{q'}} \left| \mathcal{R}_h(v) \right| \, dx   \right) \\
& +  c(p,q) \lVert  u_k - u_0 \rVert_{p^*}^\gamma \left( \lVert \mathcal{R}_h(u_k) \rVert_{p^*}^{p} +  \lVert \mathcal{R}_h(v_k) \rVert_{q^*}^q  \right)  + o_k(1).  \\
\end{align*}

\noindent
Arguing in a similar way on the term

$$ \displaystyle \into  |v_k|^{q^*\frac{p^{*}-1}{p^{*}}} \left| \mathcal{R}_h(u_k) \right| \, dx,$$

\medskip
\noindent
and denoting with 

$$ \displaystyle \delta := q^*\frac{p^{*}-1}{p^{*}} - \frac{q}{p'}=\frac{q^* \, p^* - q^* - \displaystyle \frac{q}{p'}p^*}{p^*}, $$

\medskip
\noindent
we get
 
\begin{align*}
\displaystyle
\into  |v_k|^{q^*\frac{p^{*}-1}{p^{*}}} \left| \mathcal{R}_h(u_k) \right| \, dx  
& \leq  c(p,q) \left(\into  |v|^{\delta} \left| \mathcal{T}_h(v) \right|^{\frac{q}{p'}} \left| \mathcal{R}_h(v) \right| \, dx  +  \into  |v_0|^{\delta} \left| \mathcal{R}_h(v) \right|^{\frac{q}{p'}} \left| \mathcal{R}_h(u) \right| \, dx   \right) \\
& +   c(p,q) \lVert  v_k - v_0 \rVert_{q^*}^\delta \left( \lVert \mathcal{R}_h(v_k) \rVert_{q^*}^{q} +  \lVert \mathcal{R}_h(u_k) \rVert_{p^*}^p  \right)  + o_k(1).  \\
\end{align*}

\medskip
\noindent
Altogether, denoting with $K:=K(p,q,C)$ the greatest positive constant in previous inequalities, by \eqref{calcoloPalaisSmale3} we deduce 

\begin{align*}
\displaystyle
\limsup_{k \to \infty}  & \left[  \frac{S(p,N)}{p} \lVert \mathcal{R}_h(u_k) \rVert_{p^*}^{p}  - K \lVert u_k - u_0 \rVert_{p^*}^{p^* - p} \lVert \mathcal{R}_h(u_k) \rVert_{p^*}^p   \right. \\
& + \frac{S(q,N)}{q} \lVert \mathcal{R}_h(v_k) \rVert_{q^*}^{q} - K\lVert v_k - v_0 \rVert_{q^*}^{q^* - q} \lVert \mathcal{R}_h(v_k) \rVert_{q^*}^q  \\
&  -K \lVert  u_k - u_0 \rVert_{p^*}^\gamma \left( \lVert \mathcal{R}_h(u_k) \rVert_{p^*}^{p} +  \lVert \mathcal{R}_h(v_k) \rVert_{q^*}^q  \right)  \\
& \left. -K\lVert  v_k - v_0 \rVert_{q^*}^\delta \left( \lVert \mathcal{R}_h(v_k) \rVert_{q^*}^{q} +  \lVert \mathcal{R}_h(u_k) \rVert_{p^*}^p  \right)\right] \\
\leq & \quad K \into |u|^{p^*-1} \left| \mathcal{R}_h(u) \right| \, dx + K \into |v|^{q^{*}\frac{p^*-1}{p^{*}}} \left| \mathcal{R}_h(u) \right| \, dx \\
& + K \into |u|^{p^*\frac{q^{*}-1}{q^{*}}} \left| \mathcal{R}_h(v) \right| \, dx + K \into |v|^{q^*-1} \left| \mathcal{R}_h(v) \right| \, dx \\
& + K \into \left| \mathcal{R}_h(u) \right| \, dx  + K \into \left| \mathcal{R}_h(v) \right| \, dx + o_h(1) \\
& + K \left( \into |u|^{p^*-p} \left| \mathcal{T}_h(u)  \right|^{p-1} \left| \mathcal{R}_h(u) \right| dx  + \into |u_0|^{p^*-p} \left| \mathcal{R}_h(u)  \right|^{p}   dx \right) \\
& + K \left( \into |v|^{q^*-q} \left| \mathcal{T}_h(v)  \right|^{q-1} \left| \mathcal{R}_h(v) \right| dx  + \into |v_0|^{q^*-q} \left| \mathcal{R}_h(v)  \right|^{q}  dx \right) \\
& + K \left(\into  |u|^{\gamma} \left| \mathcal{T}_h(u) \right|^{\frac{p}{q'}} \left| \mathcal{R}_h(v) \right| \, dx  +  \into  |u_0|^{\gamma} \left| \mathcal{R}_h(u) \right|^{\frac{p}{q'}} \left| \mathcal{R}_h(v) \right| \, dx   \right)        \\
& + K \left(\into  |v|^{\delta} \left| \mathcal{T}_h(v) \right|^{\frac{q}{p'}} \left| \mathcal{R}_h(v) \right| \, dx  + \into  |v_0|^{\delta} \left| \mathcal{R}_h(v) \right|^{\frac{q}{p'}} \left| \mathcal{R}_h(u) \right| \, dx   \right)       .\\
\end{align*}

\noindent
Considering the left side of previous inequality, we have

\begin{align*}
\displaystyle
& \quad \frac{S(p,N)}{p} \lVert \mathcal{R}_h(u_k) \rVert_{p^*}^{p}  - K \lVert u_k - u_0 \rVert_{p^*}^{p^* - p} \lVert \mathcal{R}_h(u_k) \rVert_{p^*}^p    \\
& + \frac{S(q,N)}{q} \lVert \mathcal{R}_h(v_k) \rVert_{q^*}^{q} - K\lVert v_k - v_0 \rVert_{q^*}^{q^* - q} \lVert \mathcal{R}_h(v_k) \rVert_{q^*}^q  \\
& - K \lVert  u_k - u_0 \rVert_{p^*}^\gamma \left( \lVert \mathcal{R}_h(u_k) \rVert_{p^*}^{p} +  \lVert \mathcal{R}_h(v_k) \rVert_{q^*}^q  \right)  \\
& -  K\lVert  v_k - v_0 \rVert_{q^*}^\delta \left( \lVert \mathcal{R}_h(v_k) \rVert_{q^*}^{q} +  \lVert \mathcal{R}_h(u_k) \rVert_{p^*}^p  \right)\\
= & \quad \left[ \frac{S(p,N)}{p} - K \left( \lVert u_k - u_0 \rVert_{p^*}^{p^* - p}  + \lVert  u_k - u_0 \rVert_{p^*}^\gamma + \lVert  v_k - v_0 \rVert_{q^*}^\delta \right)    \right] \lVert \mathcal{R}_h(u_k) \rVert_{p^*}^{p} \\
& +  \left[ \frac{S(q,N)}{q} - K \left( \lVert v_k - v_0 \rVert_{q^*}^{q^* - q}  + \lVert  v_k - v_0 \rVert_{q^*}^\delta + \lVert  u_k - u_0 \rVert_{p^*}^\gamma \right)   \right] \lVert \mathcal{R}_h(v_k) \rVert_{q^*}^{q}.
\end{align*}

\medskip
\noindent
Now, let us consider

$$ \displaystyle \overline{B_R(z_0)}:=\left\lbrace z \in X \; : \; \lVert z - z_0 \rVert \leq r \right\rbrace,$$

\medskip
\noindent
where we recall $R>0$ is a general radius.

\medskip
\noindent
Since $\{z_k\}_k=\{(u_k,v_k)\}_k \subset \overline{B_R(z_0)} $, we can choose $R$ small enough such that for any $k$ we have

$$ \displaystyle  \frac{S(p,N)}{p} - K \left( \lVert u_k - u_0 \rVert_{p^*}^{p^* - p}  + \lVert  u_k - u_0 \rVert_{p^*}^\gamma + \lVert  v_k - v_0 \rVert_{q^*}^\delta \right) \geq M $$

\medskip
\noindent
and 

$$ \displaystyle \frac{S(q,N)}{q} - K \left( \lVert v_k - v_0 \rVert_{q^*}^{q^* - q}  + \lVert  v_k - v_0 \rVert_{q^*}^\delta + \lVert  u_k - u_0 \rVert_{p^*}^\gamma \right)\geq M $$

\medskip
\noindent
for some $M >0$. Therefore,

\begin{align*}
\displaystyle
 M \limsup_{k \to \infty}  \left(\lVert \mathcal{R}_h(u_k) \rVert_{p^*}^p \right. & + \left. \lVert \mathcal{R}_h(v_k) \rVert_{q^*}^q \right) \\
& \leq  K \into |u|^{p^*-1} \left| \mathcal{R}_h(u) \right| \, dx + K \into |v|^{q^{*}\frac{p^*-1}{p^{*}}} \left| \mathcal{R}_h(u) \right| \, dx \\
& + K \into |u|^{p^*\frac{q^{*}-1}{q^{*}}} \left| \mathcal{R}_h(v) \right| \, dx + K \into |v|^{q^*-1} \left| \mathcal{R}_h(v) \right| \, dx \\
& + K \into \left| \mathcal{R}_h(u) \right| \, dx  + K \into \left| \mathcal{R}_h(v) \right| \, dx + o_h(1) \\
& + K \left( \into |u|^{p^*-p} \left| \mathcal{T}_h(u)  \right|^{p-1} \left| \mathcal{R}_h(u) \right| dx  + \into |u_0|^{p^*-p} \left| \mathcal{R}_h(u)  \right|^{p}   dx \right) \\
& + K \left( \into |v|^{q^*-q} \left| \mathcal{T}_h(v)  \right|^{q-1} \left| \mathcal{R}_h(v) \right| dx  + \into |v_0|^{q^*-q} \left| \mathcal{R}_h(v)  \right|^{q}  dx \right) \\
& + K \left(\into  |u|^{\gamma} \left| \mathcal{T}_h(u) \right|^{\frac{p}{q'}} \left| \mathcal{R}_h(v) \right| \, dx  +  \into  |u_0|^{\gamma} \left| \mathcal{R}_h(u) \right|^{\frac{p}{q'}} \left| \mathcal{R}_h(v) \right| \, dx   \right)        \\
& + K \left(\into  |v|^{\delta} \left| \mathcal{T}_h(v) \right|^{\frac{q}{p'}} \left| \mathcal{R}_h(v) \right| \, dx  + \into  |v_0|^{\delta} \left| \mathcal{R}_h(v) \right|^{\frac{q}{p'}} \left| \mathcal{R}_h(u) \right| \, dx   \right)       .
\end{align*}

\medskip
\noindent
Taking into account Remark \ref{TroncataRresto}, we can apply dominated convergence theorem (with respect to $h$) to the right side of previous inequality, obtaining that

$$ \displaystyle \lim_{h \to \infty} \, \limsup_{k \to \infty}  \left(\lVert \mathcal{R}_h(u_k) \rVert_{p^*}^p + \lVert \mathcal{R}_h(v_k) \rVert_{q^*}^q \right) =0,$$

\medskip
\noindent
by which we deduce

$$ \displaystyle \lim_{h \to \infty} \, \limsup_{k \to \infty}  \; \lVert \mathcal{R}_h(u_k) \rVert_{p^*}  = 0 \qquad \text{and} \qquad \lim_{h \to \infty} \, \limsup_{k \to \infty}   \; \lVert \mathcal{R}_h(v_k) \rVert_{q^*} = 0.  $$ 

\medskip
\noindent
Now, observing that 

$$ \displaystyle \lVert u_k - u \rVert_{p^*} \leq \lVert \mathcal{T}_h(u_k) - \mathcal{T}_h(u) \rVert_{p^*} + \lVert \mathcal{R}_h(u_k) \rVert_{p^*} + \lVert \mathcal{R}_h(u) \rVert_{p^*},$$

\medskip
\noindent
we have

$$ \displaystyle \limsup_{k \to \infty} \, \lVert u_k - u \rVert_{p^*} \leq  \limsup_{k \to \infty} \, \lVert \mathcal{R}_h(u_k) \rVert_{p^*} + \lVert \mathcal{R}_h(u) \rVert_{p^*}.$$

\medskip
\noindent
By dominated convergence theorem, we have 

$$ \displaystyle \lim_{h \to \infty}  \, \lVert \mathcal{R}_h(u) \rVert_{p^*}=0, $$

\medskip
\noindent
hence we deduce 

$$ \displaystyle \lVert u_k - u \rVert_{p^*} \to 0 ,$$

\medskip
\noindent
and similarly

$$ \displaystyle \lVert v_k - v \rVert_{q^*} \to 0. $$

\medskip
\noindent
Finally, we infer that

$$ \displaystyle \into \left[ F_s(u_k,v_k) (u_k - u) + F_t(u_k,v_k)(v_k - v) \right] \, dx \to 0,$$

\medskip
\noindent
and coming back to \eqref{ipotesiHapbq}, we deduce  

\begin{equation*}
\displaystyle \limsup_{k\to \infty} \, \langle H_{a,p,b,q}(z_k),z_k-z    \rangle 
\leq 0,
\end{equation*}

\medskip
\noindent
whence, by applying Proposition \ref{HapbqS+}, we conclude that $\lVert z_k - z \rVert \to 0$.

\qed

\section{A Poincar\'e-Hopf formula for systems}\label{finitezzagruppicritici}

\medskip
\noindent
In this Section we prove Theorem \ref{TheorFinGrupCrit}, that is each isolated critical point of the functional $J_{a,b}$ defined in \eqref{functionalPoincareHopf} and associated to system \eqref{pqPoincareHopf}, has critical groups of finite type and the  Poincaré-Hopf formula holds. We retrace the same ideas contained in \cite{CD1}, where the interested reader can find all proofs and details.

We will use the following result (see \cite[Proposition 4.1]{CD1}), which showes that a Poincaré-Hopf formula holds in finite dimension.

\begin{proposition}\label{prop4.1CingolaniDegiovanni}

Let $Y$ be a finite dimensional space, $W$ an open subset of $Y$ and $f:W \to \mathbb{R}$ a functional of class $C^1$. Let $\bar{w} \in W$ and let $r, \sigma >0$ be such that $\overline{B_{4r}(\bar{w})} \subseteq W$ and\\

\begin{itemize}
\item[$(a)$] $f'(w) \neq 0$ for every $ w \in \overline{B_{4r}(\bar{w})} $ with $\left|f(w) - f(\bar{w}) \right| \geq 3 \sigma r /5;$\\
\item[$(b)$] 
$$ \displaystyle \inf \left\lbrace \lVert f'(w) \rVert \, : \quad w \in \overline{B_{4r}(\bar{w})} \setminus B_r(\bar{w}) \right\rbrace > \sigma.$$
\end{itemize}

\medskip
\noindent
If we set

\begin{align*}
\displaystyle
\beta(w) & := f(\bar{w}) + \sigma \left( r  - ( \lVert w - \bar{w} \rVert - r)^+    \right),\\
\Sigma   & :=  \left\lbrace w \in Y \, : \quad f(\bar{w}) - 2 \sigma r \leq f(w) \leq \beta(w) \right\rbrace ,\\
A        & := \left\lbrace w \in \Sigma \, : \quad f(w) \leq f(\bar{w}) - \sigma r \right\rbrace,
\end{align*}

\medskip
\noindent
then $H_{*} ( \Sigma, A)$ is of finite type and we have

\begin{equation*}
\mathrm{deg}(f',B_r(\bar{w}),0) = \sum_{m\geq 0}
(-1)^m\,\rnk{H_m(\Sigma,A)}.
\end{equation*}

\end{proposition}

\medskip
\noindent
Let us recall also the following result, that gives information on the critical groups.

\begin{proposition}\label{proposizioneGruppiCritici}
Let $Y$ be a normed vector space and $f:Y \to \mathbb{R}$ a $C^1$ function. Let  $\bar{w} \in Y$ and $r, \sigma >0$ such that $\overline{B_{4r}(\bar{w})}$ is complete and

\begin{itemize}

\item[$(i)$] $f$ satisfies $(PS)$ on $\overline{B_{4r}(\bar{w})}$;

\item[$(ii)$] $f'(w) \neq 0$ for every $ w \in \overline{B_{4r}(\bar{w})} \setminus \{ \bar{w} \}$;

\item[$(iii)$] $$ \displaystyle \inf \left\lbrace \lVert f'(w) \rVert \, : \quad w \in \overline{B_{4r}(\bar{w})} \setminus B_r(\bar{w}) \right\rbrace > \sigma.$$
\end{itemize}

\medskip
\noindent
If we set

\begin{align*}
\displaystyle
\beta(w) & := f(\bar{w}) + \sigma \left( r  - ( \lVert w - \bar{w} \rVert - r)^+    \right),\\
\Sigma   & :=  \left\lbrace w \in Y \, : \quad f(\bar{w}) - 2 \sigma r \leq f(w) \leq \beta(w) \right\rbrace ,\\
A        & := \left\lbrace w \in \Sigma \, : \quad f(w) \leq f(\bar{w}) - \sigma r \right\rbrace,
\end{align*}

\medskip
\noindent
then we have $C_m(f,\bar{w}) \approx H_m(\Sigma, A)$ for every $m \geq 0$.
\end{proposition}

\medskip
\noindent
{\mbox {\it Proof of Theorem~\ref{TheorFinGrupCrit}.~}} 
Let $\bar{z}$ be an isolated critical point of the functional $J_{a,b}$ defined in \eqref{functionalPoincareHopf}. Moreover, taking into account that $J_{a,b}: X \to \mathbb{R}$ is of class $C^1$ and by means of Theorem \ref{JabS+} its derivative is locally of class $(S)_+$, there exist $L,r>0$ such that $J_{a,b}$ is Lipschitz continuous of constant $L$ on $ \overline{B_{4r}(\bar{z})}$, $J_{a,b}'$ is of class $(S)_+$ on  $\overline{B_{4r}(\bar{z})}$ and $J_{a,b}'(z) \neq 0$ for every $ z \in  \overline{B_{4r}(\bar{z})} \setminus \{ \bar{z} \}$. In particular (see \cite[Lemma 4.1, Lemma 4.2 and Lemma 4.3]{CD1}), there exist $\sigma \in [0,L]$ and a finite dimensional subspace $Y$ of $X$, with $\bar{z} \in Y$, such that, setting

\begin{align*}
\displaystyle
\beta(z) & := J_{a,b}(\bar{z}) + \sigma \left( r  - ( \lVert z - \bar{z} \rVert - r)^+    \right),\\
\Sigma   & :=  \left\lbrace z \in X \, : \quad J_{a,b}(\bar{z}) - 2 \sigma r \leq J_{a,b}(z) \leq \beta(z) \right\rbrace ,\\
A        & := \left\lbrace z \in \Sigma \, : \quad J_{a,b}(z) \leq J_{a,b}(\bar{z}) - \sigma r \right\rbrace,
\end{align*}

\medskip
\noindent
the following facts hold:

\begin{itemize}
\item[•] $H_{*} ( \Sigma, A)$ is of finite type;

\item[•] for any  $m \geq 0$, we have

\begin{equation}\label{contoTheorgruppicritici1}
\displaystyle H_{m} ( \Sigma, A) \approx H_{m} ( \Sigma \cap Y, A \cap Y) \qquad \text{ for any } m \geq 0.
\end{equation}

\end{itemize}

\medskip
\noindent
Since assumptions of Proposition \ref{prop4.1CingolaniDegiovanni} and Proposition \ref{proposizioneGruppiCritici} are satisfied, we also deduce 

\begin{equation}\label{contoTheorgruppicritici2}
\mathrm{deg}\left((J_{a,b}\bigl|_{U \cap Y} )',B_r(\bar{z}) \cap Y,0 \right) = \sum_{m\geq 0}
(-1)^m\,\rnk{H_{m} ( \Sigma \cap Y, A \cap Y)}
\end{equation}

\medskip
\noindent
and

\begin{equation}\label{contoTheorgruppicritici3}
\displaystyle C_m(J_{a,b},\bar{z}) \approx H_m(\Sigma, A) \qquad \text{ for any } m \geq 0.
\end{equation}

\medskip
\noindent
Moreover, the definition and properties of degree for maps of class $(S)_+$ (see \cite{browder1983}) give

\begin{equation}\label{contoTheorgruppicritici4}
\mathrm{deg}\left(J_{a,b}',B_r(\bar{z}),0 \right)=\mathrm{deg}\left((J_{a,b}\bigl|_{U \cap Y} )',B_r(\bar{z}) \cap Y,0 \right).
\end{equation}

\medskip
\noindent
Altogether, by \eqref{contoTheorgruppicritici1}, \eqref{contoTheorgruppicritici2}, \eqref{contoTheorgruppicritici3} and \eqref{contoTheorgruppicritici4}, we infer that

\begin{equation}
\mathrm{deg}\left(J_{a,b}',B_r(\bar{z}),0 \right) = \sum_{m\geq 0}
(-1)^m\,\rnk{ C_m(J_{a,b},\bar{z})  }.
\end{equation}

\qed

\section*{Appendix A}
\addcontentsline{toc}{section}{Appendix A}
\label{app:A}

\begin{lemma}\label{smc} For any  $c\geq 0$  and for any $r > 1$, the map $d_{c,r}:\re^N\to \re^N$ defined as
 
\begin{equation}\label{definizionemappadcr}
\displaystyle  d_{c,r}(\xi)=(c+|\xi|^2)^{\frac{r-2}{2}}\xi  
\end{equation} 

\noindent
is strictly monotone, that is 

\begin{equation}\label{strictmonotncondition}
	  \left[ d_{c,r}(\xi)-d_{c,r}(\eta) \right] \cdot (\xi-\eta)>0
\end{equation}

\noindent
for any $ \xi, \eta \in \re^N,$ with $\eta\neq \xi$.

\end{lemma}

\begin{proof}
Let $ \xi, \eta \in \re^N,$ with $\eta\neq \xi$. If $|\xi|=|\eta|$, then

\[\bigl(d_{c,r}(\xi)- d_{c,r}(\eta)\bigr)\cdot \bigl(\xi-\eta\bigr)=(c+|\eta|^2)^{\frac{r-2}{2}}|\xi-\eta|^2>0.\]

\medskip
\noindent
Let now $|\xi|\neq|\eta|$. First of all, we observe that

\begin{align*}
\displaystyle
 \left[ d_{c,r}(\xi)-d_{c,r}(\eta) \right] \cdot (\xi-\eta) 
 = & \left[ (c+|\xi|^2)^{\frac{r-2}{2}}\xi - (c+|\eta|^2)^{\frac{r-2}{2}}\eta \right] \cdot (\xi-\eta) \\
 = & (c+|\xi|^2)^{\frac{r-2}{2}}|\xi|^2 + (c+|\eta|^2)^{\frac{r-2}{2}}|\eta|^2 -  \xi \cdot \eta \left[ (c+|\xi|^2)^{\frac{r-2}{2}} +  (c+|\eta|^2)^{\frac{r-2}{2}}    \right],
\end{align*}

\medskip
\noindent
hence \eqref{strictmonotncondition} is equivalent to 

\begin{align}\label{equivalentformstrictmonotoncondition}
\displaystyle 
\xi \cdot \eta \left[ (c+|\xi|^2)^{\frac{r-2}{2}} +  (c+|\eta|^2)^{\frac{r-2}{2}}    \right] < |\xi|^2 (c+|\xi|^2)^{\frac{r-2}{2}} + |\eta|^2(c+|\eta|^2)^{\frac{r-2}{2}}.
\end{align}

\medskip
\noindent
Let us consider the real function $f:[0,+ \infty) \to [0, + \infty)$ defined as

$$f(t) = \left(c+t^2\right)^{\frac{r-2}{2}} t.$$

\noindent
If $c = 0$, then $f(t)=t^{r-1}$ is strictly increasing, hence strictly monotone.\\

\noindent
If $ c >0$, for any $t \geq 0$ we have 

\begin{align*}
\displaystyle
 f'(t) = \left(c+t^2\right)^{\frac{r-2}{2}} + t^2(r-2)\left(c+t^2\right)^{\frac{r-4}{2}} = \left(c+t^2\right)^{\frac{r-4}{2}}[ c + (r-1)t^2] >0,
\end{align*} 

\noindent
hence $f(t)$ is strictly increasing, so strictly monotone.

\medskip
\noindent
Therefore

$$ \left(f(s) - f(t)  \right) (s - t) > 0, $$

\noindent
that is

\begin{align*}
\displaystyle
 ts [\left(c+t^2\right)^{\frac{r-2}{2}} + \left(c+s^2\right)^{\frac{r-2}{2}}] < t^2 \left(c+t^2\right)^{\frac{r-2}{2}} + s^2 \left(c+s^2\right)^{\frac{r-2}{2}}
\end{align*}

\noindent
for any $s,t \in \mathbb{R}$ with $s \neq t$.

\medskip
\noindent
By applying previous one with $t=|\xi|$ and $s=|\eta|$, we obtain

\begin{align*}
|\xi||\eta|\left({(c+|\xi|^2)^{\frac{r-2}{2}}}+ {(c+|\eta|^2)^{\frac{r-2}{2}}}\right)   < |\xi|^2{(c+|\xi|^2)^{\frac{r-2}{2}}}+|\eta|^2{(c+|\eta|^2)^{\frac{r-2}{2}}}
\end{align*}

\noindent
by which, using also that $\xi \cdot \eta \leq \left| \xi \cdot \eta \right| \leq |\xi| |\eta|,$ we deduce \eqref{equivalentformstrictmonotoncondition} and so \eqref{strictmonotncondition}.

\end{proof}

\begin{lemma}\label{convexity}
Let $c \geq 0$ and $r>1$. The functional $f:W_0^{1,r}(\Omega) \to \mathbb{R}$ defined as

\begin{equation}\label{funzionaleff}
\displaystyle f(w):= \frac{1}{r} \into \left( c + | \nabla w |^2 \right)^{\frac{r}{2}} \, dx,
\end{equation}

\noindent
is convex.
\end{lemma}

\begin{proof}
Let us begin by proving the convexity of the function  $g:[0,+ \infty) \to [0,+\infty)$ defined as 

$$g(s):=\left(c +s^2\right)^{\frac r 2}.$$

\noindent
If $c=0,$ it's clear.\\
Now let $c > 0$. For any $s \geq 0$ we have

$$ \displaystyle g'(s)= r \; s \left( c  + s^2 \right)^{\frac{r}{2}-1} \geq 0 \qquad \text{and} \qquad  g''(s)= r \left( c + s^2 \right)^{\frac{r}{2} - 2 } \left( c +  (r-1) s^2 \right) \geq 0,$$

\noindent
by which we deduce that $g$ is increasing and convex.

\medskip
\noindent
Let us now consider the function $\psi: \mathbb{R}^N \to \mathbb{R}$ defined as

$$ \displaystyle \psi(\xi):= g(|\xi|)= \left(c +| \xi| ^2\right)^{\frac r 2}.$$

\medskip
\noindent
Considering the triangular inequality of $\left| \; \cdot \;  \right|_{\mathbb{R}^N}$ and the fact that $g$ is increasing and convex, we deduce that $\psi$ in convex, in fact for every $\xi_1, \xi_2 \in \mathbb{R}^N$ and for every  $t \in [0,1]$, we have 

\begin{align*}
\displaystyle
\psi(t\, \xi_1 + (1-t) \, \xi_2) 
\leq g (       t \, \left| \xi_1 \right|   + (1-t) \, \left|  \xi_2 \right| ) 
\leq t \, g( \left| \xi_1 \right| ) + (1-t) \, g( \left| \xi_2 \right|) 
= t \, \psi( \xi_1 ) + (1-t) \, \psi( \xi_2 ).
\end{align*}

\medskip
\noindent
Finally, taking into account the increasing monotonicity and the linearity of the Lebesgue integral, we deduce that the functional $f:W_0^{1,r}(\Omega) \to \mathbb{R}$  defined in \eqref{funzionaleff} is convex.

\end{proof}

\chapter{Anisotropic Quasilinear Systems and Regularity Results}\label{SezReg}

\bigskip
\noindent
In this Chapter, we present new results concerning $L^{\infty}$-estimates for solutions to non-autonomous quasilinear systems involving operators in divergence form and nonlinearities up to critical growth. These operators take the form $-\text{div} [\nabla \Psi (\nabla u)]$, where the function $\Psi: \mathbb{R}^N \to \mathbb{R}$ satisfies a general assumption that involves several nonlinear operators, such as the $p$-Laplacian and $p$-area for $1 < p < N$, arising in various nonlinear physical phenomena.

To prove $L^{\infty}$-estimates for this class of systems, the classical Stampacchia's Lemma \cite{S} could not be directly applied due to the singular cases occurring when $1 < p < 2$. We overcome this difficulty through a careful application of a generalization of Stampacchia's Lemma, established by Mammoliti in \cite{BDO}. In the context of quasilinear elliptic systems of two equations, where the framework is defined by the product space $X := W_0^{1,p}(\Omega) \times W_0^{1,q}(\Omega)$, an additional challenge arises from the potential coupling of $u$ and $v$ within the nonlinearity. Furthermore, since the nonlinearity may exhibit critical growth, we must proceed rigorously by establishing preliminary results up to the critical threshold. Finally, the nonlinearity may also depend on the independent variable $x \in \Omega$ and a parameter $\delta$ in a given interval $ I \subseteq \mathbb{R}$.

We establish the $L^{\infty}$-boundedness of all solutions; furthermore, we show that these solutions are bounded uniformly in $\delta$ on any sufficiently small ball of $X$. In particular, our results can also be applied to establish the $L^{\infty}$-boundedness of solutions to eigenvalue problems for quasilinear elliptic systems, such as the de Thélin eigenvalue problem \cite{DT}, which will be addressed in Chapter \ref{CapitolodeThelin}.

\medskip
This Chapter is based on the paper \cite{BCV3}.

\bigskip

\section{Assumptions and Main results}

\bigskip
\noindent
Let us consider the following non-autonomous quasilinear  system

\begin{equation}\label{pqLinfinitoboundedness}
	\begin{cases}
		\begin{array}{ll}
			-\text{\rm div}  \; [ \nabla \Psi_1 \left( \nabla u  \right)] = H_s(\delta,x, u,v) & \hbox{in} \ \Omega,
			\medskip \\
			-\text{\rm div} \; [ \nabla \Psi_2 \left( \nabla v  \right)]= H_t(\delta,x,u,v) & \hbox{in} \ \Omega, \medskip\\
			u=v=0  & \hbox{on} \ \partial\Omega.
		\end{array}
	\end{cases}
\end{equation}

\medskip
\noindent
Denoting with $I \subseteq \mathbb{R}$ an interval, the  function $H:I\times \Omega \times \re^2 \to \mathbb{R}$  is such that $H(\delta,x,\cdot,\cdot) \in C^1(\re^2, \re)$ for any $\delta \in I$ and for a.e. $x\in \Omega$; $H_s(\delta, \cdot,s,t)$, $H_t(\delta, \cdot,s,t)$ are measurable for every $(\delta,s,t) \in I\times \re^2$, and $H(\delta,x,0,0) \in  L^{1}(\Omega)$ for any $\delta \in I$. We assume:

\begin{itemize}
\item[$(\Psi)$] $\Psi_{1}, \Psi_{2}: \mathbb{R}^N \to \mathbb{R}$ of class $C^1$ with $\Psi_1(0)=0$, $\nabla \Psi_1(0)=0$  and $\Psi_2(0)=0$, $\nabla \Psi_2(0)=0$. Moreover, given $c \geq 0,$  $r >1$ and denoting with $\Psi_{c,r}: \mathbb{R}^N \to \mathbb{R}$ the function defined as

\begin{equation}\label{rArea}
	\displaystyle \Psi_{c,r}(\xi):= \frac{1}{r} \biggl[ \left( c +
	\left| \xi \right|^2\right)^{ \frac{r}{2}} - c^{\frac{r}{2}} \biggr] ,
\end{equation}	

\medskip
\noindent
we assume that there exist $a \geq 0$, $1<p<N$ and $ \frac{1}{p} < \nu_1 \leq C_1$ such that $ \displaystyle \left( \Psi_1 - \nu_1 \Psi_{a,p}    \right)$  and  $ \displaystyle \left( C_1 \Psi_{a,p} - \Psi_1   \right)$  are both convex, and there exist $b \geq 0$, $1<q<N$ and $ \frac{1}{q} < \nu_2 \leq C_2$ such that $ \displaystyle \left( \Psi_2 - \nu_2 \Psi_{b,q}    \right)$ and $ \displaystyle \left( C_2 \Psi_{b,q}   - \Psi_2 \right)$ are both convex.\\
	
\item[$(\mathcal{H})$]
	there exists $C_0>0$ such that for a.e. $x\in \Omega$ and every $(\delta,s,t)\in I\times \re^2$ we have
	
	\[|H_s(\delta,x,s,t)|\leq C_0 \left( 1 + |s|^{p^*-1} +|t|^{q^*\frac{p^*-1}{p^*}}   \right),
	\]
	
	\[|H_t(\delta,x,s,t)|\leq C_0 \left( 1 + |s|^{p^*\frac{q^*-1}{q^*}}+ |t|^{q^*-1}  \right).
	\]

\end{itemize}

Condition $(\Psi)$ was introduced in \cite{CDV}, where Cingolani Degiovanni and Vannella obtained regularity results for the scalar version of system \eqref{pqLinfinitoboundedness}, considering a nonlinearity that is independent of the parameter $\delta$ and is allowed to exhibit critical growth. (See also \cite{dibenedetto1983, GV, lieberman,tolksdorf1983, tolksdorf1984} and \cite{CY,HO} for recent results).  In particular, assumption $(\Psi)$ involves various anisotropic operators. For instance, we mention the interesting papers \cite{ACCFM,ACF} in which, letting $ B(t):= t^p/p \text{ for } t>0,$  Antonini, Cianchi, Ciraolo, Farina and Maz'ya considered  

$$ \displaystyle  \Psi_1=B \circ \mathscr{H},$$

\medskip
\noindent
where the norm $ \mathscr{H}: \mathbb{R}^N \to \mathbb{R}$ (different from the Euclidean norm) is of class $C^2$ in $\mathbb{R}^N \setminus \{O \} $, and  its anisotropic unit ball is uniformly convex (see for instance the paper \cite{CFV} of Cozzi, Farina and Valdinoci). Notice that if $\mathscr{H}$ is the Euclidean norm $\mathscr{H}(\xi)=|\xi|$, then we are in the isotropic case, and choosing $ B(t):= t^p/p \text{ for } t>0,$  we get $-\text{\rm div}  \; [ \nabla \Psi_1 \left( \nabla u  \right)]= -\Delta_p u$, hence we recover the $p$-Laplacian.  From a mathematical point of view, the study of anisotropic operators is interesting since the lack of isotropy produces a richer geometric structure. From a physical point of view, anisotropic operators and setting arise in different contexts, such as surface energy (see e.g. \cite{GIGA}) crystallography (see e.g. \cite{TAYLOR1,TAYLOR2,WULFF}), thermodynamics (see e.g. \cite{GURTIN}), noise-removal procedures in digital image processing (see e.g. \cite{ESE}).

Quasilinear elliptic systems were studied by Marino and Winkert in \cite{MW}, where $L^{\infty}$-estimates results are obtained under general conditions on the data (see assumption $(\tilde{H})$ and Theorem 3.4. in \cite{MW}) that allowed the application of Moser’s iteration technique (see also \cite{CSS,denapolimarani,GPZ,MarPer} for existence results and \cite{BDMS,CM,DeFP,MINGIONE,MMSV} for regularity results).\\
Here we aim to establish $L^\infty$-estimates for systems of the form \eqref{pqLinfinitoboundedness} under the broad assumptions $(\Psi)$ and $(\mathcal{H})$. It is immediate to see that, for any $1<p<N,$ $ a \geq 0$ and $1<q<N,$  $ b \geq 0$, the functions $\Psi_1 =\Psi_{a,p}$ and $\Psi_2 =\Psi_{b,q}$ satisfy assumption $(\Psi)$ (it is enough to take $\nu_1=C_1=1$ and $\nu_2=C_2=1$ respectively). Therefore, taking into account that 

\begin{equation}\label{rAreaGradiente}
\displaystyle  
\nabla \Psi_{c,r}(\xi):=
\begin{cases}
\left( c +
	\left| \xi \right|^2\right)^{ \frac{r-2}{2}} \xi  & \text{ if } \xi \neq 0, \bigskip \\
0 & \text{ if } \xi=0,
\end{cases}
\end{equation}

\medskip
\noindent
system \eqref{pqLinfinitoboundedness} generalizes models involving $(p,q)$-Laplacian and $(p,q)$-area type operators.

In fact, if  we take $\Psi_1 =\Psi_{0,p}$ and $\Psi_2 =\Psi_{0,q}$, and moreover we take a function $H$ that does not depend on the variable $\delta$, system \eqref{pqLinfinitoboundedness} becomes

\begin{equation}\label{pqGENERICO}
	\begin{cases}
		\begin{array}{ll}
			- \Delta_p u = H_s(x, u,v) & \hbox{in} \ \Omega, \medskip
			\\
			- \Delta_q v = H_t(x,u,v) & \hbox{in} \ \Omega, \medskip\\
			u=v=0  & \hbox{on} \ \partial\Omega.
		\end{array}
	\end{cases}
\end{equation}

In the case $p=q=2$, a priori bounds for semilinear elliptic systems similar to \eqref{pqGENERICO} have been established in \cite{CDM1,CDM2,DeY,QS} by employing different approaches.

Now, if we  consider $\Psi_1 =\Psi_{a,p}$ and $\Psi_2 =\Psi_{b,q}$ with $a,b \geq 0$, by \eqref{rAreaGradiente} problem \eqref{pqLinfinitoboundedness} becomes

\begin{equation}\label{pqcasospecifico}
	\begin{cases}
		\begin{array}{ll}
			-\text{\rm div} \left( (a+|\nabla u|^2)^{\frac{p-2}{2}}\nabla u\right) = H_s(\delta,x, u,v) & \hbox{in} \ \Omega;
			 \medskip \\
			-\text{\rm div} \left( (b+|\nabla v|^{2})^{\frac{q-2}{2}}\nabla v\right)= H_t(\delta,x,u,v) & \hbox{in} \ \Omega; \medskip \\
			u=v=0  & \hbox{on} \ \partial\Omega.
		\end{array}
	\end{cases}
\end{equation}

Obtaining $L^{\infty}$-estimates for weak solutions of \eqref{pqcasospecifico} is challenging due to the coupling between $u$ and $v$ in $(\mathcal{H})$. In \cite{CCMV}, Carmona, Cingolani, Mart\'{i}nez-Aparicio and Vannella proved such estimates assuming $a=b$, $2 \leq p, q < N$ and a subcritical nonlinearity $H: \mathbb{R}^2 \to \mathbb{R}$.

Under assumptions $(\Psi)$ and $(\mathcal{H})$, we extend the results established for system \eqref{pqcasospecifico} in several directions. Specifically, we consider the range $1<p,q<N$ and allow the constants $a$ and $b$ to differ. Furthermore, the nonlinearity $H$ may exhibit critical growth and can depend explicitly on both the spatial variable $x$ and a parameter $\delta$. Finally, we treat more general operators involving functions $\Psi_1$ and $\Psi_2$ satisfying condition $(\Psi)$.

To include the singular case (i.e. when $1<p<2$ or $1<q<2$), the classical Stampacchia's Lemma \cite[Lemma 4.1]{S} is no longer applicable. Instead, we employ \cite[Lemma 2.1]{BDO} (stated in Section \ref{SezMainResuLinfinito}), the proof of which relies on the following generalization of Stampacchia's Lemma established in \cite[Lemma A.1]{BDO}:

\begin{lemma} Let $\varphi: \mathbb{R}^+ \to \mathbb{R}^+$ be a non increasing function such that 

$$ \displaystyle \varphi(h) \leq \frac{c_0}{(h-k)^{\rho}}k^{\theta \rho} \, [\varphi(k)]^{1+ \lambda} \qquad \forall \, h>k \geq k_0,$$

\medskip
\noindent
where $c_0>0$, $k_0 \geq 0$, $\rho > 0$, $0 \leq \theta < 1$ and $\lambda>0$. Then there exists $k^*>0$ such that $\varphi(k^*)=0$.
\end{lemma}

Lastly, given that $H$ may grow critically, we generalize the technical lemma from \cite[Lemma 2.1]{Vann1} (see Section \ref{SezPrelResultLinfinito} for the statement), provided the radius $R$ is chosen to be sufficiently small:

\begin{lemma}\label{V0}
Let $r\in (1,N)$ and denote by $r^{*}$  the conjugate Sobolev exponent of~$r$.
For any $\varepsilon>0$ and $w_0\in W^{1,r}_0(\Omega)$,  there exist $\sigma>0$ and $R>0$ such that

\[\int_{\{|w(x)|\,\geq \sigma\}} \hspace{-1mm}|w(x)|^{r^*}\,  dx<\varepsilon \]
	
\medskip 
\noindent 
for any $w\in \overline{B_R(w_0)}=\{ w \in  W^{1,r}_0(\Omega)\ :\ \|w-w_0\|_{1,r}\leq R \}$.
\end{lemma}

\medskip
\noindent
Weak solutions of \eqref{pqLinfinitoboundedness} correspond to the critical points of the Euler functional $I_{ \delta, \Psi_1, \Psi_2}: X \to \mathbb{R}$  defined for any $z=(u,v) \in X$ as

\begin{align*}
	I_{\delta,\Psi_1,\Psi_2}(z) = \into  \Psi_1(\nabla u) \, dx + \into \Psi_2(\nabla v )  \, dx   - \into H(\delta,x,u,v) \,dx.
\end{align*}

\medskip
\noindent
We can now state the main result of this Chapter: 

\begin{theorem}\label{key} 
If $(u,v)$ is a solution of \eqref{pqLinfinitoboundedness} and $(\mathcal{H})$ and $(\Psi)$ hold, then $(u,v)\in \left(L^\infty(\Omega)\right)^2$.
 Moreover, for any fixed $z_0=(u_0,v_0)\in X$  there exists a  suitably small $R >0$ such that for any $\delta \in I$ and  denoting by 

	\[ 
	D_{R, \delta}(z_0):=\{
	z \in 
	X\ :\  \| z-z_0\|\leq R, \   I'_{\delta,\Psi_1,\Psi_2}(z)=0
	\},
	\]

\medskip
\noindent
there exists $C>0$ such that

	\[ \|u\|_\infty,\, \|v\|_\infty\, \leq C
	\qquad  \forall \, z=(u,v)\in D_{R, \delta}(z_0),
	\]

\medskip	
\noindent
 where both $R$ and $C$ depend on $z_0, C_0,\left| \Omega \right|, p,q, N,$ $ a, b, \nu_1$ and $\nu_2$, but not on $\delta$.
		
\end{theorem}

So in the present Chapter we show carefully that, for any arbitrary $z_0 \in X$  there exists $R>0$  such that  the $\left(L^\infty(\Omega)\right)^2$-norm of any weak solution to \eqref{pqLinfinitoboundedness} belonging to $B_R(z_0)$ is bounded by a constant $C$ that is independent of $\delta \in I $.

As a consequence of Theorem \ref{key}, we obtain $L^{\infty}$-boundedness of solutions also for system \eqref{pqGENERICO}. System of the form \eqref{pqGENERICO} include also different types of eigenvalue problem, such as the eigenvalue problem introduced by de Thélin in \cite{DT}:

\begin{equation}\label{EigenDeThelin00}
	\begin{cases}
		\begin{array}{ll}
			- \Delta_p u = \lambda |u|^{\alpha-1} |v|^{\beta+1} u
			&  \text{ in }\Omega,
			\medskip \\
			- \Delta_q v  = \lambda  |u|^{\alpha+1} |v|^{\beta - 1}v 
			&  \text{ in } \Omega, \medskip \\
			u=v=0  & \text{ on }  \partial\Omega,
		\end{array}
	\end{cases}
\end{equation}

\medskip
\noindent
where $1 <p,q < N$ and $\alpha, \beta \geq 0$ satisfy $(\alpha +1)/ p + (\beta +1)/q =1$.

By our main result we deduce that if $(u,v)$ is an eigenfunction associated to an eigenvalue $\bar{\lambda}$ of \eqref{EigenDeThelin00}, then both $ u \in L^\infty(\Omega)$ and $v \in L^{ \infty}(\Omega)$, and their norms in $L^\infty(\Omega)$ are bounded by a constant that depends on $\bar{\lambda}$.

\section{$L^\tau$-Integrability for any $1 < \tau < \infty$}\label{SezPrelResultLinfinito}

\bigskip
\noindent
Let $\Psi_1$ be a function as in assumption $(\Psi)$. Then $\Psi_1$ is strictly convex and satisfies

\begin{equation}\label{disuguaglianza1}
\displaystyle \nu_1 \Psi_{a,p}( \xi) \leq \Psi_1(\xi) \leq C_1 \Psi_{a,p}( \xi) \qquad \forall \, \xi \in \mathbb{R}^N,
\end{equation}

\begin{equation}\label{contoconvessità}
\displaystyle \nu_1( a + \left| \xi \right|^2   )^{\frac{p-2}{2}} |\xi|^2 \leq  \nabla \Psi_1 ( \xi ) \cdot \xi  \leq  C_1( a + \left| \xi \right|^2   )^{\frac{p-2}{2}} |\xi|^2 \qquad \forall \, \xi \in \mathbb{R}^N,
\end{equation}

\medskip
\noindent
and the number $p$ for which $\Psi_1$ satisfies assumption $(\Psi)$ is unique.

\medskip
\noindent
In fact, by convexity of $(\Psi_1 - \nu_1 \Psi_{a,p})$ and $ \displaystyle \left( C_1 \Psi_{a,p} - \Psi_1   \right)$, and by using the fact that $\Psi_1(0)=0$, $\Psi_{a,p}(0)=0$,  $\nabla \Psi_1(0)=0$  and $\nabla \Psi_{a,p}(0)=0$, we easily get \eqref{disuguaglianza1}.
Observe now that by the strict convexity of $\Psi_{a,p}$ (see the Proof of Lemma \ref{convexity}), we get

\begin{align*}
\displaystyle 
 0 < \nu_1 \left( \nabla \Psi_{a,p} (\xi) - \nabla \Psi_{a,p} (\eta)   \right) \cdot \left( \xi - \eta   \right) 
 & \leq \left( \nabla \Psi_1 (\xi) - \nabla \Psi_1 (\eta)   \right) \cdot \left( \xi - \eta   \right) \\
 & \leq C_1 \left( \nabla \Psi_{a,p} (\xi) - \nabla \Psi_{a,p} (\eta)   \right) \cdot \left( \xi - \eta   \right)
\end{align*}

\medskip
\noindent
for any $ \xi, \eta \in \mathbb{R}^N$ with $ \xi \neq \eta$, and so we deduce that $\Psi_1$ is strictly convex and, evaluating in $\eta=0$, we also infer \eqref{contoconvessità}. Moreover, by \eqref{disuguaglianza1}  we also get the uniqueness of $p$.

\medskip
\noindent
Similarly, we deduce that $\Psi_2$ considered in assumption $(\Psi)$ is strictly convex and satisfies

\begin{equation*}
\displaystyle \nu_2 \Psi_{b,q}( \xi) \leq \Psi_2(\xi) \leq C_2 \Psi_{b,q}( \xi) \qquad \forall \, \xi \in \mathbb{R}^N,
\end{equation*}

\begin{equation*}
\displaystyle \nu_2( b + \left| \xi \right|^2   )^{\frac{p-2}{2}} |\xi|^2 \leq  \nabla \Psi_2 ( \xi ) \cdot \xi  \leq  C_2( b + \left| \xi \right|^2   )^{\frac{p-2}{2}} |\xi|^2 \qquad \forall \, \xi \in \mathbb{R}^N,
\end{equation*}

\medskip
\noindent
and the number $q$ for which $\Psi_2$ satisfies assumption $(\Psi)$ is unique.

\medskip
\noindent
Let us now consider the Euler functional $I_{ \delta, \Psi_1, \Psi_2}: X \to \mathbb{R}$ associated to system \eqref{pqLinfinitoboundedness}, defined for any $z=(u,v) \in X$ as

\begin{align}\label{funzionale}
	I_{\delta,\Psi_1,\Psi_2}(z) = \into  \Psi_1(\nabla u) dx + \into \Psi_2(\nabla v ) dx   - \into H(\delta,x,u,v) \,dx. 
\end{align}

\medskip
\noindent
By assumption  $(\mathcal{H})$ and Nemytskii operator theory it follows that the part of \eqref{funzionale} involving the nonlinearity $H$ is of class $C^1$ on $X$. By Lemma \ref{contogradconvess} and Proposition \ref{funzionalediclasseC1} also the part of \eqref{funzionale} involving the anisotropic-type operators is of class $C^1$ on $X$. Hence the functional $I_{\delta, \Psi_1, \Psi_2 }$ is of class $C^1$ on $X$ and for any $z_0=(u_0,v_0)$, $z=(u,v) \in X$ we have

\begin{align*}
	\langle I'_{\delta,\Psi_1,\Psi_2}(z_0), z \rangle 
&= \displaystyle  \into \nabla \Psi_1 \left( \nabla u_0 \right) \cdot \nabla u \, dx + \into\nabla \Psi_2 \left( \nabla v_0  \right) \cdot \nabla v \,dx \\
& - \into   \bigl(H_s(\delta,x,u_0,v_0) u +  H_t(\delta,x,u_0,v_0) v\bigr)\, dx.
\end{align*}

\medskip
\noindent
Since in $(\mathcal{H})$ we are assuming that the function $H$ could grow critically, in order to prove Theorem \ref{key} we need to proceed carefully. Let us recall the following result (see \cite[Lemma 2.1]{Vann1}):

\medskip

\begin{lemma}\label{V2023}
Let $r\in (1,N)$ and denote by $r^{*}$ the conjugate Sobolev exponent of~$r$, namely ${r^*=rN/(N-r)}$.
If $R,\varepsilon>0$, $w_0\in W^{1,r}_0(\Omega)$ and $ \tilde{r} \in [1,r^*)$, there is $\sigma>0$ such that

	\[\int_{\{|w(x)|\,\geq \sigma\}} \hspace{-1mm}
	|w(x)|^{\tilde{r}}\,  dx<\varepsilon
	\]

\medskip 
\noindent
for any $w\in \overline{B_R(w_0)}:=\{w \in  W^{1,r}_0(\Omega)\ :\ \|w-w_0\|_{1,r}\leq R \}$.

\end{lemma}

We need to improve previous lemma up to critical growth, and we can do this provided we take a suitably small radius $R$.

\begin{lemma}\label{V}
Let $r\in (1,N)$ and denote by $r^{*}$  the conjugate Sobolev exponent of~$r$.
For any $\varepsilon>0$ and $w_0\in W^{1,r}_0(\Omega)$,  there exist $\sigma>0$ and $R>0$ such that

\[\int_{\{|w(x)|\,\geq \sigma\}} \hspace{-1mm}|w(x)|^{r^*}\,  dx<\varepsilon \]
	
\medskip 
\noindent 
for any $w \in \overline{B_R(w_0)}$.
\end{lemma}

\begin{proof}
By contradiction, assume that there are $\varepsilon>0$ and $w_0\in W^{1,r}_0(\Omega)$ such that, for any $\sigma>0$ and $R>0$ there exists $w_{\sigma,R} \in \overline{B_R(w_0)}$ such that

$$
\displaystyle \int_{\{|w_{\sigma,R}(x)|\,\geq \sigma\}} \hspace{-1mm} |w_{\sigma,R}(x)|^{r^*}\,  dx \geq \varepsilon.
$$

\medskip
\noindent
Hence, we deduce that for any $R>0$ there exists a sequence $\{ w_k \}_k \subset \overline{B_R(w_0)}$ such that

\begin{equation}\label{calcoloLinfinitoboundedness}
\displaystyle \qquad \qquad \quad \int_{\{|w_{k}(x)|\,\geq k\}} \hspace{-1mm} |w_{k}(x)|^{r^*}\,  dx \geq \varepsilon \qquad \forall \, k \in \mathbb{N}.
\end{equation}

\medskip
\noindent
Up to subsequence, there exists $\bar{w} \in W_0^{1,r}(\Omega)$ such that $w_k$ converges to $\bar{w}$ weakly in $W_0^{1,r}(\Omega)$ and strongly in $L^t(\Omega)$ for any $1 \leq t < r^*$.\\
Denoting with $E_k:= \{ x \in \Omega \, : \, \left| w_{k}(x) \right |\,\geq k \}$, by Chebyshev inequality we know that

$$ \displaystyle \left| E_k \right| \leq \frac{1}{k} \into |w_k| \, dx, $$

\noindent
by which we deduce that $\left| E_k \right| \to 0$ as $k \to + \infty$.\\
Observe now that

\begin{align}\label{calcolo2}
\displaystyle 
\int_{\{|w_k(x)|\,\geq k \}} \hspace{-1mm}|w_k(x)|^{r^*}\,  dx \nonumber
& = \int_{E_k} \hspace{-1mm}|w_k(x)- w_0(x) + w_0(x)|^{r^*}\,  dx \\
& \leq 2^{r^*} \int_{E_k} |w_k(x)- w_0(x)|^{r^*}\,  dx + 2^{r^*} \int_{E_k} |w_0(x)|^{r^*}\,  dx \nonumber\\
& \leq 2^{r^*} \into |w_k(x)- w_0(x)|^{r^*}\,  dx + 2^{r^*} \int_{E_k} |w_0(x)|^{r^*}\,  dx.
\end{align}

\medskip
\noindent
Since $\left| E_k \right| \to 0$ as $k \to + \infty$ and $w_0 \in L^{r^*}(\Omega)$, by absolute continuity of the Lebesgue integral we deduce 

$$ \displaystyle \int_{E_k} |w_0(x)|^{r^*}\,  dx \to 0 \qquad \text{as } k \to + \infty.$$
\medskip
\noindent
Taking into account Sobolev-Gagliardo-Nirenberg inequality, there exists a constant $C:=C(r,N)$ such that

$$ \displaystyle \into |w_k(x)- w_0(x)|^{r^*}\,  dx \leq C \|w_k-w_0\|_{1,r}^{r^*} \leq C R^{r^*} \qquad \forall k \in \mathbb{N}.$$

\medskip
\noindent
By choosing  $R= \frac{1}{2C^{\frac{1}{r^*}}} \bar{\varepsilon}^{\frac{1}{r^*}}$ with $0 < \bar{\varepsilon} < \varepsilon$ and passing to $ \displaystyle \limsup_{k \to  \infty}$ in \eqref{calcolo2}, we obtain 

$$ \displaystyle \limsup_{k \to  \infty} \int_{\{|w_k(x)|\,\geq n \}} \hspace{-1mm}|w_k(x)|^{r^*}\,  dx \leq \bar{\varepsilon},$$

\medskip
\noindent
that is a contradiction with \eqref{calcoloLinfinitoboundedness}. 
\end{proof}

\medskip
\noindent
To prove Theorem \ref{key}, inspired by \cite{CCMV}, \cite{GV} and \cite{Vann1},  we firstly prove the following result.

\begin{proposition}\label{integrability}
 For any fixed $z_0=(u_0,v_0)\in X$  there exists a  suitably small $R >0$ such that for any $\delta \in I$,  denoting by

	\[ 
	D_{R, \delta}(z_0)=\{
	z\in 
	X\ :\  \| z-z_0\|\leq R, \   I'_{\delta,\Psi_1,\Psi_2}(z)=0
	\},
	\]

\medskip
\noindent
and for any $\gamma>1$ there exists $C>0$ such that 

\begin{equation*}
	\into |\bar u|^{\gamma  p^*} \, dx \leq C \quad \text{and} \quad \into |\bar v|^{\gamma  q^*}\, dx  \leq  C  
	\qquad \qquad \forall \bar{z}=(\bar u, \bar v) \in D_{R, \delta}(z_0),
\end{equation*}

\medskip	
\noindent
 where both  $R$ and $C$ depend on $z_0, C_0,\left| \Omega \right|, p,q, N,$ $ a, b, \nu_1$, $\nu_2$ and $\gamma$ but not on $\delta$.

\end{proposition}

\begin{proof} For simplicity, we will omit $dx$. For every $\gamma,\, t,\, k>1$ we define

\begin{align*}
	h_{k,\gamma}(s):=&
	\begin{cases}
		s|s|^{\gamma-1} & \hbox{if }|s|\leq k,
		 \medskip \\
		\gamma k^{\gamma -1} s + \text{sign}(s)(1-\gamma) k^{\gamma} & \hbox{if } |s|>k,
	\end{cases}
	\\
	\\
	\Phi_{k,t,\gamma}(s) :=&\int_{0}^{s}\left| h'_{k,\gamma}(r)\right|^{\frac
		t\gamma}dr.
\end{align*}

\medskip
\noindent
Observe that $h_{k,\gamma}$ and $\Phi_{k,t,\gamma}$ are
$C^{1}$-functions with bounded derivative, depending on $\gamma, t$ and $k$.\\
In particular, if $w\in \sobr$  then $h_{k,\gamma}(w)$ and $\Phi_{k,t,\gamma}(w)$ are in $W_{0}^{1,r}(\Omega)$.

\medskip
\noindent
We see that 
\begin{equation}\label{hnew}
	|h'_{k,\gamma}(s)| \leq \gamma|h_{k,\gamma}(s)|^{\frac{\gamma-1}{\gamma}},
\end{equation}

\medskip
\noindent
and for every $t\geq \gamma$  there exists a positive constant $C$, depending on $\gamma$ and $t$ but independent of $k$, such that

\begin{equation}  \label{desgv}
	|s|^{\frac{t}{\gamma}-1}|\Phi_{k,t,\gamma}(s)|\leq C
	|h_{k,\gamma}(s)|^{\frac t\gamma}
\end{equation}

\begin{equation}\label{desgv2}
	|\Phi_{k,t,\gamma}(s)|\leq
	C |h_{k,\gamma}(s)|^{\frac1\gamma\left(1+t\frac {\gamma-1}\gamma \right)}
\end{equation}

\medskip
\noindent
and

\begin{equation}\label{desgv3}
	\left|h_{k,\gamma}\left(|s|^{\frac{q^*}{p^*}}\right)\right|^{p^*}
	\leq C\left|h_{k^{\frac{p^*}{q^*}},\gamma}\left(s\right)\right|^{q^*}.
\end{equation}

\medskip
\noindent
Moreover, for any $u \in  W_0^{1,p}(\Omega)$ we have

\begin{align}\label{gradhkgamma}
\displaystyle
 \nabla h_{k,\gamma}(u) = h'_{k,\gamma} (u) \nabla u,
\end{align}

\medskip
\noindent
and

\begin{align}\label{gradPhikgamma}
\displaystyle
\nabla \Phi_{k,\gamma p,\gamma }(u) = \Phi'_{k,\gamma p,\gamma } (u) \nabla \bar u = \left|h'_{k,\gamma} (u)\right|^p  \nabla u. 
\end{align}

\medskip
\noindent
For any $p>1$, $a \geq 0$ and $s \geq 0$, we have

\begin{equation}\label{disr}
	s^p\leq \left(a +s^2\right)^{\frac{p-2}{2}}s^2  
	+ a^{\frac{p}{2}}.
\end{equation}

\medskip
\noindent
The inequality is obvious if $p\geq 2$ or $s=0$.\\
Otherwise, if $p\in(1,2)$ and $s \neq 0$, we have

$$ \displaystyle s^p \leq \left(a +s^2\right)^{\frac{p}{2}}=\left(a +s^2\right) \left(a +s^2\right)^{\frac{p-2}{2}}
\leq
\left(a +s^2\right)^{\frac{p-2}{2}}s^2 + a^{\frac{p}{2}}.
$$

\medskip
\noindent
Let us fix $R>0$. For any  $\delta \in I$, let us consider  $\bar z=(\bar u, \bar v) \in D_{R,  \delta}(z_0)$. 

\medskip
\noindent
In particular,  

$$ \left\langle I'_{\delta,\Psi_1,\Psi_2}(\bar z), \left( \Phi_{k,\gamma p,\gamma}(\bar u),0 \right)\right\rangle=0 \qquad \text{ for any } k, \gamma >1.$$

\medskip
\noindent
By Sobolev-Gagliardo-Nirenberg inequality, there exists $c:=c(p, N)>0$ such that

$$ \displaystyle \left( \into |u|^{p^*} \right)^{\frac{p}{p^*}} \leq c \into | \nabla u|^p   \qquad \text{ for any } u \in W_0^{1,p}(\Omega).$$

\medskip
\noindent
Furthermore, taking into account \eqref{contoconvessità}, \eqref{gradhkgamma}, \eqref{gradPhikgamma} and (\ref{disr}), we have 

\begin{align*}
\left(\into |h_{k,\gamma}(\bar u)|^{p^*}\right)^{\frac
p{p^*}}
& \leq c \into |\nabla \bar
	u|^{p}|h'_{k,\gamma}(\bar u)|^{p} \nonumber
	\\
	&
	\leq c \into (a + |\nabla \bar u|^2)^{\frac{p-2}{2}}|\nabla \bar
	u|^{2} \ |h'_{k,\gamma}(\bar u)|^{p} + c \, a^{\frac{p}{2}} \into |h'_{k,\gamma}(\bar u)|^{p} \nonumber\\
	&
\leq \frac{c}{\nu_1} \into \nabla \Psi_1 ( \nabla \bar{u} ) \cdot \nabla \bar{u} \ |h'_{k,\gamma}(\bar u)|^{p} + c \, a^{\frac{p}{2}} \into |h'_{k,\gamma}(\bar u)|^{p}\nonumber\\
	&
= \frac{c}{\nu_1} \into \nabla \Psi_1 ( \nabla \bar{u} ) \cdot \nabla \Phi_{k,\gamma p,\gamma }(\bar u)   + c \, a^{\frac{p}{2}} \into |h'_{k,\gamma}(\bar u)|^{p} \nonumber\\
	& = \frac{c}{\nu_1} \into H_s(\delta,x,\bar u,\bar v) \Phi_{k,\gamma p,\gamma }(\bar u)
	+ c \, a^{\frac{p}{2}} \into |h'_{k,\gamma}(\bar u)|^{p}  \\
		& \leq \tilde{c}\into \left| H_s(\delta,x,\bar u,\bar v) \right| \left| \Phi_{k,\gamma p,\gamma }(\bar u) \right| 
	+\tilde{c} \into |h'_{k,\gamma}(\bar u)|^{p} .
\end{align*}

\medskip
\noindent
Considering also $(\mathcal{H})$, we get

\begin{align}\label{integ1}
\left(\into |h_{k,\gamma}(\bar u)|^{p^*}\right)^{\frac	p{p^*}} 
\leq \hat{c} \left(\into (|\bar u|^{p^*-1}+1)|\Phi_{k,\gamma p,\gamma }(\bar u)|
	+ \into |\bar v|^{q^*\frac{p^*-1}{p^*}}|\Phi_{k,\gamma p,\gamma }(\bar u)| 
	\right) \
+ \hat{c} \, \into |h'_{k,\gamma}(\bar u)|^{p}.
\end{align}

\medskip
\noindent
By (\ref{hnew}) and H\"older's inequality with $\frac{p^*\gamma}{p(\gamma-1)}$ and $ \left( \frac{p^*\gamma}{p(\gamma-1)} \right)'$, we get

\begin{equation}\label{hne}
	\hat{c} \, \into |h'_{k,\gamma}(\bar u)|^{p}\leq \hat{c} \, \gamma^p \,\into |h_{k,\gamma}(\bar u)|^{p \frac{\gamma - 1}{\gamma}} \leq 
	c_0 \left(\into |h_{k,\gamma}(\bar u)|^{p^*}\right)^{\frac{p}{p^*}\frac{\gamma-1}{\gamma}},
\end{equation}

\medskip
\noindent
where $c_0>0$ depends on  $C_0, \left| \Omega \right|, p, N, a, \nu_1$ and $\gamma$  but not on $k$ and $\delta$.

\medskip
\noindent
For any $\sigma>1$, $r>1$ and $w$ in $\sobr$, we denote by

\[\Omega_{\sigma,w}:=\{x\in \Omega\ : \ |w(x)| > \sigma\}.\]

\medskip
\noindent
From now on, when necessary, we redefine the positive constant $c_1$,  depending on $C_0,\left| \Omega \right|, p,q, N, a, \nu_1$   and $\gamma$  but not on $k$ and $\delta$. Therefore, using \eqref{desgv}, \eqref{desgv2} and H\"older inequality with conjugated exponents $\frac{p^* \gamma}{p \gamma + 1 - p}$ and $\left( \frac{p^* \gamma}{p \gamma + 1 - p} \right)'$, and with conjugated exponents $\frac{p^*}{p^*-p}$ and  $\frac{p^*}{p}$ respectively, we have

\begin{align}\label{integ2}
\displaystyle
\into (|\bar u|^{p^*-1}+1)|\Phi_{k,\gamma p,\gamma }(\bar u)| 
& = \into |\Phi_{k,\gamma p,\gamma }(\bar u)| + \int_{\Omega \setminus \Omega_{\sigma,\bar{u}}} |\bar u|^{p^*-1}|\Phi_{k,\gamma p,\gamma }(\bar u)| + \int_{\Omega_{\sigma,\bar{u}}} |\bar u|^{p^*-1}|\Phi_{k,\gamma p,\gamma }(\bar u)| \nonumber \\
& \leq (\sigma^{p^*-1}+1)\into|\Phi_{k,\gamma p,\gamma }(\bar u)| +\int_{\Omega_{\sigma,\bar u}} |\bar u|^{p^*-p} |\bar u|^{p-1} |\Phi_{k,\gamma p,\gamma }(\bar u)| \nonumber\\
&\leq 2\sigma^{p^*-1}\into|\Phi_{k,\gamma p,\gamma }(\bar u)|  + c_1\int_{\Omega_{\sigma,\bar u}} |\bar u|^{p^*-p} |h_{k,\gamma}(\bar u)|^p  \\
&\leq c_1\sigma^{p^*-1}\into |h_{k,\gamma}(\bar u)|^{\frac{p\gamma +1 -p}{\gamma}} +c_1\int_{\Omega_{\sigma,\bar u}} |\bar u|^{p^*-p} |h_{k,\gamma}(\bar u)|^p. \nonumber \\
&\leq c_1 \sigma^{p^*-1}\left(\into |h_{k,\gamma}(\bar u)|^{p^*}\right)^{\frac{p}{p^*}\frac{\gamma p+1-p}{\gamma p}} + c_1\|\bar u  \|^{p^*-p}_{L^{p^*}(\Omega_{\sigma,\bar u})}  \left(\into
	\left|h_{k,\gamma }(\bar u)\right|^{p^*}\right)^{\frac p {p^*}}. \nonumber 
\end{align}

\noindent
Now, we deal with the integral $\into |\bar v|^{q^*\frac{p^*-1}{p^*}}|\Phi_{k,\gamma p,\gamma }(\bar u)| $ in \eqref{integ1}.
By using \eqref{desgv}, \eqref{desgv2}, \eqref{desgv3}, the fact that $\Phi_{k,\gamma p,\gamma}(s)$ is nondecreasing for $s \geq 0$ and H\"older inequality as before, we obtain

\begin{align*} 
\into |\bar v|^{q^*\frac{p^*-1}{p^*}}|\Phi_{k,\gamma p,\gamma }(\bar u)|
	& = \int_{ \Omega \setminus \Omega_{\sigma,\bar v}}  |\bar v|^{q^*\frac{p^*-1}{p^*}}|\Phi_{k,\gamma p,\gamma }(\bar u)| + \int_{ \Omega_{\sigma,\bar v}} |\bar v|^{\frac{q^*}{p^*}(p^*-p)}|\bar v|^{\frac{q^*}{p^*}(p-1)}|\Phi_{k,\gamma p,\gamma }(\bar u)| \\
	&\leq  \sigma^{q^*\frac{p^*-1}{p^*}}\into |\Phi_{k,\gamma p,\gamma }(\bar u)| 
	+\int_{\Omega_{\sigma,\bar v}\cap \{|\bar v|^\frac{q^*}{p^*}\leq |\bar u|\}} |\bar v|^{\frac{q^*}{p^*}(p^*-p)} |\bar u|^{p-1} |\Phi_{k,\gamma p,\gamma }(\bar u)| 
	\\
	& \quad +\int_{\Omega_{\sigma,\bar v}\cap \{|\bar v|^\frac{q^*}{p^*}\geq |\bar u|\}} |\bar v|^{\frac{q^*}{p^*}(p^*-p)} \bigl(|\bar v|^{\frac{q^*}{p^*}}\bigr)^{p-1} |\Phi_{k,\gamma p,\gamma }(|\bar v|^{\frac{q^*}{p^*}})| 
	\\
	&\leq c_1\sigma^{q^*\frac{p^*-1}{p^*}}\into |h_{k,\gamma}(\bar u)|^{\frac{p\gamma +1 -p}{\gamma}}
	+c_1 \int_{\Omega_{\sigma,\bar v}\cap \{|\bar v|^\frac{q^*}{p^*}\leq |\bar u|\}} |\bar v|^{\frac{q^*}{p^*}(p^*-p)}|h_{k,\gamma}(\bar u)|^p
	\\
	&\quad +c_1 \int_{\Omega_{\sigma,\bar v}\cap \{|\bar v|^\frac{q^*}{p^*}\geq |\bar u|\}} |\bar v|^{\frac{q^*}{p^*}(p^*-p)}
	|h_{k,\gamma}\bigl(|\bar v|^{\frac{q^*}{p^*}}\bigr)|^p
	\\
	&\leq c_1\sigma^{q^*\frac{p^*-1}{p^*}} 
	\Bigl(\into |h_{k,\gamma
	}(\bar u)|^{p^*}\Bigr)^{\frac{p}{p^*}\frac{\gamma p+1-p}{\gamma p}}\\
	&\quad +c_1 \|\bar v \|^{\frac{q^*}{p^*}(p^*-p)}_{L^{q^*}(\Omega_{\sigma,\bar v})} 
	\left(	
	\left(\into
	\left|h_{k,\gamma }(\bar u)\right|^{p^*}\right)^{\frac p {p^*}}
	+\left(\int_{\Omega}
	\bigg| h_{k^{\frac{p^*}{q^*}},\gamma }\left(\bar v\right)\bigg|^{q^*}\right)^{\frac p {p^*}}
	\right).
	\nonumber
\end{align*}

\medskip
\noindent
Combining with  \eqref{integ1}, \eqref{hne} and \eqref{integ2}, we get

\begin{align*}
\left(\into |h_{k,\gamma}(\bar u)|^{p^*}\right)^{\frac
		p{p^*}} 
	&\leq c_1
	\sigma^{p^*-1}\left(\into |h_{k,\gamma
	}(\bar u)|^{p^*}\right)^{\frac{p}{p^*}\frac{\gamma p+1-p}{\gamma p}}
	+ c_1\|\bar u  \|^{p^*-p}_{L^{p^*}(\Omega_{\sigma,\bar u})}  \left(\into
	\left|h_{k,\gamma }(\bar u)\right|^{p^*}\right)^{\frac p {p^*}}
	\\
	& \quad + c_1\sigma^{\frac{q^*}{p^*}(p^*-1)} 
	\left( \into |h_{k,\gamma
	}(\bar u)|^{p^*}\right)^{\frac{p}{p^*}\frac{\gamma p+1-p}{\gamma p}}\\
	&\quad +c_1 \|\bar v \|^{\frac{q^*}{p^*}(p^*-p)}_{L^{q^*}(\Omega_{\sigma,\bar v})} 
	\left(	
	\left(\into
	\left|h_{k,\gamma }(\bar u)\right|^{p^*}\right)^{\frac p {p^*}}
	+\left(\into
	\bigg| h_{k^{\frac{p^*}{q^*}},\gamma }\left(\bar v\right) \bigg|^{q^*}\right)^{\frac p {p^*}}
	\right)\\
	&\quad +c_1 \left(\into |h_{k,\gamma}(\bar u)|^{p^*}\right)^{\frac{p}{p^*}\frac{\gamma-1}{\gamma}},
\end{align*}

\medskip
\noindent
where the constant $c_1$ depends on $C_0,\left| \Omega \right|, p,q, N, a, \nu_1$  and $\gamma$  but not on $k$ and $\delta$.

\bigskip
\noindent
By using Lemma \ref{V} with $\varepsilon=(4c_1)^{\frac{p^*}{p-p^*}}$ and $u_0 \in W_0^{1,p}(\Omega)$, and  with $\varepsilon=(4 c_1)^{\frac{p^*}{p-p^*}}$ and $v_0 \in W_0^{1,q}(\Omega)$ respectively, we infer that there exist  $\sigma_1>1$ and $R_1>0$ such that, for any $\sigma\geq \sigma_1$ and for any $(\bar u,\bar v)\in D_{R_1, \delta}(z_0),$

$$ \displaystyle   c_1\|\bar u  \|^{p^*-p}_{L^{p^*}(\Omega_{\sigma,\bar u})} + c_1 \|\bar v \|^{\frac{q^*}{p^*}(p^*-p)}_{L^{q^*}(\Omega_{\sigma,\bar v})} \leq \frac{1}{2}.$$

\medskip
\noindent
Hence, we infer that 

\begin{align}\label{k0}
\displaystyle
        \frac{1}{2} \left(\into |h_{k,\gamma}(\bar u)|^{p^*}\right)^{\frac{p}{p^*}} 
& \leq  c_1 \,  \Bigl(\sigma^{p^*-1}+\sigma^{\frac{q^*}{p^*}(p^*-1)} \Bigr) \left(\into |h_{k,\gamma}(\bar u)|^{p^*}\right)^{\frac{p}{p^*}\frac{\gamma p+1-p}{\gamma p}} \\
&	\quad +  c_1  \,  \| \bar v \|^{\frac{q^*}{p^*}(p^*-p)}_{L^{q^*}(\Omega_{\sigma,\bar v})} \left(\into
	\bigg| h_{k^{\frac{p^*}{q^*}},\gamma }\left(\bar v\right) \bigg|^{q^*}\right)^{\frac p {p^*}}\nonumber\\
&	\quad +  c_1  \,  \left(\into |h_{k,\gamma}(\bar u)|^{p^*}\right)^{\frac{p}{p^*}\frac{\gamma-1}{\gamma}} . \nonumber
\end{align}

\medskip
\noindent
If $\eta \in (0,1)$, by using Young inequality with conjugated exponents $\frac{1}{\eta}$ and $ \frac{1}{1-\eta}$, we obtain that

\begin{equation}\label{inyou}
\displaystyle	ax^{\eta} \leq \left( \frac{x}{8} \right)^{\eta} 8a \leq \frac x 8 + 
	(8a)^{1/(1-\eta)} \qquad \forall \, a,x \geq 0.
\end{equation}

\bigskip
\noindent
Therefore,  since  $ \frac{\gamma p+1-p}{\gamma p}<1$ and $ \frac{\gamma-1}{\gamma}<1$, we infer

\begin{align*}
	c_1 \Bigl(\sigma^{p^*-1}+\sigma^{\frac{q^*}{p^*}(p^*-1)}\Bigr)
\left(\into |h_{k,\gamma
	}(\bar u)|^{p^*}\right)^{\frac{p}{p^*}\frac{\gamma p+1-p}{\gamma p}}
	\hspace{-4mm}\leq \ \frac 1 8 \left(\into |h_{k,\gamma
	}(\bar u)|^{p^*}\right)^{\frac{p}{p^*}}+c_1\, \Bigl(\sigma^{p^*-1}+\sigma^{\frac{q^*}{p^*}(p^*-1)}\Bigr)^{\frac{\gamma p}{p-1}}
\end{align*}

\medskip
\noindent
and 

\[c_1 \left(\into |h_{k,\gamma}(\bar u)|^{p^*}\right)^{\frac{p}{p^*}\frac{\gamma-1}{\gamma}}
\leq \frac 1 8 \left(\into |h_{k,\gamma
}(\bar u)|^{p^*}\right)^{\frac{p}{p^*}}+c_1,
\]

\medskip
\noindent
so that \eqref{k0} becomes

\begin{align*} 
\frac 1 4\left(\into |h_{k,\gamma}(\bar u)|^{p^*}\right)^{\frac
		p{p^*}} \leq 
	c_1 \left[ \Bigl(\sigma^{p^*-1}+\sigma^{\frac{q^*}{p^*}(p^*-1)}\Bigr)^{\frac{\gamma p}{p-1}} + 1 \right]
	+  c_1 \, \|\bar v \|^{\frac{q^*}{p^*}(p^*-p)}_{L^{q^*}(\Omega_{\sigma,\bar v})} 
	\left(\into
	\bigg|h_{k^{\frac{p^*}{q^*}},\gamma }\left(\bar v\right)\bigg|^{q^*}\right)^{\frac p {p^*}}\hspace{-3mm}.
\end{align*}

\bigskip
\noindent
Thus,  we have shown that for any 
$k,\gamma >1$ there exists a positive constant $c_1$ depending on $C_0,\left| \Omega \right|, p,q, N, a, \nu_1$  and $\gamma$  but not on $k$ and $\delta$, and there are  $\sigma_1>1$ and a suitably small $R_1>0$, depending on $c_1$ and $(u_0,v_0)$,  such that for any $(\bar u, \bar v) \in D_{R_1,  \delta}(z_0)$ and any $\sigma\geq \sigma_1$ we have

\begin{align}\label{integ4}
\into |h_{k,\gamma}(\bar u)|^{p^*}\leq 
	c_1\left[\Bigl(\sigma^{p^*-1}+\sigma^{\frac{q^*}{p^*}(p^*-1)}\Bigr)^{\frac{\gamma p^*}{p-1}}
	 + 1 \right] +  c_1\, \|\bar v \|^{\frac{q^*}{p}(p^*-p)}_{L^{q^*}(\Omega_{\sigma,\bar v})} 
	\into
	\bigg|h_{k^{\frac{p^*}{q^*}},\gamma }\left(\bar v\right)\bigg|^{q^*}. 
\end{align}

\medskip
\noindent
Reasoning in a similar way, exploiting that 

$$\left\langle I'_{\delta,\Psi_1,\Psi_2}(\bar z), \left(0,\Phi_{\tilde k,\gamma q,\gamma}(\bar v) \right)\right\rangle=0 \qquad \text{ for any } \tilde k, \gamma>1,$$

\medskip
\noindent
we infer that for any $\tilde{k},\gamma >1$ there exists a positive constant $c_2$ depending on $C_0,\left| \Omega \right|, p,q, N, b, \nu_2$  and $\gamma$  but not on $\tilde{k}$ and $\delta$, and there are  $\sigma_2>1$ and a suitably small $R_2>0$, depending on $c_2$ and $(u_0,v_0)$,  such that for any $(\bar u, \bar v) \in D_{R_2,  \delta}(z_0)$ and any $\sigma\geq \sigma_2$ we have

\begin{align}\label{integ5}
\into |h_{\tilde k,\gamma}(\bar v)|^{q^*}
\leq c_2 \left[ \Bigl(\sigma^{q^*-1}+\sigma^{\frac{p^*}{q^*}(q^*-1)} \Bigr)^{\frac{\gamma q^*}{q-1}}
	+1 \right] +  c_2 \, \|\bar u \|^{\frac{p^*}{q}(q^*-q)}_{L^{p^*}(\Omega_{\sigma,\bar u})} 
	\into
	\bigg|h_{\tilde k^{\frac{q^*}{p^*}},\gamma }\left(\bar u\right)\bigg|^{p^*}.
\end{align}

\bigskip
\noindent
Let us consider a suitable $\bar{\sigma}>1$ and a suitably small $\bar{R}>0$ such that both \eqref{integ4} and \eqref{integ5} hold for any $\sigma \geq \bar{\sigma}$ and for any $0<R \leq \bar{R}$.\\
Setting $\tilde{k}=k^{\frac{p^*}{q^*}}$ in (\ref{integ5}) and substituting in (\ref{integ4})  we obtain

\begin{align*}
\displaystyle
    \into |h_{k,\gamma}(\bar u)|^{p^*} 
&	\leq c_1 \left[ \Bigl(\sigma^{p^*-1}+\sigma^{\frac{q^*}{p^*}(p^*-1)} \Bigr)^{\frac{\gamma p^*}{p-1}} + 1 \right]\\
&	\quad +  c_1 \, c_2 \, \|\bar v \|^{\frac{q^*}{p}(p^*-p)}_{L^{q^*}(\Omega_{\sigma,\bar v})} 
	 \left[ \Bigl(\sigma^{q^*-1}+\sigma^{\frac{p^*}{q^*}(q^*-1)} \Bigr)^{\frac{\gamma q^*}{q-1}}
	+ 1 \right] \nonumber  \\
& \quad  + c_1 \, c_2 \, \|\bar v \|^{\frac{q^*}{p}(p^*-p)}_{L^{q^*}(\Omega_{\sigma,\bar v})}  \, \|\bar u \|^{\frac{p^*}{q}(q^*-q)}_{L^{p^*}(\Omega_{\sigma,\bar u})} \into|h_{k,\gamma }\left(\bar u\right)|^{p^*}.
\end{align*}

\medskip
\noindent
By using again Lemma \ref{V} with $\varepsilon=2^{\frac{q}{q-q^*}}$ and $u_0 \in W_0^{1,p}(\Omega)$, and  with $\varepsilon=(c_1c_2)^{\frac{p}{p-p^*}}$ and $v_0 \in W_0^{1,q}(\Omega)$ respectively, we infer that there exist  $\tilde{\sigma} \geq \bar{\sigma}$ and $ \tilde{R}  \in (0,\bar{R}]$ such that, for any $\sigma\geq \tilde{\sigma}$ and for any $(\bar u,\bar v)\in D_{\tilde{R}, \delta}(z_0),$

$$ \displaystyle  c_1 c_2\, \|\bar v \|^{\frac{q^*}{p}(p^*-p)}_{L^{q^*}(\Omega_{\sigma,\bar v})} \leq 1 \quad \text{and} \quad \|\bar u \|^{\frac{p^*}{q}(q^*-q)}_{L^{p^*}(\Omega_{\sigma,\bar u})} \leq \frac{1}{2}.$$

\medskip
\noindent
Hence, we finally get that for any $(\bar u, \bar v) \in D_{R, \delta}(z_0)$

\[\into |h_{k,\gamma}(\bar u)|^{p^*} \leq c_3  \qquad \text{ for any } k,\gamma>1, \]

\noindent
and similarly

\[\into |h_{k,\gamma}(\bar v)|^{q^*} \leq c_3 \qquad \text{ for any } k,\gamma>1, \]

\medskip
\noindent
where the positive constant $c_3$ depends on $C_0,\left| \Omega \right|, p,q, N, a, \nu_1, b, \nu_2$  and $\gamma$  but not on $k$ and $\delta$, and in addition $R>0$ is a suitable radius.

\medskip
\noindent
Thus we can apply Fatou Lemma, and passing to limit for $k \to \infty$ we get 

\begin{equation*}
	\into |\bar u|^{\gamma  p^*} \leq C \quad \text{and} \quad \into |\bar v|^{\gamma  q^*} \leq  C 
	\qquad \qquad \forall (\bar u, \bar v) \in D_{R, \delta}(z_0),
\end{equation*}

\noindent
where the positive constant $C$ still depends on $C_0,\left| \Omega \right|, p,q, N, a, \nu_1, b, \nu_2$  and $\gamma$ but not on $\delta$.

\end{proof}

\section{Uniformly locally $L^{\infty}$-estimate}\label{SezMainResuLinfinito}

\medskip
\noindent
In this Section we prove Theorem~\ref{key}, and we will use a result proved in \cite[Appendix]{BDO}. In particular, in \cite[Lemma A.1]{BDO} it was shown the following generalization of Stampacchia Lemma:

\begin{lemma}\label{GeneralizzazioneStampacchiaLemma}
Let $\varphi: \mathbb{R}^+ \to \mathbb{R}^+$ be a nonincreasing function such that 

$$ \displaystyle \varphi(h) \leq \frac{c_0}{(h-k)^{\rho}}k^{\theta \rho} \, [\varphi(k)]^{1+ \lambda} \qquad \forall \, h>k \geq k_0,$$

\medskip
\noindent
where $c_0>0$, $k_0 \geq 0$, $\rho > 0$, $0 \leq \theta < 1$ and $\lambda>0$. Then there exists $k^*>0$ such that $\varphi(k^*)=0$.
\end{lemma}

\medskip
\noindent
Thanks to previous result, in \cite[Lemma 2.1]{BDO} it was proved the following Theorem:

\begin{theorem}\label{Mammoliti}
Let $r>1$ and let $w$ be a function in $ W_0^{1,r}(\Omega)$ such that, for $k$ greater then some $k_0$, 

\begin{equation}\label{inequalityBDO}
\displaystyle
\int_{A_k} \left| \nabla w  \right|^r \, dx \leq c k^{\theta r} \left|  A_k \right|^{\frac{r}{r^*} + \varepsilon},
\end{equation}

\medskip
\noindent
where $\varepsilon >0$, $0 \leq \theta < 1$, $r^*=\frac{Nr}{N-r}$ and 
$$A_k= \left\{ x\in \Omega\ :\ \left|w (x) \right| > k\right\}.$$

\medskip
\noindent
Then, the norm of $w$ in $L^{\infty}(\Omega)$ is bounded by a constant that depends  on  $c,\theta,r, N, \varepsilon, k_0$ and $ \left| \Omega \right|$. 
\end{theorem}

\medskip
\noindent
{\mbox {\it Proof of Theorem~\ref{key}.~}} From now on, we will omit $dx$ and we will redefine the positive constant $C$ when necessary. Fixing $m>N$, by Proposition \ref{integrability} and assumption  $(\mathcal{H})$ we infer that there are $C, R >0$ depending on $C_0,\left| \Omega \right|, p,q, N, a, \nu_1, b, \nu_2$ and $m$  but not on $\delta$, such that for any $(\bar u, \bar v) \in D_{R, \delta}(u_0,v_0)$ we have $H_s(\delta,x,\bar u,\bar v) \in L^m(\Omega)$, $ H_t(\delta,x,\bar u,\bar v)\in L^m(\Omega)$ and

\begin{equation}\label{dis1}
\displaystyle \|H_s(\delta,x,\bar u,\bar v)\|_m^{p'}, \quad  \|H_t(\delta,x,\bar u,\bar v)\|_m^{q'} \leq C.
\end{equation}

\medskip
\noindent 
Fixed $(\bar u, \bar v) \in D_{R, \delta}(z_0)$, for any $k \geq 1$ we denote by  

$$A^{\bar{u}}_k= \left\{ x\in \Omega\ :\ | \bar{u} (x)| > k\right\},$$

\medskip
\noindent
and we consider the map $\mathcal{R}_k: \mathbb{R} \to \mathbb{R}$ defined as 

$$ \displaystyle
\mathcal{R}_k(r):=
\begin{cases}
0 & \text{ if } |r| \leq k, \bigskip \\
 \text{sign}(r) \Bigl[ \, |r| - k \, \Bigr] & \text{ if } |r| > k. \\
\end{cases}
$$

\noindent
Since $\mathcal{R}_k(\bar{u}) \in W_0^{1,p}(\Omega)$ and $\langle I'_{\delta,\Psi_1,\Psi_2}(\bar{u},\bar{v}),\left(\mathcal{R}_k(\bar u),0\right)\rangle=0,$  by using \eqref{contoconvessità}, \eqref{disr},   and  denoting by $\bar{f}_s:=H_s(\delta,x, \bar{u}, \bar{v}) $, we have

\begin{align}\label{calc1}
\displaystyle
 \int_{A^{\bar{u}}_k} \left| \nabla \bar{u} \right|^p & = \int_{A^{\bar{u}}_k} \left| \nabla \mathcal{R}_k( \bar{u}) \right|^p =  \int_{A^{\bar{u}}_k} \left| \nabla \bar{u} \right|^p \mathcal{R}_k'(\bar{u})^p = \int_{A^{\bar{u}}_k} \left| \nabla \bar{u} \right|^p \mathcal{R}_k'(\bar{u}) \\
& \leq \int_{A^{\bar{u}}_k} \left( a + \left| \nabla \bar{u} \right|^2  \right)^{\frac{p-2}{2}} \left| \nabla \bar{u}  \right|^2 \mathcal{R}_k'(\bar{u}) + a^{\frac{p}{2}} \int_{A^{\bar{u}}_k} \mathcal{R}_k'(\bar{u})  \nonumber \\
& \leq \into \left( a + \left| \nabla \bar{u} \right|^2  \right)^{\frac{p-2}{2}} \left| \nabla \bar{u}  \right|^2 \mathcal{R}_k'(\bar{u})   + a^{\frac{p}{2}} \left| A^{\bar{u}}_k \right|  \nonumber \\
& \leq \frac{1}{\nu_1} \into \nabla \Psi_1 ( \nabla \bar{u} ) \cdot \nabla \bar{u} \, \mathcal{R}_k'(\bar{u})  + a^{\frac{p}{2}} \left| A^{\bar{u}}_k \right| \nonumber \\
& = \frac{1}{\nu_1} \into \nabla \Psi_1 ( \nabla \bar{u} ) \cdot \nabla \mathcal{R}_k(\bar{u})  + a^{\frac{p}{2}} \left| A^{\bar{u}}_k \right| \nonumber \\
& = \frac{1}{\nu_1} \into \bar{f}_s \, \mathcal{R}_k( \bar{u}) + a^{\frac{p}{2}} \left| A^{\bar{u}}_k \right| \nonumber \\
& = \frac{1}{\nu_1} \int_{A^{\bar{u}}_k} \bar{f}_s \, \mathcal{R}_k( \bar{u}) + a^{\frac{p}{2}} \left| A^{\bar{u}}_k \right| \nonumber \\
& \leq  \frac{1}{\nu_1} \int_{A^{\bar{u}}_k} \left| \bar{f}_s \right| \left| \mathcal{R}_k(\bar{u})  \right|  + a^{\frac{p}{2}} \left| A^{\bar{u}}_k \right|. \nonumber  
\end{align}

\medskip
\noindent
By Sobolev-Gagliardo-Nirenberg inequality,  H\"older inequality with conjugated exponents  $({p^*})'$ and  $p^*$, and Young inequality with conjugated exponents $p'$ and $p$, we have 

\begin{align*}
\displaystyle 
\int_{A^{\bar{u}}_k} \left| \bar{f}_s \right| \left| \mathcal{R}_k(\bar{u})  \right|  
& \leq \left(  \int_{A^{\bar{u}}_k} \left| \bar{f}_s \right|^{{p^*}'}  \right)^{\frac{1}{{p^*}'}} \left( \int_{A^{\bar{u}}_k} \left| \mathcal{R}_k(\bar{u})  \right|^{p^*}  \right)^{\frac{1}{p^*}}   \\
& = \left(  \int_{A^{\bar{u}}_k} \left| \bar{f}_s \right|^{{p^*}'}  \right)^{\frac{1}{{p^*}'}} \left( \into \left| \mathcal{R}_k(\bar{u})  \right|^{p^*}  \right)^{\frac{1}{p^*}}  \nonumber \\
& \leq  C \left(  \int_{A^{\bar{u}}_k} \left| \bar{f}_s \right|^{{p^*}'}  \right)^{\frac{1}{{p^*}'}} \left( \into \left|  \nabla \mathcal{R}_k(\bar{u})  \right|^{p}  \right)^{\frac{1}{p}} \nonumber \\
& =  C \left( \int_{A^{\bar{u}}_k}\left| \bar{f}_s \right|^{{p^*}'} \right)^{\frac{1}{{p^*}'}} \left( \int_{A^{\bar{u}}_k} \left|  \nabla \bar{u}  \right|^{p}   \right)^{\frac{1}{p}} \nonumber  \\
& \leq \frac{1}{p'} C^{p'} \left( \int_{A^{\bar{u}}_k}\left| \bar{f}_s \right|^{{p^*}'}  \right)^{\frac{p'}{{p^*}'}} + \frac{1}{p}  \int_{A^{\bar{u}}_k} \left|  \nabla \bar{u}  \right|^{p} .  \nonumber
\end{align*}

\medskip
\noindent
Therefore, we have 

\begin{equation}\label{nuovoconto}
\displaystyle \frac{1}{\nu_1} \int_{A^{\bar{u}}_k} \left| \bar{f}_s \right| \left| \mathcal{R}_k(\bar{u})  \right| \, \leq \frac{1}{ \nu_1 \, p'} C^{p'} \left( \int_{A^{\bar{u}}_k}\left| \bar{f}_s \right|^{{p^*}'}  \right)^{\frac{p'}{{p^*}'}} + \frac{1}{ \nu_1 \, p}  \int_{A^{\bar{u}}_k} \left|  \nabla \bar{u}  \right|^{p} . 
\end{equation}

\medskip
\noindent
Taking into account \eqref{calc1} and \eqref{nuovoconto}, we have proved

\begin{equation}\label{calc3}
\displaystyle \left( 1 - \frac{1}{\nu_1 \, p }  \right)\int_{A^{\bar{u}}_k} \left| \nabla \bar{u} \right|^p   \leq \frac{C^{p'}}{\nu_1 \, p'} \left(  \into \left| \bar{f}_s \right|^{{p^*}'} \right)^{\frac{p'}{{p^*}'}} +  a^{\frac{p}{2}} \left| A^{\bar{u}}_k \right|,
\end{equation}

\medskip
\noindent
by which, considering that $\nu_1 > \frac{1}{p}$, we get

\begin{equation}\label{calcnu1}
\displaystyle\int_{A^{\bar{u}}_k} \left| \nabla \bar{u} \right|^p   \leq \frac{p}{p'(\nu_1 \, p - 1)}C^{p'}\left(  \into \left| \bar{f}_s \right|^{{p^*}'}  \right)^{\frac{p'}{{p^*}'}} +  \frac{\nu_1 \, p}{\nu_1 \, p -1}a^{\frac{p}{2}} \left| A^{\bar{u}}_k \right|,
\end{equation}

\medskip
\noindent
Since $m>N$, we deduce that 

$$ \displaystyle m > {p^*}'=\frac{Np}{Np-N+p} \quad \text{and} \quad p' \left( \frac{1}{{p^*}'} - \frac{1}{m}  \right) >1.$$

\medskip
\noindent
Now, by  (\ref{dis1}) and  H\"older inequality with conjugated exponents $\displaystyle \frac{m}{{p^*}'}$ and $\displaystyle \frac{m}{{m-p^*}'}$, we have

\begin{align}\label{calc4}
\displaystyle
\left( \int_{A^{\bar{u}}_k} \left| \bar{f}_s \right|^{{p^*}'}  \right)^{\frac{p'}{{p^*}'}} 
& \leq  
  \left( \int_{A^{\bar{u}}_k} |\bar{f}_s|^m \right)^\frac{p'}{m} \, |A^{\bar{u}}_k|^{p' \left( \frac{1}{{p^*}'} - \frac{1}{m} \right)}  \\
& =  \|H_s(\delta,x,\bar u,\bar v)\|_m^{p'} \, |A^{\bar{u}}_k|^{p' \left( \frac{1}{{p^*}'} - \frac{1}{m} \right)} \nonumber \bigskip\\
& \leq  C \, |A^{\bar{u}}_k|^{p' \left( \frac{1}{{p^*}'} - \frac{1}{m} \right)} .   \nonumber
\end{align}

\medskip
\noindent
Considering \eqref{calcnu1} and  \eqref{calc4}, we have proved that there exists $C>0$ such that

\begin{equation}\label{calc5}
\displaystyle \int_{A^{\bar{u}}_k} \left| \nabla \bar{u} \right|^p   \leq  C \left( \, |A^{\bar{u}}_k|^{p' \left( \frac{1}{{p^*}'} - \frac{1}{m} \right)} +   \,  \left| A^{\bar{u}}_k \right| \right).
\end{equation}

\medskip
\noindent
Now, by using  Lemma \ref{V2023}, there exists $k_0 \geq 1$ such that for any $(\bar u, \bar v) \in D_{R, \delta}(z_0)$ and for any $k$ greater then $k_0$, we get

$$ \displaystyle \left| A^{\bar{u}}_k  \right| < 1. $$

\medskip
\noindent  
In fact,  by using Lemma \ref{V2023} with $u_0 \in W_0^{1,p}(\Omega)$, $r=1$ and $\varepsilon=1$, we deduce that there exists $\sigma_1>0$ such that

\begin{equation*}
\displaystyle
\int_{\{|u(x)|\,\geq \sigma_1\}} \hspace{-1mm}|u(x)|\,  dx< 1
\end{equation*}

\medskip 
\noindent
for any $u\in \overline{B_{R}(u_0)}=\{u\in  W^{1,p}_0(\Omega)\ :\ \|u-u_0\|_{1,p}\leq R \}$.\\
Let us denote with $k_1:= \displaystyle \max \{ 1 , \sigma_1 \}.$\\
Observe that

\begin{equation*}
\displaystyle
1 > \int_{\{|u(x)|\,\geq \sigma_1\}} \hspace{-1mm}|u(x)|\,  dx \geq \int_{A^{u}_{k_1}} \hspace{-1mm}|u(x)|\,  dx > \int_{A^{u}_{k_1}} \hspace{-1mm} k_1\,  dx = k_1 \left|  A^{u}_{k_1} \right|,
\end{equation*}

\medskip
\noindent
by which

$$ \displaystyle \left|  A^{u}_{k} \right| \leq  \left|  A^{u}_{k_1} \right| < \frac{1}{k_1} \leq 1 $$

\medskip
\noindent
for any $u\in \overline{B_{R}(u_0)}$ and for any $k \geq k_1$.\\

\medskip
\noindent
Similarly, we deduce that there exists $k_2 \geq 1$ such that 

$$ \displaystyle \left|  A^{v}_{k} \right| < 1 $$

\medskip
\noindent
for any $v\in \overline{B_{R}(v_0)}=\{v\in  W^{1,q}_0(\Omega)\ :\ \|v-v_0\|_{1,q}\leq R \}$ and for any $k \geq k_2$.\\

\medskip
\noindent
Denoting with $\displaystyle k_0:= \max\{k_1,k_2\}$, for any $(\bar u, \bar v) \in D_{R, \delta}(u_0,v_0)$ and for any $k \geq k_0$, we get

$$ \displaystyle \left| A^{\bar{u}}_k  \right| < 1. $$

\medskip
\noindent
Considering also $ \displaystyle p' \left( \frac{1}{{p^*}'} - \frac{1}{m}  \right) >1,$ by \eqref{calc5} we deduce that for any $k \geq k_0$ we have

\begin{equation}\label{calc6}
\displaystyle \int_{A^{\bar{u}}_k} \left| \nabla \bar{u} \right|^p   \leq 2 \bar{C}  \,  \left| A^{\bar{u}}_k \right|.
\end{equation}

\medskip
\noindent
Therefore, by applying Theorem \ref{Mammoliti} with $\varepsilon= 1 - \frac{p}{p^*}$ and $\theta=0$, we deduce that $\bar u\in L^\infty(\Omega)$. Similarly, we deduce that  $\bar v\in L^\infty(\Omega)$.\\
Finally, we have proved that there exist a suitable small $R>0$ and a constant $C>0$ depending on $C_0,\left| \Omega \right|, p,q, N, a, \nu_1, b$ and $\nu_2$  but not on $\delta$, such that for any $(\bar u, \bar v) \in D_{R,\delta}(z_0)$,
$\bar u$ and $\bar v$ are in $L^\infty(\Omega)$ and
\[\|\bar u\|_{\infty}, \ \|\bar v\|_{\infty}\leq C.
\]

\qed

\begin{remark}
We point out that if $ p \geq 2$ the constant $C$ does not depend on $a$. In fact, in the proof of Theorem \ref{key}  we can use instead of \eqref{disr} the inequality 

$$ \displaystyle 	s^p \leq \left(a +s^2\right)^{\frac{p-2}{2}}s^2, $$

\medskip
\noindent 
that holds for any $p \geq 2$ and for any $s \geq 0$.\\
Similarly, if $q \geq 2$  the constant $C$ does not depend on $b$.
\end{remark}

\bigskip


\section*{Appendix B}
\addcontentsline{toc}{section}{Appendix B}
\label{app:B}

\begin{lemma}\label{contogradconvess} Let  $\Psi: \mathbb{R}^N \to \mathbb{R}$ be a convex $C^1$-function. Assume that there exist $c>0, d>0$ and $r>1$ such that

\begin{equation}\label{proprconvPsi}
\displaystyle \left| \Psi(\xi) \right| \leq c\left|\xi \right|^r + d  \qquad \forall \, \xi \in \mathbb{R}^N.
\end{equation}

\medskip
\noindent
Then there exists  a constant $\tilde{c}$ such that

\begin{equation}\label{proprconvgradPsi}
\displaystyle \left| \nabla \Psi(\xi) \right| \leq \tilde{c} \left( \left|\xi \right|^{r-1}  + 1 \right) \qquad \forall \, \xi \in \mathbb{R}^N.
\end{equation}
\end{lemma}

\begin{proof} By convexity of $\Psi$, we know that

\begin{equation*}
 \displaystyle \Psi(\xi + h) - \Psi(\xi) \geq \nabla \Psi(\xi) \cdot h \qquad \forall \, \xi, h \in \mathbb{R}^N.
\end{equation*}

\medskip
\noindent
Let $\xi \in \mathbb{R}^N$ with $\xi \neq 0$. Taking into account previous relation and \eqref{proprconvPsi} , we get

\begin{align*}
\displaystyle 
\left| \nabla \Psi(\xi) \right|
= \sup_{ |\eta| \leq 1 } \; \nabla \Psi (\xi) \cdot \eta  
& = \sup_{ |\eta| \leq 1 } \; \frac{\nabla \Psi (\xi) \cdot \left|\xi \right| \, \eta}{\left| \xi \right|}\\
& \leq \sup_{ |\eta| \leq 1 } \frac{\Psi(\xi + \left|\xi\right|  \eta) - \Psi(\xi)}{\left| \xi \right|} 
 \leq \sup_{ |\eta| \leq 1 } \frac{\left| \Psi(\xi + \left|\xi\right|  \eta) - \Psi(\xi) \right|}{\left| \xi \right|}
 \leq \bar{c} \left| \xi \right|^{r-1} + \frac{\tilde{d}}{\left| \xi \right|},
\end{align*}

\medskip
\noindent
where $\bar{c}$ and $\tilde{d}$ are suitable constants.

\medskip
\noindent
Observe now that 

$$ \displaystyle 
\lim_{ \left| \xi \right| \to + \infty} \frac{\bar{c} \left| \xi \right|^{r-1} + \frac{\tilde{d}}{\left| \xi \right|}}{\bar{c} \left| \xi \right|^{r-1} }=1,$$

\medskip
\noindent
hence there exists a constant $R>0$ such that

$$ \displaystyle \bar{c} \left| \xi \right|^{r-1} + \frac{\tilde{d}}{\left| \xi \right|} \leq 2 \bar{c} \left| \xi \right|^{r-1} \qquad \text{ if } \left| \xi \right| > R,$$

\medskip
\noindent
by which

$$ \displaystyle  \left| \nabla \Psi(\xi) \right|  \leq 2 \bar{c} \left| \xi \right|^{r-1} \qquad \text{ if } \left| \xi \right| > R.$$

\medskip
\noindent
Now, since $\Psi$ is a function of class $C^1$ in $\mathbb{R}^N$, we deduce that there exists $M>0$ such that

$$ \displaystyle  \left| \nabla \Psi(\xi) \right|   \leq M \qquad  \text{ if } \left| \xi \right| \leq R.$$

\medskip
\noindent
Finally, denoting with $\tilde{c}:= \max \{ 2 \bar{c}, M   \},$ we get the thesis. 

\end{proof}

\begin{proposition}\label{funzionalediclasseC1}
Let $\Psi: \mathbb{R}^N \to \mathbb{R}$ be a convex $C^1$-function. Let $r>1$ and let us assume that there exist $c>0$ and $d>0$  such that

\begin{equation}\label{proprconvPsiprop}
\displaystyle \left| \Psi(\xi) \right| \leq c\left|\xi \right|^r + d  \qquad \forall \, \xi \in \mathbb{R}^N.
\end{equation}

\medskip
\noindent
Let $\Omega \subset \mathbb{R}^N$ be a bounded domain and let us define $f: W_0^{1,r}(\Omega) \to \mathbb{R}$ as 

\begin{equation}\label{funzionalef}
 \displaystyle f(w):= \into \Psi \left( \nabla w(x) \right ) \, dx.
\end{equation}

\medskip
\noindent
Then $f$ is a $C^1$-functional on $W_0^{1,r}(\Omega)$ and for any $w_0, w \in W_0^{1,r}(\Omega)$ we have

\begin{equation}\label{Gateauxdiff}
 \displaystyle \langle f'(w_0), w \rangle =   \into \nabla \Psi \left( \nabla w_0(x) \right) \cdot \nabla w(x) \, dx.
\end{equation}

\end{proposition}

\begin{proof}
The functional $f$ is well defined, in fact for any $ w \in W_0^{1,r}(\Omega)$ inequality \eqref{proprconvPsiprop} implies

$$ \displaystyle \left| f(w) \right| \leq \into \left| \Psi \left( \nabla w(x) \right ) \right|\, dx \leq c \lVert \nabla w \rVert_r^r + d \left| \Omega \right| < + \infty.$$

\medskip
\noindent
Let us now prove that $f$ is Gateaux differentiable in any $w_0 \in W_0^{1,r}(\Omega)$, where the Gateaux differential in $w_0$ is the map $f'(w_0)$ whose action is defined in \eqref{Gateauxdiff}.\\
Hence, we need to show that for any $w_0, w \in W_0^{1,r}(\Omega)$ we have

\begin{equation}\label{Gateauxdaprovare}
\displaystyle \lim_{t \to 0} \frac{ f(w_0 + tw) - f(w_0)}{t} = \into \nabla \Psi \left( \nabla w_0(x) \right) \cdot \nabla w(x) \, dx.
\end{equation}

\medskip
\noindent
Observe that 

$$ \displaystyle \frac{ f(w_0 + tw) - f(w_0)}{t} = \into \frac{ \Psi ( \nabla w_0(x) + t \nabla w(x)) - \Psi ( \nabla w_0(x))}{t} \, dx.$$

\medskip
\noindent
Let $\xi, \eta \in \mathbb{R}^N$ and let us define the function $h:[-1,1] \to \mathbb{R}$ as $ \displaystyle h(t):= \Psi( \xi + t \, \eta).$\\
The function $h$ is of class $C^1$ with $ \displaystyle h'(t)= \nabla \Psi( \xi + t\, \eta) \cdot \eta $  and $ h'(0)= \nabla \Psi( \xi) \cdot \eta,$ whence

$$ \displaystyle \lim_{t \to 0} \frac{ \Psi ( \nabla w_0(x) + t \nabla w(x)) - \Psi ( \nabla w_0(x))}{t} = \nabla \Psi( \nabla w_0(x) ) \cdot \nabla w(x) \qquad \text{ a.e. in } \Omega.$$

\medskip
\noindent
By Lemma \ref{contogradconvess}, we know that

\begin{equation}\label{proprconvgradPsidimostrazione}
\displaystyle \left| \nabla \Psi(\xi) \right| \leq \tilde{c} \left( \left|\xi \right|^{r-1}  + 1 \right) \qquad \forall \, \xi \in \mathbb{R}^N.
\end{equation}

\medskip
\noindent
By Lagrange Theorem applied to the function $h$, there exists $ \theta \in \mathbb{R}$ with $\left| \theta \right| \leq |t| \leq 1$ such that 

\begin{align*}
\displaystyle \left| \frac{ \Psi ( \nabla w_0(x) + t \nabla u(x)) - \Psi ( \nabla w_0(x))}{t}    \right| 
& = \left| \nabla \Psi( \nabla w_0(x) + \theta \, \nabla w(x)) \cdot \nabla w(x)  \right| \\
& \leq \left| \nabla \Psi( \nabla w_0(x) + \theta \, \nabla w(x)) \right| \left| \nabla w(x)  \right|.
\end{align*} 

\medskip
\noindent
By \eqref{proprconvgradPsidimostrazione}, we get

\begin{align*}
\displaystyle 
\left| \nabla \Psi( \nabla w_0(x) + \theta \, \nabla w(x)) \right| \left| \nabla w(x)  \right|
& \leq  \tilde{c} \left( \left|\nabla w_0(x) + \theta \, \nabla w(x)\right|^{r-1}  + 1 \right)\left| \nabla w(x)  \right| \\
& \leq \tilde{c} \left( [ \, \left|\nabla w_0(x) \right| + \left| \nabla w(x)\right| \, ]^{r-1}  + 1 \right)\left| \nabla w(x)  \right| \\
& \leq C  \left|\nabla w_0(x) \right|^{r-1} \left| \nabla w(x)  \right| + C\left| \nabla w(x)\right|^{r} + \tilde{c} \left| \nabla w(x)  \right|,
\end{align*}

\medskip
\noindent
with

$$ \displaystyle C  \left|\nabla w_0(x) \right|^{r-1} \left| \nabla w(x)  \right| + C\left| \nabla w(x)\right|^{r} + \tilde{c} \left| \nabla w(x)  \right| \in L^1(\Omega).$$

\medskip
\noindent
Hence, by applying dominated convergence theorem we get \eqref{Gateauxdaprovare}.

\medskip
\noindent
Let us prove that the Gateaux derivative $f': W_0^{1,r}(\Omega) \to W^{-1,r'}(\Omega)$ is continuous.\\
Let $\{w_k\}_k \subset W_0^{1,r}(\Omega)$  be a sequence such that $w_k \to w$ in $W_0^{1,r}(\Omega)$. Up to subsequence, we have:

\begin{itemize}
\item[$i)$] $\nabla w_k \to \nabla w$ in $(L^r(\Omega))^N$ as $k \to \infty$;
\item[$ii)$] $\nabla w_k(x) \to \nabla w(x)$ a.e. $x \in \Omega$ as $k \to \infty$;
\item[$iii)$] there exists $g \in L^1(\Omega)$ such that $\left| \nabla w_k(x) \right|^r \leq g(x)$ a.e. $x \in \Omega$, and for all $k \in \mathbb{N}$.
\end{itemize}

\medskip
\noindent
We have

$$ \displaystyle \langle f'(w) - f'(w_k), \bar{w}   \rangle= \into  \left( \, \nabla \Psi ( \nabla w(x)) - \nabla \Psi( \nabla w_k(x)) \, \right) \cdot \nabla \bar{w}(x) \, dx $$

\medskip
\noindent
and

\begin{align*}
\displaystyle
\left| \langle f'(w) - f'(w_k), \bar{w}   \rangle \right| 
& =    \left| \into  \left( \, \nabla \Psi ( \nabla w(x)) - \nabla \Psi( \nabla w_k(x)) \, \right) \cdot \nabla \bar{w}(x) \, dx   \right| \\
& \leq \into \left| \left( \, \nabla \Psi ( \nabla w(x)) - \nabla \Psi( \nabla w_k(x)) \, \right) \cdot \nabla \bar{w}(x) \right| \, dx \\
& \leq \into \left|  \, \nabla \Psi ( \nabla w(x)) - \nabla \Psi( \nabla w_k(x)) \,  \right| \, \left| \nabla \bar{w}(x) \right| \, dx \\
& \leq \left( \into \left|  \, \nabla \Psi ( \nabla w(x)) - \nabla \Psi( \nabla w_k(x)) \,  \right|^{\frac{r}{r-1}} \, dx  \right)^{\frac{r-1}{r}}   \lVert \nabla \bar{w} \rVert_r 
\end{align*}

\medskip
\noindent
Therefore, we deduce

$$ \displaystyle \lVert f'(w) - f'(w_k) \rVert_{W_0^{-1,r}(\Omega)} \leq  \left( \into \left|  \, \nabla \Psi ( \nabla w(x)) - \nabla \Psi( \nabla w_k(x)) \,  \right|^{\frac{r}{r-1}} \, dx  \right)^{\frac{r-1}{r}} .$$

\medskip
\noindent
Taking into account $ii)$ and the fact that $\Psi$ is a $C^1$ function, we easily deduce that 

$$ \displaystyle \nabla \Psi( \nabla w_k(x) ) \to \nabla \Psi ( \nabla w(x)) \qquad \text{ a.e. in } \Omega.$$

\medskip
\noindent
Considering \eqref{proprconvgradPsidimostrazione} and $iii)$, we have

\begin{align*}
\displaystyle
\left|  \, \nabla \Psi ( \nabla w(x)) - \nabla \Psi( \nabla w_k(x)) \,  \right|^{\frac{r}{r-1}} 
& \leq  \left( \, \left| \nabla \Psi ( \nabla w(x)) \right| + \left| \nabla \Psi( \nabla w_k(x))  \right|    \, \right)^{\frac{r}{r-1}} \\
& \leq C \left(   \, \left| \nabla \Psi ( \nabla w(x)) \right|^{\frac{r}{r-1}} +  \left| \nabla \Psi ( \nabla w_k(x)) \right|^{\frac{r}{r-1}} \,    \right) \\
& \leq  \tilde{C} \left(   \, \left( \, \left| \nabla w(x) \right|^{r-1} + 1 \right)^{\frac{r}{r-1}} +  \left(\, \left| \nabla w_k(x) \right|^{r-1} + 1 \right)^{\frac{r}{r-1}} \,    \right) \\
& \leq \tilde{\tilde{C}} \left(  \, \left| \nabla w(x) \right|^r +   \left| \nabla w_k(x) \right|^r + 1 \,    \right) \\
& \leq 2\tilde{\tilde{C}}   g(x)  + \tilde{\tilde{C}} 
\end{align*}

\medskip
\noindent
a.e. $x \in \Omega$ and for all $k \in \mathbb{N}$, where $2\tilde{\tilde{C}}   g(x)  + \tilde{\tilde{C}} \in L^1(\Omega)$.\\
Therefore, we can apply dominated convergence theorem to say that 

$$ \displaystyle \left( \into \left|  \, \nabla \Psi ( \nabla w(x)) - \nabla \Psi( \nabla w_k(x)) \,  \right|^{\frac{r}{r-1}} \, dx  \right)^{\frac{r-1}{r}} \to 0 \qquad \text{ as } k \to \infty,$$

\medskip
\noindent
whence

$$ \displaystyle \lVert f'(w) - f'(w_k) \rVert  \to 0 \text{ as } k \to \infty.$$

\medskip
\noindent
This holds for a subsequence of the original sequence $\{  w_k \}_k$. By using the Lemma of sub-sub-sequence we obtain that $f'$ is continuous, and so $f$ is of class $C^1$.

\end{proof}

\bigskip

\chapter{Amann-Zehnder type results for Resonance Quasilinear Elliptic Systems}\label{CAPAMANNZEHNDER}

\bigskip
\noindent
The interaction between a nonlinearity and the spectrum of $-\Delta$ at $0$ and infinity has been employed in the celebrated paper \cite{amann_zehnder} of Amann and Zehnder, under the assumption of nonresonance at infinity. The quasilinear counterpart of the Amann-Zehnder existence result has been obtained by Cingolani, Degiovanni and Vannella in \cite{CDV} (see also \cite{CD,DrabekRobinson}) for a class of quasilinear elliptic equations, introducing a Saddle Theorem where linear subspaces are substituited by symmetric cones (see Theorem 7.1 in \cite{CDV}).  

In this Chapter, we present Amann-Zehnder type results for a class of quasilinear elliptic systems whose principal parts consist of $(p,q)$-Laplace or $(p,q)$-area type operators. The nonlinearities are assumed to exhibit $(p,q)$-linear growth at infinity, covering both the asymptotic resonant and nonresonant cases. The space setting is given by the Banach product space $X := W_0^{1,p}(\Omega) \times W_0^{1,q}(\Omega)$, hence we exploit recent Morse-theoretic techniques in Banach spaces (see \cite{CV,CV2,CV4,U}). First,  we establish a Saddle Point Theorem specifically tailored for $(p,q)$-homogeneous sets, since we lose the structure of symmetric cones. Second, we establish a relationship between the local behaviour of the functional near its critical points and differential notions like the Morse index. Moreover, since the nonlocal coupling induces a lack of regularity in the functional, we overcome this difficulty by using a penalization method and uniform regularity estimates shown in Chapter \ref{SezReg}.

\medskip
The present Chapter is based on the paper \cite{BCV}.

\bigskip

\section{Amann-Zehnder type results}\label{SezioneIntroduttivaNontrivialSolution}

\medskip
\noindent
In 1980, Amann and Zehnder \cite{amann_zehnder} studied the asymptotically linear elliptic problem

\begin{equation}
	\label{semiIntro} 
	\left\{
	\begin{array}{ll}
		- \, \Delta  u = g(u)
		 \qquad&\text{in $\Omega$}, \\
		\noalign{\medskip}
	\quad \quad 	u=0
		\qquad&\text{on $\partial\Omega$},
	\end{array}
	\right.
\end{equation}

\medskip
\noindent
where  $\Omega$ is a bounded domain in $\mathbb{R}^N$ with smooth boundary, $g: \mathbb{R} \to \mathbb{R}$ is a $C^1$-function such that $g(0)=0$ and there exists $\bar \lambda \in \re$ such that

$$ \displaystyle \lim_{|s|\to\infty} g'(s) = \bar \lambda. $$

\medskip
\noindent
They proved that problem \eqref{semiIntro} admits a nontrivial solution $u$, supposing that $\bar{\lambda}$ is not an eigenvalue of $- \Delta$, the so-called nonresonance condition at infinity, and that there exists
some eigenvalue of $- \Delta$ between $\bar{\lambda}$ and $g'(0)$. Using Morse theory for manifolds with boundary, the same result has been found by Chang in \cite{chang1981-cpam}.
Successively,  in \cite{lazer_solimini1988} Lazer and Solimini recognized that such nontrivial solution can be detected  combining mini-max characterization of the critical point and Morse index estimates (see also \cite{MasielloPisani}). More precisely, the basic idea in \cite{lazer_solimini1988} is to recognize that the energy functional associated with the asymptotically linear problem \eqref{semiIntro} has a saddle geometry, which implies that a suitable Poincaré polynomial is not trivial, and also to show that a certain critical group at zero is trivial, to ensure the existence of a solution $u \neq 0$ of \eqref{semiIntro}.

Recently the quasilinear counterpart of the Amann-Zehnder existence result has been obtained in \cite{CDV} (see also \cite{CD,DrabekRobinson}) 
for a class of quasilinear elliptic equations, involving either the $p$-Laplace operator or the operator related to the $p$-area functional, and a nonlinearity  with $p$-linear growth at infinity.\\ 
More precisely, in \cite{CDV} Cingolani, Degiovanni and Vannella  studied the following problem

\begin{equation}\label{problemaCingolaniDegiovanniVannella}
	\left\{
	\begin{array}{ll}
		- \text{div} \left[ \left( k^2 + | \nabla u |^2   \right)^{\frac{p-2}{2}} \nabla u  \right] = g(u),
		 \qquad&\text{in $\Omega$}, \medskip \\
	 u=0,
		\qquad&\text{on $\partial\Omega$},
	\end{array}
	\right. 
\end{equation}

\medskip
\noindent
where $\Omega$ is a bounded open subset of $\mathbb{R}^N$, $N \geq 1$, with $\partial \Omega$ of class $C^{1,\alpha}$ for some $\alpha \in ]0,1]$, while $k \geq 0,$ $p>1$  are real numbers and $g: \mathbb{R} \to \mathbb{R}$ is a $C^1$-function such that $g(0)=0$ and there exists $\bar \lambda \in \re$ such that

$$ \displaystyle \lim_{|s|\to\infty} \frac{g(s)}{|s|^{p-2}s} = \bar \lambda. $$

\medskip
\noindent
Problem \eqref{problemaCingolaniDegiovanniVannella} gives rise to new difficulties when $p \neq 2$. In particular, as regards the eigenvalue problem $-\Delta_p u = \lambda |u|^{p-2}u$, the spectral properties  are not yet well understood.  Denoted $\sigma(-\Delta_p)$ the set of such eigenvalues $\lambda$,
one can define in at least three different ways a diverging sequence $(\lambda_m)$ of eigenvalues of $-\Delta_p$ (see~\cite{CD, DrabekRobinson}), but it is not known if they agree for $m \geq 3$ and if the whole set $\sigma(-\Delta_p)$ is covered. Hence it is more difficult to recognize a saddle structure with a related information on a suitable Poincaré polynomial. On the other hand, for functionals defined on Banach spaces serious difficulties arise in extending Morse theory, and so in providing an estimate of the critical groups at zero by some Hessian type notion. Therefore problem \eqref{problemaCingolaniDegiovanniVannella} was addressed in \cite{CDV} exploiting new tecniques of Morse theory in Banach spaces  (see \cite{CV,CV2,CV4,U}) and introducing a Saddle Theorem where linear subspaces are substituited by symmetric cones (see Theorem 7.1 in \cite{CDV}).  
 
In this Chapter we present Amann-Zehnder type results for a class of quasilinear elliptic systems whose principal parts consist of $(p,q)$-Laplace or $(p,q)$-area type operators. 
Precisely we will seek for nontrivial solutions for the following autonomous quasilinear system

\begin{equation}\label{pqNontrivialSolutions}
	\begin{cases}
		\begin{array}{ll}
			-\text{\rm div} \left((a+|\nabla u|^2)^{\frac{p-2}{2}}\nabla u \right) = G_s(u,v) & x\in\Omega,
			\medskip \\
			-\text{\rm div} \left((a+|\nabla v|^{2})^{\frac{q-2}{2}}\nabla v \right)= G_t(u,v) & x\in \Omega, \medskip	\\
			u=v=0  & x\in \partial\Omega,
		\end{array}
	\end{cases}
\end{equation}

\medskip
\noindent
where $2\leq p<N,\ 2\leq q<N$, $a\geq 0$,  and $G \in C^{1}(\re^2, \re)$ satisfies the following conditions:  
\medskip
\noindent
\begin{itemize}
	\item[${(a_1)}$] 
	$ \nabla G(0,0)=(0,0)$ and there exists $\bar\lambda\in \re$ such that
	\[
	\lim_{|(s,t)|\to  \infty} \frac{G_s(s,t)-\bar{\lambda} F_s(s,t) 
	}{|s|^{p-1}+|t|^{q\frac{p-1}{p}} } =0 \qquad \text{and} \qquad 
	\lim_{|(s,t)|\to  \infty} \frac{G_t(s,t)-\bar \lambda F_t(s,t)}{|s|^{p\frac{q-1}{q}}+|t|^{q-1} }
	=0,\,
	\]	
	where $F \in C^1(\re^2, \re)$ is defined by  
	\begin{align*}
		F(s,t)= 	  \frac{1}{p}  | s|^p + \frac{1}{q}	| t|^q + \frac{1}{(\alpha +1)(\beta +1)}
		|s|^\alpha |t|^\beta s t
	\end{align*}
	with $\alpha, \beta > 0$ such that $(\alpha +1)/ p + (\beta +1)/q =1$;

	\smallskip
	
	\item[$(a_2)$] there is a suitable neighborhood $U$ of $(0,0)$ such that $G\in C^2(U,\re)$.
	\end{itemize}

\medskip
\noindent
Weak solutions of problem \eqref{pqNontrivialSolutions} correspond to critical points of the Euler functional $J_a:X \to \mathbb{R}$ defined for any $z=(u,v) \in X$ by setting

\begin{align}\label{functionalNontrivialsolutions}
 J_a(z) =  \frac{1}{p} \into  \left(a+|\nabla u|^2\right)^{\frac p 2} \ dx + \frac{1}{q}\into \left(a+|\nabla v|^2\right)^{\frac q 2} \ dx - \into G(u,v) \ dx.
\end{align}

\medskip
\noindent
Observe now that for any $(s,t) \in \mathbb{R} \times \mathbb{R}$ we have

$$ \displaystyle F_s(s,t)= |s|^{p-2}s + \frac{1}{\beta+1} |s|^{\alpha} |t|^{\beta} t \qquad \text{and} \qquad  F_t(s,t)= |t|^{q-2}t + \frac{1}{\alpha+1} |s|^{\alpha} |t|^{\beta} s.$$

\medskip
\noindent
By $(\alpha +1)/ p + (\beta +1)/q =1,$ we deduce $\frac{\alpha}{p-1} + \frac{p(\beta+1)}{q(p-1)}=1$ and      $\frac{\beta}{q-1} + \frac{q(\alpha+1)}{p(q-1)}=1$, and so by Young's inequality we infer that there exists a constant $\bar{C}$ such that for any $(s,t) \in \mathbb{R} \times \mathbb{R}$ we have

$$ \displaystyle \left| F_s(s,t) \right| \leq \bar{C} \left( |s|^{p-1} + |t|^{q\frac{p-1}{p}}   \right) \qquad \text{and} \qquad \left| F_t(s,t) \right| \leq \bar{C} \left( |t|^{q-1} + |s|^{p\frac{q-1}{q}}   \right). $$

\medskip
\noindent
Considering now the limits in assumption $(a_1)$, we infer that there is $C >0$ such that for any $(s,t) \in \mathbb{R} \times \mathbb{R}$ we have

\begin{equation}\label{condizionidicrescitasuFsFt}
\left| G_s(s,t)\right| \leq C \left(  1 +  |s|^{p-1}+ |t|^{q\frac{p-1}{p}}\right)   \qquad \text{and} \qquad \left| G_t(s,t) \right| \leq C \left( 1+ |t|^{q-1} + |s|^{p\frac{q-1}{q}}   \right).
\end{equation}

\medskip
\noindent
Therefore the functional $J_a$ is of class $C^1$ on $X$ and for any $z_0=(u_0,v_0),z=(u, v) \in X$ it results
	
\begin{align*}
		\langle J'_a(z_0), z \rangle = & \displaystyle  \into (a+|\nabla
		u_0|^2)^{\frac {p-2}{2}}\nabla u_0 \cdot \nabla u \ dx +\into
		(a+|\nabla v_0|^2)^{\frac {q-2}{2}}\nabla v_0 \cdot \nabla v \ dx  \\
		&-\displaystyle\into (  G_s(u_0,v_0) u +  G_t(u_0,v_0) v)
		dx.
\end{align*}

\medskip
\noindent
Existence, nonexistence and regularity results for  quasilinear elliptic systems  are obtained by various authors, see for instance \cite{boccardodefiguerido, BM, candela, morais-souto, DT, dingxiao, perera, velin-dethelin}. In this Chapter we are interested to find nontrivial solution of \eqref{pqNontrivialSolutions}  in both the asymptotic resonant and nonresonant cases, that is we consider the following nonlinear eigenvalue problem

\begin{equation}\label{eigen}
	\begin{cases}
		\begin{array}{ll}
			- \Delta_p u  = \lambda |u|^{p-2} u +
			\frac{\lambda}{\beta +1}|u|^{\alpha} |v|^{\beta} v & x\in\Omega,
			 \medskip \\
			- \Delta_q v = \lambda |v|^{q-2} v +
			\frac{\lambda}{\alpha +1}|u|^{\alpha} |v|^{\beta} u, & x\in\Omega, \medskip \\
			u=v=0,  & x\in \partial\Omega,
		\end{array}
	\end{cases}
\end{equation}

\medskip
\noindent
where $\ 2 \leq p <N$, $\ 2 \leq q < N$,  
 $\alpha > 0$ and $ \beta> 0$ are real numbers satisfying
$(\alpha +1)/ p + (\beta +1)/q =1$.

\medskip
\noindent
The real number $\lambda$ is called an eigenvalue of $(\ref{eigen})$ if there exists a nontrivial solution $(u,v)$ of
$(\ref{eigen})$. 

\medskip

In \cite{Tang},  the existence of an unbounded sequence of minimax eigenvalues was proved by using
$\ze_2$-cohomological index of Fadell and Rabinowitz \cite{fadell_rabinowitz1977,fadell_rabinowitz1978}.

\medskip
\noindent
Precisely, let $\Phi : X \to \re$, $\Psi:X \to \re$ be the following functionals 

\begin{align*}
	\Phi(u,v) =&  \frac{1}{p} \int_{\Omega}
	|\nabla u(x)|^p \ dx + \frac{1}{q}\int_{\Omega}
	|\nabla v(x)|^q \ dx,
\end{align*}
\begin{align*}
	\Psi(u,v) =  \int_{\Omega} F(u,v) \ dx. 
\end{align*}

\medskip
\noindent
By Young's inequality we infer that

\begin{equation}\label{disYo}
	\frac 1 p \frac{\beta}{\beta+1}|s|^p+\frac 1 q\frac{\alpha}{\alpha+1} |t|^q\leq F(s,t)\leq \frac 1 p \frac{\beta+2}{\beta+1}|s|^p+\frac 1 q \frac{\alpha+2}{\alpha+1}|t|^q.
\end{equation}

\medskip
\noindent
Note that $\Phi$ and $\Psi$ are $(p,q)$-homogeneous, i.e.

\begin{equation}\label{pqh}
	\Phi(t^{\frac 1 p}u,t^{\frac 1 q}v)=t\Phi(u,v), \quad \Psi(t^{\frac 1 p}u,t^{\frac 1 q}v)=t\Psi(u,v) 
\end{equation}
for all $t>0$ and $(u,v)\in X$.

\medskip
\noindent
Moreover $\lambda$ is an eigenvalue if and only if there is $(u,v)\in X\setminus \{(0,0)\}$
such that

\[\Phi'(u,v)=\lambda \Psi'(u,v). \]

\medskip
\noindent
In \cite{Tang} it is proved that
$(\ref{eigen})$ has a nondecreasing and unbounded sequence of eigenvalues with the variational characterization

\begin{equation}\label{lambdak}
	\lambda_k = \inf_{A \in \Sigma_k} \sup_{(u,v) \in A} \Phi(u,v),
\end{equation}

\medskip
\noindent
where 
$\Sigma$ is the $C^1$ manifold

$$
\Sigma := \{(u,v) \in X  \ | \Psi(u,v)=1 \}$$

\medskip
\noindent
and

\medskip
\noindent
$$
\Sigma_k := \{ A \subset \Sigma  \ | \, A \ is\  symmetric,\ compact\ and\  i(A) \geq k \}.
$$

\medskip

\medskip
\noindent
Here $i(A)$ denotes the $\ze_2$-cohomological index of $A$, introduced by Fadell-Rabinowitz \cite{fadell_rabinowitz1977,fadell_rabinowitz1978}. 
For a matter of convenience, we also put $\la_0=-\infty$.

In \cite{SZ} it is shown that the first eigenvalue $\lambda_1$ is simple and isolated with a first strictly positive eigenfunction $\psi =(\psi_1, \psi_2)$, i.e.  $\psi_1>0, \psi_2>0$ in $\Omega$ (see also \cite{boccardodefiguerido, denapolimarani,DT, drabek}).
However it is not clear if the set of the eigenvalues described by $(\ref{lambdak})$ contains all the eigenvalues of problem \eqref{eigen}.

Furthermore we remark that in the context of quasilinear systems, with  $p \neq q$, the level set
$E_{\lambda_k} = \{(u,v) \in X \ | \Phi(u,v)= \lambda_k \Psi(u,v) \}$
is not a linear space or a cone, but by $(\ref{pqh})$ we infer that $E_{\lambda_k}$ is a $(p,q)$-homogeneous set, i.e.  
$(t^{\frac 1 p}u,t^{\frac 1 q}v) \in E_{\lambda_k}$ for any $t \geq 0$ and 
$(u,v) \in E_{\lambda_k}$.
This makes  delicate to 
recognize a saddle type geometry for the energy functional associated to $\eqref{pqNontrivialSolutions}$.

In Section \ref{SezioneLinkingComologico} we derive a quite general abstract Linking Theorem (see Theorem \ref{th key}), that extends Theorem 7.1 in \cite{CDV} which deals with symmetric cones.  This abstract result can be successfully applied in the framework of systems and allows to detect the existence of a minimax critical point for $J_a$ with a suitable nontrivial critical group.
Finally, in order to prove that the solution we found is not trivial, we need to compute the critical groups at the origin in terms of differential notions. However the applicability of Morse arguments in Banach spaces present severe difficulties since classical Morse Lemma and generalized Morse lemma  \cite{gromoll_meyer1969-t} are so far known, also due to the lack of Fredholm properties of the second derivative of the functionals. Moreover the energy functional $J_a$ is not $C^2$, so that we cannot define the classical notion of Morse index. 
In order to overcome this lack of regularity, the first idea is to consider the quadratic form $Q : X\to \re$ defined as

\begin{align}\label{m0}
	Q(z)= Q(u,v)	 =   \, a^{\frac{p-2}{2}}\, \bigl(\int_{\Omega} |\nabla u|^2 \, dx\bigr)  + a^{\frac{q-2}{2}}\bigl(\int_{\Omega} |\nabla v|^2 \, dx \bigr)\nonumber
		\\ 
		- \int_\Omega \bigl(  G_{ss}(0,0) u^2  + 
		G_{tt}(0,0) v^2 + 2  G_{st}(0,0) u v  \bigr) \ dx \\
		 =
		a^{\frac{p-2}{2}}\, \|u \|_{1,2}^2+ a^{\frac{q-2}{2}}\, \|v \|_{1,2}^2 - \int_\Omega H_G(0,0) [z]^2 \, dx,\nonumber
\end{align}

\medskip
\noindent
where $H_G(0,0)$ is the Hessian matrix of $G$ at $(0,0)$. Then we consider $m_0$ and $m_0^*$ defined as 
the supremum of the dimensions of the subspaces of $X$ where~$Q$ is negative definite and negative semidefinite respectively.  These objects will play the role of Morse index and large Morse index.

The second problem becomes to relate such differential notions $m_0$ and $m_0^*$ to the behaviour of the functional $J_a$ near the origin.
This delicate issue will be obtained by means of a penalized functional of $J_a$, which is  $C^2$ just in a ball centered at the origin. 
We stress that even if $p \geq 2$ and $q \geq 2$,  the functional $F$ is just $C^1$ on $X$, differently from the case of a single equation. Roughly speaking, we say that the eigenvalue problem  $(\ref{eigen})$ is intrinsically less regular than the eigenvalue problem for the single $p$-Laplace equation, due to the presence of the right-hand side coupled terms.

\medskip
\noindent
The main goals of this Chapter are the following:

\begin{theorem}\label{key1} Suppose that $2 \leq p <N, 2 \leq q <N$, $a \geq 0$ and  $(a_1)-(a_2)$ hold. Assume that $\bar \lambda$ is not an eigenvalue for system $(\ref{eigen})$ and denote by $m_{\infty}$ the integer such that
	\[	\lambda_{m_{\infty}} < \bar\lambda < \lambda_{m_{\infty}+1}.
	\]
	If
	\[m_{\infty} \not\in [m_0,m_0^*],
	\]
	then there exists a nontrivial weak solution $z=(u,v)$ of \eqref{pqNontrivialSolutions}.
\end{theorem}

We remark that in Theorem \ref{key1} the parameter $\bar \lambda$ is not an eigenvalue of  \eqref{eigen} and  $\lambda_m < \bar \lambda < \lambda_{m+1}$ for some $m \in \N$. 
This condition seems to be possible, for instance, if $m=1$ taking into account that  $\lambda_1$ is isolated and it is the only eigenvalue of \eqref{eigen} to which corresponds a componentwise positive eigenfunction  \cite{SZ}.
However such nonresonant restriction would be quite severe and  it becomes relevant to face with  systems at resonance.

In the following theorems the value $\bar \lambda$ is allowed to be an eigenvalue of \eqref{eigen},
under an additional condition of $G$ at infinity.

\smallskip
\begin{theorem}\label{key2} Suppose that $2 \leq p <N, 2 \leq q <N$, $a\geq 0$ and  $(a_1)-(a_2)$ hold.  
	Denote by $m_\infty$ the integer such that
		\[
		\lambda_{m_{\infty}} < \bar\lambda \leq
		\lambda_{m_{\infty}+1}.
		\]
		If 
		\[
		m_{\infty} \not\in [m_0,m_0^*]
		\]
		and 
		\[ (b_-) \qquad
		\lim_{|(s,t)|\to  \infty} \left[ G(s,t) - \frac{1}{p} G_s(s,t)s - \frac{1}{q} G_t(s,t)t\right] = -\infty,\,
		\]
	then there exists a nontrivial weak solution $z=(u,v)$ of \eqref{pqNontrivialSolutions}.
\end{theorem}

\begin{theorem}\label{key3} Suppose that $ 2 \leq p < N, 2 \leq q <N $, $a= 0$ and  $(a_1)-(a_2)$ hold.  
Denote by $m_\infty$ the integer such that
\[
\lambda_{m_{\infty}} \leq \bar\lambda <
\lambda_{m_{\infty}+1}.
\]
	If	
	\[
	m_{\infty} \not\in [m_0,m_0^*]
	\]
	and
	\[ (b_+)\qquad
	\lim_{|(s,t)|\to  \infty} \left[ G(s,t) - \frac{1}{p} G_s(s,t)s - \frac{1}{q} G_t(s,t)t\right] = +\infty,\,
	\]
	then there exists a nontrivial weak solution $z=(u,v)$ of \eqref{pqNontrivialSolutions}.
\end{theorem}

\medskip
We remark that the nontrivial solutions of \eqref{pqNontrivialSolutions}, including those found through Theorems \ref{key1}, \ref{key2} and \ref{key3},  cannot be semitrivial if we assume that there exists $\varepsilon>0$ such that $G_t(s,0)\neq 0$ for any $s \in (-\varepsilon, \varepsilon)\setminus \{0\}$ and  $G_s(0,t)\neq 0$ for any $t \in (-\varepsilon, \varepsilon)\setminus \{0\}$ (see Proposition \ref{c1} and also \cite[Lemma 3.1]{chang-wang}).

\bigskip
Finally, we notice that if $\alpha=\beta=0$ then $p=q=2$. In this case, minimax arguments and Morse theory are applied in  \cite{furtadodepaiva} for deriving existence of nontrivial solutions of a semilinear elliptic system with a $C^2$ nonlinear function $G$.
 
 \bigskip

This Chapter is organized as follows: in Section \ref{SezionePSCPSNontrivialSolution}, we show that the functional $J_{a}$ associated with system \eqref{pqNontrivialSolutions} satisfies the PS condition or the CPS condition, depending on whether resonance conditions occur. In Section \ref{SezioneLinkingComologico}, we prove a Saddle Theorem based on cohomological linking and provide related results concerning $(p,q)$-homogeneous sets. Consequently, in Section \ref{sect:z3} we show that suitable $(p,q)$-homogeneous sets satisfy the geometric assumptions of the Saddle Theorem, thereby obtaining a solution. Finally, in Section \ref{SezioneGruppiCriticiInZero} we present estimates of critical groups at zero, which allow us to ensure that the obtained solution is nontrivial.

\section{PS and CPS conditions}\label{SezionePSCPSNontrivialSolution}

\begin{proposition}\label{PSlim}
	If $\{z_n\}_n$ is a bounded sequence in $X$ such that $\ \|J_\alpha'(z_n)\|\rightarrow~0$, then
	$\{z_n\}_n$ has a convergent subsequence in $X$.
\end{proposition}

\begin{proof}
Let $\{z_n\}_n=\{(u_n,v_n)\}_n \subset X$ be a bounded sequence  such that $\ \|J_a'(z_n)\|\rightarrow~0$. Since  $\{z_n\}_n$ is bounded in $X$, there exists a subsequence, still denoted by $\lbrace{z_n\rbrace}_n,$ that converges to some $z=(u,v)$ weakly in $X$ and strongly in $L^p(\Omega)\times L^q(\Omega)$.\\
In particular, by $\ \|J_a'(z_n)\|\rightarrow~0$ we get

$$ \displaystyle \langle J_a'(z_n),(u_n-u,0)\rangle =   \into (a+|\nabla
		u_n|^2)^{\frac {p-2}{2}}\nabla u_n \cdot \nabla (u_n-u) \, dx + \into  G_s(u_n,v_n)(u_n-u) \, dx \to 0.$$ 
		
\medskip
\noindent
Moreover, since $u_n \to u$ in $L^p(\Omega)$, by first inequality in \eqref{condizionidicrescitasuFsFt} we deduce

$$ \displaystyle \into  G_s(u_n,v_n)(u_n-u) \, dx \to 0.$$

\medskip
\noindent
Therefore we obtain 

$$ \displaystyle \langle H_{a,p}(u_n), u_n - u \rangle = \into (a+|\nabla u_n|^2)^{\frac {p-2}{2}}\nabla u_n \cdot \nabla (u_n-u) \, dx \to 0.$$

\medskip
\noindent
By Lemma \ref{s+PoincareHopf} we know that $H_{a, p}$ is of class $(S)_+$, so that $u_n \to u$ strongly in $W_0^{1,p}(\Omega)$.\\
In the same way, we also get that  $v_n \to v$ strongly in $\sobq$.

\end{proof}	
	
\medskip

\begin{proposition}\label{PS}
If $\bar{\lambda}$ is not an eigenvalue of \eqref{eigen}, then $J_a$ satisfies the $(PS)$ condition at any level $c \in \re$, namely if  $\{z_n\}_n=\{(u_n,v_n)\}_n$ is a sequence in $X$ such that
	$J_a(z_n)\rightarrow c$  and $\|J_a'(z_n)\|\rightarrow 0$, then $\{z_n\}_n$ has a convergent subsequence in $X$.
\end{proposition}

\begin{proof}
Let $\{z_n\}_n=\{(u_n,v_n)\}_n$ be a sequence in $X$ such that
$J_a(z_n)$ is bounded and $\|J_a'(z_n)\|~\rightarrow~0$. Due to the previous proposition, it suffices to prove that the sequence $\{z_n \}_n$ is bounded in $X$.\\
By contradiction, there exists a subsequence, still denoted with $\lbrace{ z_n \rbrace}_n,$  such that $\|(u_n,v_n)\|~\rightarrow~\infty$.
In particular

$$r_n:=\|u_n\|_{1,p}^{p}+ \|v_n\|_{1,q}^q \to + \infty.$$

\medskip
\noindent
Let us set

$$ \displaystyle \bu_n= \frac{u_n}{r_n^{1/p}} \qquad \text{and} \qquad  \bv_n=\frac{v_n}{r_n^{1/q}},$$

\medskip
\noindent
and observe that $\bz_n:=(\bu_n,\bv_n)$ satisfies

$$ \displaystyle \lVert \bu_n \rVert_{1,p}^p + \lVert \bv_n \rVert_{1,q}^q = \frac{\|u_n\|_{1,p}^{p}+ \|v_n\|_{1,q}^q }{r_n}= 1 \qquad \text{ for any } n \in \mathbb{N}.$$

\medskip
\noindent
Hence $ \lVert \bz_n \rVert \leq 2$ for any $n \in \mathbb{N}$, so there exists a subsequence, still denoted by $\lbrace{\bz_n\rbrace}_n,$ that converges to some $\bz=(\bu,\bv)$ weakly in $X$ and strongly in $L^p(\Omega)\times L^q(\Omega)$.

\medskip
\noindent
In particular

\begin{equation}\label{PALSMA1}
	\frac{1}{\ r_n^{\frac{p-1}{p}}}\ \bigl\langle J_a'(z_n),(\bu_n-\bu,0)\bigr\rangle \rightarrow 0
\end{equation}

\medskip
\noindent
and
	
\begin{equation}\label{PALSMA2}
\frac{1}{\ r_n^{\frac{q-1}{q}}}\ \bigl\langle J_a'(z_n),(0,\bv_n -\bv)\bigr\rangle \rightarrow 0.
\end{equation}

\medskip
\noindent
Considering  first inequality in \eqref{condizionidicrescitasuFsFt}, we deduce

\begin{align*}
\frac{1}{\ r_n^{\frac{p-1}{p}}}\ \left| \into G_s(z_n) (\bu_n-\bu)dx \right| \leq   C \left(  \into  |\bu_n-\bu|dx +  \into |\bu_n|^{p-1} |\bu_n-\bu|dx + \into  |\bv_n|^{\frac{q(p-1)} p} |\bu_n-\bu|dx\right).
\end{align*}

\medskip
\noindent
Hence, since $\bu_n \to \bu$ in $L^p(\Omega)$, by H\"{o}lder inequality we obtain

\[ \lim_{n \to \infty} \frac{1}{\ r_n^{\frac{p-1}{p}}}\  \into G_s(z_n)
(\bu_n-\bu)dx =0
\]

\medskip
\noindent
which, combined with (\ref{ps1}), gives

\[ \into\left(\frac{a}{r_n^{2/p}}+|\nabla \bu_n|^2\right)^{\frac{p-2}{2}}\nabla \bu_n \cdot \nabla  \left(\bu_n-\bu \right) dx \ \rightarrow 0.\]

\medskip
\noindent
Using now the convexity (see Lemma \ref{convexity})  of the functional $f:W_0^{1,p}(\Omega) \to \mathbb{R}$ defined as

\begin{equation*}
\displaystyle f(u):= \into \left( \frac{a}{r_n^{2/p}} + | \nabla u |^2 \right)^{\frac{p}{2}} \, dx,
\end{equation*}

\medskip
\noindent
we get

\begin{align*}
\|\bu\|_{1,p}^p & \leq \liminf_{n \to \infty} \|\bu_n\|_{1,p}^p
\leq \limsup_{n \to \infty
} \|\bu_n\|_{1,p}^p \leq
\limsup_{n \to \infty
}\into\left(\frac{a}{r_n^{2/p}}+|\nabla \bu_n |^2\right)^{\frac{p}{2}} \, dx \\
& \leq \limsup_{n \to \infty
}\left[ \into\left(\frac{a}{r_n^{2/p}}+|\nabla \bu|^2\right)^{\frac{p}{2}} + p \into\left(\frac{a}{r_n^{2/p}}+|\nabla \bu_n |^2\right)^{\frac{p-2}{2}} \nabla \bu_n  \cdot \nabla  \left(\bu_n - \bu \right) \, dx \right] \leq \|\bu \|_{1,p}^p ,
\end{align*}

\medskip
\noindent
i.e. $\|\bu_n\|_{1,p}^p \to \|\bu\|_{1,p}^p$, and  by uniform convexity of $W_0^{1,p}(\Omega)$, we deduce $\bu_n\rightarrow \bu$ in $\sob$.
Analogously, by \eqref{PALSMA2}, we get that $\bv_n \rightarrow \bv $ strongly in $\sobq$.

\medskip
\noindent
In particular

\[\|\bu\|_{1,p}^p + \|\bv \|_{1,q}^q =\lim_{n \to \infty} \left\lbrace \|\bu_n \|_{1,p}^p +\|\bv_n \|_{1,q}^q \right\rbrace =1,\]

\medskip
\noindent
so that $\ (\bu,\bv)\neq (0,0)$.

\medskip
\noindent
Since $\|J_a'(z_n)\|~\rightarrow~0$, for any $\tilde{z}=(\tilde{u},\tilde{v}) \in X$ we have

\begin{equation}\label{convergenzatozeroPS}
\displaystyle \left\langle J_a'(z_n) ,   \left( \frac{\tilde{u}}{r_n^{\frac{p-1}{p}}}   , \frac{\tilde{v}}{r_n^{\frac{q-1}{q}}}   \right) \right\rangle \to 0,
\end{equation}

\medskip
\noindent
where 

\begin{align*}
\displaystyle 
\left\langle J_a'(z_n) ,   \left( \frac{\tilde{u}}{r_n^{\frac{p-1}{p}}}   , \frac{\tilde{v}}{r_n^{\frac{q-1}{q}}}   \right) \right \rangle = & \into\left(\frac{a}{r_n^{2/p}}+|\nabla \bu_n|^2\right)^{\frac{p-2}{2}}\nabla \bu_n \cdot \nabla  \tilde{u} \, dx \\
+ & \into\left(\frac{a}{r_n^{2/q}}+|\nabla \bv_n|^2\right)^{\frac{q-2}{2}}\nabla \bv_n \cdot \nabla \tilde{v} \, dx \\
- & \bar{\lambda} \left[ \frac{1}{r_n^{\frac{p-1}{p}}}\into F_s(z_n) \tilde{u} \, dx + \frac{1}{r_n^{\frac{q-1}{q}}}\into F_t(z_n) \tilde{v} \, dx  \right] \\
- & \bar{\lambda} \left[  \frac{1}{r_n^{\frac{p-1}{p}}}\into r_1(z_n) \tilde{u} \, dx + \frac{1}{r_n^{\frac{q-1}{q}}}\into r_2(z_n) \tilde{v} \, dx \right].
\end{align*}

\medskip
\noindent
By assumption $(a_1)$, we get

\[G_s(s,t)=\bar{\lambda} F_s(s,t) +r_1(s,t) \qquad \text{and} \qquad G_t(s,t)=\bar{\lambda} F_t(s,t) +r_2(s,t),
\]

\medskip
\noindent
where $r_1: \mathbb{R} \times \mathbb{R} \to \mathbb{R}$ and $r_2: \mathbb{R} \times \mathbb{R} \to \mathbb{R}$ are continuous functions such that

\begin{equation}\label{r1}
	\lim_{|(s,t)|\to \infty}\  \frac{r_1(s,t)}{|s|^{p-1}+|t|^{q\frac{p-1}{p}} }=0 \qquad \text{and} \qquad \lim_{|(s,t)|\to \infty}\  \frac{r_2(s,t)}{|t|^{q-1}+|s|^{p\frac{q-1}{q}} }=0.
\end{equation}

\medskip
\noindent
Hence there is $c>0$ such that 

\begin{equation}\label{r11}
\left |r_1(s,t)\right | \leq c \left( |s|^{p-1}+|t|^{q\frac{p-1}{p}} +1\right) \qquad \text{and} \qquad \left |r_2(s,t)\right | \leq c \left( |t|^{q-1}+|s|^{p\frac{q-1}{q}} +1\right)
\end{equation}

\medskip
\noindent
for any $(s,t)\in \re^2$. We want to prove that

\begin{equation*}
\lim_{n \to \infty}   \frac{1}{r_n^{\frac{p-1}{p}}}	\into |r_1(z_n(x))|\ |\tilde{u}(x)|dx =0 \qquad \text{and} \qquad \lim_{n \to \infty}   \frac{1}{r_n^{\frac{q-1}{q}}}	\into |r_2(z_n(x))|\ |\tilde{v}(x)|dx =0.
\end{equation*}

\medskip
\noindent
Recalling that $\|\bu_n \|_{1,p},\ \|\bv_n\|_{1,q} \leq 1$, let $\bar{c}>0$ be such that

\begin{equation}\label{barc}
\| \bu_n\|_p^{p-1}+ 
\| \bv_n\|_q^{q\frac{p-1}{p}}
\leq \bar{c} \qquad {\rm for\ any\ } n \in \na.
\end{equation}

\medskip
\noindent
Let us fix $\varepsilon >0$. By (\ref{r1}), there is $\delta_\varepsilon >0$ such that

\begin{equation}\label{deleps}
 |(s,t)|>\delta_\varepsilon \quad \Rightarrow  \quad 
\left| r_1(s,t)\right |< \frac{\varepsilon}{2 \bar{c}}\left(|s|^{p-1}+|t|^{q\frac{p-1}{p}} \right).
\end{equation}

\medskip
\noindent
Denoting by
$\Omega ^-_{n,\varepsilon}=\left\{ x \in \Omega \ : \ |z_n(x)|\leq \delta_\varepsilon \right\}\ $ and $\ \Omega ^+_{n,\varepsilon}=\Omega \setminus \Omega ^-_{n,\varepsilon}$,

\[x\in \Omega ^-_{n,\varepsilon} \quad  \Rightarrow \quad |u_n(x)|, \,
|v_n(x)|\leq  \delta_\varepsilon 
\]

\medskip
\noindent
so, by (\ref{r11}),

\[
  \frac{1}{r_n^{\frac{p-1}{p}}}\int_{\Omega ^-_{n,\varepsilon}} |r_1(z_n(x))|\ |\tilde{u}(x)|dx
\leq \frac{c}{r_n^{\frac{p-1}{p}}}\,\left(|\delta_\varepsilon|^{p-1}
+|\delta_\varepsilon|^{q\frac{p-1}{p}}+1
\right)\, 
\|\tilde{u}\|_1\]

\medskip
\noindent
and, choosing $n$  big enough, 

\begin{equation}\label{r-}
  \frac{1}{r_n^{\frac{p-1}{p}}}	\int_{\Omega ^-_{n,\varepsilon}} |r_1(z_n(x))|\ |\tilde{u}(x)|dx <\frac{\varepsilon}{2}\,\|\tilde{u}\|_p\, . \qquad 
\end{equation}

\medskip
\noindent
On the other hand, by (\ref{deleps}), H\"{o}lder inequality and (\ref{barc}),

\begin{align*}
&  \frac{1}{r_n^{\frac{p-1}{p}}}	\int_{\Omega ^+_{n,\varepsilon}} |r_1(z_n(x))|\ |\tilde{u}(x)|dx\\
&\leq  \frac{\varepsilon}{2 \bar{c}}  \frac{1}{r_n^{\frac{p-1}{p}}}   \int_{\Omega ^+_{n,\varepsilon}}\left(|u_n(x)|^{p-1}+|v_n(x)|^{q\frac{p-1}{p}}\right)|\tilde{u}(x)|dx\\
& \leq \frac{\varepsilon}{2 \bar{c}} \into \left(|\bu_n(x)|^{p-1}+|\bv_n(x)|^{q\frac{p-1}{p}}\right)|\tilde{u}(x)|dx\\
&\leq \frac{\varepsilon}{2} \ \|\tilde{u}\|_p
\end{align*}

\medskip
\noindent
which, together with (\ref{r-}), gives

$$\lim_{n \to \infty}   \frac{1}{r_n^{\frac{p-1}{p}}}	\into |r_1(z_n(x))|\ |\tilde{u}(x)|dx =0.$$

\medskip
\noindent
Similarly, we also get 
$$\lim_{n \to \infty}   \frac{1}{r_n^{\frac{q-1}{q}}}	\into |r_2(z_n(x))|\ |\tilde{v}(x)|dx =0.$$

\medskip
\noindent
Since $\bu_n \to \bu$ in $W_0^{1,p}(\Omega)$ and $\bv_n \to \bv$ in $W_0^{1,q}(\Omega)$, we get

$$ \displaystyle \lim_{n \to \infty} \into\left(\frac{a}{r_n^{2/p}}+|\nabla \bu_n|^2\right)^{\frac{p-2}{2}}\nabla \bu_n \cdot \nabla  \tilde{u} \, dx = \into |\nabla \bu|^{p-2} \nabla \bu \cdot \nabla  \tilde{u} \, dx,$$

$$ \displaystyle \lim_{n \to \infty} \into\left(\frac{a}{r_n^{2/q}}+|\nabla \bv_n|^2\right)^{\frac{q-2}{2}}\nabla \bv_n \cdot \nabla  \tilde{v} \, dx = \into |\nabla \bv|^{q-2} \nabla \bv \cdot \nabla  \tilde{v} \, dx,$$

\medskip
\noindent
and

$$ \lim_{n \to \infty} \frac{1}{r_n^{\frac{p-1}{p}}} \into F_s(u_n(x) ,v_n(x) ) \tilde{u}(x)\,dx =\into F_s(\bu(x) ,\bv(x) ) \tilde{u}(x)\,dx, $$

$$ \lim_{n \to \infty} \frac{1}{r_n^{\frac{q-1}{q}}} \into F_t(u_n(x) ,v_n(x) ) \tilde{v}(x)\,dx =\into F_t(\bu(x) ,\bv(x) ) \tilde{v}(x)\,dx. $$

\medskip
\noindent
Considering these convergences in \eqref{convergenzatozeroPS}, we have proved

\begin{align*}
\displaystyle \into |\nabla \bu|^{p-2} \nabla \bu \cdot \nabla  \tilde{u} \, dx 
& = \bar{\lambda} \into F_s(\bu,\bv) \tilde{u} \, dx \\
\into |\nabla \bv|^{q-2} \nabla \bv \cdot \nabla  \tilde{v} \, dx 
& = \bar{\lambda} \into F_t(\bu,\bv) \tilde{v} \, dx
\end{align*}

\medskip
\noindent
for any  $\tilde{z}=(\tilde{u},\tilde{v}) \in X$. In other words, $(\bu,\bv)\neq (0,0)$ solves (\ref{eigen}) with $\lambda=\bar{\lambda}$, so the contradiction proves that $\{z_n\}_n$ is bounded in $X$.
\end{proof}

\bigskip
\begin{proposition}\label{CPS}
	If $\bar{\lambda}$ is an eigenvalue of (\ref{eigen}) and $(b_-)$ holds, so  that
	\[
			\lim_{|(s,t)|\to  \infty} \left[ G(s,t) - \frac{1}{p} G_s(s,t)s - \frac{1}{q}G_t(s,t)t\right] = -\infty, 
	\]
	then $J_a$ satisfies the $(CPS)$ condition at any level $c \in \re$, namely any sequence $\{z_n\}_n=\{(u_n,v_n)\}_n$  in $X$ such that
	$J_a(z_n)\rightarrow c$  and $\|J_a'(z_n)\|(1+\|z_n\|)\rightarrow 0$ has a convergent subsequence in $X$.
	
	\medskip
	
	Moreover, if $\bar{\lambda}$ is an eigenvalue of (\ref{eigen}) and $(b_+)$ holds, so  that
	\[
	\lim_{|(s,t)|\to  \infty} \left[ G(s,t) - \frac{1}{p} G_s(s,t)s - \frac{1}{q}G_t(s,t)t\right] = +\infty, 
	\]
	then $J_0$ satisfies the $(CPS)$ condition at any level $c\in \re$.
\end{proposition}

\begin{proof}
Let $\{z_n\}_n=\{(u_n,v_n)\}_n$ be a sequence in $X$ such that
$J_a(z_n)$ is bounded and

$$\|J_a'(z_n)\|(1+\|z_n\|)\rightarrow 0.$$

\medskip
\noindent
Firstly we notice that previous condition implies $\|J_a'(z_n)\|\rightarrow 0,$ hence by Proposition \ref{PSlim} it suffices to prove that the sequence $\{z_n \}_n$ is bounded in $X$.\\
By contradiction, there exists a subsequence, still denoted with $\lbrace{ z_n \rbrace}_n,$  such that $\|(u_n,v_n)\|~\rightarrow~\infty$.
As in  the proof of Proposition~\ref{PS}, we set 

$$r_n:=\|u_n\|_{1,p}^{p}+ \|v_n\|_{1,q}^q \to + \infty,$$

\medskip
\noindent
and

$$ \displaystyle \bu_n= \frac{u_n}{r_n^{1/p}} \qquad \text{and} \qquad  \bv_n=\frac{v_n}{r_n^{1/q}},$$

\medskip
\noindent
so that $\{ \bz_n \}_n=\{(\bu_n,\bv_n)\}_n$ is bounded in $X$ and there exists a subsequence, still denoted by $\lbrace{\bz_n\rbrace}_n,$ that converges to some $\bz=(\bu,\bv)$ weakly in $X$ and strongly in $L^p(\Omega)\times L^q(\Omega)$. Arguing as in the proof of  Proposition~\ref{PS}, we deduce that $\lbrace{\bz_n\rbrace}_n$ converges to $\bz=(\bu,\bv)$ strongly in $X$, with $\bz \neq (0,0)$.\\
Let us show that $\{z_n\}_n$ is bounded in $X$ when  $(b_-)$ holds and $a\geq 0$.\\
Denoting by 

$$K(s,t)=G(s,t) - \frac{1}{p} G_s(s,t)s - \frac{1}{q}G_t(s,t)t,$$

\medskip
\noindent
by $(CPS)$ condition there exists a constant $C \in \re$ such that
	
\begin{align}\label{calcoloperReverseFatouLemma}
C
& \geq J_a(z_n) - \left\langle J_a' (z_n), \left( \frac{u_n}{p}, \frac{v_n}{q} \right) \right\rangle \nonumber \\
& = \frac{1}{p} \into  \left(a+|\nabla u_n|^2\right)^{\frac p 2} \ dx - \frac{1}{p} \into  \left(a+|\nabla u_n|^2\right)^{\frac{p-2}{2}} \left| \nabla u_n \right|^2 \, dx \nonumber \\
& \quad + \frac{1}{q}\into \left(a+|\nabla v_n|^2\right)^{\frac{q}{2}} \ dx - \frac{1}{q} \into  \left(a+|\nabla v_n|^2\right)^{\frac{q-2}{2}} \left| \nabla v_n \right|^2 \, dx \nonumber \\
& \quad - \left[ \into  \left( G(u_n(x),v_n(x))  -  \frac{1}{p} G_s(u_n(x),v_n(x)) u_n(x)   -   \frac{1}{q} G_t(u_n(x),v_n(x)) v_n(x) \right) \, dx  \right] \nonumber \\
&=  \frac{a}{p} \into \left(a+|\nabla u_n|^2\right)^{\frac{p-2}{2}} \, dx  + \frac{a}{q} \into \left(a+|\nabla v_n|^2\right)^{\frac{q-2}{2}}\, dx  -\into K(z_n(x)) \, dx \nonumber \\
& \geq -\into K(z_n(x))\, dx.
\end{align}
    
\medskip
\noindent
From $(b_-)$, there is $c_1 \in \re$ such that 

\[ K(s,t)\leq c_1 \qquad \forall (s,t) \in \re^2, \]

\medskip
\noindent
hence, by reverse Fatou's lemma we get

$$ \displaystyle \limsup_{n \to \infty} \into K(z_n(x)) \, dx \leq \into \limsup_{n \to \infty} K(z_n(x)) \, dx$$

\medskip
\noindent
and taking into account inequality \eqref{calcoloperReverseFatouLemma} we deduce

\begin{equation}\label{lsk}
		-C \leq \into \limsup_{n\to \infty}K\bigl(z_n(x)\bigr)\, dx.
\end{equation}

\medskip
\noindent
As we are assuming $r_n \to +\infty$, by $(b_-)$ we get that, for almost every $x\in \Omega$
	\begin{align*}
& \bz(x)\neq (0,0) \ \Rightarrow \   |z_n(x)|=\left|\bigl(r_n^{1/p}\bu_n(x),r_n^{1/q}\bv_n(x)\bigr)\right|\to \infty
	\ \Rightarrow \ K\bigl(z_n(x)\bigr) \to -\infty.
\end{align*}
	So, from (\ref{lsk}), we infer that $\bz(x)=(0,0)$ almost everywhere in $\Omega$, which contradicts $\bz\neq (0,0)$.
	Consequently $\{z_n\}_n$ is bounded and, by Proposition~\ref{PSlim}, has a convergent subsequence in $X$. 
	
	\bigskip
	
	If instead $(b_+)$ holds, $J_0(z_n)\rightarrow c\in \re\ $  and $\|J_0'(z_n)\|(1+\|z_n\|)\rightarrow 0$, then
	\[\lim_{n\to \infty}\into K(z_n(x))\, dx= -c.
	\]
	
	\smallskip
	
	Then, reasoning as before, we infer again that $\{z_n\}_n$ has a convergent subsequence in $X$. 
\end{proof}

\bigskip

\section{Cohomological linking}\label{SezioneLinkingComologico}
\label{sect:z2}

Throughout this section, $Y$ denotes a Banach space and 
$f:Y\to\re$ a $C^1$ function.
We also denote by $H^*$ the
Alexander-Spanier cohomology \cite{spa} with coefficients in $\ze_2=\{-1,1\}$.
\par
Let us recall the following definitions (see \cite{CCMV,chang,D,DL}).

\begin{definition}
	\label{df:coho}
	Let $D,S,A$ be three subsets of $Y$, $m$ a nonnegative integer.
	
	We say that {\em $(D,S)$ links $A$ cohomologically in dimension} $m$
	(over $\ze_2$) if $S\subseteq D$, $S \cap A =\emptyset$ and
	the restriction homomorphism
	$H^m(Y,Y\setminus A)\to H^m(D,S)$ is not identically zero.
\end{definition}

\smallskip
If $(D,S)$ links $A$ cohomologically in some dimension, then $(D, S)$ links $A$  (c.f. \cite[Definition 5.1]{CD}). Moreover if $(D,S)$ links $A$, then $D \cap A \not = \emptyset$.

\smallskip

\begin{definition}
	\label{criticalgroup} Let $G$ be an abelian group.
	Let the $m$-th critical group of $f$ at $z \in Y$ with coefficient in $G$ be defined by 
	\[
	C_m(f;z) = H^m\left(f^{c}, f^{c}\setminus \{z\}\right) \,,
	\]
	where $c=f(z)$, $f^c =\{ u \in Y \ | \ f(u) \leq c\}$.
\end{definition}

\bigskip
In general it can happen that $C_m(f,z)$ is not finitely generated for some $m$ and that $C_m(f,z) \neq 0$ for infinitely many $m$'s.

Now we assume that  $f$ is Gauteax differentiable and we denote by $U$ an open subset of $X$.	
If $z$ is an isolated critical point of $f$ and $f':U \to X'$ is a demicontinuous function (namely it is continuous from the strong topology to the weak topology) and of class $(S)_+$ in a neighborhood of $u$, then  $C_*(f,z)$ is of finite type (see \cite[Theorem 1.1]{CD1} and \cite[Theorem 3.4]{ad}).

\bigskip
Let us introduce a general result concerning a $C^1$-functional in a Banach space
which extends Theorem 3.2 in \cite{CD} to the case in which $f$ satisfies only the (CPS) condition.

\begin{theorem}\label{th key}
	Let $Y$ be a Banach space, $f\in C^1(Y,\re)$, $D,S,A$ be three subsets of $Y$, $m \in \na$. 
	\par Assume that
	
	$(D,S)$ links $A$ cohomologically in dimension $m$ over $\ze_2$,
	\[ \sup_S f<\inf_A f=a,  \qquad b=\sup_D f<+\infty,
	\]
	$f$ verifies the (CPS)  condition at any level $c\in [a,b]$ and $f^{-1}([a,b])$ contains a finite number of critical points.
	\newline
	Then there exists a critical point $z$ of $f$ with
	$a\leq f(z)\leq b$ and $C_m(f;z)\neq~ \!\{0\}$.
\end{theorem}

\bigskip

\begin{proof} 
	Since $D \cap A \not = \emptyset$, then $a \leq b$.
	We aim to apply \cite[Theorem 5.2]{D}.
	If $(E,d)$ is a metric space, $f:E\to \re$ is a continuous function 
	and $z\in E$, let us consider the notion of {\sl weak slope} $|df|(z)$ 
	as defined in \cite[Definition~2.1]{DM} (see also \cite[Definition~5.1]{K}).
	Consequently we say that
	\newline$\bullet$ $z$ is a critical point if \ $|df|(z)=0$; 
	\newline$\bullet$ $f$ satisfies the $(PS)$ condition at a level $c \in \re$ if    any sequence  $z_n$ in $E$ such that
	$f(z_n)\rightarrow c$  and $|df|(z_n)\rightarrow 0$
	has a convergent subsequence in $E$. 
	\par
	In particular, if $E$ is also a Banach space and $f$ is $C^1$, we have
	\begin{equation*}\label{we}
		|df|(z)=\|f'(z)\| \qquad \forall\, z \in E.
	\end{equation*}
	Due to \cite[Theorem 4.1,  Remark 4.4]{C}, there exists a metric $\tilde{d}$ on $Y$, 
	topologically equivalent to the metric induced by $\| \  \|_Y$, such that $(Y,\tilde{d})$ 
	is complete and, denoting by $|\tilde{d}f|$ the weak slope of $f$ with respect to the metric $\tilde{d}$, 
	\begin{equation*}\label{dtilde}
		|\tilde{d}f|(z)=\|f'(z)\|_{Y'}(1+\|z\|_Y) \qquad \forall \, z\in Y.
	\end{equation*}
	Therefore, as $f$ satisfies the $(CPS)$ condition at any level $c \in [a,b]$,     
	then $f$ satisfies the $(PS)$ condition at any level $c \in [a,b]$, with respect to $\tilde{d}$.
	\par
	By \cite[Theorem 7.5]{D} the assumptions of \cite[Theorem 5.2]{D} are satisfied.
	Then the assertion follows, taking into account \cite[Proposition 7.3 and Remark 5.3]{D}.
	 
\end{proof}

\bigskip
\bigskip

\bigskip
Denoting by $\mathcal{A}$ the class of symmetric subsets of $Y$, Fadell and Rabinowitz (see \cite{fadell_rabinowitz1977,fadell_rabinowitz1978}) constructed an index theory $i:\mathcal{A}\to \na \cup \{\infty\}$ with the following properties:

\bigskip

\begin{itemize}
	\item[(i)] Definiteness: $i(A)\geq 0,\ i(A)=0 \Leftrightarrow A=\emptyset$;
	\item[(ii)] Monotonicity: If there is an odd continuous map $A\to A'$, then
	\[i(A)\leq i(A'),\]
	in particular, equality holds if $A$ and $A'$ are homeomorphic;
	\item[(iii)] Subadditivity:
	\[i(A\cup A')\leq i(A)+i(A');
	\]
	\item[(iv)] Continuity: If $A$ is closed, there is a neighborhood $U\in \mathcal{A}$ of $A$ such that
	\[i(U)=i(A);
	\]
	\item[(v)] Neighborhood of zero: If $U$ is a bounded symmetric neighborhood of $0$ in $Y$,
	\[i(\partial U)=\dim Y;\]
	\item[(vi)] Stability: If $A$ is closed and $A *\ze_2$ is the join of $A$ with $\ze_2$, realized in $Y \oplus \re$,
	\[i(A *\ze_2)=i(A)+1\]
	where $A *\ze_2$ is the union of all line segments in $Y \oplus \re$, joining $\{1\}$ and $\{-1\}$ to points of $A$;
	\item[(vii)] Piercing property: Assume that $A,\ A_0,\ A_1$
	are closed and
	\[\varphi:A\times [0,1]\to A_0\cup A_1 \]
	is an odd continuous map such that $\varphi(A\times [0,1])$ is closed, and
	\[ \varphi(A\times \{0\})\subset A_0, \qquad \varphi(A\times \{1\})\subset A_1.\]
	Then
	\[i\left(\varphi(A\times [0,1])\cap A_0 \cap A_1\right)\geq i(A).\]	
\end{itemize}

\bigskip

We can state the next theorem whose proof follows  from \cite[Theorem 3.6]{CD} and 
\cite[Theorem 2.7]{DL}.  For reader's convenience, we give a sketch of the proof.

\begin{theorem}\label{thmCD}
	Let $S,A$ be two symmetric subsets of $Y$ with $S\cap A=\emptyset$.
	Assume that $0\in A$ and $i(S) = i(Y\setminus A)=m<\infty$.
	\par
	Then $(Y,S)$ links $A$ cohomologically in dimension $m$ over $\ze_2$.
\end{theorem}

\begin{proof}
	We want to prove that 
	$H^m(Y,Y\setminus A)\to H^m(Y,S)$ is not identically zero. So, considering the exact sequence
	\[H^m(Y,Y\setminus A)\to H^m(Y,S)\to H^m(Y\setminus A,S)
	\]
	it suffices to show that
	\begin{equation}\label{gno}
		H^m(Y,S)\to H^m(Y\setminus A,S) \qquad {\rm is \ not\ injective.}
	\end{equation}
	
	If $m=0$, then $S=\emptyset$ and $A=Y$, so $H^0(Y,S)\neq \{0\}$ while $H^0(Y\setminus A,S)=\{0\}$, whence (\ref{gno}) is proved.
	
	\bigskip
	
	Let us denote by $j:H^{m-1}(Y\setminus A)\to H^{m-1}(S)$ and $i:H^{m-1}(Y)\to H^{m-1}(S)$ the homomorphisms induced by inclusions.
		We recall that $H^q(Y)=\{0\}$ if $q\neq 0$, while the dimension of $H^0(Y)$ is 1.

	By Lemma 3.4 in \cite{CD}, we infer that
	\begin{itemize}
		\item[$(i)$]	if $m=1$, the dimension of $im(j)$ is at least 2, while the dimension of $im(i)$  is at most 1.
		\item[$(ii)$]
		if $m\geq 2$, the dimension of $im(j)$ is at least 1, while $im(i)=\{0\}$.
		\end{itemize}
\smallskip
Hence $im(i)$ is a proper subset of $im(j)$, for any $m\geq 1$.

Consider now the commutative diagram
\[
\begin{CD}
H^{m-1}(Y) @>{i}>> H^{m-1}(S) @>{\delta}>>
H^m(Y,S) \\
@V{}VV @V{\mathrm{Id}}VV @V{}VV \\
H^{m-1}(Y\setminus A) @>{j}>> H^{m-1}(S) @>{{\bar \delta}}>>
H^m(Y\setminus A,S)
\end{CD}
\]
where the horizontal rows are exact sequences and the vertical rows are induced by
inclusions. \newline Let $\omega \in im(j)\setminus im(i)$,
then $\delta\omega\neq 0$
in $H^m(Y,S)$, while ${\bar \delta}\omega=\left(\delta
\omega\right)_{|(Y\setminus A,S)}=0$. \newline
Hence $(\ref{gno})$ is proved also for any $m\ge 1$, as $H^m(Y,S)\to H^m(Y\setminus A,S)$ is not injective.
\end{proof}

\bigskip
Now we consider the Banach space $X = W_0^{1,p}(\Omega)\times
W_0^{1,q}(\Omega)$.
We recall that the sequence $(\lambda_m)$ is defined by (\ref{lambdak}), i.e. 
\[ 	
\lambda_m = \inf_{A \in \Sigma_m} \sup_{(u,v) \in A} \Phi(u,v)
\]
where $ \Sigma_m = \{ A\subset \Sigma  \ | \, A \ is\  symmetric,\ compact\ and\  i(A) \geq m \}.$
Let us recall the following Lemma, proved in \cite{DL}.

\begin{lemma}\label{lemma Bm}
	For every symmetric and open subset $A$ of $\Sigma$
	\[i(A)=\sup \{ i(K)\, | \, K \ is \ compact\ and\  symmetric\  with\  K\subseteq A \}.
	\]
\end{lemma}

\begin{remark}\label{Bm}
	Denoting by
	$\mathcal{B}_m = \{ A\subset \Sigma  \ | \, A \ is\  
	symmetric\ and\  i(A) \geq m \}$,
	the previous Lemma assures that
		\[ 	\lambda_m = \inf_{A \in \mathcal{B}_m} \sup_{(u,v) \in A} \Phi(u,v).\]
	
\end{remark}

\bigskip

\begin{proposition}\label{indicivari}
Let $\alpha \in \re$ and $r>0$ and set
\[ X^\alpha=\{z \in X\ | \ \Phi(z)\leq \alpha\Psi(z)\}, \qquad 
\Sigma^\alpha=\Sigma \cap  X^\alpha,\]
\[ S^\alpha_r=\{z \in X^\alpha\ | \ \Phi(z)=r\},
\]
\[ \stackrel{\hspace{-2mm}\circ}{X^\alpha}=\{z \in X\ | \ \Phi(z)< \alpha\Psi(z)\},
 \qquad 
\stackrel{\hspace{-2mm}\circ}{\Sigma^\alpha}=\Sigma \, \cap \stackrel{\hspace{-2mm}\circ}{X^\alpha}.
\]

\smallskip
	
Then
\begin{itemize} 
	\item[$(a)$] $i(\Sigma^{\alpha})\geq m$, for any $\alpha \geq \lambda_m$;
	\item[$(b)$] $\lambda_m< \lambda_{m +1}\ \Rightarrow \ i(\Sigma^{\alpha})= m$, for any $\alpha \in [\lambda_m, \lambda_{m +1})$;
	\item[$(c)$]  $\lambda_m< \lambda_{m +1}\ \Rightarrow \ i\bigl(\stackrel{\hspace{-7mm}\circ}{\Sigma^{\lambda_{m +1}}}\bigr)=m$;
	\item[$(d)$] $\Sigma^{\alpha}$ is an odd strong deformation retract of $X^\alpha\setminus \{0\}$;
	\item[$(e)$] $\stackrel{\hspace{-2mm}\circ}{\Sigma^\alpha}$ is an odd strong  deformation retract of $\stackrel{\hspace{-2mm}\circ}{X^\alpha}\setminus \{0\}$;
	\item[$(f)$] $S^\alpha_r$ is an odd strong  deformation retract of $X^\alpha\setminus \{0\}$.
	
\end{itemize}

\end{proposition}

\begin{proof}
	Let us consider $\alpha >\lambda_1$, otherwise the proof is even simpler.
	\newline By contradiction, assume that $i(\Sigma^{\lambda_m})\leq m-1$. 
	Continuity assures that there is a closed neighborhood $N$ of $\Sigma^{\lambda_m}$ such that
	$i(N)=i(\Sigma^{\lambda_m})$. As $N$ is also a neighborhood of the critical set $ \{z 
	\ | \ \Phi_{|\Sigma}(z)=\lambda_m\}$, by the equivariant deformation theorem, there exist $\varepsilon >0$ and an odd continuous map $\eta:\Sigma ^{\lambda_m+ \varepsilon}\hspace{-1mm}\to \Sigma ^{\lambda_m- \varepsilon}\cup N=N$. So, due to monotonicity,
	$i(\Sigma ^{\lambda_m+ \varepsilon})\leq i(N)=i(\Sigma^{\lambda_m})\leq m-1$.
	By definition of $\lambda_m$, there is $\bar A \in \Sigma_m$ such that \ $\sup \Phi(\bar A)<\lambda_m+ \varepsilon$. Hence $\bar A \subset \Sigma ^{\lambda_m+ \varepsilon} $ and we have the contradiction
	$m\leq i(\bar A)\leq i(\Sigma ^{\lambda_m+ \varepsilon})\leq m-1$, which proves $(a)$.
	
	\smallskip
	
	 Now, if $\lambda_m< \lambda_{m +1}$ and $\alpha \in [\lambda_m, \lambda_{m +1})$, by $(a)$ we know that $i(\Sigma^{\alpha})\geq m$.  
	By contradiction, assume that $i(\Sigma^{\alpha})\geq m+1$. Recalling Remark \ref{Bm}, $\Sigma^{\alpha} \in \mathcal{B}_{m+1}$ and 
	\[\lambda_{m +1}=\inf_{A\in \mathcal{B}_{m+1}} \sup_A \Phi\leq \sup_{\Sigma^{\alpha} } \Phi=\alpha<\lambda_{m +1}
	\]
	so we conclude that $i(\Sigma^{\alpha})= m$.
	\newline Moreover, still due to $(a)$,  $i\bigl(\stackrel{\hspace{-7mm}\circ}{\Sigma^{\lambda_{m +1}}}\bigr)\geq m$. Let $K$ be a symmetric and compact subset of $\stackrel{\hspace{-7mm}\circ}{\Sigma^{\lambda_{m +1}}}$, so that $\max \Phi(K)\leq \alpha$ for a suitable $\alpha \in [\lambda_m, \lambda_{m +1})$.
	 $(b)$ gives that $i(K)\leq i(\Sigma^{\alpha})=m$, hence $(c)$ comes from
	 Lemma \ref{lemma Bm}.
	 
	 \smallskip
	 
	 For any $s\in [0,1]$ and $z=(u,v)\in X\setminus \{0\}$, we define
	 $\gamma_{s,z}=1-s+\frac{s}{\Psi(z)}$ 
	 and $\eta (s,(u,v))=\left(\gamma_{s,z}^{1/p}u,\gamma_{s,z}^{1/q}v\right)$.
	 
	 $\eta$ is clearly continuous and odd. In addition, due to the $(p,q)$-homogeneity of $\Phi$ and $\Psi$, showed in (\ref{pqh}),
	 \[\Phi(\eta(s,z))=\gamma_{s,z}\Phi(z) \quad {\rm and} \quad \Psi(\eta(s,z)) =
	 \gamma_{s,z}\Psi(z). \]
	 Hence $(d)$ and $(e)$ are proved as
	 \[z\in X^\alpha \setminus \{0\}\quad \Rightarrow \quad 
	 \eta(0,z)=z, \ \eta(1,z) \in \Sigma^\alpha \ {\rm and}\ \
	 \eta(s,z) \in X^\alpha \setminus \{0\} \ \  \forall s\in [0,1]
	 \]
	 $z\in \Sigma \quad \Rightarrow\quad \eta(s,z)=z$ for any $s\in [0,1].$
	 
	 Analogously we obtain $(f)$ setting 
	 
	 $\delta_{s,z}=1-s+s\frac{r}{\Phi(z)}$ \quad and \quad
	 $\bar \eta(s,(u,v))=\left(\delta_{s,z}^{1/p}u,\delta_{s,z}^{1/q}v\right)$.
	 \end{proof}

\bigskip

Let $m \in \na\ $ be such that $\la_m<\la_{m+1}$.

If $m\geq 1$ we set
\[\xm := \{ z= (u,v) \in X  \ |    \ \Phi(u,v) \leq \la_m 
\Psi(u,v)\}
\]
\[\xp := \{ z= (u,v) \in X  \ |   \ \Phi(u,v) \geq \la_{m+1} \Psi(u,v) \}\]

and if $m=0$
\[X_-^0 := \{0\} \qquad X_+^0 :=X.\]

\bigskip

We underline that the sets  $\xm$ and $\xp$ are symmetric and
$(p,q)$-homogeneous, i.e., 
\[ (u,v) \in X^m_{\pm} \Rightarrow (\delta^{\frac 1 p}u,\delta^{\frac 1 q}v) \in X^m_{\pm} \qquad \qquad \hbox{ for any } \delta >0.
\]
In particular, if $p\neq q$ and $m\geq 1$, the sets $\xp$ and $\xm$ are not cones.

\bigskip

\begin{theorem}\label{th dr}
	
	Setting 
	
	\[ D_r=\{z\in \xm\ |\ \Phi(z)\leq r\} \qquad  S_r=\{z\in \xm\ |\ \Phi(z)= r\},
	\]
	 
	 \bigskip
	
	$(D_r,S_r)$ links $\xp$ cohomologically in dimension $m$ over $\ze_2$, for any $r>0$.
\end{theorem}

\begin{proof}
	From Proposition \ref{indicivari} we obtain that
	
	\[ i(S_r)=i(X \setminus \xp)=m.
	\]
	As $S_r$ and $\xp$ are symmetric, $0\in \xp$ and $S_r\cap \xp=\emptyset$, Theorem \ref{thmCD} gives that
	\begin{equation}\label{prefive}
	(X,S_r)\  {\rm links}\ \xp \hbox{ cohomologically in dimension } m \hbox{ over } \ze_2.
	\end{equation}
	
	Consider the commutative diagram
	
	\begin{align*}
		H^{m-1}&(X)\to H^{m-1}(\xm \!\setminus\!\{0\})\to H^{m}(X,\xm\!\setminus\! \{0\})
		\to H^m (X) \to  H^m (\xm \!\setminus\!\{0\})\\
		&\downarrow \gamma_1 \hspace{15mm} \downarrow \gamma_2 \hspace{20mm}	\downarrow \gamma_3 \hspace{20mm}	\downarrow \gamma_4 \hspace{16mm}	\downarrow \gamma_5
		\\
		H^{m-1}&(D_r)\ \longrightarrow\ H^{m-1}(S_r)\ \longrightarrow \ 
		H^m( D_r,S_r)\ \,\longrightarrow \ \, H^{m}(D_r)\ \longrightarrow \ H^{m}(S_r)\\
		\end{align*}
	
	where the horizontal rows are exact sequences and the vertical rows are induced by
	inclusions. 
	
	In particular, $ \gamma_1,\ \gamma_2,\ \gamma_4$ and $\gamma_5$ are isomorphisms, as $D_r$ is homotopic to $X$ and $S_r$
	is homotopic to $\xm\setminus \{0\}$.
	Hence also $\gamma_3$ is an isomorphism, which, combined with (\ref{prefive}), completes the proof.
\end{proof}

{\vskip 1cm}

\section{Critical Groups Estimates at zero}\label{SezioneGruppiCriticiInZero}

From now on we assume that  $p,\,q \in[2,N)$, $a \geq 0$, $H \in C^1(\re^2,\re)$ and there exist $C>0$, $p'\in(p,p^*)$,  $q'\in (q,q^*)$ 
such that  

\begin{equation}\label{primaNontrivialSolutions}
	\begin{array}{l}
		|H_s(s,t)|\leq C
		\left(|s|^{p'-1}+|t|^{q'\frac{p'-1}{p'}}+1\right),\\
		\ \\
		| H_t(s,t)| \leq C
		\left(|s|^{p'\frac{q'-1}{q'}}+|t|^{q'-1}+ 1 \right).
	\end{array}
\end{equation}

\medskip
\noindent
We denote by $I_{a,H}$ the $C^1$-functional defined for any $z=(u,v) \in X$ by

\begin{equation}\label{IH}
I_{a,H}(z)  
=  \frac{1}{p} \into  \left(a+
|\nabla u|^2\right)^{\frac{p}{2}} \ dx + \frac{1}{q} \into \left(a+
|\nabla v|^2\right)^{\frac{q}{2}} \ dx  - \into
H(u,v) \ dx.
\end{equation}

\begin{theorem}\label{Hs}
	Let  $p,\,q \in[2,N)$, $a \geq 0$. 
		Assume that $H \in C^1(\re^2,\re)$ satisfies $\eqref{primaNontrivialSolutions}$ and $\nabla H(0,0)=(0,0)$.  
			
	\medskip
	
	If 
	$\mathbf{0}=(0,0)$ is an isolated critical point for $I_{a,H}$, there are $\overline {H\,} \in C^1(\re^2,\re)$ and $\bar\eta>0$ such that:
	
	\bigskip
	
	\begin{itemize}
		\item 	$\overline {H\,}(s,t)=H(s,t)$ when $|s|,|t|\leq \bar\eta$;
		\item $\mathbf{0}$ is an isolated critical point for $I_{a,\overline {H\,}}$ and $C_m(I_{a,H},\,\mathbf{0})=C_m(I_{a,\overline {H\,}},\,\mathbf{0})$, for any  $m\in \na$;
		\item $\overline {H\,}(\re^2)=\overline {H\,}([-2\bar\eta,2\bar\eta]^2)\subset H\left([-2\bar\eta,2\bar\eta]^2\right)$. 
	\end{itemize}
	In particular, if $(s,t)\notin [-2\bar\eta,2\bar\eta]^2$, then
	
	\medskip
	
	$\overline {H\,}(s,t)=\begin{cases}
	\begin{array}{ll}
	H(0,0)& {if\ }|s|\geq 2\bar\eta,\ |t|\geq 2\bar\eta,\\
	\overline {H\,}(0,t)& {if\ }|s|\geq 2\bar\eta,\ |t|\leq 2\bar\eta,\\
	\overline {H\,}(s,0)& {if\ }|s|\leq 2\bar\eta,\ |t|\geq 2\bar\eta.
	\end{array}
	\end{cases}$

\medskip
\noindent
In addition, if $H$ is  $C^2$ at least in an open set  containing $[-2\bar{\eta},2\bar{\eta}]^2$, then  $\overline {H\,}\in C^2(\re^2,\re)$.

\end{theorem}

\begin{proof}
Consider a $C^{\infty}$-function $\theta:\re \to  [0,1]$ such that $\theta(\tau)=1$ for $|\tau|\leq 1$ and $\theta(\tau)=0$ for $|\tau|\geq 2$. For every $\de\in [0,1]$, we define 

$$H^{\de}(s,t)=H\Bigl(\te(\de s)s,\, \te(\de t)t\Bigr),$$

\medskip
\noindent
so that
		
\begin{equation}\label{secondaNontrivialSolutions}
		\begin{array}{l}
		H_s^\de(s,t)=H_s\Bigl(
			\te(\de s)s,\, \te(\de t)t
			\Bigr)
			\bigl(\te'(\de s)\de s+\te(\de s)
			\bigr)
			\\
			\ 
			\\
			H_t^\de(s,t)=H_s\Bigl(
			\te(\de s)s,\, \te(\de t)t
			\Bigr)
			\bigl(\te'(\de t)\de t+\te(\de t)
			\bigr).
		\end{array}
\end{equation} 

\medskip
\noindent
As the function $\tau \mapsto |\te'(\tau)\tau +\te(\tau)|$ is bounded in $\re$, by $\eqref{primaNontrivialSolutions}$ and  $\eqref{secondaNontrivialSolutions}$ there is a  constant $C_0>0$, independent from $\de$, such that
	
\begin{equation}\label{terza}
	\begin{array}{l}
	| H_s^\de(s,t)|\leq C_0
	\left(|s|^{p'-1}+|t|^{q'\frac{p'-1}{p'}}+1\right),\\
	\ \\
	|  H_t^\de(s,t)| \leq C_0
	\left(|s|^{p'\frac{q'-1}{q'}}+|t|^{q'-1}+ 1 \right).
	\end{array}
\end{equation}

\medskip
\noindent
Hence the functional $I_{a,H^\de} :X\to \re$ is $C^1$ for any $a \geq 0$ and $\de \in [0,1]$, 
and for any $z_0=(u_0,v_0) \in X,$ $z=(u, v) \in X$ it results

\begin{align*}
\langle I'_{a,H^\de}(z_0), z \rangle = & \displaystyle  \into (a+|\nabla
	u_0|^2)^{\frac {p-2}{2}}\nabla u_0\nabla u  \ dx +\into
	(a+|\nabla v_0|^2)^{\frac {q-2}{2}}\nabla v_0 \nabla v \ dx \\
	&-\displaystyle\into (  H_s^\de(u_0,v_0) u +  H_t^\de(u_0,v_0) v)
	dx.
\end{align*}

\medskip
\noindent
Let $r>0$ be such that $\mathbf{0}$ is the unique critical point of $I_{\a,H^0}=I_{\a,H}$ in

\[	D_r =\left\{z=(u,v)\in X:\,\,\|\nabla u\|_p + \|\nabla v\|_q\leq r\right\}. \]

\medskip
\noindent
The map $\de \mapsto I_{a,H^\de}$ is continuous from $[0,1]$ to $C^1_b(D_r,\re)$. Moreover, reasoning as in the proof of Proposition~\ref{PSlim}, we see that $I_{a,H^\de}$ satisfies the Palais-Smale condition in $D_r$, for any $a \geq 0$ and $\de \in [0,1]$.\\
We claim that there exists $\bar \de \in (0,1]$ such that $\mathbf{0}$ is the unique critical point of $I_{a,H^{\de}}$ in $D_r$, for any $\de \in [0,\bar \de]$.\\
Assume, by contradiction, that $\de_j \to 0$ and $z_j=(u_j,v_j)\ \in \ D_r\setminus \{\mathbf{0}\}$ is a critical point of $I_{\a,H^{\de_j}}$. Taking into account (\ref{terza}), by Theorem \ref{key1} in Chapter \ref{SezReg} (see also \cite[Theorem 1.1]{Vann1}) we infer that every 
	$u_j$ and $v_j$ are in  $L^\infty(\Omega)$ and there is $C_1>0$, depending on $\Omega,\ p, \ q, \ p', \ q'$ and $r$ but independent from $j$, such that
	\[\|u_j\|_\infty, \ \|v_j\|_\infty\leq C_1,
	\] 
	thus $|\de_j u_j(x))|< 1$ and $ |\de_j v_j(x))|< 1$ in $\Omega$, when $j$ is big enough. \newline
	Therefore, by $\eqref{secondaNontrivialSolutions}$, $z_j=(u_j,v_j)$ is a critical point of $I_{a,H}$ and a contradiction
	follows.
	\par
	Setting $\overline {H\,}=H_{\bar \de}$, from~\cite[Theorem~5.2]{CH} we deduce 
	that	$C_m(I_{a,H},\,\mathbf{0})=C_m(I_{a,\overline {H\,}},\,\mathbf{0})$ (for related results, see also \cite[Theorem~I.5.6]{chang}) and the assertion follows.

\end{proof}

\bigskip
\noindent
Assume that $H\in C^2(\re^2,\re)$ and there exist $C>0$, $p'\in (p,p^*),\  q'\in (q,q^*)$ such that

\begin{equation}\label{sesta}
\begin{array}
c| H_{ss}(s,t)|\leq C
\left(|s|^{p'-2}+|t|^{q'\frac{p'-2}{p'}}+1\right),\\
\\
| H_{st}(s,t)|, | H_{ts}(s,t)|\leq C
\left(|s|^{p'-1-\frac{p'}{q'}}+|t|^{q'-1-\frac{q'}{p'}}+ 1 \right), \\
\\
| H_{tt}(s,t)|\leq C
\left(|s|^{p'\frac{q'-2}{q'}}+|t|^{q'-2}+1\right).
\end{array}
\end{equation}

\medskip
\noindent
The functional $I_{a,H}$ defined in \eqref{IH} is of class $C^2$ and  for any $z_0=(u_0,v_0), \ z_1=(u_1,v_1), \ z_2=(u_2, v_2) \in X$, we have

\begin{align*}
\langle I_{a,H}''(z_0) z_1, z_2 \rangle 
& = \int_\Omega
\Bigl((a+|\nabla u_0|^2)^{\frac{p-2}{2}} (\nabla u_1 |\nabla u_2) + (p-2) (a+|\nabla u_0|^2)^{\frac{p-4}{2}} (\nabla u_0|\nabla u_1) (\nabla u_0| \nabla u_2)\Bigr) \, dx \\
& + \int_\Omega \Bigl((a+|\nabla v_0|^2)^{\frac{q-2}{2}} (\nabla v_1 |\nabla v_2) + (q-2) (a+|\nabla v_0|^2)^{\frac{q-4}{2}} (\nabla v_0|\nabla v_1) (\nabla v_0 | \nabla v_2)\Bigr) \, dx \\
& - \int_\Omega \bigl( H_{ss}(u_0,v_0) u_1 u_2 + H_{tt}(u_0,v_0) v_1 v_2 +  H_{st}(u_0,v_0) u_1  v_2 +
H_{ts}(u_0,v_0) u_2  v_1 \bigr) \ dx.
\end{align*}

\bigskip
\noindent
We recall the following notion.

\begin{definition}
	If $Y$ is a Banach space, $I\in C^2(Y,\re)$ and $z_0$ a critical point of $I$,  the Morse index $m(I,z_0)$ of $I$ at $z_0$ is the supremum of the dimensions of the subspaces of $Y$ where $I''(z_0)$ is negative definite. The large
	Morse index $m^*(I, z_0)$ is the sum of $m(I,z_0)$
	and the dimension of the kernel of $I''(z_0)$.
\end{definition}

\bigskip
Arguing as in \cite[Theorem 1.4]{CCMV} we can establish the following 
critical groups estimates for the $C^2$-functionals $I_a$ associated to systems of  $(p,q)$-area equations. Precisely, we have

\begin{theorem}\label{critgroupcomp}
	Let $a >0$ and $H\in C^2(\re^2,\re)$ satisfying assumptions (\ref{sesta}). If $z_0$ is a critical point of the functional
	$I_{a,H}$, then $m(I_{a,H},z_0)$ and
	$m^*(I_{a,H},z_0)$ are finite and
	\[
	C_m(I_{a,H},z_0) = \{0 \}
	\]
	whenever $m < m(I_{a,H},z_0)$ or
	$m > m^*(I_{a,H},z_0)$.
\end{theorem}

\bigskip

In the following theorem we refer to $m_0$ and $m^*_0$ which have been defined through (\ref{m0}).

\begin{theorem}
	\label{thm:general}
	Let  $a \geq 0$ and $p,\,q \in[2,N)$. Let $G \in C^{1}(\re^2, \re)$ a function that satisfies the conditions $(a_1)$ and $(a_2)$.
	If $\mathbf{0}$ is an isolated critical point of the functional $J_a$ defined in $\eqref{functionalNontrivialsolutions}$, then 
	
	\[
	C_m(J_a, \mathbf{0}) = \{0\}\qquad\text{
		whenever $m < m_0$ or $m > m^*_0$}\,.
	\]
	
\end{theorem}

\begin{proof}
\par\noindent 
From (\ref{IH}), $J_a=I_{a,G}$. Due to assumption $(a_1)$ and Theorem~\ref{Hs}, there is $\overline {G\,} \in C^1(\re^2,\re)$ and
$\bar \eta >0$ such that

\begin{equation}\label{e}
	|s|,\,|t| \leq \bar{\eta}\quad \Rightarrow \quad	\overline {G\,}(s,t)=G(s,t);
\end{equation}

\begin{equation}\label{et}
	|s|,\,|t| \geq 2\bar{\eta}\quad \Rightarrow \quad	\overline {G\,}(s,t)=G(0,0);
\end{equation}

\begin{equation}\label{eta}
	C_m(J_{a},\,\mathbf{0})=C_m(I_{a,\overline {G\,}},\,\mathbf{0}), \quad \text{for any}\ m\in \na.
\end{equation}

\bigskip
\noindent
Taking into account assumption $(a_2)$, where $U$ is an open set containing $[-2\bar{\eta},2\bar{\eta}]\times \nobreak[-2\bar{\eta},2\bar{\eta}]$, Theorem~\ref{Hs} assures that $\overline {G\,} \in C^2(\re^2,\re)$. Furthermore by (\ref{et}) $\overline {G\,}$ satisfies (\ref{sesta}), so that $I_{a,\overline{G\,}}\in C^2(X,\re)$. Moreover (\ref{e}) and (\ref{m0}) imply $\langle I_{a, \overline{G\,}}''(\mathbf{0})z,z\rangle=Q(z)$ for any $z$ in $X$, thus
 $m(I_{a,\overline {G\,}},\,\mathbf{0})=m_0$ and $m^*(I_{a,\overline {G\,}},\,\mathbf{0})=m^*_0$.
 
\medskip
\noindent
If ${a>0}$, then Theorem~\ref{critgroupcomp} holds, thus $C_m(I_{a,\overline {G\,}},\,\mathbf{0})=0$ whenever $m<m(I_{a,\overline {G\,}},\,\mathbf{0})$ or $m>m^*(I_{a,\overline {G\,}},\,\mathbf{0})$, so by (\ref{eta}) the assertion follows.

\medskip
\noindent
If instead $a=0$, we have to make a further distinction.

\smallskip

If $H_G(0,0)$ is null or negative semidefinite we have $m_0=0$, $m_0^*=+\infty$
and the assertion is obvious.

\smallskip

If $H_G(0,0)$ is negative definite, we have $m_0=m_0^*=0$ and there is $\mu>0$ such that
\begin{equation}\label{negdef}
	H_G(0,0)\,[\xi]^2\leq -\mu |\xi|^2 \quad \text{for any }
	\xi \in \re^2.
\end{equation}

Moreover, due to assumption $(a_1)$, there are $p'>p,\ q'>q,\ C>0$ such that $W^{1,p}_0(\Omega)\subset L^{p'}(\Omega)$, $W^{1,q}_0(\Omega)\subset L^{q'}(\Omega)$ and 

\[|G(s,t)|\leq C\bigl(|s|^{p'}+|t|^{q'}+1\bigr).
\]

\noindent
Combining with \eqref{negdef} and redefining $C>0$, we get
\[-G(s,t)\geq -G(0,0)-C\bigl(|s|^{p'}+|t|^{q'}\bigr) \qquad \text{for any }
(s,t) \in \re^2,
\]

\noindent
thus

\[J_0(u,v)-J_0(\mathbf{0})\geq \frac 1 p\, {\|u\|_p}^p+\frac 1 q\, {\|v\|_q}^q -C\bigl({|u|_{p'}}^{p'}+{|v|_{q'}}^{q'}\bigr)
\]
\[\geq {\|u\|_p}^p\left(\frac 1 p -C{\|u\|_{p}}^{p'-p}\right)+{\|v\|_q}^q\left(\frac 1 q -C{\|v\|_{q}}^{q'-q}\right)\geq 0
\]

\noindent
if $\|(u,v)\|$ is small enough.

Hence $\mathbf{0}$ is a local minimum for $J_0$ and,
by the excision property, we have

\[
C_m(J_0,\mathbf{0}) \approx
H^m \left(\{0\},\emptyset\right),
\]

\noindent
so the assertion follows.

\smallskip

If $H_G(0,0)$ is positive semidefinite or indefinite, we have $m_0=m_0^*=+\infty$.

As $I_{a,\overline{G\,}}\in C^2(X,\re)$, from \cite[Theorem~3.1]{6} we infer that
$C_{m}(I_{a,\overline {G\,}},\,\mathbf{0})=\{0\}$ for any $m$, so by (\ref{eta})
the assertion follows.

\end{proof}

{\vskip 1cm} 

\section{Geometry of $J_a$ and Proof of the main results}
\label{sect:z3}

Let us recall that  $F:\re^2\to \re$  is the $C^1$-function defined by 

\[F(s,t)=\frac 1 p |s|^p+\frac 1 q |t|^q +\frac{1}{(\alpha+1)(\beta+1)}|s|^{\alpha} |t|^{\beta} s t. 
\]

\bigskip
\noindent
Assuming that $G(0,0)=0$, from assumption $(a_1)$ we see that 

\[G(s,t)=\bar{\lambda}F(s,t)+R(s,t),
\]

\medskip
\noindent
where $R\in C^1(\re^2,\re)$, $R(0,0)=0$ and 
 
\begin{equation}\label{eq rf}
 	\lim_{|(s,t)|\to \infty}\frac{R_s(s,t)}{|s|^{p-1}+ |t|^{q \frac{p-1}p}}=0, \qquad 	\lim_{|(s,t)|\to \infty}\frac{R_t(s,t)}{|s|^{p \frac{q-1}q}+|t|^{q-1}}=0.
\end{equation}

\begin{theorem}\label{thm rf}
 	Following the previous notations, we have 
 	
 	\[\lim_{|(s,t)|\to \infty}\frac{R(s,t)}{F(s,t)}=0.
 	\]
\end{theorem}
 
\begin{proof}
 	From Young's inequality and  (\ref{disYo}), there exists  $c>0$ be such that
 	\begin{equation*}
 	\frac 1 p |s|^p+|s||t|^{q \frac{p-1}p}+\frac 1 q |t|^q	\leq c F(s,t)\qquad \forall (s,t)\in \re^2.
 	\end{equation*}
 	By (\ref{eq rf}), for every $\varepsilon >0$ there are $\eta_\varepsilon, \, c_\varepsilon, \, \delta_\varepsilon>0$ such that
 	\begin{eqnarray*}
 	|(s,t)|>\eta_\varepsilon \ \Rightarrow  |R_s(s,t)|\leq 
 		\frac{\varepsilon}{2c}\left(|s|^{p-1}+ |t|^{q \frac{p-1}p}\right), \\ 
 |R_t(s,t)|\leq 
 		\frac{\varepsilon}{2c}\left(|s|^{p \frac{q-1}q}+|t|^{q-1}\right);\\
 		\\
 		|(s,t)|\leq\eta_\varepsilon \Rightarrow   |R_s(s,t)|\leq c_\varepsilon, \quad |R_t(s,t)|\leq c_\varepsilon;\\
 		\\
 		|(s,t)|>\delta_\varepsilon\ \Rightarrow  {\displaystyle 
 			F(s,t)>
 			\frac{4 c_\varepsilon \eta_\varepsilon}{\varepsilon }}\, .
 	\end{eqnarray*}
 	 	For any $r\in\re$ let us denote by $r^+=\max\{0,r\}$ \ and $r^-=\min\{0,r\}$.	 	
 	As 
 	 	\[R(s,t)=\int_{0}^{s}R_s(\sigma,t)d\sigma + \int_{0}^{t}R_t(0,\tau)d\tau, 
 	\]
 	 we infer that
 	\[
 	|R(s,t)|\leq\int_{s^-}^{s^+}|R_s(\sigma,t)|d\sigma + \int_{t^-}^{t^+}|R_t(0,\tau)|d\tau\]
 	\[
 	=
 	\int_{ 
 		\{\sigma \in [s^-,s^+] \ :\ |(\sigma,t)|\leq \eta_\varepsilon\}}
 	\hspace{-20mm}|R_s(\sigma,t)|d\sigma \ + \int_{ \{\sigma \in [s^-,s^+]  \ :\ 
 		|(\sigma,t)|> \eta_\varepsilon\}}\hspace{-20mm}|R_s(\sigma,t)|d\sigma\ \]
 	\[
 	+\int_{ 
 		\{\tau \in [t^-,t^+] \ :\ |\tau|\leq \eta_\varepsilon\}}
 	\hspace{-14mm}|R_t(0,\tau)|d\tau + \int_{ 
 		\{\tau \in [t^-,t^+] \ :\ |\tau|> \eta_\varepsilon\}}
 	\hspace{-12mm}|R_t(0,\tau)|d\tau\]	
 		\[ 
 	\leq  2 c_\varepsilon \eta_\varepsilon + \frac \varepsilon 2 F(s,t)
 	\]
 	therefore
 	\[|(s,t)|>\delta_\varepsilon\ \Rightarrow \quad \frac{|R(s,t)|}{F(s,t)}
 	< \varepsilon.
 	\]	
 	
 \end{proof}
 	
 	\bigskip

\begin{proposition}\label{c1}
If $(a_1)$ holds and $(u,v)$ is a weak solution to problem \eqref{pqNontrivialSolutions}, then 
 	
\begin{itemize}
 	\item[(i)]  $(u,v) \in \left(C^{1,\eta}(\overline{\Omega})\right)^2$ for some $\eta \in (0,1)$;
 	\item[(ii)] for every bounded set $A\subset X$, there exists $\eta \in (0,1)$ and $M>0$ such that $\|u\|_{C^{1,\eta}(\overline{\Omega})}\leq M$ and $\|v\|_{C^{1,\eta}(\overline{\Omega})}\leq M$ for every $(u,v)\in A$ solving \eqref{pqNontrivialSolutions};
 	\item[(iii)] if there is $\varepsilon >0$ such that
 	
\begin{align*}\label{a3}
 		 G_s(0,t)\neq 0 \qquad \forall \, t \in (-\varepsilon,\varepsilon)\setminus \{0\}\nonumber
 		\\
 	(*)	\hspace{6.5cm}\\
 		G_t(s,0)\neq 0 \qquad \forall \, s \in (-\varepsilon,\varepsilon)\setminus \{0\},\nonumber
\end{align*}

 	then every nontrivial solution $(u,v)$ of \eqref{pqNontrivialSolutions} is also not semitrivial, i.e. $u\neq 0$ and $v\neq 0$.
\end{itemize}
 
\end{proposition}
 	
\begin{proof}
 	If  $(a_1)$ holds and $A$ is a  bounded subset of $X$, by Theorem 1.1 in \cite{Vann1} all the weak solutions $(u,v)\in A$ are in $\left(L^\infty(\Omega)\right)^2$ and there is $M_0>0$ such that
 		$\|u\|_\infty \leq M_0$ and $\|v\|_\infty \leq M_0$.
 		Hence, due to \cite{lieberman}, we infer (i) and (ii).
 		\newline Let us prove (iii) arguing by contradiction. Suppose that $(a_1)$ and $(*)$ hold and that there is an $(u,0)\neq (0,0)$ solving \eqref{pqNontrivialSolutions}.
 		By (i) we get that $u\in C(\overline{\Omega})\setminus \{0\}$, $M=\max_{\overline{\Omega}} |u|>0$ and there is $\bar x \in \Omega$ such that
 		$0<|u(\bar x)|<\min \{\varepsilon, M\}$. Hence, by $(*)$,  $G_t(u(\bar x),0)\neq 0$, in contradiction with \eqref{pqNontrivialSolutions}. 
 		\end{proof} 
 	
 	\bigskip
 	
\begin{lemma}\label{lem 1}
If $r\geq 2$, $a \geq 0$ and $\varepsilon >0$, there is $c(r,a,  \varepsilon)>0$ such that

\[ \frac{1}{r} (a +t^2)^{\frac{r}{2}} \leq \frac 1 r (1+\varepsilon)|t|^r+ c(r,a, \varepsilon) \]

\medskip
\noindent
for any $t \in \re$.
\end{lemma}

\begin{proof}
 The assertion is trivial if $r=2$, so let us consider the case $r>2$.
 		\newline Setting $f_a(t)=\frac 1 r(a +t^2)^{\frac r 2}-\frac 1 r (1+\varepsilon)|t|^r$, we see that  
 ${\displaystyle \lim_{|t| \to +\infty}}f_a(t)=-\infty$. 
 		
Moreover $f'_a(t)=t|t|^{r-2} \left((\frac{a}{t^2}+1)^{\frac{r-2}2}-1-\varepsilon\right)$ and, denoting by $t_a=\left(\frac{a}{(1+\varepsilon)^{\frac{2}{r-2}}-1}\right)^{1/2}$\!\!\!\!, we get

\[ 
\max_{t\in \re} f_a(t)= f_a(\pm t_a)=
\frac{a^{r/2}}{r}  
\frac{1+\varepsilon}{\left((1+\varepsilon)^{\frac{2}{r-2}}-1\right)^{\frac{r-2}2}}.
\]

\end{proof}

\medskip
\noindent
Setting

	\[L(s,t)=\frac 1 p|s|^p+\frac 1 q|t|^q \qquad 
	L_a(s,t)=\frac{1}{p} (a +s^2)^{\frac{p}{2}} + \frac{1}{q} (a +t^2)^{\frac{q}{2}}
	\] 

\medskip
\noindent
and recalling Theorem \ref{thm rf} we obtain the following result.

\begin{corollary}\label{cork}
For any $\varepsilon >0$ there is $c_\varepsilon >0$ such that
\begin{equation}\label{eq la}
L(s,t)\leq L_a(s,t)\leq (1+\varepsilon) L(s,t)+c_\varepsilon	
\end{equation}
\begin{equation}\label{eq fg}
-(\bar{\lambda}+\varepsilon)F(s,t)-c_\varepsilon\leq - G(s,t)\leq 
-(\bar{\lambda}-\varepsilon)F(s,t)+c_\varepsilon
\end{equation}
for any $(s,t) \in \re^2$ .
\end{corollary}

\begin{lemma}\label{lem-infty}
Let $\Gamma:\re^2\to \re$ be a $C^1$ function such that

\begin{equation}\label{0key}
	\lim_{|(s,t)|\to \infty}\frac{\Gamma(s,t)}{F(s,t)}=0
\end{equation}

\medskip
\noindent
and

\begin{equation}\label{eq-infty}
	\lim_{|(s,t)|\to \infty} \Gamma(s,t)-\frac 1 p\Gamma_s(s,t)s-\frac 1 q\Gamma_t(s,t)t=-\infty.
\end{equation}

\medskip
\noindent
Then

\begin{equation}\label{eq2-infty}
\lim_{|(s,t)|\to \infty} \Gamma(s,t)=-\infty.
\end{equation}

\medskip
\noindent
In the same way, if $ \ \Gamma$ satisfies condition (\ref{0key}) and

\begin{equation}\label{eq-infty+}
\lim_{|(s,t)|\to \infty} \Gamma(s,t)-\frac{1}{p} \Gamma_s(s,t)s-\frac{1}{q} \Gamma_t(s,t)t=+\infty,
\end{equation}

\medskip
\noindent
then

\begin{equation}\label{eq2-infty+}
\lim_{|(s,t)|\to \infty} \Gamma(s,t)=+\infty.
\end{equation}

\end{lemma}

\begin{proof}
	Let observe  $F(x^{1/p}s,x^{1/q}t)=xF(s,t)$, for any $x\geq 0$ and $(s,t)\in \re^2$.
	
	Hence if $(s,t)\in \re^2\setminus \{(0,0)\}$ then, by (\ref{0key}),
	\begin{equation}\label{eq0}
	\lim_{x\to +\infty}\frac{\Gamma(x^{1/p}s,x^{1/q}t)}{x}=0.
	\end{equation}
Let us fix $M>0$. By (\ref{eq-infty}) there exists $\delta >0$ such that
	\newline$\Gamma(s,t)-\frac 1 p\Gamma_s(s,t)s-\frac 1 q\Gamma_t(s,t)t<-M$ 
	\ when $|(s,t)|>\delta$.
	
	\bigskip
	
	If $(s',t')\in \re^2$ is such that $|(s',t')|>\delta$, then
	
	 \begin{align*}
	 	&\frac{d}{dx}\frac{\Gamma(x^{\frac 1 p}s',x^{\frac 1 q}t')+M}{x}\\
	 	=&\frac {\frac 1 p\Gamma_s(x^{\frac 1 p}s',x^{\frac 1 q}t')x^{\frac 1 p}s'+ \frac 1 q\Gamma_t(x^{\frac 1 p}s',x^{\frac 1 q}t')x^{\frac 1 q}t'-\Gamma(x^{\frac 1 p}s',x^{\frac 1 q}t') -M
	 	}{x^2}>0 \qquad {\rm if\ } x\geq 1
	 \end{align*}
	 
	so that
	\[\frac{\Gamma(x^{\frac 1 p}s',x^{\frac 1 q}t')+M}{x}\leq \frac{\Gamma(\bar x^{\frac 1 p}s',\bar x^{\frac 1 q}t')+M}{\bar x} \qquad {\rm whenever}\ 1\leq x\leq \bar x.
	\]
	By (\ref{eq0}), passing to the limit for $\bar x \to +\infty$, we get
	\[\frac{\Gamma(x^{\frac 1 p}s',x^{\frac 1 q}t')+M}{x}\leq 0 \qquad \forall \ x\geq 1
	\]
	hence \ $\Gamma(s',t')\leq -M$, which proves (\ref{eq2-infty}).
	
	\bigskip
	
	Finally (\ref{eq2-infty+}) can  be proved applying  (\ref{eq2-infty}) to the function $\bar \Gamma=-\Gamma$.
	\end{proof}

\medskip
\noindent
Now, we can prove Theorem \ref{key1}, Theorem \ref{key2} and Theorem \ref{key3}.

\bigskip
\noindent
\emph{Proof of Theorem~\ref{key1}.}

As $\bar\lambda$ is not an eigenvalue for system $(\ref{eigen})$, let $m_{\infty}\in \na$ be such that  

$\lambda_{m_{\infty}} < \bar\lambda < \lambda_{m_{\infty}+1}$ and assume 
$m_{\infty} \not\in [m_0,m_0^*]$.

Setting
\begin{gather*}
	X_- 
	=\{ z  \in X  \ |    \ \Phi(z) \leq \lambda_{m_{\infty}}
		 \Psi(z)\}\\
		 X_+=\{ z \in X  \ |   \ \Phi(z) \geq \lambda_{m_{\infty}+1} \Psi(z) \}\\
		  D_r=\{z\in 	X_- \ |\ \Phi(z)\leq r\} \qquad  S_r=\{z\in 	X_-\ |\ \Phi(z)= r\},	
\end{gather*}

Theorem \ref{th dr} assures that
	$(D_r,S_r)$ links $X_+$ cohomologically in dimension $m_{\infty}$ over $\ze_2$, for any $r>0$.
	
	\medskip
	
	We take into consideration
 only the case $m_{\infty}\geq 1$, as when $m_{\infty}=0$ the proof is similiar and simpler.

	Let $\alpha',\alpha'', \beta'$ be such that
	\[\lambda_{m_{\infty}} < \alpha'<\alpha''< \bar\lambda<\beta'
	< \lambda_{m_{\infty}+1}.
	\]
	
	By Corollary \ref{cork}, there is $C>0$ such that 
	\begin{gather*}
		-\beta'F(s,t)-C\leq -G(s,t)\leq -\alpha''F(s,t)+C\\
		L(s,t)\leq L_a(s,t)\leq \frac{\alpha'}{\lambda_{m_{\infty}}} L(s,t)+C.
	\end{gather*}
	
	Therefore, for any $(u,v)\in X_+$,
	\begin{eqnarray*}
		&J_a(u,v)=\int_{\Omega}L_a(|\nabla u(x)|,|\nabla v(x)|)dx-\int_{\Omega}G(u(x),v(x))dx\\
		&\geq  \int_{\Omega}L(|\nabla u(x)|,|\nabla v(x)|)dx
		-\int_{\Omega}\left(\beta'F(u(x),v(x)) +C\right) dx\\
		&\geq \left( \lambda_{m_{\infty}+1}-\beta'\right)\int_{\Omega}F(u(x),v(x))dx-C|\Omega|
		\\
		& 
		\geq -C|\Omega|
	\end{eqnarray*}
	so that
	
	\[\inf_{ z\in X_+} J_a(z)=\bar{a}>-\infty. \]

 On the other hand, for any $(u,v)\in X_-$,
\begin{eqnarray*}
	&J_a(u,v)=\int_{\Omega}L_a(|\nabla u(x)|,|\nabla v(x)|)dx-\int_{\Omega}G(u(x),v(x))dx\\
	&\leq \frac{\alpha'}{\lambda_{m_{\infty}}} \int_{\Omega}L(|\nabla u(x)|,|\nabla v(x)|)dx  + C|\Omega|
	-\alpha''\int_{\Omega}F(u(x),v(x))dx  + C|\Omega| \\
	&\leq  2C|\Omega| +\left(\frac{\alpha'-\alpha''}{\lambda_{m_{\infty}}}
	\right) \int_{\Omega}L(|\nabla u(x)|,|\nabla v(x)|)dx
\end{eqnarray*}

so that 

\[\lim_{\substack{\|z\|\to \infty\\ z\in X_-}} J_a(z)=-\infty\quad \Rightarrow \quad \exists \, 
r>0\ :
\ \sup_{S_r} J_a < \bar{a}=\inf_{X_+} J_a.
\]

\medskip
\noindent
Moreover, as $D_r$ is compact, $\bar{b}={\displaystyle \max_{D_r}}\,\ja<+\infty$.

If  $J_a^{-1}([\bar{a},\bar{b}])$ contains an infinite number of critical points, the theorem is obviously proved.
Otherwise, applying Theorem \ref{th key} combined with Proposition \ref{PS}, 
there exists a critical point $z_0$ of $J_a$ with
$\bar{a} \leq J_a(z_0)\leq \bar{b}$ and 

\[
C_{m_{\infty}}(J_a,z_0)\neq \{0\}.	
\]

\medskip
\noindent
As $m_{\infty} \not\in [m_0,m_0^*]$, Theorem~\ref{thm:general} gives that

\[C_{m_{\infty}}(J_{a},\,\mathbf{0})
= \{0\},
\] 

\medskip
\noindent
hence $z_0\neq \mathbf{0}$.

\qed

\medskip
\noindent
\emph{Proof of Theorem~\ref{key2}.}
	As \[
	\lim_{|(s,t)|\to  \infty} \left[ G(s,t) - \frac{1}{p} G_s(s,t)s - \frac{1}{q} G_t(s,t)t\right] = -\infty\, 
	\]
	and 
	$F(s,t) - \frac{1}{p} F_s(s,t)s - \frac{1}{q}F_t(s,t)t=0$ for any $(s,t)\in \re^2$,
	by Theorem \ref{thm rf} and Lemma \ref{lem-infty} we get
	\[\lim_{|(s,t)|\to  \infty}  G(s,t) - \bar \lambda  F(s,t)= -\infty,\
	\]
	hence $G(s,t) - \bar \lambda  F(s,t)$ is bounded from above.
	
	Therefore, recalling Corollary \ref{cork} and letting $\alpha',\ \alpha''$ be such that
	\[\lambda_{m_\infty}<\alpha'<\alpha''<\bar \lambda\leq \lambda_{m_\infty+1},
	\]
	there is $C>0$ such that 
	\begin{gather*}
	-\bar \lambda F(s,t)-C\leq -G(s,t)\leq -\alpha''F(s,t)+C\\
	L(s,t)\leq L_a(s,t)\leq \frac{\alpha'}{\lambda_{m_{\infty}}} L(s,t)+C.
	\end{gather*}

	So we can conclude as in the proof of Theorem \ref{key1}, defining $X_-$ and $X_+$ in the same way. 
	\qed
	
	\bigskip
	
\bigskip

\noindent\emph{Proof of Theorem~\ref{key3}.}

Firstly, reasoning as in the previous proof, by assumption $(b_+)$ we get

\begin{equation}\label{b+}
	\lim_{|(s,t)|\to  \infty}  G(s,t) - \bar \lambda  F(s,t)= +\infty.
\end{equation}

Let $\beta'\in (\bar \lambda, \lambda_{m_{\infty+1}})$. By Corollary~\ref{cork} there is  $C>0$ such that
 
\[	-\beta'F(s,t)-C\leq -G(s,t).
\]

We define $X_-$ and $X_+$ as in the proof of Theorem \ref{key1} and, reasoning in the same way, we see that

\[\inf_{ z\in X_+}J_0(z)=\bar{a} >-\infty.
\]

So it suffices to prove that, for any $M>0$, we can choose a suitably big $r>0$ such that

\[\sup_{ z\in S_r}J_0(z)<-M. \]

\bigskip
\noindent
By contradiction, there is $M>0$ and a sequence $\{z_n\}_n=\{(u_n,v_n)\}_n$ in $X_-$ such that $\Phi(z_n)= n$ and

\begin{equation}\label{M}
	J_0(z_n)\geq -M.
\end{equation}

As $z_n\in X_-$ we see that

	\begin{equation}\label{ast}
	J_0(z_n)=\Phi(z_n)-\int_{\Omega}G(z_n(x))\, dx
		\leq
		\int_{\Omega}\bar \lambda F(z_n(x))-G(z_n(x))\, dx.
	\end{equation}

	Let \quad $\delta_n=1/\Phi(z_n)$, $\bu_n=\delta_n^{1/p}u_n$, $\bv_n =\delta_n^{1/q}v_n$,
	$\bz_n =(\bu_n,\bv_n)$.
	
\medskip
\noindent
As $\Phi(\bz_n)=1$, up to subsequence  $\{ \bz_n \}_n$ converges to some $\bz=(\bu,\bv)$ weakly in $X$ and strongly in $L^p(\Omega)\times L^q(\Omega)$.
	
	It is easy to see that every $\bz_n \in X_-$, as well as $\bz \in X_-$.
	
	Moreover
	
\[\Psi(\bz)=\lim_{n\to \infty}\Psi(\bz_n)\geq  \lim_{n\to \infty}\frac{1}{\lambda_{m_{\infty}}}\ \Phi(\bz_n)=\frac{1}{\ \lambda_{m_{\infty}}}>0,\]
	
	hence
	
\begin{equation}\label{neq0}
		\bz \neq 0.
\end{equation}
	
Since
	\[|u_n(x)|^p+|v_n(x)|^q=\Phi(z_n)\left(|\delta_n^{1/p}u_n(x)|^p+|\delta_n^{1/q}v_n(x)|^q\right)
	\geq n \left(|\bu_n(x)|^p+|\bv_n(x)|^q\right)\]
	
	and 
	
	\[|\bu_n (x)|^p+|\bv_n(x)|^q \longrightarrow |\bu(x)|^p+|\bv(x)|^q \qquad \hbox {a.e. in }
	 \Omega,\]
	 
\medskip
\noindent
taking into account (\ref{b+}), (\ref{ast})  and (\ref{neq0}),  up to subsequence we get
	 
	 \[\lim_{n\to \infty}J_0(z_n)=-\infty,
	 \]
	 
\medskip
\noindent
which contradicts (\ref{M}).

 \qed \\

\bigskip

\chapter[On the de Thélin eigenvalue problem and Landesman-Lazer conditions]{On the de Thélin eigenvalue problem and Landesman-Lazer conditions for quasilinear systems}\label{CapitolodeThelin}

\bigskip
\noindent
In this Chapter, we investigate the eigenvalue problem for quasilinear elliptic systems originally introduced by de Thélin in \cite{DT}. We prove that the first eigenvalue $\lambda_1$ is not only simple (in a suitable sense), but also isolated. Moreover, we provide a variational characterization for a sequence of eigenvalues $\{\lambda_k\}_k$, employing a deformation lemma for $C^1$-submanifolds established by Bonnet in \cite{BON}. Furthermore, in the spirit of the work of Landesman and Lazer \cite{LANDLAZ}, we prove the existence of a weak solution for a quasilinear elliptic system at resonance with $\lambda_1$. We introduce new sufficient Landesman-Lazer type conditions, extending the results obtained for quasilinear elliptic equations by Arcoya and Orsina in \cite{AO}. These conditions ensure that the Euler functional satisfies the PS condition and, under specific assumptions, is coercive. In the other cases, since the set of eigenfunctions associated with $\lambda_1$ for quasilinear systems is not necessarily a vector subspace of $X := W_0^{1,p}(\Omega) \times W_0^{1,q}(\Omega)$, the standard splitting argument used in \cite{AO} is not applicable. To overcome this difficulty, we exploit the isolation of $\lambda_1$ and the variational characterization of $\lambda_2$ to obtain a solution.

\medskip
This Chapter is based on the paper \cite{ABC}.

\bigskip
\section{de Thélin eigenvalue problem}

\medskip
\noindent
In this Chapter, we study the following eigenvalue problem for quasilinear elliptic systems, originally introduced by de Thélin \cite{DT}:

\begin{equation}\label{EigenDeThélin0}
	\begin{cases}
		\begin{array}{ll}
			- \Delta_p u = \lambda |u|^{\alpha-1} |v|^{\beta+1} u
			&  \text{ in }\Omega,
			\medskip \\
			- \Delta_q v  = \lambda  |u|^{\alpha+1} |v|^{\beta - 1}v 
			&  \text{ in } \Omega, \medskip \\
			u=v=0  & \text{ on }  \partial\Omega,
		\end{array}
	\end{cases}
\end{equation}

\medskip
\noindent
where  $\lambda \in \R$ and $p, q,\alpha, \beta$ are real numbers satisfying:

\begin{equation}\label{condizionepqalphabeta}
	\displaystyle  1 < p < N , \quad 1 < q < N, \quad \alpha >0, \quad \beta >0, \quad \text{and}  \quad \frac{\alpha + 1}{p} + \frac{\beta + 1}{q} = 1.
\end{equation}

\medskip
\noindent
A real number $\lambda$ is said to be an eigenvalue of \eqref{EigenDeThélin0} if there exists $(\varphi, \psi) \in X$, $(\varphi, \psi) \neq (0,0)$ such that for any $(u, v) \in X$ we have

\begin{align}\label{significatoautofunzione}
\displaystyle
& (\alpha+1) \int_{\Omega} |\nabla \varphi |^{p-2} \nabla \varphi \cdot \nabla u \, dx 
+ (\beta+1) \int_{\Omega} |\nabla \psi|^{q-2} \nabla \psi \cdot \nabla v \, dx \nonumber \\
= &  (\alpha + 1) \into \left| \varphi \right|^{\alpha-1}  \left| \psi \right|^{\beta+1} \varphi u \, dx 
+ (\beta + 1) \into \left| \varphi \right|^{\alpha+1}  \left| \psi \right|^{\beta-1} \psi v  \, dx.
\end{align}

\medskip
\noindent
The nontrivial weak solution $(\varphi, \psi)$  is then called eigenfunction associated with the eigenvalue $\lambda$.\\
Notice that if $(\varphi,\psi)$ is a nontrivial (i.e. $(\varphi, \psi) \neq (0,0)$) eigenfunction of \eqref{EigenDeThélin0} associated with an eigenvalue $\lambda \in \mathbb{R}$, then it cannot be semitrivial (i.e. $ \varphi \neq 0 \neq \psi$). In fact, if for example $\psi \equiv 0$ and we test \eqref{significatoautofunzione} in $(u,v)=(\varphi,0),$ we get

$$ \displaystyle \into \left| \nabla \varphi \right|^p \, dx = 0,$$

\medskip
\noindent
by which $\varphi \equiv 0$, that is a contradiction with the nontriviality of $(\varphi,\psi).$

\medskip
\noindent
Let us introduce the functionals $\Phi : X \to \re$ and $\Psi:X \to \re$ defined as

\begin{align}\label{funzionalePhi}
\displaystyle
\Phi (u,v) =  \frac{\alpha+1}{p} \int_{\Omega} |\nabla u|^p \, dx + \frac{\beta+1}{q}\int_{\Omega} |\nabla v|^q \, dx
\end{align}

\medskip
\noindent
and

\begin{align}\label{funzionalePsi}
	\Psi(u,v) =  \into \left| u \right|^{\alpha+1}  \left| v \right|^{\beta+1} \, dx
\end{align}

\medskip
\noindent
for any $(u,v) \in X$, and let us denote with $\Sigma$  the  $C^1$ manifold

\begin{align*}
\displaystyle 
\Sigma:=\Biggl\{ z=(u,v) \in X \, : \,  \Psi(z)=1   \Biggr\}.
\end{align*}

In \cite[Theorem 1]{DT} de Thélin proved that problem \eqref{EigenDeThélin0} admits a smallest eigenvalue $\lambda_1>0$ defined as

\begin{equation}\label{definizioneprimoautovalore}
 \displaystyle \lambda_1 := \inf_{z \in \Sigma} \Phi(z).
\end{equation}

In \cite[Theorem 2]{DT}  he also proved, by means of Tolksdorf regularity results \cite{tolksdorf1983} and Vazquez's maximum principle \cite{VAZ},  that there exists a nontrivial eigenfunction  $(\varphi_0,\psi_0)$ associated with $\lambda_1$ such that $\varphi_0 \in C^{1,\eta}(\overline{\Omega})$, $\psi_0 \in C^{1,\xi}(\overline{\Omega})$ and $\varphi_0 > 0,$ $\psi_0 > 0$ in $\Omega$. In addition, by using results of Diaz and Saa \cite{DIAZSAA} and Anane \cite{AL}, in \cite[Theorem 3]{DT} it was proved that if $(\varphi_0,\psi_0)$ is an eigenfunction associated with $\lambda_1$  satisfying $\varphi_0 > 0,$ $\psi_0> 0$ in $\Omega$ and $(\varphi_0,\psi_0) \in \Sigma$, then is unique.\\

In this Chapter we prove that the eigenvalue $\lambda_1$ is isolated, and it is simple in the sense stated below.\\
Firstly, if $\lambda$ is an eigenvalue with a nontrivial eigenfunction $(\varphi, \psi)$, by using the relation $(\alpha+1)/p + (\beta+1)/q=1$ it turns out that 

$$ \displaystyle E:= \bigg\{ ( \left|  \theta \right|^{\frac{1}{p}} \varphi ,  \left|  \theta \right|^{\frac{1}{q}} \psi ) \, \text{sgn} (\theta) \, : \, \theta \in \mathbb{R} \bigg\} \,  \bigcup \, \bigg\{ (- \left|  \theta \right|^{\frac{1}{p}} \varphi ,  \left|  \theta \right|^{\frac{1}{q}} \psi ) \, \text{sgn} (\theta) \, : \, \theta \in \mathbb{R} \bigg\} $$

\noindent
is a set of eigenfunctions associated with $\lambda$. However, if $(\bar{\varphi},\bar{\psi)}$ is another eigenfunction associated with the same eigenvalue $\lambda$, it could  be the case that  $(\bar{\varphi},\bar{\psi)} \not \in E$. \cite[Theorem~3]{DT} proves that this possibility is excluded for the eigenvalue $\lambda_1$ (simplicity of  $\lambda_1$).

\begin{theorem}\label{simplicityINTRODUZIONE}
The set of all eigenfunctions associated with $\lambda_1$ is given by

$$ \displaystyle E_1:= \bigg\{ ( \left|  \theta \right|^{\frac{1}{p}} \varphi_0 ,  \left|  \theta \right|^{\frac{1}{q}} \psi_0 ) \, \text{sgn} (\theta) \, : \, \theta \in \mathbb{R} \bigg\} \,  \bigcup \, \bigg\{ (- \left|  \theta \right|^{\frac{1}{p}} \varphi_0 ,  \left|  \theta \right|^{\frac{1}{q}} \psi_0 ) \, \text{sgn} (\theta) \, : \, \theta \in \mathbb{R} \bigg\} $$

\noindent
where $(\varphi_0,\psi_0)$ is the unique nontrivial eigenfunction satisfying $(\varphi_0,\psi_0) \in \Sigma$,    $\varphi_0 >0$ and $ \psi_0>0$  in $\Omega$.
\end{theorem}

Previous result and Picone's identity (see \cite[Theorem 1]{AH} and \cite[Lemma 24]{AG}), allow us to prove that if $(\varphi,\psi)$ is an eigenfunction of \eqref{EigenDeThélin0} associated with the eigenvalue $\lambda$ with $\lambda > \lambda_1$, then both $\varphi$ and $\psi$ change sign (see Proposition \ref{funzionicambianosegno}).\\
By means of this result, we are able to prove that:

\begin{theorem}\label{lambda1isolatoINTRODUZIONE}
The eigenvalue $\lambda_1$ is isolated, that is there exists $a > \lambda_1$  such that $\lambda_1$ is the only eigenvalue of \eqref{EigenDeThélin0} in  $[0,a]$. 
\end{theorem}

Moreover, inspired by \cite{DrabekRobinson}, we also construct a sequence $\{\lambda_k\}_k$ of eigenvalues of \eqref{EigenDeThélin0} with a suitable variational characterization involving the unit sphere $S^{k-1}$ of $\mathbb{R}^k$ for any $k \geq 1$. Such a sequence will be constructed on the $C^1$ manifold 

\begin{align*}
\displaystyle 
\mathcal{M}:=\Biggl\{ z=(u,v) \in X \, : \,  \Phi(z)=1   \Biggr\}.
\end{align*}

\medskip
\noindent
In particular, setting

$$
\mathcal{M}_k := \{ A=\alpha(S^{k-1}) \subset \mathcal{M}  \ | \,  \alpha: S^{k-1} \to \mathcal{M} \text{ is odd and continuous}  \},
$$

\medskip
\noindent
and

\begin{equation*}
	c_k := \sup_{A \in \mathcal{M}_k} \min_{(u,v) \in A} \Psi(u,v),
\end{equation*}

\medskip
\noindent
we prove in Proposition \ref{prop 2.9} that the sequence $\{\lambda_k\}_k \subset \mathbb{R}$ where

\begin{equation}\label{lambdakINTRODUZIONE}
	\lambda_k := \frac{1}{c_k},
\end{equation}

\medskip
\noindent
is a sequence of eigenvalues of \eqref{EigenDeThélin0}, and $\lambda_1$ is exactly the first eigenvalue defined in \eqref{definizioneprimoautovalore} (see Remark \ref{definizionicoincidono}).

In order to construct a sequence of eigenvalues of \eqref{EigenDeThélin0}, we notice that the manifold $\mathcal{M}$ is not of class $C^{1,1}$ when $ 1 <p < 2$ or $1<q<2$. We will apply a deformation lemma for $C^1$ submanifolds of a Banach space (see Theorem~\ref{BON}) based in \cite[Theorem 2.5]{BON}. This theorem requires a coupled Palais-Smale condition.

It is unknown if the set of the eigenvalues described by \eqref{lambdakINTRODUZIONE} contains all the eigenvalues of problem \eqref{EigenDeThélin0}. Actually, the eigenvalue structure in the context of quasilinear elliptic systems is a challenging open problem. 

Really, also for the single scalar eigenvalue problem $-\Delta_p u = \eta |u|^{p-2}u$, the spectral properties of $-\Delta_p$ are not yet well understood. We mention \cite{AL} and \cite{ANANETSOULI,DrabekRobinson,PERERA} for results concerning the first eigenvalue and the second one, (see also \cite{CD} and \cite{CV3}). Moreover it is also possible to define, in at least three different ways, a diverging sequence $\{\eta_k\}_k$ of eigenvalues (see \cite{CUESTA,DrabekRobinson,PERERA,SZU}), but it is unknown whether these definitions are equivalent for $k \geq 3$ or whether the whole set of eigenvalues is covered.

In the context of quasilinear elliptic systems, we also recall that there are several ways to define an eigenvalue problem (see for instance \cite{SZ} for an other type of eigenvalue problem different from \eqref{EigenDeThélin0}, and also \cite{boccardodefiguerido} for a more general class of eigenvalue problems involving $(p,q)$-homogeneous functions).

Then we consider the eigenvalue problem \eqref{EigenDeThélin0}, by proving that $\lambda_1$ is isolated and simple, and admits the sequence $\{\lambda_k \}_k$ of eigenvalues defined in \eqref{lambdakINTRODUZIONE}, with the idea of proving the same results for other eigenvalue problems in a forthcoming paper.

Here we also concern with sufficient Landesman-Lazer type conditions for a quasilinear elliptic systems in resonance around $\lambda_1$. When talking about Landesman-Lazer conditions, one generally refers to the pioneering work \cite{LANDLAZ}, where Landesman and Lazer considered the boundary value problem 

\begin{equation}\label{LandLaz}
\begin{cases}
- \Delta u = \lambda u + f(u) - h(x)  & x \in \Omega, \medskip\\
u=0 & x \in \partial \Omega,
\end{cases}
\end{equation}

\medskip
\noindent
where $\Omega$ is a bounded domain in $\mathbb{R}^N$, $N \geq 1$, $h$ is a real function in $L^2(\Omega)$, $\lambda$ is a positive number and $f$ is a real valued function which is bounded and continuous on $\mathbb{R}$. Moreover, there exist the limits

$$ \displaystyle f_{- \infty}:=  \lim_{s \to - \infty} f(s)  \qquad \text{and} \qquad f^{+ \infty}:=  \lim_{s \to + \infty} f(s).$$

They studied the existence of solutions of \eqref{LandLaz} in the resonance case, i.e.  they took into account the linear eigenvalue problem 

\begin{equation*}
\begin{cases}
- \Delta u = \lambda u & x \in \Omega, \medskip\\
u=0 & x \in \partial \Omega,
\end{cases}
\end{equation*}

\noindent
and considered the case in which $\lambda$ is a simple eigenvalue, that is the eigenspace is one-dimensional and spanned by the vector $w \neq 0$.

Denoting with $\Omega^+:= \{ x \in \Omega \, : \, w(x) >0 \}$, $\Omega^-:= \{ x \in \Omega \, : \, w(x) <0 \}$ and assuming that $f_{- \infty} < f(s) < f^{+ \infty}$ for any $s \in \mathbb{R}$, Landesman and Lazer proved that the inequalities 

\begin{align*}
\displaystyle
f_{-\infty} \int_{\Omega^+} |w| \, dx - f^{+\infty} \int_{\Omega^-} |w| \, dx < \into h w \, dx  < f_{+\infty} \int_{\Omega^+} |w| \, dx - f^{-\infty}  \int_{\Omega^-} |w| \, dx
\end{align*} 

\medskip
\noindent
are sufficient for the existence of a solution of the boundary value problem \eqref{LandLaz} (see also \cite{HESS} for an improvement of the proof in which $\lambda$ is also isolated, and \cite{AMBROSETTIMANCINI} for the existence of multiple solutions).

An extension of previous result was obtained in \cite{ALP}, where Lazer et al. proved existence of solutions of \eqref{LandLaz} without assuming $f_{- \infty} < f(s) < f^{+ \infty}$ for any $s \in \mathbb{R}$, and showing that also the reverse inequalities 

\begin{align*}
\displaystyle
f_{+\infty} \int_{\Omega^+} |w| \, dx - f^{-\infty} \int_{\Omega^-} |w| \, dx < \into h w \, dx  < f_{-\infty} \int_{\Omega^+} |w| \, dx - f^{+\infty}  \int_{\Omega^-} |w| \, dx
\end{align*} 

\medskip
\noindent
are sufficient for the existence of a solution of \eqref{LandLaz}.

The quasilinear counterpart, that is quasilinear elliptic equations involving $p$-Laplacian,  was studied by several authors in different contexts (see for instance \cite{ANGOS,AO,BDK}).\\

Here we also concern with  the existence of a weak solution for a quasilinear elliptic systems in resonance around $\lambda_1$, under new sufficient Landesman-Lazer type conditions, 
extending the results by Arcoya and Orsina \cite{AO}.\\

\medskip
\noindent
We then consider the following system:

\begin{equation}\label{Syst0}
	\begin{cases}
		\begin{array}{ll}
			-\Delta_p u = \lambda_1 |u|^{\alpha-1} |v|^{\beta+1} u + \frac{1}{\alpha + 1} [F_s(x,u,v) - h_1(x) ]& x\in\Omega,
			\bigskip
			\\
			-\Delta_q v= \lambda_1  |u|^{\alpha+1} |v|^{\beta - 1}v  + \frac{1}{\beta + 1} [ F_t(x,u,v) - h_2(x)] & x\in \Omega,
			\bigskip
			\\
			u=v=0  & x\in \partial\Omega.
		\end{array}
	\end{cases}
\end{equation}

\medskip
\noindent
Here the nonlinearity  $F:\Omega \times \mathbb{R}^2 \to \mathbb{R}$ is a $C^1$-Carathéodory function (that is, $F(x, s, t)$ is measurable with respect to $ x \in \Omega$ for every $(s,t) \in \mathbb{R}^2$, and of class $C^1$ with respect to $(s,t) \in \mathbb{R}^2$ for a.e. $ x \in \Omega$), while $h_1 \in L^{p'}(\Omega)$ and  $h_2 \in L^{q'}(\Omega)$, where $p':= p/(p-1)$ and $q':=q/(q-1)$.

\medskip
\noindent
We assume that $F$ satisfies also the following conditions:\\

\begin{itemize}
\item[\textbf{(F1)}] $F(x,0,0) \in L^1(\Omega)$ and there exists $M>0$ such that
$$ \displaystyle |F_s(x,s,t)| \leq M \quad \text{and} \quad |F_t(x,s,t)| \leq M \qquad \text{ for all } (x,s,t) \in \Omega \times \mathbb{R}^2.$$\\

\item[\textbf{(F2)}] For almost every $x \in \Omega$, there exist
\begin{align*}
\displaystyle 
& F_s^{++}(x)= \lim_{\substack{s \to + \infty \\ t  \to  + \infty}} F_t(x,s,t), & F_t^{++}(x)= \lim_{\substack{s \to + \infty \\ t  \to  + \infty}} F_s(x,s,t), \\
& F_s^{+-}(x)= \lim_{\substack{s \to + \infty \\ t  \to  - \infty}} F_t(x,s,t), & F_t^{+-}(x)= \lim_{\substack{s \to + \infty \\ t  \to  - \infty}} F_s(x,s,t), \\
& F_s^{-+}(x)= \lim_{\substack{s \to - \infty \\ t  \to  + \infty}} F_t(x,s,t), & F_t^{-+}(x)= \lim_{\substack{s \to - \infty \\ t  \to  + \infty}} F_t(x,s,t), \\
& F_s^{--}(x)= \lim_{\substack{s \to - \infty \\ t  \to  - \infty}} F_t(x,s,t), & F_t^{--}(x)= \lim_{\substack{s \to - \infty \\ t  \to  - \infty}} F_t(x,s,t).
\end{align*}
\end{itemize}

\medskip
\noindent
Moreover, setting

$$ \displaystyle (\varphi_1,\psi_1):=\left( \frac{\varphi_0}{ \left( \lVert \varphi_0 \rVert_{1,p}^p + \lVert \psi_0 \rVert_{1,q}^q \right)^{\frac{1}{p}}}  ,\frac{\psi_0}{ \left( \lVert \varphi_0 \rVert_{1,p}^p + \lVert \psi_0 \rVert_{1,q}^q \right)^{\frac{1}{q}}} \right),$$

\medskip
\noindent
we infer that  $(\varphi_1,\psi_1)$ is an eigenfunction associated with $\lambda_1$ with 

$$\displaystyle \lVert \varphi_1 \rVert^p_{1,p} + \lVert \psi_1 \rVert^q_{1,q}=1.$$

\medskip
\noindent
In particular, by $(p,q)$-homogeneity we can write $E_1$ given by Theorem \ref{simplicityINTRODUZIONE} also as

\begin{equation}\label{InsiemeE1riscrittoINTRODUZIONE}
 \displaystyle E_1= \bigg\{ ( \left|  \theta \right|^{\frac{1}{p}} \varphi_1 ,  \left|  \theta \right|^{\frac{1}{q}} \psi_1 ) \, \text{sgn} (\theta) \, : \, \theta \in \mathbb{R} \bigg\} \,  \bigcup \, \bigg\{ (- \left|  \theta \right|^{\frac{1}{p}} \varphi_1 ,  \left|  \theta \right|^{\frac{1}{q}} \psi_1 ) \, \text{sgn} (\theta) \, : \, \theta \in \mathbb{R} \bigg\}.
\end{equation}

\noindent
We then introduce the following sufficient Landesman-Lazer type conditions, involving $p,q,F,h_1,h_2$ and the eigenfunction $(\varphi_1,\psi_1)$ (we omit $dx$ for clarity of reading):\\

\newpage

\noindent
if $p < q$,  either

\begin{equation}\label{LLpminoreq1}
\begin{cases}
\displaystyle \int_{\Omega} F_s^{++}\varphi_1 \,  < \int_{\Omega} h_1 \varphi_1 \,  < \int_{\Omega} F_s^{--} \varphi_1 \bigskip\\
\text{and} \bigskip\\
\displaystyle
 \int_{\Omega} F_s^{+-}\varphi_1 \,   < \int_{\Omega} h_1 \varphi_1 \,  < \int_{\Omega} F_s^{-+} \varphi_1,
 \end{cases}
\end{equation}

\medskip
\noindent
or

\begin{equation}\label{LLpminoreq2}
\begin{cases}
\displaystyle
 \int_{\Omega} F_s^{--}\varphi_1 \,   < \int_{\Omega} h_1 \varphi_1 \,  < \int_{\Omega} F_s^{++} \varphi_1 \,  \bigskip\\
\text{and} \bigskip\\
\displaystyle
 \int_{\Omega} F_s^{-+}\varphi_1 \,  < \int_{\Omega} h_1 \varphi_1 \,  < \int_{\Omega} F_s^{+-} \varphi_1;
 \end{cases}
\end{equation}

\medskip
\noindent
if $p=q$,  either

\begin{equation}\label{LLpp1}
\begin{cases}
\displaystyle
 \into F_s^{++}\varphi_1  + F_t^{++}\psi_1 \,  <  \into h_1 \varphi_1 +  h_2 \psi_1  \,   <  \into F_s^{--} \varphi_1 +  F_t^{--} \psi_1 \, \bigskip\\
\text{and} \bigskip\\
\displaystyle
    \into F_s^{+-} \varphi_1 - F_t^{+-} \psi_1 \, <  \into h_1 \varphi_1  - h_2 \psi_1 \,  < \into F_s^{-+} \varphi_1- F_t^{-+} \psi_1,
    \end{cases}
\end{equation}

\medskip
\noindent
or

\begin{equation}\label{LLpp2}
\begin{cases}
\displaystyle
    \into F_s^{--} \varphi_1 + F_t^{--} \psi_1 \, <  \into h_1 \varphi_1  + h_2 \psi_1 \, < \into F_s^{++} \varphi_1+ F_t^{++} \psi_1 \,  \bigskip\\
\text{and} \bigskip\\
\displaystyle
    \into F_s^{-+} \varphi_1 - F_t^{-+} \psi_1 \, <  \into h_1 \varphi_1  - h_2 \psi_1 \,  < \into F_s^{+-} \varphi_1 - F_t^{+-} \psi_1;
\end{cases}
\end{equation}

\medskip
\noindent
if $p>q$,  either

\begin{equation}\label{LLpmaggioreq1}
\begin{cases}
\displaystyle
 \int_{\Omega} F_t^{++} \psi_1 \,  < \int_{\Omega} h_2 \psi_1 \,  < \int_{\Omega} F_t^{--} \psi_1  \bigskip \\
\text{and} \bigskip \\
\displaystyle
 \int_{\Omega} F_t^{-+} \psi_1\,  < \int_{\Omega} h_2 \psi_1\,  < \int_{\Omega} F_t^{+-} \psi_1, 
\end{cases} 
\end{equation}

\medskip
\noindent
or

\begin{equation}\label{LLpmaggioreq2}
\begin{cases}
\displaystyle
 \int_{\Omega} F_t^{--} \psi_1\,  < \int_{\Omega} h_2 \psi_1\,  < \int_{\Omega} F_t^{++} \psi_1  \bigskip \\
\text{and} \bigskip\\
\displaystyle
 \int_{\Omega} F_t^{+-} \psi_1\, < \int_{\Omega} h_2 \psi_1\,  < \int_{\Omega} F_t^{-+} \psi_1.
 \end{cases} 
\end{equation}

\medskip
\noindent
Our main result is the following.

\begin{theorem}\label{mainresult}
Let $F:\Omega \times \mathbb{R}^2 \to \mathbb{R}$ be a $C^1$-Carathéodory function satisfying \textbf{(F1)} and \textbf{(F2)}. Consider $h_1 \in L^{p'}(\Omega)$ and  $h_2 \in L^{q'}(\Omega)$. If we assume that 
\begin{itemize}

\item[•]  in the case $p<q$,  either condition \eqref{LLpminoreq1} or  condition \eqref{LLpminoreq2} holds true;

\item[•]  in the case $p=q$,  either condition \eqref{LLpp1} or  condition \eqref{LLpp2} holds true;

\item[•]  in the case $p>q$,  either condition \eqref{LLpmaggioreq1} or  condition \eqref{LLpmaggioreq2} holds true;
\end{itemize}
then there exists at least a weak solution $(u,v) \in X$ for problem \eqref{Syst0}. 
\end{theorem}

In the proof of Theorem \ref{mainresult}, we show that in any case the Euler functional $J$ associated with \eqref{Syst0} satisfies the compactness condition of Palais-Smale.

In addition, under the assumptions \eqref{LLpminoreq1} if $p<q$, \eqref{LLpp1} if $p=q$ and \eqref{LLpmaggioreq1} if $p>q$, the functional $J$ is coercive on $X$, and so \eqref{Syst0} has a solution. 

On the other hand, that is under the assumptions  \eqref{LLpminoreq2} if $p<q$, \eqref{LLpp2} if $p=q$ and \eqref{LLpmaggioreq2} if $p>q$, the extension of the results contained in \cite{AO} to problem \eqref{Syst0} is not straightforward, since the set $E_1$  given by \eqref{InsiemeE1riscrittoINTRODUZIONE} is not a vector subspace of $X$.

Indeed, in \cite{AO} it was possible to apply the saddle point theorem of Rabinowitz (see  \cite{RAB}), exploiting the fact that whenever $W_0^{1,p}(\Omega)$ is splitted as the direct sum of $\langle \varphi \rangle$ and $V$, where $\varphi$ is a nontrivial eigenfunction associated with the first eigenvalue $\mu_1$ of $-\Delta_p$ on $\Omega$ with respect to Dirichlet boundary condition, then from standard argument it follows that there exists $\bar{\lambda} > \mu_1$ such that $\into \left| \nabla u \right|^p \, dx \geq \bar{\lambda} \into |u|^p \, dx $ for any $u \in V$. In particular, this argument requires just the fact that the eigenspace $\langle  \varphi \rangle$ associated with $\mu_1$ is a finite-dimensional vector subspace of $W_0^{1,p}(\Omega)$, without using any other information on the spectrum of the $p$-Laplacian, even the fact that $\mu_1$ is isolated.

Here we cannot apply the same argument, since the set $E_1$ given by \eqref{InsiemeE1riscrittoINTRODUZIONE} is never a vector subspace of $X$, even in the case $p=q$, in which case we have the union of two vector subspaces. We overcome this difficulty  by using the fact that $\lambda_1$ is isolated and the variational characterization of $\lambda_2$, therefore we are able to prove that two suitable sets satisfy the geometry of saddle point and are linked in a way that allows the application of standard minimax theorems (see \cite[Theorem 3.4]{STRUWE}).\\

This Chapter is organized as follows: in Section \ref{simplicitylambda1} we prove the simplicity of $\lambda_1$ in the sense stated in Theorem \ref{simplicityINTRODUZIONE}; Section \ref{isolationoflambda1} is devoted to the proof of Theorem \ref{lambda1isolatoINTRODUZIONE}, that is the isolation of $\lambda_1$; in Section \ref{sezionesequencaautovalori} we show that the sequence $\{\lambda_k\}_k$ defined in \eqref{lambdakINTRODUZIONE} is a sequence of eigenvalues of \eqref{EigenDeThélin0}.  Finally, in Section \ref{sezionePalaisSmale} we show that the Euler functional $J$ associated with problem \eqref{Syst0} satisfies the compactness condition of Palais-Smale in under any of the previously introduced Landesman-Lazer sufficient conditions, while Section \ref{GeometryofJ} is devoted to the geometry of $J$, so that Theorem \ref{mainresult} is proved.

\bigskip

\section{Simplicity of $\lambda_1$}\label{simplicitylambda1}

\bigskip
\noindent
The functionals $\Phi : X \to \re$ and $\Psi:X \to \re$ defined in \eqref{funzionalePhi} and \eqref{funzionalePsi} are of class $C^1$ and for any $(u_0,v_0), (u,v) \in X$ we have

\begin{align*}
\displaystyle \langle \Phi '(u_0,v_0), (u,v)  \rangle = (\alpha+1) \int_{\Omega} |\nabla u_0|^{p-2} \nabla u_0 \cdot \nabla u \, dx + (\beta+1) \int_{\Omega} |\nabla v_0|^{q-2} \nabla v_0 \cdot \nabla v \, dx
\end{align*}

\medskip
\noindent
and

\begin{align*}
\displaystyle \langle \Psi '(u_0,v_0), (u,v)  \rangle = (\alpha + 1) \into \left| u_0 \right|^{\alpha-1}  \left| v_0 \right|^{\beta+1} u_0 u \, dx + (\beta + 1) \into \left| u_0 \right|^{\alpha+1}  \left| v_0 \right|^{\beta-1} v_0 v \, dx.
\end{align*}

\medskip
\noindent
By \eqref{condizionepqalphabeta} we get that $\Phi$ and $\Psi$ are $(p,q)$-homogeneous functionals, i.e.

\begin{equation*}
	\Phi(\theta^{\frac 1 p}u,\theta^{\frac 1 q}v)=\theta\Phi(u,v), \quad \Psi(\theta^{\frac 1 p}u,\theta^{\frac 1 q}v)=\theta\Psi(u,v) \qquad \forall \theta>0,\ (u,v)\in X.
\end{equation*}

Notice that $\Psi(u,v)=0$ whenever $u(x)v(x)=0$ a.e. $x \in \Omega$. Hence, denoting with $\tilde{X}$ the set of $(u,v) \in X$ such that the product $uv$ is different from the zero function on $\Omega$, i.e.
  
$$\tilde{X}:= \{ (u,v) \in X \, : \, u v \not \equiv 0\},$$

\medskip
\noindent
we can characterize $\displaystyle \lambda_1 := \inf_{z \in \Sigma} \Phi(z)$ also as

$$ \displaystyle \lambda_1 = \inf_{z \in \tilde{X}} \frac{\Phi(z)}{\Psi(z)}.$$

\medskip
\noindent
In fact, it suffices to show that $A=B$ where

$$ \displaystyle A:=\Biggl\{ \frac{\Phi(u,v)}{\Psi(u,v)}  ; \: \; (u,v) \in \tilde{X} \Biggr\} \quad
\text{and} \quad B:=\Biggl\{ \Phi(u,v) ; \: \; (u,v) \in  \Sigma \Biggr\}.$$

\medskip
\noindent
We have $B \subseteq A$ (in fact $\Sigma \subset \tilde{X}$ and $\Psi(u,v)=1$ when $(u,v) \in \Sigma$).\\
Let now $z \in A$, that is $\displaystyle z=\frac{\Phi(u,v)}{\Psi(u,v)} $ for some $(u,v) \in \tilde{X}$.\\
Since $\displaystyle \Psi(u,v) \neq 0$,  we can define

$$\displaystyle (\bar{u},\bar{v}):= \left( \frac{u}{\Psi(u,v)^{\frac{1}{p}}}   ,  \frac{v}{\Psi(u,v)^{\frac{1}{q}}}  \right).$$

\medskip
\noindent
In particular, we have $(\bar{u},\bar{v}) \in \Sigma$, by which  $\Phi(\bar{u},\bar{v}) \in B$. Moreover we get

\begin{align*}
\displaystyle
 \Phi(\bar{u},\bar{v})= \frac{\Phi(u,v)}{\Psi(u,v)} = z,
\end{align*}
and then $z \in B$. We have proved $A \subseteq B$ and we conclude that $A=B$.

In \cite[Theorem 2]{DT} de Thélin showed through Tolksdorf regularity results \cite{tolksdorf1983} that every nontrivial eigenfunction  $(\varphi,\psi)$ satisfies $\varphi \in C^{1,\eta}(\overline{\Omega})$, $\psi \in C^{1,\xi}(\overline{\Omega})$. In addition, through Vazquez's maximum principle \cite{VAZ}, for every nontrivial eigenfunction $(\varphi_0, \psi_0)$ associated with $\lambda_1$ with   $\varphi_0 \geq 0$  and $\psi_0 \geq 0$, we have  necessarily $\varphi_0 > 0,$ $\psi_0 > 0$ in $\Omega$ and  $\frac{\partial \varphi_0}{\partial \nu}, \frac{\partial \psi_0}{\partial \nu} < 0$ on $\partial \Omega$, where $\nu$ denotes the outward normal to $\partial \Omega$. This result is easily extended to any eigenvalue $\lambda$ of \eqref{EigenDeThélin0}.

\begin{lemma}\label{autofunzionipositive}
If $(\varphi,\psi)$ is a nontrivial eigenfunction associated with an eigenvalue $\lambda \geq \lambda_1$ of \eqref{EigenDeThélin0}, then $\varphi \in C^{1,\eta}(\overline{\Omega})$ and $\psi \in C^{1,\xi}(\overline{\Omega})$ for some $\eta,\xi \in (0,1)$.\\
If $\varphi \geq 0$ and $\psi \geq 0$, then $\varphi >0, \psi >0$ in $\Omega$ and $\frac{\partial \varphi}{\partial \nu}, \frac{\partial \psi}{\partial \nu} < 0$ on $\partial \Omega$.
\end{lemma}

\begin{proof}
Let $(\varphi,\psi)$ be a nontrivial eigenfunction associated with an eigenvalue $\lambda \geq \lambda_1$.
By \cite[Theorem 1.3]{BCV} we know that both $\varphi \in L^{\infty}(\Omega)$ and $\psi \in L^{\infty}(\Omega)$.\\ Since assumptions of \cite[Theorem 1]{lieberman} are satisfied, we infer that $\varphi \in C^{1,\eta}(\overline{\Omega})$ and $\psi \in C^{1,\xi}(\overline{\Omega})$ for some $\eta,\xi \in (0,1)$.

Now, let us assume $\varphi \geq 0$ and $\psi \geq 0$ and let us consider the first equation of \eqref{EigenDeThélin0}:
$$ \displaystyle  \Delta_p \varphi = -\lambda |\varphi|^{\alpha-1} |\psi|^{\beta+1} \psi.$$
Since assumptions  of \cite[Theorem 5]{VAZ} are satisfied, we get $\varphi > 0$ in $\Omega$ and  $\frac{\partial \varphi}{\partial \nu} < 0$ on $\partial \Omega$. Similarly we obtain  $\psi >0$ in $\Omega$ and $\frac{\partial \psi}{\partial \nu} < 0$ on $\partial \Omega$.

\end{proof}

\medskip
\noindent
{\mbox {\it Proof of Theorem~\ref{simplicityINTRODUZIONE}.~}  It is contained in \cite{DT}, but for the reader's convenience, we provide a sketch here.
First of all, we observe that any $(\varphi, \psi) \in E_1$ is an eigenfunction of \eqref{EigenDeThélin0}  associated with the eigenvalue $\lambda_1$. Conversely, let $(\varphi,\psi)$ be an eigenfunction associated with $\lambda_1$. Since

$$ \displaystyle \frac{ \Phi (|\varphi|,|\psi|) }{\Psi (|\varphi|,|\psi|)}=\frac{ \Phi (\varphi,\psi) }{\Psi (\varphi,\psi)} = \lambda_1,   $$

\medskip
\noindent
we deduce that $(\left|\varphi\right|, \left|\psi \right|)$ is an eigenfunction associated with $\lambda_1$.\\We know that $|\varphi| \in C^{1,\eta}(\overline{\Omega})$, $|\psi| \in C^{1,\xi}(\overline{\Omega})$ and $|\varphi| > 0,$ $|\psi| > 0$ in $\Omega$. Therefore both $\varphi$ and $\psi$ cannot change sign. In fact, if by contradiction they change sign, since $\Omega$ is a connected set we would have $x_0 \in \Omega$ such that $\varphi(x_0)=0$, and similarly for $\psi$, but this is not possible.\\
Let us assume for example $\varphi > 0$ and $\psi > 0$.
Denoting with

$$ \displaystyle \theta := \into \left| \varphi \right|^{\alpha+1}  \left| \psi \right|^{\beta+1} \, dx, $$

\noindent
we have that 

$$ \displaystyle (\bar{\varphi}, \bar{\psi}):= \left(  \frac{\varphi}{\theta^{\frac{1}{p}}} ,  \frac{\psi}{\theta^{\frac{1}{q}}} \right) $$

\medskip
\noindent
is an eigenfunction associated with $\lambda_1$ such that $\bar{\varphi} > 0$, $\bar{\psi} > 0$ and $(\bar{\varphi}, \bar{\psi}) \in \Sigma$. By \cite[Theorem 3]{DT} 
there is a unique eigenfunction $(\varphi_0,\psi_0)$ associated with $\lambda_1$  satisfying $(\varphi_0,\psi_0) \in \Sigma$,    $\varphi_0 >0$ and $ \psi_0>0$  in $\Omega$ and
we deduce $(\bar{\varphi}, \bar{\psi})= (\varphi_0,\psi_0)$, hence $(\varphi,\psi) \in E_1$. We have proved the result assuming $\varphi > 0$ and $\psi > 0$.

The same result holds if we argue as before on the vector $(\varphi, - \psi)$ if $\varphi > 0$ and $\psi <0$, $(-\varphi,  \psi)$ if $\varphi <  0$ and $\psi > 0$, and $(-\varphi, - \psi)$ if $\varphi < 0$ and $\psi < 0$.

\qed

\bigskip

\section{Isolation of $\lambda_1$}\label{isolationoflambda1}

\bigskip
\noindent
To prove that $\lambda_1$ is isolated, we follow the ideas of \cite[Theorem 1]{AH}  by using the Picone's identity (see \cite[Lemma 24]{AG}).

\begin{theorem}\label{Picone}
Let $u \geq 0$ and $v >0$ be differentiable functions, and let $r >1$. Denote

\begin{align*}
\displaystyle
L_r(u,v) & =  \left| \nabla u \right|^r + (r-1) \left( \frac{u}{v} \right)^r \left| \nabla v \right|^r - r \left( \frac{u}{v} \right)^{r-1}\left| \nabla v \right|^{r-2} \nabla v \cdot \nabla u, \\
R_r(u,v) & =  \left| \nabla u \right|^r - \left| \nabla v \right|^{r-2} \nabla v \cdot  \nabla \left(    \frac{u^{r}}{v^{r-1}} \right).
\end{align*}

\medskip
\noindent
Then 

$$ \displaystyle L_r(u,v) = R_r(u,v) \geq 0.$$

\medskip
\noindent
Moreover, $L_r(u,v)=0$ a.e. in $\Omega$ if and only if $u=kv$ for some constant $k$. 
\end{theorem}

\begin{proposition}\label{funzionicambianosegno}
Let $(\varphi,\psi)$ be an eigenfunction of \eqref{EigenDeThélin0} associated with the eigenvalue $\lambda$ with $\lambda > \lambda_1$. Then both $\varphi$ and $\psi$ change sign. 
\end{proposition}

\begin{proof} 
Let $(\varphi,\psi)$ be an eigenfunction of \eqref{EigenDeThélin0} associated with the eigenvalue $\lambda$ with $\lambda > \lambda_1$. By contradiction, assume that both $\varphi$ and $\psi$ have constant sign, for example $\varphi \geq 0$ and $\psi \geq 0$. We know that 
 $\varphi \in C^{1,\eta}(\overline{\Omega})$ and $\psi \in C^{1,\xi}(\overline{\Omega})$ for some $\eta, \xi \in (0,1)$, and that $\varphi > 0, \psi >0$ in $\Omega$.

Using that  $(\varphi_0, \psi_0)$ is  an eigenfunction associated with $\lambda_1$ and $(\varphi,\psi)$ is an eigenfunction  associated with $\lambda$, for any $(u,v) \in X$ we have

\begin{align}\label{autofunzionephi0psi0}
\displaystyle 
& \into \left| \nabla \varphi_0  \right|^{p-2} \nabla \varphi_0 \cdot  \nabla u \, dx + \into \left| \nabla \psi_0   \right|^{q-2} \nabla \psi_0 \cdot  \nabla v \, dx \\
=  & \lambda_1 \left[ \into |\varphi_0|^{\alpha-1} |\psi_0|^{\beta+1} \varphi_0 u \, dx   + \into |\varphi_0|^{\alpha+1} |\psi_0|^{\beta-1} \psi_0 v \, dx     \right], \nonumber
\end{align}

\medskip
\noindent
and

\begin{align*}
\displaystyle 
 \into \left| \nabla \varphi  \right|^{p-2} \nabla \varphi \cdot  \nabla u \, dx + \into \left| \nabla \psi   \right|^{q-2} \nabla \psi \cdot  \nabla v \, dx 
=  \lambda \left[ \into |\varphi|^{\alpha-1} |\psi|^{\beta+1} \varphi u \, dx   + \into |\varphi|^{\alpha+1} |\psi|^{\beta-1} \psi v \, dx     \right]. \nonumber
\end{align*}

\medskip
\noindent
By testing \eqref{autofunzionephi0psi0} in $ \displaystyle (u,v)= \left( \frac{\alpha + 1}{p}  \varphi_0  , \frac{\beta +1}{q}  \psi_0\right), $  we get

\begin{align*}
\displaystyle
\frac{\alpha + 1}{p} \into \left| \nabla \varphi_0 \right|^p \, dx  + \frac{\beta + 1}{q} \into \left| \nabla \varphi_0 \right|^p \, dx  = \lambda_1 \into \varphi_0^{\alpha + 1} \psi_0^{\beta +1}    \, dx.
\end{align*}

\medskip
\noindent
By testing second equation in
$ (u,v)= \displaystyle \left(  \frac{\alpha + 1}{p} \frac{\varphi_0^{p}}{\varphi^{p-1}} ,\frac{\beta + 1}{q} \frac{\psi_0^{q}}{\psi^{q-1}}    \right),$ we get

\begin{align*}
\displaystyle
& \frac{\alpha + 1}{p} \into \left| \nabla \varphi  \right|^{p-2} \nabla \varphi \cdot  \nabla  \left( \frac{\varphi_0^{p}}{\varphi^{p-1}}    \right)  \, dx  +  \frac{\beta + 1}{q} \into \left| \nabla \psi  \right|^{q-2} \nabla \psi \cdot  \nabla  \left( \frac{\psi_0^{q}}{\psi^{q-1}}    \right)  \, dx \\
= &  \lambda \left[  \frac{\alpha + 1}{p} \into \varphi^{\alpha+1-p} \psi^{\beta+1} \varphi_0^{p} \, dx  +  \frac{\beta + 1}{q} \into \varphi^{\alpha+1} \psi^{\beta+1-q} \psi_0^{q} \, dx     \right].
\end{align*}

\medskip
\noindent
By Theorem \ref{Picone}, and recalling that $\lambda_1 < \lambda,$   we have

\begin{align}\label{dimfunznoncambsegn}
\displaystyle 
0 
& \leq \frac{\alpha+1}{p} \into L_p(\varphi_0, \varphi) \, dx + \frac{\beta+1}{q} \into L_q(\psi_0, \psi) \, dx \\
& = \frac{\alpha+1}{p} \into R_p(\varphi_0, \varphi) \, dx + \frac{\beta+1}{q} \into R_q(\psi_0, \psi) \, dx   \nonumber\\
& = \frac{\alpha+1}{p} \into \left| \nabla \varphi_0  \right|^p \, dx - \frac{\alpha+1}{p} \into \left| \nabla \varphi \right|^{p-2} \nabla \varphi \cdot  \nabla \left(    \frac{\varphi_0^{p}}{\varphi^{p-1}} \right) \, dx \nonumber \\
& \quad + \frac{\beta+1}{q} \into \left| \nabla \psi_0  \right|^q \, dx  \nonumber  - \frac{\beta+1}{q} \into \left| \nabla \psi \right|^{q-2} \nabla \psi \cdot  \nabla \left(    \frac{\psi_0^{q}}{\psi^{q-1}} \right) \, dx \nonumber \\
& =  \lambda_1 \into \varphi_0^{\alpha + 1} \psi_0^{\beta +1}    \, dx         \nonumber \\
& \quad - \lambda \left[  \frac{\alpha + 1}{p} \into \varphi^{\alpha+1-p} \psi^{\beta+1} \varphi_0^{p} \, dx  +  \frac{\beta + 1}{q} \into \varphi^{\alpha+1} \psi^{\beta+1-q} \psi_0^{q} \, dx     \right]\nonumber \\
& \leq \lambda  \left[ \into ( \, \varphi_0^{\alpha + 1} \psi_0^{\beta +1}     -  \frac{\alpha + 1}{p} \varphi^{\alpha+1-p} \psi^{\beta+1} \varphi_0^{p}   - \frac{\beta + 1}{q} \varphi^{\alpha+1} \psi^{\beta+1-q} \psi_0^{q} \, ) \, dx    \right].  \nonumber
\end{align}

\medskip
\noindent
However, by Young inequality with conjugated exponents $p/(\alpha + 1)$ and $q/(\beta +1)$, we also have

\begin{align*}
\displaystyle
\varphi_0^{\alpha + 1} \psi_0^{\beta +1}  
& = \left( \varphi_0^{\alpha + 1} \frac{\psi^{\frac{(\alpha+1)(\beta+1)}{p}}}{\varphi^{\frac{(\alpha+1)(\beta+1)}{q}}}  \right) \, \left(  \psi_0^{\beta +1} \frac{\varphi^{\frac{(\alpha+1)(\beta+1)}{q}}}{\psi^{\frac{(\alpha+1)(\beta+1)}{p}}} \right) \\
& \leq \frac{\alpha+1}{p} \varphi_0^p \psi^{\beta+1} \varphi^{\alpha+1-p} + \frac{\beta+1}{q} \psi_0^q \varphi^{\alpha+1}  \psi^{\beta+1-q},
\end{align*}

\medskip
\noindent
whence, considering \eqref{dimfunznoncambsegn}, we get 

$$ \displaystyle L_p(\varphi_0, \varphi) = 0 \qquad \text{and} \qquad L_q(\psi_0, \psi)=0.$$

\medskip
\noindent
Applying again Theorem \ref{Picone}, we deduce $\varphi_0 = c_1 \varphi$ and $\psi_0= c_2 \psi$ for some positive constant $c_1$ and $c_2$. Observe now that by \eqref{autofunzionephi0psi0} we infer

$$ \displaystyle \into \left| \nabla \varphi_0 \right|^p \, dx = \lambda_1 \into \varphi_0^{\alpha + 1} \psi_0^{\beta +1}   \, dx = \into \left| \nabla \psi_0 \right|^q \, dx.$$

\medskip
\noindent
Substituting $\varphi_0 = c_1 \varphi$ and $\psi_0= c_2 \psi$, we get $c_1^p=c_2^q$. Denoting by $\theta:= 
c_1^p$, we have shown that $(\varphi_0, \psi_0) = (\theta^{\frac{1}{p}} \varphi, \theta^{\frac{1}{q}} \psi)$ by which $(\varphi_0, \psi_0)$ is an eigenfunction associated with $\lambda$, hence a contradiction.

The proof is similar in the other sign cases: either $\psi \leq 0 \leq \varphi,$ or $\psi,\varphi \leq 0$ or 
$\varphi \leq 0 \leq \psi$.

\end{proof}

\medskip
\noindent
{\mbox {\it Proof of Theorem~\ref{lambda1isolatoINTRODUZIONE}.~}
By the characterization of $\lambda_1$, for every eigenvalue $\lambda$ of \eqref{EigenDeThélin0} we have $\lambda \geq \lambda_1$. We prove the theorem by contradiction, assuming that $\{\mu_n\}_n$ is a sequence of eigenvalues of \eqref{EigenDeThélin0} such that $\mu_n > \lambda_1$ for any $n \geq 2$ and  $\mu_n \to \lambda_1$ as $n \to + \infty$. For every $n \geq 2$ consider an eigenfunction $z_n=(\varphi_n,\psi_n) \in X$ associated with $\mu_n$ with 

$$ \|\varphi_n\|_{1,p}^{p}+ \|\psi_n\|_{1,q}^q=1.$$

\medskip
\noindent
Hence the sequence $\{z_n\}_n$ is bounded and there exists a subsequence, still denoted by $\lbrace{z_n\rbrace}_n,$ that converges to some $z=(\varphi,\psi)$ weakly in $X$ and strongly in $L^p(\Omega)\times L^q(\Omega)$. Since $\varphi_n$ converges to $\varphi$  strongly in $L^p(\Omega)$ and $\psi_n$ converges to $\psi$  strongly in $L^q(\Omega)$, we have $\varphi_n(x) \to \varphi(x)$ and $\psi_n(x) \to \psi(x)$ a.e. in $\Omega$ as $n \to + \infty$. Moreover, passing to a subsequence, we can assume that  there exists $\bar{\varphi} \in L^{p}(\Omega)$ and $\bar{\psi} \in L^{q}(\Omega)$ such that $|\varphi_n(x)|\leq \bar{\varphi}(x)$ and $|\psi_n(x)|\leq \bar{\psi}(x)$ a.e. in $\Omega$ and for all $n \geq 1$.\\
Since $z_n$ is an eigenfunction associated with $\mu_n$, for any $(u,v) \in X$ it satisfies

\begin{align}\label{zneigenfunction}
\displaystyle 
& \into \left| \nabla \varphi_n  \right|^{p-2} \nabla \varphi_n \cdot  \nabla u \, dx + \into \left| \nabla \psi_n   \right|^{q-2} \nabla \psi_n \cdot  \nabla v \, dx \\
= & \mu_n \left[ \into |\varphi_n|^{\alpha-1} |\psi_n|^{\beta+1} \varphi_n u \, dx   + \into |\varphi_n|^{\alpha+1} |\psi_n|^{\beta-1} \psi_n v \, dx     \right]. \nonumber
\end{align}

\medskip
\noindent
By testing previous equation in $(u,v)=(\varphi_n - \varphi,0)$, we get

\begin{align}\label{compisol1}
\displaystyle 
\into \left| \nabla \varphi_n  \right|^{p-2} \nabla \varphi_n \cdot  \nabla ( \varphi_n - \varphi) \, dx = \mu_n \into |\varphi_n|^{\alpha-1} |\psi_n|^{\beta+1} \varphi_n (\varphi_n - \varphi) \, dx. 
\end{align}

\medskip
\noindent
By Young inequality with conjugated exponents $p/(\alpha+1)$ and $q/(\beta+1)$, we get

\begin{align*}
\displaystyle 
&\left|  \, |\varphi_n(x)|^{\alpha-1} |\psi_n(x)|^{\beta+1} \varphi_n(x) (\varphi_n(x) - \varphi(x) )  \,  \right|\\
\leq &  \, |\varphi_n(x)|^{\alpha} |\psi_n(x)|^{\beta+1}(  \left| \varphi_n(x) \right| + \left| \varphi(x) \right| ) \\
\leq & 2 \bar{\varphi}(x)^{\alpha+1} \bar{\psi}(x)^{\beta+1} \\
\leq & 2 \left( \frac{\alpha+1}{p}  \bar{\varphi}(x)^p   + \frac{\beta + 1}{q}   \bar{\psi}(x)^q \right) \in L^1(\Omega).
\end{align*} 

\medskip
\noindent
Hence, by dominated convergence theorem we have 

$$ \displaystyle \lim_{n \to \infty} \into |\varphi_n|^{\alpha-1} |\psi_n|^{\beta+1} \varphi_n (\varphi_n - \varphi) \, dx =0.$$

\medskip
\noindent
Considering also that $\mu_n \to \lambda_1>0$ and passing to limit as $n \to + \infty$ in \eqref{compisol1}, we get 

$$ \displaystyle \lim_{n \to \infty} \into \left| \nabla \varphi_n  \right|^{p-2} \nabla \varphi_n \cdot  \nabla ( \varphi_n - \varphi) \, dx =0.$$

\medskip
\noindent
Now, the weak lower semicontinuity and the convexity of  $\| \cdot \|^p_{1,p}$ imply

\begin{align*}
\|\varphi\|_{1,p}^p & \leq \liminf_{n \to \infty} \|\varphi_n\|_{1,p}^p
\leq \limsup_{n \to \infty} \|\varphi_n\|_{1,p}^p = \limsup_{n \to \infty} \into |\nabla \varphi_n|^p \, dx\\
& \leq \limsup_{n \to \infty} \left[ \into |\nabla \varphi|^p \, dx + p \into |\nabla \varphi_n|^{p-2} \nabla \varphi_n \cdot \nabla  \left(\varphi_n - \varphi \right) \, dx \right] = \|\varphi\|_{1,p}^p,
\end{align*}

\medskip
\noindent
i.e. $\|\varphi_n\|_{1,p} \to \|\varphi\|_{1,p},$ and by uniform convexity of $W_0^{1,p}(\Omega)$ we deduce $\varphi_n \rightarrow \varphi$ strongly in $\sob$. Analogously, we get that $\psi_n \rightarrow \psi$ strongly in $\sobq$.\\
In particular,

\[\|\varphi\|_{1,p}^p + \|\psi \|_{1,q}^q =
\lim_{n \to \infty} \left( \|\varphi_n\|_{1,p}^p +\|\psi_n\|_{1,q}^q \right)
=1\]

\medskip
\noindent
and $z=(\varphi,\psi) \neq 0.$
By testing \eqref{zneigenfunction} in $(u,v) =\left( \frac{\alpha+1}{p} \varphi_n, \frac{\beta + 1}{q} \psi_n \right)$, we get 

$$ \frac{\alpha+1}{p} \into \left| \nabla \varphi_n \right|^p \, dx + \frac{\beta+1}{q} \into \left| \nabla \psi_n \right|^q \, dx = \mu_n \into \left| \varphi_n \right|^{\alpha+1} \left| \psi_n \right|^{\beta+1} \, dx,$$

\medskip
\noindent
and passing to limit as $n \to + \infty$ we deduce

$$ \frac{\alpha+1}{p} \into \left| \nabla \varphi \right|^p \, dx + \frac{\beta+1}{q} \into \left| \nabla \psi \right|^q \, dx = \lambda_1 \into \left| \varphi \right|^{\alpha+1} \left| \psi \right|^{\beta+1} \, dx.$$

\medskip
\noindent
By last identity and $z=(\varphi,\psi) \neq (0,0),$ we get $z \in \tilde{X}$ and

$$ \displaystyle \lambda_1= \frac{\displaystyle \frac{\alpha+1}{p} \into \left| \nabla \varphi \right|^p \, dx + \frac{\beta+1}{q} \into \left| \nabla \psi \right|^q \, dx}{ \displaystyle \into \left| \varphi \right|^{\alpha+1} \left| \psi \right|^{\beta+1} \, dx  } = \frac{\Phi(z)}{\Psi(z)},$$

\medskip
\noindent
i.e. $z=(\varphi,\psi)$ is an eigenfunction associated with $\lambda_1$.\\
By Theorem \ref{simplicityINTRODUZIONE} and \eqref{InsiemeE1riscrittoINTRODUZIONE} we conclude that only one of the following possibilities occurs:

$$ \displaystyle (\varphi,\psi)= (\varphi_1,\psi_1), \qquad (\varphi,\psi)= (-\varphi_1,-\psi_1), \qquad (\varphi,\psi)= (-\varphi_1,\psi_1), \qquad (\varphi,\psi)= (\varphi_1,-\psi_1).$$

\medskip
\noindent
Let us assume $(\varphi,\psi)= (\varphi_1,\psi_1)$.\\
By Proposition \ref{funzionicambianosegno} we know that both $\varphi_n$ and $\psi_n$ change sign for any $n \geq 2$, hence both $\varphi_n^-:= -\min\{ \varphi_n, 0    \}$ and  $\psi_n^-:= -\min\{ \psi_n, 0    \}$ are different from zero.\\
By testing \eqref{zneigenfunction} in $(u,v)= \left(  \frac{\alpha + 1}{p} \varphi_n^-, \frac{\beta + 1}{q} \psi_n^- \right)$, we get

\begin{align*}
 \displaystyle \frac{\alpha+1}{p} \int_{\Omega} \left| \nabla \varphi_n^- \right|^p \, dx + \frac{\beta+1}{q} \int_{\Omega} \left| \nabla \psi_n^- \right|^q \, dx 
& =  \, \mu_n \int_{\Omega}  \left|\varphi_n^-\right|^{\alpha+1} \left| \psi_n^-\right|^{\beta+1} \, dx,\\
& =  \, \mu_n \int_{\Omega_n^-}  \left|\varphi_n^-\right|^{\alpha+1} \left| \psi_n^-\right|^{\beta+1} \, dx
\end{align*}

\medskip
\noindent
where $\Omega_n^- := \{ x \in \Omega \; : \; \varphi_n(x) < 0 \; \text{ or }  \;  \psi_n(x) < 0 \}$.\\
Applying Young and  H\"older inequalities we obtain

\begin{align*}
\displaystyle \int_{\Omega_n^-}  \left|\varphi_n^-\right|^{\alpha+1} \left| \psi_n^-\right|^{\beta+1} \, dx
& \leq   \left[ \frac{\alpha+1}{p} \int_{\Omega_n^-} \left|\varphi_n^-\right|^{p} \, dx   + \frac{\beta+1}{q} \int_{\Omega_n^-} \left| \psi_n^-\right|^{q} \, dx     \right] \\
& \leq   \left[ \frac{\alpha+1}{p} \left| \Omega_n^- \right|^{\frac{p}{N}} \lVert \varphi_n^-  \rVert_{p^*}^p  +  \frac{\beta+1}{q} \left| \Omega_n^- \right|^{\frac{q}{N}} \lVert \psi_n^-  \rVert_{q^*}^q \right].
\end{align*}

\medskip
\noindent
For $p^*=Np/(N-p)$ and $q^*=Nq/(N-q)$  we use now the Sobolev embeddings $\lVert u \rVert_{p^*} \leq \mathcal{S}(p,N) \lVert \nabla u \rVert_{p}$ and $\lVert v \rVert_{q^*} \leq \mathcal{S}(q,N) \lVert \nabla v \rVert_{q}$ to get from the above identity that 

\begin{align*}
\displaystyle
&\frac{\alpha+1}{p}  \int_{\Omega} \left| \nabla \varphi_n^- \right|^p \, dx + \frac{\beta+1}{q} \int_{\Omega} \left| \nabla \psi_n^- \right|^q \, dx  \\
 & \leq 
 \mu_n C \max \left\lbrace \left| \Omega_n^- \right|^{\frac{p}{N}}, \left| \Omega_n^- \right|^{\frac{q}{N}}   \right\rbrace \left[ \frac{\alpha+1}{p} \int_{\Omega} \left| \nabla \varphi_n^- \right|^p \, dx + \frac{\beta+1}{q} \int_{\Omega} \left| \nabla \psi_n^- \right|^q \, dx   \right],
\end{align*}

\medskip
\noindent
where $ \displaystyle C:= \max \left\lbrace \mathcal{S}(p,N)^p, \mathcal{S}(q,N)^q \right\rbrace.$ This implies that 

$$ \displaystyle \frac{1}{\mu_n C} \leq \max \left\lbrace \left| \Omega_n^- \right|^{\frac{p}{N}}, \left| \Omega_n^- \right|^{\frac{q}{N}}   \right\rbrace ,$$

\medskip
\noindent
whence

\begin{equation*}
 \displaystyle 0 < \frac{1}{C \lambda_1} \leq \liminf_{n \to \infty}  \, \max \left\lbrace \left| \Omega_n^- \right|^{\frac{p}{N}}, \left| \Omega_n^- \right|^{\frac{q}{N}}   \right\rbrace.
\end{equation*}

\medskip
\noindent
Since, up to subsequence,  $\varphi_n(x) \to \varphi_1(x)$, $\nabla \varphi_n(x) \to \nabla \varphi_1(x)$, $\psi_n(x) \to \psi_1(x)$, and $\nabla \psi_n(x) \to \nabla \psi_1(x)$ a.e. in $\Omega$, we can apply Severini-Egorov theorem to say that in the complement of a set of arbitrarily small measure we have $\varphi_n$ converges uniformly to $\varphi_1$, $ \nabla \varphi_n$ converges uniformly to $ \nabla \varphi_1$, $\psi_n$ converges uniformly to $\psi_1$ and $\nabla \psi_n$ converges uniformly to $\nabla \psi_1$. Since  $\varphi_1>,$  $\psi_1>0$ in $\Omega$ and $\frac{\partial \varphi_1}{\partial \nu},\frac{\partial \varphi_1}{\partial \nu}<0$ on $\partial \Omega$, we conclude that there exists a subset of $\Omega$ of arbitrarily small measure such that in its complement $\varphi_n$ and $\psi_n$ are positive  for $n$ sufficiently large, that is a contradiction with the last inequality.

Let us recall that in the previous argument we have assumed $(\varphi,\psi)= (\varphi_1,\psi_1)$.\\
In the rest of cases we can argue similarly. For instance, if $(\varphi,\psi)= (\varphi_1,-\psi_1)$, then we arrive to a contradiction by considering 
$\varphi_n^-:= -\min\{ \varphi_n, 0    \}$,  $\psi_n^+:= \max\{ \psi_n, 0    \}$, and
 
\begin{align*}
\displaystyle
\Omega_n^{-,+} := \{ x \in \Omega \; : \; \varphi_n(x) < 0 \; \text{ or }  \;  \psi_n(x) > 0  \}.
\end{align*}
 
\qed

\bigskip

\section{A sequence of eigenvalues}\label{sezionesequencaautovalori}

\bigskip
\noindent
Let us introduce the $C^1$ functional $ Q :  X \setminus \{(0,0)\} \to \mathbb{R}$ defined by 

\begin{equation}\label{RQ}
 \displaystyle Q(z):=\frac{ \Psi (z)}{ \Phi(z) },
\end{equation}

\medskip
\noindent
where $\Phi$ and $\Psi$ are defined in \eqref{funzionalePhi} and \eqref{funzionalePsi} respectively.\\
Observe that

\begin{equation}\label{derivataQ}
 \displaystyle Q'(z) = \frac{1}{\Phi(z)} \left[ \, \Psi' (z) - Q(z) \, \Phi' (z)    \,  \right] \qquad \forall z=(u,v) \in X \setminus \{(0,0) \}.
\end{equation}

\medskip
\noindent
We note that a real number $\lambda \neq 0$ is an eigenvalue of \eqref{EigenDeThélin0} if and only if there exists $\tilde{z} \in \tilde{X}$ such that

$$ \displaystyle \Phi'(\tilde{z})= \lambda \Psi'(\tilde{z}),$$

\medskip
\noindent
or equivalently, if $1/\lambda$ is a critical value for $Q$.\\
Let us denote with $\mathcal{M}$  the set

\begin{align}\label{nuovomanifold}
\displaystyle 
\mathcal{M}:=\Biggl\{ z=(u,v) \in X \, : \,  \Phi(z)=1   \Biggr\}.
\end{align}

$\mathcal{M}$  is a $C^1$ manifold of codimension $1$ in $X$ (see for instance \cite[Example 27.2]{DEM}) and for any $z \in \mathcal{M}$ the tangent space $T_z\mathcal{M}$ of $\mathcal{M}$ at $z$ can be identified with

$$ \displaystyle \text{Ker }( \Phi'(z) )= \{ h \in X \, : \, \langle \Phi'(z),h \, \rangle=0    \}.$$

\medskip
\noindent
Observe that by definition of Q and \eqref{derivataQ} the restriction $Q_{\mathcal{M}}:= Q\big|_{ \mathcal{M}}$ of $Q$ to $\mathcal{M}$ satisfies 

$$ \displaystyle   Q_{ \mathcal{M}}(z)=\Psi(z) \quad \text{ and } \quad  Q'(z) =  \Psi' (z) - \Psi(z) \, \Phi' (z) \qquad \forall z=(u,v) \in \mathcal{M}.$$

\medskip
\noindent
Moreover, (see \cite[Example 27.3]{DEM}), the derivative $Q_{\mathcal{M}}'$ of $Q_{\mathcal{M}}$ is defined on the tangent bundle $T\mathcal{M}$, and in particular for any $z \in \mathcal{M}$ the map $Q_{\mathcal{M}}'(z)$, defined on $T_z \mathcal{M}$, can be identified  as

\begin{equation}\label{derivataQsuM}
\displaystyle Q_{ \mathcal{M}}'(z)= [ Q' (z) -  \langle \, Q'(z) , e(z)   \, \rangle \, \Phi'(z)]\big|_{T_z \mathcal{M}}, \qquad \forall \,    z=(u,v) \in \mathcal{M},    
\end{equation}

\medskip
\noindent
where $e: \mathcal{M} \to X$ is given by 

$$ \displaystyle e(z)= \left( \frac{u}{p} , \frac{v}{q} \right) \qquad \forall \, z=(u,v) \in \mathcal{M}.$$

\medskip
\noindent
Notice that

$$ \displaystyle \langle \, \Psi'(z) , e(z)   \, \rangle = \Psi(z) \qquad \text{and} \qquad \langle \, \Phi'(z) , e(z)   \, \rangle = \Phi(z)=1 \qquad \forall z=(u,v) \in \mathcal{M},$$

\medskip
\noindent
hence by \eqref{derivataQsuM} we deduce

\begin{equation}\label{derivataQsuM3}
\displaystyle Q_{ \mathcal{M}}'(z)=  [\Psi' (z) -  \Psi(z) \, \Phi'(z)]\big|_{T_z \mathcal{M}} =Q'(z)\big|_{T_z \mathcal{M}} \qquad \forall z=(u,v) \in \mathcal{M}.
\end{equation}

We recall that a real number $c$ is a critical value of $Q_{ \mathcal{M}}$ if there exists $z \in \mathcal{M}$ such that $Q_{ \mathcal{M}}(z)=c$ and $ Q_{ \mathcal{M}}'(z) \equiv 0$. It follows from standard arguments that critical values of $Q_{ \mathcal{M}}$ correspond to critical values of $Q$. In particular,  we deduce that $\lambda \neq 0$ is an eigenvalue of \eqref{EigenDeThélin0} if and only if $1/\lambda$ is a critical value for $Q_{ \mathcal{M}}$.

Given $z \in \mathcal{M}$, we have $Q_{ \mathcal{M}}'(z) \in T^*_z \mathcal{M},$ where $ T^*_z \mathcal{M}$ denotes the dual space of $ T_z \mathcal{M}$ endowed with the norm $ \lVert \cdot \rVert_z^*$. Denoting with $\lVert \cdot \rVert^*$ the norm of $X^*$, by \cite[Proposition 3.54]{PAO} it follows that

\begin{align}\label{normaspaziotangente}
\displaystyle   \lVert Q_{ \mathcal{M}}'(z) \rVert^*_z= \min_{\mu \in \mathbb{R}} \lVert \,  Q'(z)  - \mu \Phi'(z) \, \rVert^*=\min_{\mu \in \mathbb{R}} \lVert \,  \Psi'(z)  -( \Psi(z) +  \mu) \Phi'(z) \, \rVert^* .
\end{align}

As it has been mentioned,  the construction of  a sequence of eigenvalues of \eqref{EigenDeThélin0} uses a deformation lemma for $C^1$ submanifolds of a Banach space (see Theorem~\ref{BON}) based in \cite[Theorem 2.5]{BON} that requires a coupled Palais-Smale condition.

Firstly we observe that $0$ is the only critical value of $\Phi$, therefore for any $\varepsilon >0$ sufficiently small, the manifold
$$ \mathcal{M}_{\varepsilon}:= \Biggl\{ z =(u,v) \in X \, : \, \Phi(z)=1+ \varepsilon \Biggr\} $$
is well defined, hence definitions and results proved for $Q$ on $\mathcal{M}$ still hold on $\mathcal{M}_{\varepsilon}$.

\begin{lemma}\label{QPhiPSaccoppiata}
For every $c>0$ the functionals $ \displaystyle Q$ and $\Phi$ satisfy the following coupled Palais-Smale condition at level $c$ on $\mathcal{M}$ (denoted by $\widehat{(PS)}_{\mathcal{M},c}$), that is  any sequence $\{z_n\}_n = \{(u_n,v_n)\}_n \subset X$ such that \medskip
\begin{itemize}
\item[$i)$] $z_n \in \mathcal{M}_{\varepsilon_n},$ where $\varepsilon_n >0$ and $\varepsilon_n \to 0$ as $n \to \infty$; \medskip
\item[$ii)$] $Q(z_n) \to c$ as $n \to \infty$; \medskip
\item[$iii)$] $\lVert Q_{\mathcal{M}_{\varepsilon_n}}'(z_n) \rVert_{z_n}^* \to 0 \text{ or } \lVert \Phi'(z_n) \rVert \to 0$ as $n \to \infty$; \medskip
\end{itemize}
has a converging subsequence in $X$.
\end{lemma}

\begin{proof} Let $\{z_n\}_n = \{(u_n,v_n)\}_n \subset X$ a sequence satisfying $i)-iii)$.\\
By condition $i)$ we deduce that 

\begin{align*}
\displaystyle \lVert u_n \rVert_{1,p}^p \leq \frac{p}{\alpha+1} \Phi(z_n) = \frac{p(1 + \varepsilon_n)}{\alpha+1} \quad \text{and} \quad \lVert v_n \rVert_{1,q}^q \leq \frac{q}{\beta+1} \Phi(z_n) = \frac{q(1 + \varepsilon_n)}{\beta+1},
\end{align*}

\medskip
\noindent
that is $ \lbrace{ z_n \rbrace}_n $ is bounded in $X$, and so there exists a subsequence, still denoted by $\lbrace{z_n\rbrace}_n,$ that converges to some $z=(u,v)$ weakly in $X$ and strongly in $L^p(\Omega)\times L^q(\Omega)$.
By condition $iii)$ we have two possibilities.\\
Firstly we assume that $\lVert \Phi'(z_n) \rVert \to 0$ as $n \to \infty$ holds true. In particular, we have

$$ \langle \Phi'(z_n) , u_n - u \rangle \to 0,$$

\medskip
\noindent
that is

$$ \displaystyle \into \left| \nabla u_n \right|^{p-2} \nabla u_n \cdot \nabla (u_n - u) \, dx \to 0.$$

\medskip
\noindent
Now, the weak lower semicontinuity and the convexity of  $\| \cdot \|^p_{1,p}$
imply

\begin{align*}
\|u\|_{1,p}^p & \leq \liminf_{n \to \infty} \|u_n\|_{1,p}^p
\leq \limsup_{n \to \infty} \|u_n\|_{1,p}^p = \limsup_{n \to \infty} \into |\nabla u_n|^p \, dx\\
& \leq \limsup_{n \to \infty} \left[ \into |\nabla u|^p \, dx + p \into |\nabla u_n|^{p-2} \nabla u_n \cdot \nabla  \left(u_n-u\right) \, dx \right] = \|u\|_{1,p}^p,
\end{align*}

\medskip
\noindent
i.e. $\|u_n\|_{1,p}^p \to \|u\|_{1,p}^p$, and  by uniform convexity of $W_0^{1,p}(\Omega)$, we deduce $u_n\rightarrow u$ in $\sob$. Analogously,  we get that $v_n\rightarrow v$ strongly in $\sobq$.

Now, let us assume the other possibility, i.e.  $\lVert Q_{\mathcal{M}_{\varepsilon_n}}'(z_n) \rVert_{z_n}^* \to 0 $ as $n \to \infty$.
Strong convergences of $u_n$ to $u$ in $L^{p}(\Omega)$ and of $v_n$ to $ v$ in $L^{q}(\Omega)$ give $ \displaystyle \Psi(z_n) \to \Psi(z).$
Moreover, by $i)$ we infer that $\Phi(z_n) \to 1$. By using also $ii)$, we get

$$ \displaystyle \Psi(z_n) \to \Psi(z)=c.$$

\medskip
\noindent
By \eqref{normaspaziotangente} applied to $z_n \in \mathcal{M}_{\varepsilon_n}$ and the convergence to zero of $\lVert Q_{\mathcal{M}_{\varepsilon_n}}'(z_n) \rVert_{z_n}^*$, there exists a sequence $\{\mu_n\}_n \subset \mathbb{R}$ such that

$$ \displaystyle \Psi'(z_n)  -( \Psi(z_n) +  \mu_n) \Phi'(z_n) \to 0 \quad \text{in } X^* \qquad \text{as } n \to \infty.$$

\medskip
\noindent
Testing on the vector $(u_n/p,v_n/q)$, we get 

$$ \displaystyle \Psi(z_n)  -( \Psi(z_n) +  \mu_n) \Phi(z_n) \to 0 \qquad \text{as } n \to \infty,$$

\medskip
\noindent
and using that $\Phi(z_n) \to 1$ and $\Psi(z_n) \to c$, we deduce $\mu_n \to 0$ as $n \to \infty$.
Testing now on the vector $(u_n - u,0),$ we get

$$ \displaystyle \into |u_n|^{\alpha-1}  |v_n|^{\beta +1 } u_n (u_n -u) \, dx - ( \Psi(z_n) +  \mu_n)\into \left| \nabla u_n \right|^{p-2} \nabla u_n \cdot \nabla (u_n - u) \, dx \to 0,  $$

\medskip
\noindent
and considering the strong convergences of $u_n$ to $u$ in $L^p(\Omega)$ and $v_n$ to $v$ in $L^q(\Omega)$ we deduce

$$ \displaystyle ( \Psi(z_n) +  \mu_n)\into \left| \nabla u_n \right|^{p-2} \nabla u_n \cdot \nabla (u_n - u) \, dx \to 0.   $$

\medskip
\noindent
Since $\Psi(z_n) +  \mu_n\to c$ with $c \neq 0$, we get

\begin{equation*}
\displaystyle 
\into \left| \nabla u_n \right|^{p-2} \nabla u_n \cdot \nabla (u_n - u) \, dx \to 0. 
\end{equation*}

\medskip
\noindent
Arguing as before we conclude that $u_n \to u$ strongly in $W_0^{1,p}(\Omega)$. Similarly we deduce that $v_n \to v$ strongly in $W_0^{1,q}(\Omega)$.

\end{proof}

We will apply the following deformation result for $C^1$ submanifolds of a Banach space, which is a consequence of \cite[Theorem 2.5]{BON}.

\begin{theorem}\label{BON}
Let $X$ be a Banach space and $\Phi \in C^1(X, \mathbb{R})$. Assume that $1$ is a regular value of $\Phi$. If $Q$ is a $C^1$ functional on a neighborhood of $\mathcal{M}= \{ x \in X \, : \, \Phi(x)=1 \}$ and $Q$ and $\Phi$ satisfy $\widehat{(PS)}_{\mathcal{M},c}$ for  a regular value $c \in \mathbb{R}$ of $Q_\mathcal{M}$, then
there exists $\hat{\varepsilon} >0$ such that for all $0 < \varepsilon < \hat{\varepsilon}$ there is an homeomorphism $\eta$ of $\mathcal{M}$ onto $\mathcal{M}$ such that: \medskip
\begin{enumerate}
\item $ \eta (z)=z$  if $Q(z) \not \in [ c - \hat{\varepsilon}, c + \hat{\varepsilon}]$; \medskip
\item $Q(\eta(z)) \geq Q(z)$ for all $z \in \mathcal{M}$; \medskip
\item $Q(\eta(z)) \geq c + \varepsilon$ for all $z \in \mathcal{M}$ such that $Q(z) \geq c - \varepsilon$; \medskip
\item If $\mathcal{M}$ is symmetric $(\mathcal{M}=-\mathcal{M})$ and if $Q$ is even, then $\eta$ is odd. \medskip
\end{enumerate}
\end{theorem}

Now, we will construct a sequence of eigenvalues for \eqref{EigenDeThélin0}.

\begin{proposition}\label{prop 2.9}
Let $k \in \mathbb{N}$ and let $S^{k-1}$ be the unit sphere in $\mathbb{R}^k$. Set
 
$$
\mathcal{M}_k := \{ A=\alpha(S^{k-1}) \subset \mathcal{M}  \ | \,  \alpha: S^{k-1} \to \mathcal{M} \text{ is odd and continuous}  \}.
$$

\medskip
\noindent
If

\begin{equation}\label{ck}
	c_k := \sup_{A \in \mathcal{M}_k} \min_{(u,v) \in A} Q_{ \mathcal{M}}(u,v),
\end{equation}

\medskip
\noindent
then

\begin{equation}\label{lambdakSfera}
	\lambda_k := \frac{1}{c_k}
\end{equation}

\medskip
\noindent
is an eigenvalue of problem \eqref{EigenDeThélin0}.

\end{proposition}

\begin{proof}
As it has been mentioned a real number $\lambda$ is an eigenvalue of \eqref{EigenDeThélin0} if and only if $1/\lambda$ is a critical value of $ Q_{ \mathcal{M}}$. By contradiction assume  that $c_k$  is a regular value of $ Q_{ \mathcal{M}}$. By Theorem \ref{BON} there exists $ \varepsilon >0$ an homeomorphism $\eta$ of $\mathcal{M}$ onto $\mathcal{M}$ satisfying  properties $(1)-(4)$. By definition of $c_k$, we deduce that there exists an odd and continuous map $\alpha: S^{k-1} \to \mathcal{M} $ such that $ \min_{(u,v) \in A} Q_{ \mathcal{M}}(u,v) > c_k - \varepsilon$ for $A= \alpha(S^{k-1})$. Hence, taking into account property $(4)$, we deduce that $\eta(A) \in \mathcal{M}_k$. However, by $(3)$ of the cited Theorem we infer that $ \min_{(u,v) \in \eta(A) } Q_{ \mathcal{M}}(u,v) \geq c_k + \varepsilon$, which contradicts the definition of $c_k$. Therefore $c_k$ is a critical value of $ Q_{ \mathcal{M}}$ and $\lambda_k$ is an eigenvalue of \eqref{EigenDeThélin0}.

\end{proof}

\begin{remark}\label{definizionicoincidono}
We point out that the first eigenvalue $\lambda_1$ defined in \eqref{lambdakSfera} coincides with the eigenvalue $\lambda_1$ defined in \eqref{lambda1isolatoINTRODUZIONE}. In fact, by definition we know that
 
$$
\mathcal{M}_1 := \{ A=\{ z,-z\}   \ | \, z\in \mathcal{M}\},
$$

\medskip
\noindent
and since $\Psi (-z)=\Psi (z)$ we also have

$$ 
\displaystyle c_1 
= \sup_{A \in \mathcal{M}_1} \min_{ \bar{z} \in A} \Psi(\bar{z})
= \sup_{ z \in \mathcal{M}}  \Psi(z) 
=  \sup_{z \in X \setminus (0,0)} \frac{ \Psi (z)}{ \Phi(z) } 
= \frac{1}{ \displaystyle \inf_{z \in \tilde{X}} \frac{\Phi(z)}{\Psi(z)}}= \frac{1}{\lambda_1}.
$$

\end{remark}

\bigskip

\section{Palais-Smale condition of $J$}\label{sezionePalaisSmale}

\bigskip
\noindent
In order to study problem \eqref{Syst0}, let us consider the functional $J:X \to \mathbb{R}$  defined for any $z=(u,v) \in X$ as

\begin{align}\label{J}
\displaystyle
J(z)
&:=  \frac{\alpha+1}{p} \into \left| \nabla u \right|^p \, dx + \frac{\beta+1}{q} \into \left| \nabla v \right|^q \, dx    - \lambda_1   \into \left| u \right|^{\alpha+1} \left| v \right|^{\beta+1} \, dx \\
& \quad - \into F(x,u,v) \, dx + \into h_1 u \, dx + \into h_2 v \, dx \nonumber \\
& = \Phi(u,v) - \lambda_1 \Psi(u,v) - \into F(x,u,v) \, dx + \into h_1 u \, dx + \into h_2 v \, dx. \nonumber
\end{align}

\medskip
\noindent
It turns out that the functional $J$ is of class $C^1$ on $X$ with 

\begin{align*}
 \langle J'(z_0), z \rangle  
= & \displaystyle (\alpha+1) \into |\nabla
		u_0|^{p-2}\nabla u_0 \cdot \nabla u \ dx + (\beta+1)\into |\nabla v_0|^{q-2}\nabla v_0 \cdot \nabla v \ dx  \\
 & - \displaystyle \lambda_1 \left[ (\alpha+1)\into | u_0|^{\alpha-1}| v_0|^{\beta+1} u_0 u \ dx + (\beta+1)\into | u_0|^{\alpha+1}| v_0|^{\beta-1} v_0  v \ dx  \right] \\
 & -\displaystyle \into \left[  F_s(x,u_0,v_0) u +  F_t(x,u_0,v_0) v \right] \, dx +  \into h_1 u \, dx + \into h_2 v \, dx 
\end{align*}

\medskip
\noindent
for any $z_0=(u_0,v_0), z=(u, v) \in X$.

\begin{lemma}
Let $F:\Omega \times \mathbb{R}^2 \to \mathbb{R}$ be a $C^1$-Carathéodory function satisfying \textbf{(F1)} and \textbf{(F2)}. Consider $h_1 \in L^{p'}(\Omega)$ and  $h_2 \in L^{q'}(\Omega)$. If we assume that\\
 
\begin{itemize}
\item[•]  in the case $p<q$,  either condition \eqref{LLpminoreq1} or  condition \eqref{LLpminoreq2} holds true; \medskip

\item[•]  in the case $p=q$,  either condition \eqref{LLpp1} or  condition \eqref{LLpp2} holds true; \medskip

\item[•]  in the case $p>q$,  either condition \eqref{LLpmaggioreq1} or  condition \eqref{LLpmaggioreq2} holds true; \medskip
\end{itemize}

\medskip
\noindent
then $J$ satisfies the (PS) condition, namely if $ \lbrace{ z_n \rbrace}_n = \lbrace{ (u_n,v_n) \rbrace}_n$ is a sequence in $X$ such that there exists a positive constant $c$ such that  

\begin{equation}\label{PS1}
\left| J(z_n)\right| \leq c  \quad \text{ for any } n \in \mathbb{N}
\end{equation}

\medskip
\noindent
and 

\begin{equation}\label{PS2}
\|J'(z_n)\|~\rightarrow~0,
\end{equation}

\medskip
\noindent
then $\lbrace{ z_n \rbrace}_n $ has a convergent subsequence in $X$.
\end{lemma}

\begin{proof}
Let $\lbrace{ z_n \rbrace}_n=\lbrace{ (u_n,v_n) \rbrace}_n \subset X$ satisfying \eqref{PS1} and \eqref{PS2}. Let us begin by proving that the sequence $\lbrace{ z_n \rbrace}_n=\lbrace{ (u_n,v_n) \rbrace}_n$  is bounded in $X$. By contradiction, there exists a subsequence, still denoted with $\lbrace{ z_n \rbrace}_n,$  such that $\|(u_n,v_n)\|~\rightarrow~\infty$.
In particular, as $n \to + \infty$ we have

$$r_n:=\|u_n\|_{1,p}^{p}+ \|v_n\|_{1,q}^q \to + \infty.$$

\medskip
\noindent
Let us set

$$ \displaystyle \bu_n= \frac{u_n}{r_n^{\frac{1}{p}}} \qquad \text{and} \qquad  \bv_n=\frac{v_n}{r_n^{\frac{1}{q}}},$$

\medskip
\noindent
and observe that $\bz_n:=(\bu_n,\bv_n)$ satisfies 

$$ \displaystyle \lVert \bu_n \rVert_{1,p}^p + \lVert \bv_n \rVert_{1,q}^q = \frac{\|u_n\|_{1,p}^{p}+ \|v_n\|_{1,q}^q }{r_n}= 1 \qquad \text{ for any } n \in \mathbb{N}.$$

\medskip
\noindent
Hence $ \lVert \bz_n \rVert \leq 2$ for any $n \in \mathbb{N}$ and there exists a subsequence, still denoted by $\lbrace{\bz_n\rbrace}_n,$ that converges to some $\bz=(\bu,\bv)$ weakly in $X$ and strongly in $L^p(\Omega)\times L^q(\Omega)$.\\
In particular, as $n \to \infty$ we have

\begin{align}\label{ps1}
\frac{1}{\ r_n^{\frac{p-1}{p}}}\ \bigl\langle J'(z_n),(\bu_n-\bu,0)\bigr\rangle
& = \displaystyle  (\alpha+1)\into |\nabla \bu_n|^{p-2}\nabla \bu_n \cdot \nabla \left( \bu_n - \bu \right) \ dx \\
& \displaystyle \quad -\lambda_1 (\alpha+1)\into | \bu_n|^{\alpha-1} | \bv_n|^{\beta+1}\bu_n \left( \bu_n - \bu \right) \ dx \nonumber \\
& \displaystyle \quad - \frac{1}{\ r_n^{\frac{p-1}{p}}} \into  F_s(x,u_n,v_n) (\bu_n - \bu)   \, dx \nonumber \\
& \displaystyle \quad + \frac{1}{\ r_n^{\frac{p-1}{p}}} \into h_1(\bu_n - \bu)  \, dx \to 0, \nonumber
\end{align}

\medskip
\noindent
and

\begin{equation}\label{palaissmale2}
\frac{1}{\ r_n^{\frac{q-1}{q}}}\ \bigl\langle J'(z_n),(0,\bv_n-\bv)\bigr\rangle \rightarrow 0.
\end{equation}

\medskip
\noindent
Since $\bu_n \to \bu$ in $L^p(\Omega)$, $h \in L^{p'}(\Omega)$ and $F_s$ is bounded  we deduce as $n \to + \infty$ that

\begin{equation*}
\displaystyle  \frac{1}{\ r_n^{\frac{p-1}{p}}} \into h_1(\bu_n - \bu)  \, dx \, \to 0 \quad \text{ and } \quad  \frac{1}{\ r_n^{\frac{p-1}{p}}} \into  F_s(x,u_n,v_n) (\bu_n - \bu)   \, dx \, \to 0.
\end{equation*}

\medskip
\noindent
Moreover, considering $(\alpha +1)/p + (\beta + 1)/q=1$ and the strong convergences $\bu_n \to \bu$ in $L^{p}(\Omega)$ and $\bv_n \to \bv$ in $L^{q}(\Omega)$, by dominated convergence theorem we get

$$ \displaystyle \lim_{n \to \infty} \into | \bu_n|^{\alpha-1} | \bv_n|^{\beta+1}\bu_n \left( \bu_n - \bu \right) \, dx = 0.$$

\medskip
\noindent
Therefore, taking into account these convergences in \eqref{ps1}, we have proved

$$ \displaystyle \lim_{n \to \infty} \into |\nabla \bu_n|^{p-2}\nabla \bu_n \cdot \nabla \left( \bu_n - \bu \right) \, dx  = 0.$$

\medskip
\noindent
How done in the proof of Theorem \ref{lambda1isolatoINTRODUZIONE}, by using also \eqref{palaissmale2} we finally obtain that $\{\bu_n,\bv_n\}_n$ converges to $\bz=(\bu,\bv)$ strongly in $X$, where $\bz \neq (0,0)$ satisfies

$$ \displaystyle \lVert \bu \rVert_{1,p}^p + \lVert \bv \rVert_{1,q}^q=1.$$

\medskip
\noindent
Now, to prove that $\bz$ is an eigenfunction associated with the eigenvalue $\lambda_1$, it suffices to show that $\bu \neq 0$, $\bv \neq 0$ and 

$$ \displaystyle \frac{\Phi(\bz)}{\Psi(\bz)}=\lambda_1.$$

\medskip
\noindent
In the following, we  assume both $\lVert u_n \rVert_{1,p} \to \infty$ and $\lVert v_n \rVert_{1,q} \to \infty$. In the case that one of both sequences is bounded,  computations are simpler and left to reader.
Considering just the boundedness form above coming from \eqref{PS1}, we have 

$$\displaystyle \Phi(u_n,v_n) - \lambda_1 \Psi(u_n,v_n) - \into F(x,u_n,v_n) \, dx + \into h_1 u_n \, dx + \into h_2 v_n \, dx \leq c.$$

\medskip
\noindent
Dividing by $r_n$ and taking into account the definition of $\bu_n$ and $\bv_n$, we have

$$\displaystyle \Phi(\bu_n,\bv_n) - \lambda_1 \Psi(\bu_n,\bv_n) -  \into \frac{F(x,u_n,v_n)}{r_n} \, dx + \into \frac{h_1 u_n}{r_n} \, dx +  \into \frac{h_2 v_n}{r_n} \, dx \leq \frac{c}{r_n},$$

\medskip
\noindent
whence

\begin{align}\label{limsup}
\displaystyle \limsup_{n \to  \infty} \Biggl[ 
\Phi(\bu_n,\bv_n) - \lambda_1 \Psi(\bu_n,\bv_n)  
 -  \into \frac{F(x,u_n,v_n)}{r_n} \, dx + \into \frac{h_1 u_n}{r_n} \, dx +  \into \frac{h_2 v_n}{r_n} \, dx \Biggr] \leq 0. 
\end{align}

\medskip
\noindent
Since $h_1 \in L^{p'}(\Omega)$ and $h_2 \in L^{q'}(\Omega)$, by applying H\"older and Poincaré inequality and the convergence to infinity of the sequence $\{r_n\}_n$, we have 

$$ \displaystyle \lim_{n \to \infty} \into  \frac{h_1 u_n}{r_n} \, dx = 0 \qquad \text{and}  \qquad \lim_{n \to \infty} \into \frac{h_2 v_n}{r_n} \, dx =0.$$

\medskip
\noindent
For almost every $(x,s,t) \in \Omega \times \mathbb{R}^2$ mean value theorem gives

$$ \displaystyle F(x,s,t)=  F(x,0,0) + \int_0^1 F_s(x,\tau s,\tau t)s \, d \tau + \int_0^1 F_t(x,\tau s,\tau t)t \, d \tau,$$

\medskip
\noindent
hence by assumption \textbf{(F1)} we deduce

$$ \displaystyle \left| F(x,s,t)   \right| \leq \left| F(x,0,0) \right| + M |s| + M|t|.$$

\medskip
\noindent
Therefore,

$$ \displaystyle \left| \into \frac{F(x,u_n,v_n)}{r_n} \, dx \right| \leq \frac{\lVert F(x,0,0) \rVert_1}{r_n} + K_1 \lVert u_n \rVert_{1,p}^{1-p} + K_2 \lVert v_n \rVert_{1,q}^{1-q},$$

\medskip
\noindent
whence
$$ \displaystyle \lim_{n \to \infty} \into \frac{F(x,u_n,v_n)}{r_n} \, dx =  0.$$

\medskip
\noindent
Considering these convergences to zero into \eqref{limsup}, we obtain

$$\displaystyle \limsup_{n \to  \infty} \Biggl[ \Phi(\bu_n,\bv_n) - \lambda_1 \Psi(\bu_n,\bv_n) \Biggr] \leq 0.$$

\medskip
\noindent
Now, since $\bu_n \, \to \bu$ in $L^p(\Omega)$ and $\bv_n \, \to \bv$ in $L^q(\Omega)$, by Lebesgue theorem we deduce

$$ \displaystyle \Psi(\bu_n,\bv_n) \, \to \Psi(\bu,\bv) \qquad   \text{ as } n \, \to \infty$$

\medskip
\noindent
and by the strong convergences $\bu_n \, \to\bu$ in $\sob$ and $\bv_n \, \to \bv$ in $\sobq$ we also get

$$ \displaystyle \Phi(\bu_n,\bv_n) \, \to \Phi(\bu,\bv) \qquad   \text{ as } n \, \to \infty.$$

\medskip
\noindent
Therefore, we have proved $\Phi(\bu,\bv) \leq \lambda_1 \Psi(\bu,\bv)$, i.e.
 
$$ \displaystyle  \frac{\alpha+1}{p} \into \left| \nabla \bu \right|^p \, dx + \frac{\beta+1}{q} \into \left| \nabla \bv \right|^q \, dx    \leq \lambda_1   \into \left| \bu \right|^{\alpha+1} \left| \bv \right|^{\beta+1} \, dx. $$

\medskip
\noindent
Since $\bz \neq (0,0)$, the left side of previous inequality is strictly greater than zero, by which $\bu \neq 0,$ $\bv \neq 0$ and

$$ \displaystyle \frac{\Phi(\bz)}{\Psi(\bz)}=\lambda_1,$$

\medskip
\noindent
that is $\bz=(\bu,\bv)$ is an eigenfunction associated with the eigenvalue $\lambda_1$.\\
Since also $\bz$ satisfies $ \displaystyle \lVert \bu \rVert_{1,p}^p + \lVert \bv \rVert_{1,q}^q=1,$ we conclude by Theorem \ref{simplicityINTRODUZIONE} and \eqref{InsiemeE1riscrittoINTRODUZIONE} that only one of the following possibilities occurs:

$$ \displaystyle (\bu,\bv)= (\varphi_1,\psi_1), \qquad (\bu,\bv)= (-\varphi_1,-\psi_1), \qquad (\bu,\bv)= (-\varphi_1,\psi_1), \qquad (\bu,\bv)= (\varphi_1,-\psi_1).$$

\medskip
\noindent
By \eqref{PS2}, we deduce that, up to subsequence, there exists a decreasing sequence $\{ \varepsilon_n \}_n$, converging to zero, such that

$$ \displaystyle \left| \left\langle J'(z_n), \left(\frac{u_n}{p},\frac{v_n}{q}\right) \right\rangle \right| \leq \varepsilon_n \lVert z_n \rVert \qquad \forall n \in \mathbb{N}, \quad \forall z=(u,v) \in X.$$

\medskip
\noindent
Hence, taking into account also \eqref{PS1},  we have

\begin{equation*}
\displaystyle  \left(1 +  \varepsilon_n \lVert z_n \rVert \right) C \geq \left| J(z_n) - \left\langle J'(z_n), \left(\frac{u_n}{p},\frac{v_n}{q}\right)\right\rangle + \into F(x,0,0) \, dx \right|,
\end{equation*}

\medskip
\noindent
which by the mean value theorem and the definition of $\bu_n$ and $\bv_n$ implies

\begin{align}\label{calcolo}
\displaystyle
\left(1 +  \varepsilon_n \lVert z_n \rVert \right) C \geq
& \left| r_n^{\frac{1}{p}} \Biggl\{ \left( 1 - \frac{1}{p}  \right) \into h_1 \bu_n \, dx  \ -  \into  \int_0^1 F_s(x,\tau u_n,\tau v_n) \bu_n \, d \tau \, dx  \nonumber \right.\\
+&  \frac{1}{p} \into F_s(x,u_n,v_n) \bu_n \,  dx \Biggr\} \\ 
+& r_n^{\frac{1}{q}} \Biggl\{ \left( 1 - \frac{1}{q}  \right) \into h_2 \bv_n \, dx  \ -  \into  \int_0^1 F_t(x,\tau u_n,\tau v_n)\bv_n\, d \tau \, dx  \nonumber\\
+&   \left. \frac{1}{q} \into F_t(x,u_n,v_n)\bv_n \,  dx \Biggr\}\right| .  \nonumber
\end{align}

\medskip
\noindent
Let us assume $(\bu,\bv)= (\varphi_1,\psi_1)$.\\
Since $ \bu_n \to \varphi_1$ in $L^p(\Omega)$ and  $ \bv_n \to \psi_1$ in $L^q(\Omega)$, we have

$$\displaystyle \lim_{n \to \infty} \into h_1 \bu_n \, dx  = \into h_1 \varphi_1 \, dx \qquad \text{and} \qquad  \lim_{n \to \infty} \into h_2 \bv_n \, dx  = \into h_1 \psi_1 \, dx.$$

\medskip
\noindent
Now, let us observe that $u_n(x)=\bu_n(x)r_n^{\frac{1}{p}}$, where $r_n^{\frac{1}{p}} \to + \infty$ and $\bu_n(x) \to \varphi_1(x)$ a.e. $x \in \Omega$. Since $\varphi_1 > 0$ in $ \Omega$, we deduce $u_n(x) \to + \infty$ a.e. $x \in \Omega$, and similarly we have $v_n(x) \to + \infty$ a.e. $x \in \Omega$. Hence, by assumption $\textbf{(F2)}$ we get

$$ \displaystyle F_s(x,u_n(x),v_n(x)) \to F_s^{++}(x) \qquad \text{a.e. } x \in \Omega,$$

\medskip
\noindent
and

$$ \displaystyle F_t(x,u_n(x),v_n(x)) \to F_t^{++}(x) \qquad \text{a.e. } x \in \Omega,$$

\medskip
\noindent
so by Lebesgue theorem we deduce

$$ \displaystyle \lim_{n \to \infty} \into F_s(x,u_n,v_n)\bu_n \, dx = \into F_s^{++} \varphi_1 \, dx $$

\medskip
\noindent
and

$$ \lim_{n \to \infty}  \into F_t(x,u_n,v_n) \bv_n \, dx = \into F_t^{++} \psi_1 \, dx. $$

\noindent
In a similar way, by using also Fubini-Tonelli's Theorem we derive

\begin{equation*}
\displaystyle \lim_{n \to \infty} \into \int_0^1 F_s(x, \tau u_n, \tau v_n)\bu_n \, d \tau \, dx = \into F_s^{++} \varphi_1 \, dx
\end{equation*}

\medskip
\noindent
and

\begin{equation*}
\displaystyle \lim_{n \to \infty} \into \int_0^1 F_t(x, \tau u_n,\tau v_n)\bv_n \, d \tau \, dx = \into F_t^{++} \psi_1 \, dx.
\end{equation*}

\medskip
\noindent
Considering these convergences in \eqref{calcolo}, and denoting with $o(1)$  any quantity that goes to zero as $n$ goes to infinity, we have proved

\begin{align}\label{computation2}
\displaystyle
\left(1 +  \varepsilon_n \lVert z_n \rVert \right) C 
\geq & \left| r_n^{\frac{1}{p}} \left( 1 - \frac{1}{p}  \right) \left[ \into h_1 \varphi_1 \, dx - \into F_s^{++} \varphi_1 \, dx + o(1) \right] \right. \\
& \left. + r_n^{\frac{1}{q}}\left( 1 - \frac{1}{q}  \right) \left[ \into h_2 \psi_1 \, dx  - \into F_t^{++} \psi_1 \, dx + o(1) \right] \right|. \nonumber
\end{align}

\medskip
\noindent
If $p<q$, dividing  by $r_n^{\frac{1}{p}} \left(1 - \frac{1}{p} \right)$ and since $\varepsilon_n \to 0$ as $n \to \infty$, passing to limit we get

$$ \displaystyle 0 \geq \left|  \into h_1 \varphi_1 \, dx - \into F_s^{++} \varphi_1 \, dx    \right|,$$

\medskip
\noindent
that is a contradiction with both \eqref{LLpminoreq1} and \eqref{LLpminoreq2}.\\
If $p=q$, in a similar way we get 

$$ \displaystyle 0 \geq \left| \into h_1 \varphi_1 \,  + h_2 \psi_1 \, dx  - \into  F_s^{++} \varphi_1 + F_t^{++} \psi_1 \, dx   \right|,$$

\medskip
\noindent
that is a contradiction with both \eqref{LLpp1} and \eqref{LLpp2}.\\
If $p > q$, dividing  by $r_n^{\frac{1}{q}} \left(1 - \frac{1}{q} \right)$ and passing to limit we obtain
 
$$ \displaystyle 0 \geq \left|  \into h_2 \psi_1 \, dx - \into F_t^{++} \psi_1 \, dx    \right|,$$

\medskip
\noindent
that is a contradiction with both \eqref{LLpmaggioreq1} and \eqref{LLpmaggioreq2}.\\
Let us recall that we assumed $(\bu,\bv)=(\varphi_1,\psi_1)$.

If it would be  $(\bu,\bv)= (-\varphi_1,-\psi_1)$, we deduce inequality \eqref{computation2} in which $F_s^{++}$ and $F_t^{++}$ are replaced by $F_s^{--}$ and $F_t^{--}$ respectively, and reasoning as before we get a contradiction under the same assumptions.

If now we consider $(\bu,\bv)= (\varphi_1,-\psi_1)$, convergences in \eqref{calcolo} give

\begin{align*}
\displaystyle
\left(1 +  \varepsilon_n \lVert z_n \rVert \right) C 
\geq & \left| r_n^{\frac{1}{p}} \left( 1 - \frac{1}{p}  \right) \left[ \into h_1 \varphi_1 \, dx - \into F_s^{+-} \varphi_1 \, dx + o(1) \right] \right. \\
& \left. + r_n^{\frac{1}{q}}\left( 1 - \frac{1}{q}  \right) \left[  \into F_t^{+-} \psi_1 \, dx -\into h_2 \psi_1 \, dx  + o(1) \right] \right|, \nonumber
\end{align*}

\medskip
\noindent
by which we infer a contradiction with \eqref{LLpminoreq1} and \eqref{LLpminoreq2} if $p<q$, with \eqref{LLpp1} and \eqref{LLpp2} if $p=q$, and with \eqref{LLpmaggioreq1} and \eqref{LLpmaggioreq2} if $p>q$.

If it would be  $(\bu,\bv)= (-\varphi_1,\psi_1)$, we deduce previous inequality in which $F_s^{+-}$ and $F_t^{+-}$ are replaced by $F_s^{-+}$ and $F_t^{-+}$ respectively, and we reach a contradiction under the same hypothesis.

We have so proved that the sequence $\{z_n\}_n=\{(u_n,v_n)\}_n$ satisfying \eqref{PS1} and \eqref{PS2} is bounded in $X$. This implies that there exists $z=(u,v) \in X$ such that, up to subsequence, $\{z_n\}_n$ converges to $z$ weakly in $X$ and strongly in $L^p(\Omega)\times L^q(\Omega)$.\\
In particular, as $n \to \infty$ we have

\begin{equation*}
\bigl\langle J'(z_n),(u_n-u,0)\bigr\rangle \rightarrow 0
\end{equation*}

\medskip
\noindent
and

\begin{equation*}
\bigl\langle J'(z_n),(0,v_n-v)\bigr\rangle \rightarrow 0.
\end{equation*}

\medskip
\noindent
Arguing as previously done to show the strong convergence of $\bu_n$ to $\bu$ in $W_0^{1,p}(\Omega)$ and the strong convergence of $\bv_n$ to $\bv$ in $W_0^{1,q}(\Omega)$, we conclude that up to subsequence $z_n \to z$ strongly in $X$.

\end{proof}

\section{Geometry of $J$}\label{GeometryofJ}

\bigskip
\noindent
\begin{proposition}
Let $F:\Omega \times \mathbb{R}^2 \to \mathbb{R}$ be a $C^1$-Carathéodory function satisfying \textbf{(F1)} and \textbf{(F2)}. Consider $h_1 \in L^{p'}(\Omega)$ and  $h_2 \in L^{q'}(\Omega)$. If we assume that\\

\begin{itemize}

\item[•]  in the case $p<q$,  condition \eqref{LLpminoreq1}   holds true; \medskip

\item[•]  in the case $p=q$,  condition \eqref{LLpp1}         holds true; \medskip

\item[•]  in the case $p>q$,  condition \eqref{LLpmaggioreq1} holds true; \medskip

\end{itemize}
then $J$ is coercive, i.e., 

$$ \displaystyle \lim_{\lVert(u,v) \rVert \to + \infty} J(u,v)=+ \infty.$$
\end{proposition}

\begin{proof} Arguing by contradicition, suppose that there exists a constant $c$ and a sequence $\{(u_n,v_n)\}_n$ in $X$ such that

\begin{equation}\label{coercconto1}
\displaystyle J(u_n,v_n) \leq c \qquad \text{for any } n \in \mathbb{N},
\end{equation} 

\medskip
\noindent
and

$$ \displaystyle \lVert (u_n,v_n) \rVert \to + \infty \qquad \text{as } n \to + \infty.$$

\medskip
\noindent
Let us set

$$ \displaystyle r_n:=\lVert u_n \rVert_{1,p}^p + \lVert v_n \rVert_{1,q}^q, \qquad \bu_n:= \frac{u_n}{r_n^{\frac{1}{p}}} \qquad \text{and} \qquad  \bv_n:=\frac{v_n}{r_n^{\frac{1}{q}}},$$

\medskip
\noindent
and observe that $\bz_n:=(\bu_n,\bv_n)$ satisfies 

$$ \displaystyle \lVert \bu_n \rVert_{1,p}^p + \lVert \bv_n \rVert_{1,q}^q = 1 \qquad \text{ for any } n \in \mathbb{N}.$$

\medskip
\noindent
Hence $ \lVert \bz_n \rVert \leq 2$ for any $n \in \mathbb{N}$ and there exists a subsequence, still denoted by $\lbrace{\bz_n\rbrace}_n,$ that converges to some $\bz=(\bu,\bv)$ weakly in $X$ and strongly in $L^p(\Omega)\times L^q(\Omega)$.\\
Dividing \eqref{coercconto1} by $r_n$, arguing as in the proof of (PS) condition and considering also the characterization of $\lambda_1$ and the weak lower semicontinuity of both $\lVert \cdot \lVert_{1,p}^p$ and $\lVert \cdot \lVert_{1,q}^q$, we get

\begin{align*}
\displaystyle 
\Phi(\bu,\bv)  \leq \liminf_{n \to \infty} \Phi(\bu_n,\bv_n) \leq \limsup_{n \to \infty} \Phi(\bu_n,\bv_n)   \leq \lambda_1 \Psi(\bu,\bv) \leq \Phi(\bu,\bv).
\end{align*}

\medskip
\noindent
In particular, since $\lVert \bu_n \rVert_{1,p}^p + \lVert \bv_n \rVert_{1,q}^q =1 $ for any $n \in \mathbb{N},$ we deduce that $(\bu,\bv)$ is a nontrivial eigenfunction associated with $\lambda_1$, hence there exists $\theta>0$ such that only one of the following possibilities occurs:

$$ \displaystyle (\bu,\bv)= (\theta^{\frac{1}{p}}\varphi_1,\theta^{\frac{1}{q}} \psi_1), \qquad (\bu,\bv)= (-\theta^{\frac{1}{p}}\varphi_1,-\theta^{\frac{1}{q}}\psi_1),$$

$$ (\bu,\bv)= (\theta^{\frac{1}{p}}\varphi_1,-\theta^{\frac{1}{q}}\psi_1), \qquad (\bu,\bv)= (-\theta^{\frac{1}{p}}\varphi_1,\theta^{\frac{1}{q}}\psi_1).$$

\medskip
\noindent
Let us assume $(\bu,\bv)= (\theta^{\frac{1}{p}}\varphi_1,\theta^{\frac{1}{q}} \psi_1)$.\\
By the mean value theorem, the characterization of $\lambda_1$ and the definition of the normalized sequence $\bz_n=(\bu_n,\bv_n)$, we infer that

\begin{align*}
\displaystyle
c 
& \geq J(u_n,v_n) \\
& \geq - \into F(x,u_n,v_n) \, dx + \into h_1 u_n \, dx + \into h_2 v_n \, dx \nonumber \\
& = r_n^{\frac{1}{p}} \left[  \into h_1 \bu_n \, dx  \ -  \into  \int_0^1 F_s(x,\tau u_n,\tau v_n)\bu_n \, d \tau \, dx \right] \nonumber \\
&+ r_n^{\frac{1}{q}} \left[ \into h_2 \bv_n \, dx  \ -  \into  \int_0^1 F_t(x,\tau u_n,\tau v_n)\bv_n \, d \tau \, dx   \right] - \into F(x,0,0) \, dx. \nonumber
\end{align*}

\medskip
\noindent
Since $(\bu_n,\bv_n) \to (\theta^{\frac{1}{p}} \varphi_1, \theta^{\frac{1}{q}} \psi_1)$  strongly in $L^p(\Omega)\times L^q(\Omega)$, we have

\begin{align*}
& \displaystyle  \lim_{n \to \infty}\into h_1 \bu_n \, dx = \theta^{\frac{1}{p}} \into h_1 \varphi_1 \, dx,\\ &\displaystyle \lim_{n \to \infty} \int_0^1 F_s(x, \tau u_n, \tau v_n)\bu_n \, d \tau \, dx = \theta^{\frac{1}{p}} \into F_s^{++} \varphi_1 \, dx,
\end{align*}

\medskip
\noindent
and

\begin{align*}
& \displaystyle \lim_{n \to \infty} \into h_2 \bv_n \, dx  = \theta^{\frac{1}{q}} \into h_1 \psi_1 \, dx,\\
& \displaystyle \lim_{n \to \infty} \into \int_0^1 F_t(x, \tau u_n,\tau v_n)\bv_n \, d \tau \, dx = \theta^{\frac{1}{q}} \into F_t^{++} \psi_1 \, dx.
\end{align*}

\medskip
\noindent
Therefore, considering these limits in the previous inequality we deduce that

\begin{align*}
\displaystyle
c 
& \geq \theta^{\frac{1}{p}}r_n^{\frac{1}{p}} \left[  \into h_1 \varphi_1 \, dx  \ -  \into  F_s^{++} \varphi_1 \, dx + o(1) \right] \\
&+ \theta^{\frac{1}{q}}r_n^{\frac{1}{q}} \left[ \into h_2 \psi_1 \, dx  \ -  \into F_t^{++} \psi_1 \, dx  + o(1) \right]  - \into F(x,0,0) \, dx.
\end{align*}

\medskip
\noindent
In the case $p<q$, dividing by  $\theta^{\frac{1}{p}}r_n^{\frac{1}{p}}$ and passing to limit as $n \to \infty$, we have

$$ \displaystyle 0 \geq \into h_1 \varphi_1 \, dx  \ -  \into  F_s^{++} \varphi_1 \, dx,$$

\medskip
\noindent
which contradicts \eqref{LLpminoreq1}.\\
When $p=q$, similarly we obtain

$$ \displaystyle 0 \geq \into h_1 \varphi_1 \, dx + h_2 \psi_1 \, dx  -  \into  F_s^{++} \varphi_1 + F_t^{++} \psi_1  \, dx,$$ 

\medskip
\noindent
contradicting \eqref{LLpp1}.\\
If $p>q$, dividing by  $\theta^{\frac{1}{q}}r_n^{\frac{1}{q}}$ and passing to limit as $n \to  \infty$ we get

$$ \displaystyle 0 \geq \into h_2 \psi_1\, dx  \ -  \into  F_t^{++} \psi_1 \, dx,$$ 

\medskip
\noindent
contradicting \eqref{LLpmaggioreq1}.\\

These contradiction show that $(\bu,\bv)$ cannot be equal to $(\theta^{\frac{1}{p}}\varphi_1,\theta^{\frac{1}{q}}\psi_1).$

In a similar way if it would be $(\bu,\bv)=(-\theta^{\frac{1}{p}}\varphi_1,-\theta^{\frac{1}{q}}\psi_1)$, one obtains a contradiction with the hypotheses \eqref{LLpminoreq1}, \eqref{LLpp1} and \eqref{LLpmaggioreq1}, involving $F_s^{--}$ and $F_t^{--}$.

If now we assume $(\bu,\bv)=(\theta^{\frac{1}{p}}\varphi_1,-\theta^{\frac{1}{q}}\psi_1),$ we obtain

\begin{align*}
\displaystyle
c 
& \geq \theta^{\frac{1}{p}}r_n^{\frac{1}{p}} \left[  \into h_1 \varphi_1 \, dx  \ -  \into  F_s^{+-} \varphi_1 \, dx + o(1) \right] \\
&+ \theta^{\frac{1}{q}}r_n^{\frac{1}{q}} \left[  \into F_t^{+-} \psi_1 \, dx - \into h_2 \psi_1 \, dx  + o(1) \right]  - \into F(x,0,0) \, dx,
\end{align*}

\medskip
\noindent
and arguing as before we get a contradiction with  \eqref{LLpminoreq1}, \eqref{LLpp1} and \eqref{LLpmaggioreq1}.

Similarly, if it would be $(\bu,\bv)=(-\theta^{\frac{1}{p}}\varphi_1,\theta^{\frac{1}{q}}\psi_1)$, one obtains a contradiction with the hypotheses \eqref{LLpminoreq1}, \eqref{LLpp1} and \eqref{LLpmaggioreq1}, involving $F_s^{-+}$ and $F_t^{-+}$.

\end{proof}

\medskip
Under the assumption  \eqref{LLpminoreq2} if $p<q$, \eqref{LLpp2} if $p=q$ and \eqref{LLpmaggioreq2} if $p>q$, we will apply the following deformation result in order to prove  that the functional $J$ has critical points of saddle point type (see \cite[Theorem 3.4]{STRUWE}).

\begin{lemma}\label{lemmaStruwe}
Let $f$ be a functional of class $C^1$ on a Banach space $V$ that satisfies (PS) condition. Let $ c \in \mathbb{R}$ be a regular value of $f$ and $ \bar{\varepsilon} >0$. Then there is $ \varepsilon \in ( 0 , \bar{\varepsilon})$ and a continuous one parameter family of homeomorphisms $ \phi: V \times [0,1] \to V$ with the following properties:\\

\begin{enumerate}
\item[(i)] $ \phi ( z,\tau)=z$ if $\tau=0$ or if $ \left|  f(z) - c \right| \geq \bar{\varepsilon} $; \medskip
\item[(ii)]$ f(\phi ( z,\tau))$ is non-increasing in $\tau$ for any $ z \in V$; \medskip
\item[(iii)] If $f(z) \leq c + \varepsilon$, then $f( \phi(z,1)) \leq c - \varepsilon$.
\end{enumerate}
\end{lemma}

\medskip
\noindent
{\mbox {\it Proof of Theorem~\ref{mainresult}.~}} Let us consider assumption  \eqref{LLpminoreq1} if $p<q$, \eqref{LLpp1} if $p=q$ and \eqref{LLpmaggioreq1} if $p>q$. In these cases we have proved that the functional $J$ is coercive.\\

Now, by weak lower semicontinuity of both $\lVert \cdot \rVert_{1,p}$ and $\lVert \cdot \rVert_{1,q}$, and by Lebesgue theorem, we deduce that $J$ is sequentually weakly lower semicontinous, in fact for any sequence $\lbrace{ z_n \rbrace}_n=\lbrace{ (u_n,v_n) \rbrace}_n$ in $X$  weakly convergent to some $ z=(u,v) \in X$, we have

\begin{align*}
\displaystyle  
 J(z)  & =  \Phi(u,v) - \lambda_1 \Psi(u,v) - \into F(x,u,v) \, dx + \into h_1 u \, dx + \into h_2 v \, dx \\
& \leq \liminf_{n \to \infty} \frac{\alpha+1}{p} \into \left| \nabla u_n \right|^p \, dx + \liminf_{n \to \infty} \frac{\beta+1}{q} \into \left| \nabla v_n \right|^q \, dx \\
&  \quad - \lambda_1 \Psi(u,v) - \into F(x,u,v) \, dx + \into h_1 u \, dx + \into h_2 v \, dx \\
& \leq \liminf_{n \to \infty} \Phi(u_n,v_n)  + \lim_{n \to \infty} \left[ - \lambda_1 \Psi(u_n,v_n) - \into F(x,u_n,v_n) \, dx + \into h_1 u_n \, dx + \into h_2 v_n \, dx   \right] \\
& = \liminf_{n \to \infty} J(z_n).
\end{align*} 

\noindent
This and the coercivity of $J$ imply by Weierstrass Theorem  that $J$ has a global minimum point, i.e. there exists a weak solution for problem \eqref{Syst0}.\\
Let us now consider assumption  \eqref{LLpminoreq2} if $p<q$, \eqref{LLpp2} if $p=q$ and \eqref{LLpmaggioreq2} if $p>q$.\\
We will show that

\begin{align}\label{limiti}
\displaystyle
& \lim_{\theta \to + \infty} J(\theta^{\frac{1}{p}} \varphi_1,\theta^{\frac{1}{q}} \psi_1)=-\infty, & \lim_{\theta \to + \infty} J(-\theta^{\frac{1}{p}} \varphi_1,-\theta^{\frac{1}{q}} \psi_1)=-\infty, \nonumber \\
&\\
& \lim_{\theta \to + \infty} J(\theta^{\frac{1}{p}} \varphi_1,-\theta^{\frac{1}{q}} \psi_1)=-\infty, & \lim_{\theta \to + \infty} J(-\theta^{\frac{1}{p}} \varphi_1,\theta^{\frac{1}{q}} \psi_1)=-\infty.  \nonumber 
\end{align}

\medskip
\noindent
By characterization of $\lambda_1$ and the mean value formula, for any $\theta >0$ we infer that

\begin{align*}
\displaystyle
J(\theta^{\frac{1}{p}} \varphi_1,\theta^{\frac{1}{q}} \psi_1)
& = \theta^{\frac{1}{p}} \into h_1 \varphi_1 \, dx + \theta^{\frac{1}{q}} \into h_2 \psi_1 \, dx - \into F(x,\theta^{\frac{1}{p}} \varphi_1,\theta^{\frac{1}{q}} \psi_1) \, dx \\
& = \theta^{\frac{1}{p}} \left[ \into h_1 \varphi_1 \, dx  -  \into \int_0^1 F_s(x,\tau \theta^{\frac{1}{p}} \varphi_1,\tau \theta^{\frac{1}{q}} \psi_1) \varphi_1 \, d \tau \, dx \right] \\
& + \theta^{\frac{1}{q}} \left[ \into h_2 \psi_1 \, dx  -  \into \int_0^1 F_t(x,\tau \theta^{\frac{1}{p}} \varphi_1,\tau \theta^{\frac{1}{q}} \psi_1) \psi_1 \, d \tau \, dx \right] \\
& - \into F(x,0,0) \, dx.
\end{align*}

\medskip
\noindent
By \textbf{(F2)}, Fubini-Tonelli and Lebesgue theorems, it follows

$$ \displaystyle \lim_{\theta \to + \infty} \into \int_0^1 F_s(x,\tau \theta^{\frac{1}{p}} \varphi_1,\tau \theta^{\frac{1}{q}} \psi_1) \varphi_1 \, d \tau \, dx = \into F_s^{++} \varphi_1  \, dx;$$

$$ \displaystyle \lim_{\theta \to + \infty} \into \int_0^1 F_t(x,\tau \theta^{\frac{1}{p}} \psi_1,\tau \theta^{\frac{1}{q}} \psi_1) \psi_1 \, d \tau \, dx = \into F_t^{++} \psi_1  \, dx.$$

\medskip
\noindent
Therefore,

\begin{align*}
\displaystyle
J(\theta^{\frac{1}{p}} \varphi_1,\theta^{\frac{1}{q}} \psi_1)
& = \theta^{\frac{1}{p}} \left[ \into h_1 \varphi_1 \, dx  - \into F_s^{++} \varphi_1  \, dx  + o(1)\right] \\
& + \theta^{\frac{1}{q}} \left[ \into h_2 \psi_1 \, dx  -  \into F_t^{++} \psi_1  \, dx + o(1) \right] \\
& - \into F(x,0,0) \, dx,
\end{align*}

\medskip
\noindent
by which

\begin{equation*}
\displaystyle 
\lim_{\theta \to + \infty} J(\theta^{\frac{1}{p}} \varphi_1,\theta^{\frac{1}{q}} \psi_1)=- \infty.
\end{equation*}

\medskip
\noindent
Arguing in a similar way, considering the other inequalities in assumptions  \eqref{LLpminoreq2}, \eqref{LLpp2} and \eqref{LLpmaggioreq2} involving $F_s^{--}$ and $F_t^{--}$,  $F_s^{+-}$ and $F_t^{+-}$, and $F_s^{-+}$ and $F_t^{-+}$ respectively, we get the other limits in \eqref{limiti}.

Let $\lambda_2$ be the eigenvalue of \eqref{EigenDeThélin0} defined in \eqref{lambdakSfera} and let us  set

$$\Lambda_2:=\{ (u,v) \in X \; : \; \Phi(u,v) \geq \lambda_2 \Psi(u,v)     \}.$$

\medskip
\noindent
By $\textbf{(F1)}$ and the mean value formula, we have

$$ \displaystyle \left| F(x,s,t) \right| \leq M|s| + M|t| + \left|F(x,0,0) \right|.$$

\medskip
\noindent
Denoting by $\lambda_{1,p}$ and $\lambda_{1,q}$ the first eigenvalue on $\Omega$ under homogeneous Dirichlet boundary condition of $\Delta_p$ and $\Delta_q$ respectively, by H\"older and Poincaré inequality and for any $(u,v) \in X$ we have

\begin{align*}
\displaystyle
& \left|  - \into F(x,u,v) \, dx + \into h_1 u \, dx + \into h_2 v \, dx \right| \\
\leq & \into  \left| F(x,0,0) \right|  \, dx + M \into \left| u  \right| \, dx +  M \into \left| v  \right| \, dx +  \into \left| h_1 \right| \left| u  \right| \, dx + \into \left| h_2 \right| \left| v \right| \, dx \nonumber \\
\leq & \into \left|F(x,0,0)\right| \, dx + \left( M \left| \Omega \right|^{\frac{1}{p'}} +  \lVert h_1 \rVert_{p'}  \right) \lVert u \rVert_p +  \left( M \left| \Omega \right|^{\frac{1}{q'}} + \lVert h_2 \rVert_{q'} \right) \lVert v \rVert_q \nonumber \\
\leq &  \into \left|F(x,0,0)\right| \, dx + \frac{\left( M \left| \Omega \right|^{\frac{1}{p'}} +  \lVert h_1 \rVert_{p'}  \right)}{ \lambda^{\frac{1}{p}}_{1,p}} \lVert u \rVert_{1,p} +  \frac{\left( M \left| \Omega \right|^{\frac{1}{q'}} + \lVert h_2 \rVert_{q'} \right)}{\lambda^{\frac{1}{q}}_{1,q}} \lVert v \rVert_{1,q},\nonumber
\end{align*}

\noindent
On the other hand, for any $(u,v) \in \Lambda_2$ we get

\begin{align*}
\displaystyle \Phi(u,v) - \lambda_1 \Psi(u,v) 
& \geq \left( 1 - \frac{\lambda_1}{\lambda_2}\right)  \Phi(u,v)  \\
& = \frac{\alpha+1}{p}\left( 1 - \frac{\lambda_1}{\lambda_2}\right)\lVert u \rVert^p_{1,p}  + \frac{\beta+1}{q}\left( 1 - \frac{\lambda_1}{\lambda_2}\right)\lVert v \rVert^q_{1,q}. \nonumber
\end{align*}

\noindent
Altogether, for any $(u,v) \in \Lambda_2$ we have

\begin{align*}
\displaystyle
J(u,v)
& \geq \frac{\alpha+1}{p}\left( 1 - \frac{\lambda_1}{\lambda_2}\right)\lVert u \rVert^p_{1,p}  - \frac{\left( M \left| \Omega \right|^{\frac{1}{p'}} +  \lVert h_1 \rVert_{p'}  \right)}{ \lambda^{\frac{1}{p}}_{1,p}} \lVert u \rVert_{1,p} \\
& + \frac{\beta+1}{q}\left( 1 - \frac{\lambda_1}{\lambda_2}\right)\lVert v \rVert^q_{1,q} - \frac{\left( M \left| \Omega \right|^{\frac{1}{q'}} + \lVert h_2 \rVert_{q'} \right)}{\lambda^{\frac{1}{q}}_{1,q}} \lVert v \rVert_{1,q}\\
& - \into \left|F(x,0,0)\right| \, dx.
\end{align*}

\noindent
Since $\lambda_1$ is isolated, we have $\lambda_2 > \lambda_1$ and so we deduce that $J$ is coercive on $\Lambda_2$, and in particular it is bounded from below.
Considering also \eqref{limiti}, we have shown that there exists $\Theta >0$ such that, setting

$$ \displaystyle E_1^{\Theta}= \bigg\{ ( \left|  \theta \right|^{\frac{1}{p}} \varphi_1 ,  \left|  \theta \right|^{\frac{1}{q}} \psi_1 )  \text{sgn} (\theta)  : \left|\theta \right| \geq \Theta \bigg\}   \bigcup  \bigg\{ (- \left|  \theta \right|^{\frac{1}{p}} \varphi_1 ,  \left|  \theta \right|^{\frac{1}{q}} \psi_1 )  \text{sgn} (\theta)  :  \left|\theta \right| \geq \Theta \bigg\},$$ 

\medskip
\noindent
we have

$$ \displaystyle \gamma:= \sup_{(u,v) \in E_1^{\Theta}} J(u,v) < \inf_{(u,v) \in \Lambda_2} J(u,v):= \delta.$$

\medskip
\noindent
Consider the family of mappings 

$$ 
\Gamma:= \{ h \in C( [-1,1] ,X) \, : \, h(1)=-h(-1)\in E_1^{\Theta}  \}.
$$

\medskip
\noindent
Notice that $\displaystyle \tilde{h}(\tau):= ( \tau \Theta^{\frac{1}{p}} \varphi_1 , \tau \Theta^{\frac{1}{q}} \psi_1 )$,  $\tau \in [-1,1]$, belongs to $\Gamma$ and thus $\Gamma\not=\emptyset$.\\
Let us show that $h(B_1) \cap \Lambda_2 \neq \emptyset$ for every $h \in \Gamma$, where $B_1:=[-1,1]$.\\ Indeed, this is immediate if $(0,0) \in h(B_1)$. When $(0,0) \not \in h(B_1)$, we consider the odd map $\pi: X \setminus (0,0) \to \mathcal{M}$ given by 

$$ \displaystyle \pi(u,v):= \left( \frac{u}{\Phi(u,v)^{\frac{1}{p}}}   ,  \frac{v}{\Phi(u,v)^{\frac{1}{q}}}     \right), \quad (u,v)\in X\setminus (0,0),$$

\medskip
\noindent
and the map $\alpha_0:S^1 \to \mathcal{M}$ defined for every $(\xi_1,\xi_2)\in S^1$ by

$$ 
\displaystyle
\alpha_0(\xi_1,\xi_2):=
\begin{cases}
\; \; \, \pi ( \, h(\xi_1) \, ) & \text{ if } \xi_2 \geq 0, \bigskip \\
- \pi  ( \, h (-\xi_1) \, ) & \text{ if } \xi_2 < 0.
\end{cases}
$$

\medskip
\noindent
Observe that by the definition of $\alpha_0$, we have $\alpha_0(-\xi_1,-\xi_2) = -\alpha_0(\xi_1,\xi_2)$ for every $(\xi_1,\xi_2)\in S^1$ with $\xi_2\not=0$. 
On the other hand, 
since $h(1)=-h(-1)$, the oddness of $\pi$ implies that  $\alpha_0(-1,0) = -\alpha_0(1,0)$ and, consequently,  $\alpha_0$ is odd in all $S^1$. In addition, using again that $h(1)=-h(-1)$, the map $\alpha_0$ is also continuous. In particular, $\alpha_0(S^{1})\in \mathcal{M}_2$, where  $\mathcal{M}_2$ is defined in Proposition \ref{prop 2.9}.\\
Recalling that

$$ \displaystyle \lambda_2:= \frac{1}{c_2} \qquad \text{where} \quad c_2:= \sup_{A \in \mathcal{M}_2} \min_{(u,v) \in A} Q_{ \mathcal{M}}(u,v)=\sup_{A \in \mathcal{M}_2} \min_{(u,v) \in A} \frac{\Psi(u,v)}{\Phi(u,v)},$$

\medskip
\noindent
we have

$$
\min_{(u,v) \in\alpha_0(S^1)} \frac{\Psi(u,v)}{\Phi(u,v)} \leq c_2.
$$

\medskip
\noindent
Let $(u_0,v_0) \in \alpha_0 (S^1)$ be such that 
 $ \displaystyle \frac{\Psi(u_0,v_0) }{\Phi (u_0,v_0)} \leq c_2$, i.e.,
 $ \displaystyle \Phi(u_0,v_0)\geq \lambda_2 \Psi(u_0,v_0)$ and 
 so $(u_0,v_0) \in \alpha_0(S^1) \cap \Lambda_2$. 
By definition of $\alpha_0$, this implies that
 there exists $\tau \in [-1,1]$ such that $\pi(h(\tau)) \in \Lambda_2$. 
 Now we observe that the $(p,q)$-homogeneity of the functionals 
 $\Phi$ and $\Psi$ implies that  $h(\tau) \in \Lambda_2$ 
and thus $h(B_1) \cap \Lambda_2 \neq \emptyset$.\\
In particular, 

$$ 
\displaystyle  c:= \inf_{h \in \Gamma}  \max_{\tau \in B_1} J (h(\tau)) \geq \delta > \gamma.
$$

\medskip
\noindent
Now, we  show that 
$c$
is a critical value of $J$.\\
Indeed, let us assume by contradiction that $c$ is a regular value of $J$. 
Fixing $0 < \bar{\varepsilon} < c - \gamma,$ by Lemma \ref{lemmaStruwe} there 
exists $ \varepsilon \in (0, \bar{\varepsilon})$ and a family of homeomorphisms $ \phi: X \times [0,1] \to X$ satisfying the following properties:\\

\begin{enumerate}
\item[$(i)$] $ \phi ( z,\tau)=z$ if $\tau=0$ or if $ \left|  J(z) - c \right| \geq \bar{\varepsilon} $; \medskip
\item[$(ii)$]$ J(\phi ( z,\tau))$ is non-increasing in $\tau$ for any $ z \in X$; \medskip
\item[$(iii)$] If $J(z) \leq c + \varepsilon$, then $J( \phi(z,1)) \leq c - \varepsilon$. 
\end{enumerate} 

\medskip
\noindent
For any $z \in E_1^{\Theta}$, we have $J(z) \leq \gamma < c - \bar{\varepsilon}$ 
and so by property $(i)$ we infer that $\phi$ leaves the set $E_1^{\Theta}$ fixed. 
By definition of $c$ there exists $h \in \Gamma$ such that 
$ \max_{\tau \in [-1,1]} J (h(\tau)) < c + \varepsilon$. Let us define 
$\hat{h}(\cdot):= \phi( h(\cdot) , 1).$ 
Since $h(\pm 1)\in E_1^{\Theta}$, 
we get  $\hat{h}(\pm 1)= \phi( h(\pm 1) , 1) = h(\pm 1)$, from which we deduce 
that $\hat{h}(-1)= -\hat{h}(1)$ and hence
 $\hat{h} \in \Gamma$. Now,  
property $(iii)$ gives
$ \max_{\tau \in [-1,1]} J (\hat{h}(\tau)) \leq c - \varepsilon$, 
which  contradicts the definition of $c$.

\qed

\bigskip

\section*{Appendix C}
\addcontentsline{toc}{section}{Appendix C}
\label{app:C}

\bigskip
\noindent
In this Section we analyze the definitions contained in the paper \cite{BON} of Bonnet.\\
Let $E$ be a Banach space and let us denote with $S$ the $C^1$ manifold defined by

$$ \displaystyle S:= \{ x \in E \, : \, g(x)=1\},$$

\medskip
\noindent
where $g \in C^1(E, \mathbb{R})$ and $1$ is a regular value for $g$.
For any $x \in S$, we denote with $T_x S$ the tangent space of $S$ at $x$.
Let $f$ be a $C^1$ functional defined on a neighborhood of $S$.
The restriction $f_{S}:= f\big|_{S}$ of $f$ to $S$ is a $C^1$ functional on $S$ and the derivative $f_S'$ of $f_S$ is a map defined on the tangent bundle $TS$ and is such that

$$ \displaystyle  f_S' (x) = f'(x)\big|_{T_x S} \qquad \forall \, x \in S.$$

\medskip
\noindent
In particular, for any $x \in S$ we have $ f_S' (x) \in T_x^*S$, where $T_x^*S$ denotes the dual space of $ T_xS$, endowed with the norm $ \lVert \cdot \rVert_x^*$. Denoting with $\lVert \cdot \rVert^*$ the norm of $E^*$, by \cite[Proposition 3.54]{PAO} it follows that

\begin{align}\label{normaspaziotangenteBON}
\displaystyle   \lVert f_{S}'(x) \rVert^*_x= \min_{\mu \in \mathbb{R}} \lVert \, f'(x)   - \mu g'(x) \, \rVert^*.
\end{align} 

\begin{definition}
We say that $f_S$ satisfies the Palais-Smale condition at the level $c$ (henceforth denoted by $(PS)_c$) if any sequence $\{x_n\}_n \subset S$ such that

\begin{align*}
\displaystyle & f_S(x_n) \to c \qquad \text{as } n \to + \infty,\\
              &\\
              & f_S'(x_n) \to 0 \qquad \text{as } n \to + \infty,
\end{align*}

\medskip
\noindent
possesses a convergent subsequence.
\end{definition}

\begin{lemma}
Suppose that $c$ is a regular value of $f_S$ and that $f_S$ satisfies $(PS)_c$. Then there exist $\delta,\varepsilon>0$ such that

\begin{equation}\label{conseguenzaPSa}
\displaystyle \forall x \in S, \qquad |f(x) - c | \leq \varepsilon \quad \Longrightarrow \quad \lVert f_S'(x) \rVert_x^* > \delta.
\end{equation}
\end{lemma}

\begin{proof}
By contradiction, $ \forall \delta >0$ and $\forall \varepsilon >0$ there exists $x:=x(\delta,\varepsilon) \in S$ such that

$$ \displaystyle  |f(x) - c | \leq \varepsilon \qquad \text{and} \qquad \lVert f_S'(x) \rVert_x^* \leq \delta.$$

\medskip
\noindent
In particular, for any $n \in \mathbb{N}$ there exists $x_n \in S$ such that

$$ \displaystyle |f(x_n) - c | \leq \frac{1}{n} \qquad \text{and} \qquad \lVert f_S'(x_n) \rVert_{x_n}^* \leq \frac{1}{n}.$$

\medskip
\noindent
Therefore, $\{x_n\}_n$ is a Palais-Smale sequence at level $c$. Since $f_S$ satisfies $(PS)_c$, there exists $\bar{x} \in S$ such that $f_S(\bar{x})=c$ and $f_S'(\bar{x})=0$, that is a contradiction with the fact that $c$ is a regular value.

\end{proof}

\bigskip
\noindent
Now, let $a>0$ and let us consider

$$ \displaystyle F_{a}:= \{ S_{b} \, : \, b \in [0,a] \, \} $$

\medskip
\noindent
and

$$ \displaystyle F^*_{a}:= \{ S_{b} \, : \, b\in [-a, 0] \, \}, $$
\medskip
\noindent
where

$$ \displaystyle S_{b}= \{ x \in E \, : \, g(x)=1+ b\}.$$

\begin{definition}
Suppose that $c \in \mathbb{R}$  is a regular value of $f_S$, and let $a>0$. A family $F_a$ is said to be admissible for $f$ at $c$ if $f$ is defined in $S_{b}$ for any $b \in [0,a]$ and if there exist $\mu,\delta,\varepsilon>0$ such that

\begin{equation}\label{conseguenzaPSafgaccoppiate}
\displaystyle \forall b \in [0,a], \forall x \in S_{b}, \quad |f(x) - c| \leq \varepsilon \quad \Rightarrow \quad ( \lVert f_{S_{b}}'(x) \rVert_x^* > \delta \; \text{and} \; \lVert g'(x) \rVert > \mu).
\end{equation}

\medskip
\noindent
Similar definitions hold if we consider $F^*_a$.

\end{definition}

\begin{definition}
Let $c \in \mathbb{R}$. We say that $f$ and $g$ satisfy a coupled Palais-Smale condition at level $c$ on $S$ (henceforth denoted by $\widehat{(PS)}_{S,c}$) if any sequence $\{x_n\}_n$ such that 

\begin{itemize}
\item[•] $x_n \in S_{\varepsilon_n}$ where $\varepsilon_n >0$ and $\varepsilon_n \to 0$ as $n \to \infty$;\\

\item[•] $f(x_n) \to c$ as $n \to \infty$;\\

\item[•] $\lVert f_{S_{\varepsilon_n}}'(x_n) \rVert_{x_n}^* \to 0 \text{ or } \lVert g'(x_n) \rVert \to 0$ as $n \to \infty$;
\end{itemize}

\medskip
\noindent
possesses a convergent subsequence.

\end{definition}

\begin{lemma}\label{LemmaPSaccoppiata}
Suppose that $c$ is a regular value of $f_S$ and that $f$ and $g$ satisfy $\widehat{(PS)}_{S,c}$.  Then there exists $a>0$ such that $F_{a}$ is admissible for $f$ at $c$.
\end{lemma}

\begin{proof}
By contradiction $ \forall a >0$ and $\forall \mu,\delta,\varepsilon>0$ there exist $b \in [0,a]$ and $\bar{x} \in S_{b}$ such that 

$$ \displaystyle \left|f(\bar{x}) - c\right| \leq \varepsilon \qquad \text{and} \qquad ( \lVert f_{S_{b}}'(\bar{x}) \rVert_{\bar{x}}^* \leq \varepsilon  \quad \text{ or } \quad \lVert g'(\bar{x}) \rVert \leq \mu).$$

\medskip
\noindent
In particular, for any $n \in \mathbb{N}$ there exist  $b_n \in \left[0, \frac{1}{n} \right]$ and $x_n \in S_{b_n}$ such that

$$ \displaystyle |f(x_n) - c | \leq \frac{1}{n} \qquad \text{and}  \qquad \left( \lVert f_{S_{b_n}}'(x_n) \rVert_{x_n}^* \leq \frac{1}{n}  \quad \text{ or } \quad \lVert g'(x_n) \rVert \leq \frac{1}{n} \right).$$

\medskip
\noindent
Notice that $\beta_n \to 0$,  $f(x_n) \to a$  and that

$$ \displaystyle \lVert f_{S_{b_n}}'(x_n) \rVert_{x_n}^* \to 0 \quad \text{ or } \quad \lVert g'(x_n) \rVert \to 0 \qquad \text{as } n \to \infty.$$

\medskip
\noindent
Since $\{x_n\}_n$ is sequence of the type of previous definition and $f$ and $g$ satisfy $\widehat{(PS)}_{S,c}$, we infer that $ \{x_n\}_n$ has a converging subsequence. Hence there exists $\bar{x} \in S$ such that $f(\bar{x})=c$ and 

$$ \displaystyle  f_{S}'(\bar{x})  = 0 \quad \text{ or } \quad  g'(\bar{x}) = 0.$$

\medskip
\noindent
If we consider $g'(\bar{x}) = 0,$ we get a contradiction with the fact that $1$ is a regular value of $g$.\\
If we consider $ f_{S}'(\bar{x})  = 0$, we get a contradiction with the fact that $c$ is a regular value of $f_S$.\\
In any case, we have a contradiction.

\end{proof}

\medskip
\noindent
Analogous definitions and results hold if we consider the family $F_a^*$.\\

\medskip
\noindent
We can now recall \cite[Theorem 2.5]{BON}.

\begin{theorem}\label{BONNETTheorem}

Let $S:= \{ x \in E \, : \, g(x)=1 \}$ be a $C^1$ submanifold of a Banach space $E$, where $g \in C^1(E, \mathbb{R})$ and $1$ is a regular value of $g$. Suppose that $f$ is a $C^1$ functional on a neighborhood of $S$, $c \in \mathbb{R}$ is a regular value of $f_S$ and that there exists $a>0$ such that $F_{a}$ or $F_{a}^*$ is admissible for $f$ at the value $c$. Then there exists $\hat{\varepsilon} >0$ such that for all $0 < \varepsilon < \hat{\varepsilon}$ there exists an homeomorphism $\eta$ of $S$ onto $S$ such that:\\

\begin{enumerate}
\item $ \eta (x)=x$  if $f(x) \not \in [ c - \hat{\varepsilon}, c + \hat{\varepsilon}]$; \medskip
\item $f(\eta(x)) \geq f(x)$ for all $x \in S$; \medskip
\item $f(\eta(x)) \geq c + \varepsilon$ for all $x$ such that $f(x) \geq c - \varepsilon$; \medskip
\item If $S$ is symmetric $(S=-S)$ and if $f$ is even, then $\eta$ is odd.
\end{enumerate}
\end{theorem}

\medskip
\noindent
As a consequence of Lemma \ref{LemmaPSaccoppiata} and previous theorem, we get the following.

\begin{theorem}

Let $S:= \{ x \in E \, : \, g(x)=1 \}$ be a $C^1$ submanifold of a Banach space $E$, where $g \in C^1(E, \mathbb{R})$ and $1$ is a regular value of $g$. Suppose that $f$ is a $C^1$ functional on a neighborhood of $S$, $c \in \mathbb{R}$ is a regular value of $f_S$ and $f$ and $g$ satisfy $\widehat{(PS)}_{S,c}$. Then there exists $\hat{\varepsilon} >0$ such that for all $0 < \varepsilon < \hat{\varepsilon}$ there exists an homeomorphism $\eta$ of $S$ onto $S$ satisfying  properties $(1)-(4)$ of Theorem \ref{BONNETTheorem}.
\end{theorem}

}

\bigskip

\chapter{$N$-dimensional Euclidean Onofri Inequality in Weighted Sobolev Space}\label{CAPITOLOONOFRI}

\medskip
\noindent
In this Chapter, we focus on the $N$-dimensional Euclidean Onofri inequality established by Del Pino and Dolbeault in \cite{DD3} for smooth functions with compact support. We prove that the $N$-Laplacian operator, $\Delta_N$, is intrinsically linked to this inequality through its connection with the Liouville equation in $\mathbb{R}^N$. Furthermore, we extend the Euclidean Onofri inequality to a suitable weighted Sobolev space. We also establish an equivalence result with the logarithmic Moser-Trudinger inequality on the balls of $\mathbb{R}^N$; in the specific case $N=2$, we highlight the connection between the $2$-dimensional Euclidean Onofri inequality and the two-dimensional unit sphere $\mathbb{S}^2$, by using the stereographic projection from the North pole.

\medskip
\noindent
This Chapter is based on the works \cite{BCM,BCM2}.

\bigskip

\section{Euclidean Onofri Inequality in $C_0^{\infty}(\mathbb{R}^N)$}

\bigskip
\noindent
If $\Omega$ is a bounded domain of $\mathbb{R}^N$ with $N \geq 2$, it is well known that the Sobolev space $W_0^{1,N}(\Omega)$ is continuously embedded into $L^q(\Omega)$ for any
$1 \leq q < \infty$. In particular, there exists a constant $c:=c(q,N,\Omega) $ such that 

$$ \displaystyle \lVert u \rVert^q_{q} \leq c \,  \lVert \nabla u \rVert^q_{N} \qquad \text{ for any } u \in W_0^{1,N}(\Omega).$$

\medskip
\noindent
Denoting by $\mathcal{H}:=\{u \in W_0^{1,N}(\Omega) \ | \  \lVert \nabla u \rVert_N \leq 1 \}$, it is equivalent to state that

$$ \displaystyle \sup_{u \in \mathcal{H}} \, \into |u|^q \, dx \leq c(q,N,\Omega) \qquad \text{ for any } 1 \leq q < \infty. $$

\medskip
\noindent
Even if $W_0^{1,N}(\Omega) \hookrightarrow L^q(\Omega)$ for any $1 \leq q < \infty$, we have $W_0^{1,N}(\Omega) \not \hookrightarrow L^{\infty}(\Omega)$. In fact if $\Omega$ is for instance the unit ball, then a counterexample is given by the function $u(x)= \ln (1 - \ln |x|)$. 

However, in \cite{Tru} Trudinger (see also  Yudovich  \cite{Yudo} and Pohozaev \cite{poho}) improved previous power type estimate obtaining an exponential type one. He proved that there exist $\alpha >0$ and a constant $c:=c(N)>1$ such that

\begin{equation*}
\displaystyle \sup_{u \in \mathcal{H}} \, \frac{1}{|\Omega|}\int_{\Omega} e^{\alpha \left|u \right|^\frac{N}{N-1}}\,dx \le c.
\end{equation*}

\medskip
\noindent
Trudinger's proof \cite{Tru} makes use of power series expansion of the exponential function and some Sobolev estimates for the individual terms carefully observing the dependence on the exponent of the expansion. An improvement on Trudinger's result was achieved by Moser in \cite{M}, who found the best number $\alpha_N$ 
such that the above statement holds for any $\alpha \leq \alpha_N$ and it is false for $\alpha > \alpha_N$.
More precisely, denoting by $\alpha_N:= N \, \omega_{N-1}^\frac{1}{N-1}$ where $\omega_{N-1}$ denotes the measure of the surface of the unit sphere of $\mathbb{R}^N$, Moser proved in \cite[Theorem 1]{M} that there exists a constant $C_N > 1$ such that

\begin{equation}\label{supsudominilimitatiINTRODUZIONE}
\sup_{u \in \mathcal{H}} \,  \frac{1}{|\Omega|} \int_{\Omega} e^{\alpha \left|u \right|^\frac{N}{N-1}}\,dx \le C_N \qquad \forall \, \alpha \leq  \alpha_N.
\end{equation}

\medskip
\noindent
The integral on the left is actually finite for any positive $\alpha$ and for any $u \in \mathcal{H}$, but if $\alpha > \alpha_N$ it can be made arbitrarily large by an appropriate choice of $u$, hence the supremum on $\mathcal{H}$ is equal to $+ \infty.$\\

Moser's proof of the above result relies on symmetrization techniques and a change of variable that reduces the $N$-dimensional problem to a $one$-dimensional one. He obtained with the same idea a corresponding inequality of \eqref{supsudominilimitatiINTRODUZIONE} also on the two-sphere $\mathbb{S}^2$ endowed with the standard metric $g_0$. In particular, Moser proved  in \cite[Theorem 2]{M} that there exists a constant $S > 1$ such that

\begin{equation}\label{supsuS2INTRODUZIONE}
\sup_{u \in \bar{\mathcal{H}} } \,  \frac{1}{4\pi}\int_{\Sm^2} e^{4\pi u^2 }\,d \nu_{g_0} \le S,
\end{equation}

where $\bar{\mathcal{H}}:= \{u \in W^{1,2}(\mathbb{S}^2,g_0) \ | \ \bar{u}=0, \ \lVert \nabla_{g_0} u \rVert_2 \leq 1 \}$ and $ \displaystyle \bar{u}:= \frac{1}{4 \pi} \int_{\Sm^2} u \,d \nu_{g_0}$.

\medskip
\noindent
Although the proof of the two previous results is based on the same techniques, as pointed out by Moser in his paper it seems impossible to deduce one theorem from the other, therefore the equivalence between \cite[Theorem 1]{M} and \cite[Theorem 2]{M} is still an open problem.

\medskip
It is also important to mention the work of Adams \cite{Adams} for functions of higher order derivatives, and the work of Fontana \cite{Fontana} for compact smooth Riemannian manifolds.

\medskip
\noindent
Starting from \eqref{supsudominilimitatiINTRODUZIONE} and \eqref{supsuS2INTRODUZIONE}, and applying Young's inequality, one derives the so called logarithmic Moser-Trudinger inequalities. Specifically, \eqref{supsudominilimitatiINTRODUZIONE} yields the following \textit{logarithmic Moser-Trudinger inequality on bounded domains of $\mathbb{R}^N$}:

\begin{equation}\label{MTdomininonottimale}
\displaystyle \ln \left( \frac{1}{|\Omega|} \int_{\Omega} e^u \, dx \right) \leq  \frac{1}{\widetilde{\omega_N}} \int_{\Omega} | \nabla u |^N \, dx + \ln C_N,
\end{equation}

\medskip
\noindent
for any $u \in W_0^{1,N}(\Omega)$, where $ \displaystyle \widetilde{\omega_N}:= N^N \left( \frac{N}{N-1} \right)^{N-1} \omega_{N-1}$ and $C_N>1$ is the constant in \eqref{supsudominilimitatiINTRODUZIONE}.\\
Based on the work of Carleson and Chang \cite{CarlesonChang}, the following improvement of the previous inequality can be established for the unit ball $B_1 \subset \mathbb{R}^N$:

\begin{equation}\label{logCCintroduzione}
\displaystyle \ln \left( \frac{1}{|B_1|} \int_{B_1} e^{u}   dx  \right) < \frac{1}{\widetilde{\omega_N}} \int_{B_1} |\nabla u|^N  dx + \sum_{k=1}^{N-1} \frac{1}{k},
\end{equation}

\medskip
\noindent
for any $u \in  W^{1,N}_0(B_1)$.

\medskip
\noindent
Moreover, inequality \eqref{logCCintroduzione} is sharp in the sense that the functional $J:W^{1,N}_{0}(B_1) \to \mathbb{R}$ defined as 

\begin{equation}\label{Jfunctional}
\displaystyle J(u) := \frac{1}{\widetilde{\omega_N}} \int_{B_1} | \nabla u |^N \, dx  - \ln \left( \frac{1}{V_N} \int_{B_1} e^u \, dx \right),
\end{equation}

\medskip
\noindent
satisfies 

$$
\inf_{W^{1,N}_{0}(B_1)} J = - \sum_{k=1}^{N-1} \frac{1}{k}.
$$

\medskip
\noindent
Regarding the \textit{logarithmic Moser–Trudinger inequality on $\mathbb{S}^2$}, it reads as follows:

\begin{equation}\label{MTS2nonottimale}
\displaystyle  \ln \left( \frac{1}{4 \pi} \int_{\mathbb{S}^2} e^u \, d \nu_{g_0} \right) \leq  \frac{1}{16 \pi} \int_{\mathbb{S}^2} | \nabla_{g_0} u |_{g_0}^2 \, d \nu_{g_0} + \frac{1}{4 \pi} \int_{\mathbb{S}^2}u \, d \nu_{g_0} + c
\end{equation}

\medskip
\noindent
for any $\displaystyle u \in W^{1,2}(\mathbb{S}^2,g_0),$ where $c:= \ln S$ and  $S>1$ is the constant appearing in \eqref{supsuS2INTRODUZIONE}.\\
In \cite{O}, Onofri proved that the best constant $c$ for which the previous inequality holds actually vanishes. This leads to the so-called \textit{logarithmic Onofri–Moser–Trudinger inequality}, or simply the \textit{Onofri inequality}:

\begin{equation}\label{bestOnofriINTRODUZIONE}
\displaystyle  \ln \left( \frac{1}{4 \pi} \int_{\mathbb{S}^2} e^u \, d \nu_{g_0} \right) \leq  \frac{1}{16 \pi} \int_{\mathbb{S}^2} | \nabla_{g_0} u |_{g_0}^2 \, d \nu_{g_0} + \frac{1}{4 \pi} \int_{\mathbb{S}^2}u \, d \nu_{g_0}
\end{equation}

\medskip
\noindent
for any $\displaystyle u \in W^{1,2}(\mathbb{S}^2,g_0)$.\\
By means of the stereographic projection from the North Pole of $\mathbb{S}^2$, the Onofri inequality \eqref{bestOnofriINTRODUZIONE} can be equivalently reformulated as the \textit{two-dimensional Euclidean Onofri inequality}:

\begin{equation}\label{EuclideanOnofri2DIntroduzione}
\displaystyle  \ln \left(\int_{\mathbb{R}^2} e^u \, d \mu_2 \right) \leq  \frac{1}{16 \pi} \int_{\mathbb{R}^2} | \nabla u |^2 \, dx + \int_{\mathbb{R}^2} u \, d \mu_2
\end{equation}

\medskip
\noindent
for any $ u \in  W_{\mu_2}(\mathbb{R}^2):= \{ u \in L^1(\mathbb{R}^2,d\mu_2) \, : \, |\nabla u| \in L^2(\mathbb{R}^2,dx) \}$, where

$$ \displaystyle \mu_2(x):= \frac{1}{ \pi \left( 1+ |x|^2 \right)^2} \qquad \text{and} \quad d \mu_2:=\mu_2(x) \, dx. $$

\medskip
\noindent
Extension of \eqref{EuclideanOnofri2DIntroduzione} to any dimension $N \geq 2$, was obtained by del Pino and Dolbeault in \cite{DD3}. Denoting by $V_N$  the volume of the unit ball in $\R^N$, they obtained the following \textit{$N$-dimensional Euclidean Onofri inequality} for any $u \in C_0^{\infty}(\mathbb{R}^N)$:

\begin{equation}\label{EuclideanOnofriNdimensional}
\displaystyle  \ln \left( \int_{\mathbb{R}^N} e^u \, d\mu_N \right) \leq  \frac{1}{\widetilde{\omega_N}} \int_{\mathbb{R}^N} H_N(u,\mu_N) \, dx + \int_{\mathbb{R}^N} u \, d\mu_N,
\end{equation}

\medskip
\noindent
where

\begin{equation*}
\widetilde{\omega_N}:= N^N \left( \frac{N}{N-1} \right)^{N-1} \omega_{N-1},
\end{equation*}
 
\begin{equation}\label{definizionemuN}
 \mu_N(x):= \frac{1}{V_N \left( 1+ |x|^{\frac{N}{N-1}} \right)^N}, \qquad  d \mu_N(x):=\mu_N(x) \, dx,
\end{equation}

\medskip
\noindent
and

\begin{equation*}
 H_N(u,\mu_N) :=| \nabla v_N + \nabla u |^N - |\nabla v_N|^N - N |\nabla v_N |^{N-2} \nabla v_N \cdot \nabla u, 
\end{equation*}

\medskip
\noindent
with 

\begin{equation}
v_N(x):=\ln \, \mu_N(x).
\end{equation}

\medskip
\noindent
Inequality \eqref{EuclideanOnofriNdimensional}  was obtained in \cite{DD3} by considering the endpoint of  optimal Gagliardo–Nirenberg interpolation inequalities, discovered in \cite{DD1} and extended in \cite{DD2}.
We point out that the Euclidean Onofri inequality \eqref{EuclideanOnofriNdimensional} was achieved with different techniques also by Agueh, Boroushaki, and Ghoussoub in \cite{ABG}, for any 
$u \in W^{1,N}(\mathbb{R}^N)$. In this Chapter we extend inequality \eqref{EuclideanOnofriNdimensional} to a suitable weighted Sobolev space. Precisely, we introduce the following weighted Sobolev space

\begin{align*}
\displaystyle 
W_{\mu_N}(\mathbb{R}^N):= \{ u \in L^1(\mathbb{R}^N,d\mu_N): \quad  |\nabla u| \in L^N(\mathbb{R}^N,dx),                     \quad |\nabla u|^2 |\nabla v_N|^{N-2} \in L^1(\mathbb{R}^N,dx) \},
\end{align*}

\medskip
\noindent
endowed with the  norm

\begin{equation}\label{norm}
\lVert u \rVert_{\mu_N}:= \int_{\mathbb{R}^N} |u| \, d\mu_N + \lVert \nabla u \rVert_N + \left( \int_{\mathbb{R}^N}  |\nabla u|^2 |\nabla v_N|^{N-2} \, dx \right)^{\frac{1}{2}}.
\end{equation}

\medskip
\noindent
Our main result is the extension of inequality \eqref{EuclideanOnofriNdimensional} to $W_{\mu_N}(\mathbb{R}^N)$. To this aim, we firstly show that smooth compactly supported functions of $\mathbb{R}^N$ are dense in $W_{\mu_N}(\mathbb{R}^N)$ with respect to the norm stated above.

\begin{theorem}\label{densita}
Assume $N \ge 2$. Then 
$$ \displaystyle W_{\mu_N}(\mathbb{R}^N)=\overline{C_0^{\infty} (\mathbb{R}^N)}^{\lVert \cdot \rVert_{\mu_N}}.$$
\end{theorem}

\medskip
\noindent
Therefore, by density we are able to prove the following.

\begin{theorem}\label{EuclideanOnofriNdim}
Assume $N \geq 2.$ For any $u \in W_{\mu_N}(\mathbb{R}^N)$ we have
 
\begin{equation}\label{OnofriEuclideaEstesa}
\displaystyle  \ln \left( \int_{\mathbb{R}^N} e^u \, d\mu_N \right) \leq  \frac{1}{\widetilde{\omega_N}} \int_{\mathbb{R}^N} H_N(u,\mu_N) \, dx + \int_{\mathbb{R}^N} u \, d\mu_N.
\end{equation}

\end{theorem}

\medskip
\noindent
Let us now consider the functional  $I:W_{\mu_N}(\R^N) \to \mathbb{R}$, associated to \eqref{EuclideanOnofriNdim}, defined as

\begin{equation}\label{Ifunctional}
\displaystyle I(u):=   \frac{1}{\widetilde{\omega_N}} \int_{\mathbb{R}^N} H_N(u,\mu_N) \, dx + \int_{\mathbb{R}^N} u \, d\mu_N- \ln \left(\int_{\mathbb{R}^N} e^u \, d\mu_N \right). 
\end{equation} 

\medskip
\noindent
Without using \eqref{logCCintroduzione} and \eqref{OnofriEuclideaEstesa}, we also prove the following identity between the infimum.

\begin{theorem}\label{equivalence} For any $N\ge 2$, we have 

$$
\inf_{W_{\mu_N}(\R^N)} I = \inf_{W^{1,N}_0(B_1)}  J +\sum_{k=1}^{N-1}\frac{1}{k}.
$$

\end{theorem}

\medskip
\noindent
Theorem \ref{equivalence} shows that the sharp inequalities \eqref{logCCintroduzione} and \eqref{OnofriEuclideaEstesa} are equivalent. In particular, we can give a simple alternative proof of \eqref{OnofriEuclideaEstesa}. On the one hand, since $I(0)=0$, Theorem \ref{equivalence} yields

$$
0 = \inf_{W_{\mu_N}(\R^N)} I = \min_{W_{\mu_N} (\R^N)} I.  
$$

\medskip
\noindent
On the other hand it follows  from \eqref{logCCintroduzione} that for any $u\in W^{1,N}_0(B_1)$ we have
  
$$
J(u) > -\sum_{k=1}^{N-1} \frac{1}{k} = \inf_{W^{1,N}_0(B_1)} J.
$$

\medskip
\noindent
Indeed, we provide a refinement of inequality \eqref{logCCintroduzione}; specifically, in the proof of Theorem \ref{equivalence} we establish that for any  $ u \in W_0^{1,N}(B_1)$ we have

\begin{equation}\label{logCCdanoimigliorata}
\displaystyle \ln \left( \frac{1}{V_N} \int_{B_1} e^{u} \,   dx  + N-1  \right) \leq  \frac{1}{\widetilde{\omega_N}} \int_{B_1} |\nabla u|^N  dx + \sum_{k=1}^{N-1} \frac{1}{k}.
\end{equation}

\medskip
\noindent
Theorem \ref{equivalence} becomes more interesting in dimension $N=2$, considering also the equivalence between \eqref{bestOnofriINTRODUZIONE} and \eqref{EuclideanOnofri2DIntroduzione}. Indeed, as  announced in \cite{BCM} (see also \cite{IM}), it shows that the sharp two-dimensional logarithmic Moser-Trudinger inequality \eqref{logCCintroduzione} by Carleson and Chang is equivalent to the  Onofri's inequality \eqref{bestOnofriINTRODUZIONE} on $\Sm^2$. 
 
As discussed, the inequalities with non-optimal constants \eqref{MTdomininonottimale} and  \eqref{MTS2nonottimale}  are  direct corollaries of the Moser-Trudinger inequality for bounded domains  \eqref{supsudominilimitatiINTRODUZIONE} and its analogue for the  sphere  \eqref{supsuS2INTRODUZIONE}. In \cite{M}, Moser proved these theorems indipendently,
although their proofs are based on similar arguments relying on symmetrization techniques. He explicitly states that it seems impossible to deduce one Theorem from the other. Although the equivalence between Moser's theorems  is still an open problem, we have proved  the equivalence between the corresponding  logarithmic inequalities with sharp constants, i.e. \eqref{logCCintroduzione} for $N=2$ and  \eqref{bestOnofriINTRODUZIONE}.

\bigskip
The Chapter is structured as follows: Section \ref{SezExtension} focuses on the weighted Sobolev space $W_{\mu_N}(\mathbb{R}^N)$ and the proofs of Theorems \ref{densita} and \ref{EuclideanOnofriNdim}. The logarithmic Moser–Trudinger inequality on bounded domains is discussed in Section \ref{sezioneLogMTball}. Section \ref{SezEquiv} details the proof of Theorem \ref{equivalence} concerning the equivalence between \eqref{logCCintroduzione} and \eqref{OnofriEuclideaEstesa}, while Section \ref{SezN=2} examines the case $N=2$ and its relation to the unit sphere $\mathbb{S}^2$.

\section{Euclidean Onofri Inequality in $W_{\mu_N}(\mathbb{R}^N)$}\label{SezExtension}

\medskip
\noindent
We notice that $\mu_N$  defined in \eqref{definizionemuN}  is a probability density on $\mathbb{R}^N$, namely

\begin{align*}
\displaystyle
     \int_{\mathbb{R}^N} \, d \mu_N =  1.
\end{align*}

\medskip
\noindent
In fact

\begin{align*}
\displaystyle
     \int_{\mathbb{R}^N} \mu_N(x) \, dx =  N \int_0^{+ \infty} \frac{1}{\left( 1+ \rho^{\frac{N}{N-1}} \right)^N} \rho^{N-1} \, d \rho 
=  (N-1) \int_{0}^{+ \infty} \frac{1}{(1+t)^2} \left( \frac{t}{1+t} \right)^{N-2} \, dt  = 1.
\end{align*} 

\medskip
\noindent
Moreover, we have

$$ \displaystyle \nabla v_N(x)=-\frac{N^2}{N-1} \, \frac{|x|^{\frac{1}{N-1}}}{1+|x|^{\frac{N}{N-1}} } \,\frac{x}{|x|}.$$

\medskip
\noindent
By computations (see Lemma \ref{computationNlaplacian}), it follows that

\begin{equation}\label{GaussGreen}
 \displaystyle  -\Delta_N v_N(x) =N^N \left( \frac{N}{N-1} \right)^{N-1} V_N \, \mu_N(x)= c_N V_N \, e^{v_N(x)},
\end{equation}

\medskip
\noindent
where $ c_N:=N^N \left( \frac{N}{N-1} \right)^{N-1}$ and   $\Delta_N:= \text{div}\left( | \nabla u |^{N-2} \nabla u  \right)$ denotes the  $N$-Laplacian operator.

\medskip
\noindent
By \eqref{GaussGreen} we can observe that the function $v_N(x)$ is connected to the following Liouville equation

\begin{equation}\label{LiouvilleEquation}
\begin{cases}
\displaystyle - \Delta_N U(x) = e^{U(x)} & x \in \mathbb{R}^N, \bigskip\\
\displaystyle \int_{\mathbb{R}^N}   e^{U(x)} \, dx  < + \infty.
\end{cases}
\end{equation}

\medskip
\noindent
This problem has an explicit solution given by

$$ \displaystyle U(x)= \ln \left(  \frac{c_N}{ \left( 1+ |x|^{\frac{N}{N-1}} \right)^N  }   \right) \qquad x \in \mathbb{R}^N.$$

\medskip
\noindent
Observe now that we have

$$ \displaystyle -\Delta_N v_N(x) = -\Delta_N U(x) = e^{U(x)} = \frac{c_N}{ \left( 1+ |x|^{\frac{N}{N-1}} \right)^N   } = c_N V_N \mu_N(x) = c_N V_N e^{v_N(x)}.$$

\medskip
\noindent
We also recall that thanks to scaling and invariance translation, a family $U_{\lambda,\bar{x}}$ of solutions of \eqref{LiouvilleEquation} is given by

$$ \displaystyle U_{\lambda,\bar{x}}(x)= \ln \left(  \frac{c_N \lambda^N}{ \left( 1+ \lambda^{\frac{N}{N-1}}|x - \bar{x}|^{\frac{N}{N-1}} \right)^N  }   \right)$$

\medskip
\noindent
for all $\lambda >0$ and for any $\bar{x} \in \mathbb{R}^N$.\\
Moreover,  it follows that 

$$ \displaystyle \int_{\mathbb{R}^N}   e^{U_{\lambda,\bar{x}}(x)} \, dx= \int_{\mathbb{R}^N}   e^{U(x)} \, dx = c_N V_N.$$

\medskip
\noindent
In fact

\begin{align*}
\displaystyle
\int_{\mathbb{R}^N}   e^{U(x)} \, dx = c_N \int_{\mathbb{R}^N}  \frac{1}{ \left( 1+|x|^{\frac{N}{N-1}} \right)^N  } \,dx = c_N \omega_{N-1}  \int_0^{+ \infty} \frac{1}{\left( 1+ \rho^{\frac{N}{N-1}} \right)^N} \rho^{N-1} \, d \rho = C_N V_N.
\end{align*}

\medskip
\noindent
In particular, in \cite[Theorem 1.1.]{esposito} it was shown that if $U$ is a solution of \eqref{LiouvilleEquation}, then it is of the form $U_{\lambda,\bar{x}}$ for some $\lambda>0$ and $\bar{x} \in \mathbb{R}^N$.

\bigskip
\noindent
Now, we show how inequality \eqref{EuclideanOnofriNdimensional} involves somehow the $N$-Laplacian operator $\Delta_N$, by means of $H_N(u,\mu_N)$, that we recall is given by

$$
\begin{aligned} 
\displaystyle   H_N(u(x),\mu_N(x))   =  | \nabla v_N(x) + \nabla u(x) |^N - |\nabla v_N(x) |^N - N |\nabla v_N(x) |^{N-2} \nabla v_N(x) \cdot \nabla u(x).
\end{aligned}
$$

\medskip
\noindent
If $u \in C_0^{\infty}(\mathbb{R}^N)$, then there exists a suitable $r>0$ such that the open ball $B_r$, centered in the origin e with radius $r$, contains the support of $u$. Since $u$ is zero on the boundary of the ball $B_r$, Gauss-Green formula gives 

$$ \displaystyle -N \int_{B_r}   |\nabla v_N (x)|^{N-2} \nabla v_N (x)\cdot \nabla u(x) \, dx = \int_{B_r} \Delta_N v_N(x) \, u(x) \, dx,$$

\medskip
\noindent
We have thus shown that the term  $-|\nabla v_N|^{N-2} \nabla v_N \cdot \nabla u$ in $H_N(u, \mu_N)$ contains $\Delta_N$.

\medskip
\noindent
Now, let us denote by 

\begin{equation}\label{definizioneRN}
\displaystyle R_N(X,Y):=|X+Y|^N - |X|^N  -N|X|^{N-2}X \cdot Y \qquad (X,Y) \in \mathbb{R}^N \times \mathbb{R}^N.
\end{equation}

\medskip
\noindent
Considering the definition of $H_N(u(x),\mu_N(x))$, we get 

$$ \displaystyle H_N(u(x),\mu_N(x)) := R_N(\nabla v_N(x), \nabla u(x)).$$

\medskip
\noindent
Let us now prove the following estimate on $R_N$.

\begin{lemma}\label{stima} For any $N\in \N$, $N\ge 2$, there exists a constant $c_N$ such that

$$
0\le R_N(X,Y) \le c_N (|Y|^N + |Y|^2 |X|^{N-2}) \quad \forall \, X,Y \in \R^N.
$$

\end{lemma}

\begin{proof}
If $N=2$, $R_2(X,Y) = |Y|^2$ does not depend on $X$ and the thesis is trivial.\\
Let $N \geq 3$ and let $X,Y \in \mathbb{R}^N.$ We consider the $C^2$ function $f:\mathbb{R} \to \mathbb{R}$ given by 

$$f(t)= \left|tY + X  \right|^N.$$

\medskip
\noindent
By  a direct computation, we have

$$ \displaystyle f'(t)= N \left|tY + X  \right|^{N-2} ( tY + X) \cdot Y, $$

\medskip
\noindent
and

$$ \displaystyle
f''(t)=N(N-2) \left|tY + X  \right|^{N-4} \left( ( tY + X) \cdot Y \right)^2 + N \left|tY + X  \right|^{N-2}|Y|^2
$$

\medskip
\noindent
if $tY + X  \neq 0$, and $f''(t)=0$ otherwise.\\
Now, we observe that if $t \in [0,1]$ we have

\begin{equation}\label{fsecondo}
\displaystyle
\begin{split}
       0 \leq f''(t) 
& \leq   2c_N \left( |Y|^N + |Y|^2  |X|^{N-2} \right), \\
\end{split}
\end{equation}

\medskip
\noindent 
where $c_N=N(N-1)2^{N-4}.$ By Taylor's theorem with Lagrange remainder, there exists $\overline{t} \in\, ]0,1[$ such that

$$ \displaystyle \frac{1}{2} f''(\overline{t})=f(1)-f(0)-f'(0)=\left|Y + X  \right|^N  - \left| X  \right|^N - N \left|X  \right|^{N-2} X \cdot Y .$$

\medskip
\noindent
Taking into account \eqref{fsecondo}, we deduce that

$$ \displaystyle 0 \leq R_N(X,Y) \leq  c_N \left( |Y|^N + |Y|^2  |X|^{N-2} \right).$$

\end{proof}

\medskip
\noindent
Now, recalling that

$$ \displaystyle H_N(u(x),\mu_N(x)) := R_N(\nabla v_N(x), \nabla u(x)),$$

\medskip
\noindent
an immediate consequence of Lemma \ref{stima} is the following integrability condition for $H_N(u,\mu_N)$.

\begin{corollary}
Let $ u \in W_{\mu_N}(\mathbb{R}^N).$ Then 
$$ \displaystyle \int_{\mathbb{R}^N} H_N(u,\mu_N) \, dx < + \infty.$$
\end{corollary}

\medskip
\noindent
We stress the fact that we have shown how $|\nabla u| \in L^N(\mathbb{R}^N,dx)$ and $ |\nabla u|^2 |\nabla v_N|^{N-2} \in L^1(\mathbb{R}^N,dx)$ in the definition of $W_{\mu_N}(\mathbb{R}^N)$ are sufficient conditions to obtain the integrability of $H_N(u,\mu_N)$. Therefore one can wonder what happens if we do not consider $ |\nabla u|^2 |\nabla v_N|^{N-2} \in L^1(\mathbb{R}^N,dx)$, i.e. it is interesting to compare $W_{\mu_N}(\mathbb{R}^N)$ with $W^{1,N}(\mathbb{R}^N)$.

\begin{lemma}\label{W1nmun}
If $\displaystyle u \in  W^{1,N}(\mathbb{R}^N)$ is such that $|\nabla u|^2 |\nabla v_N|^{N-2} \in L^1(\mathbb{R}^N,dx)$, then $u \in W_{\mu_N}(\mathbb{R}^N).$
Moreover, we have $ W_{\mu_N}(\mathbb{R}^N) \setminus W^{1,N}(\mathbb{R}^N) \neq \emptyset$. 
\end{lemma}

\begin{proof} 
Let $\displaystyle u \in  W^{1,N}(\mathbb{R}^N)$ such that $|\nabla u|^2 |\nabla v_N|^{N-2} \in L^1(\mathbb{R}^N,dx)$.

\begin{align*}
\displaystyle 
        \int_{\mathbb{R}^N} |u(x)| \, d \mu_N  =     \int_{\mathbb{R}^N} |u(x)|  \frac{1}{V_N \left( 1+ |x|^{\frac{N}{N-1}} \right)^N}  \, dx.
\end{align*}

\medskip
\noindent
By assumption $\lVert u \rVert_N < + \infty$, while

\begin{align*}
\displaystyle
    \gamma_N:=\int_{\mathbb{R}^N} \frac{1}{\left( 1+ |x|^{\frac{N}{N-1}} \right)^{\frac{N^2}{N-1}}}  \, dx < + \infty.
\end{align*}

\medskip
\noindent
By H\"older inequality we have

$$ \displaystyle \int_{\mathbb{R}^N} |u(x)| \, d \mu_N  \leq \frac{\gamma_N^{\frac{N-1}{N}}}{V_N} \lVert u \rVert_N  < + \infty.$$

\medskip
\noindent
Therefore $u \in W_{\mu_N}(\mathbb{R}^N).$\\
Finally, the constant $u \equiv 1$ is such that $u \in W_{\mu_N}(\mathbb{R}^N) \setminus W^{1,N}(\mathbb{R}^N)$.

\end{proof}

\medskip
\noindent
The next lemma shows that the condition $|\nabla u|^2 |\nabla v_N|^{N-2}\in L^1(\R^N,dx)$ is also necessary for the integrability  of  $H(u,\mu_N)$,  at least in even dimension. 

\begin{lemma}\label{HBelow}
If $N\ge 4$ is even, then 
$$
H(u,\mu_N) \ge \frac{N}{2}|\nabla u|^{2}|\nabla v_N|^{N-2}.
$$
\end{lemma}

\begin{proof}
The proof is based on the following inequalities, which hold for $a,b \in \R$ and $k\in \N$, with $a\ge 0$ and $a+b\ge0$:\\

\begin{enumerate}
\item\label{Ineq1} if $k\ge 2$, then $(a+b)^k\ge a^{k}+ k a^{k-1} b + b^k$.\medskip  
\item\label{Ineq2}  if $k\ge 3$, then $(a+b)^k\ge a^{k}+k a^{k-1} b + k a b^{k-1} + b^k$. 
\end{enumerate}

\medskip
\noindent
Both $(1)$ and $(2)$ can be easily proved by induction on $k$.
Let $X,Y\in \R^N$. If $\frac{N}{2}$ is even, using inequality \eqref{Ineq1} with $k = \frac{N}{2}$, $a = |X|^2$ and $b = 2 X\cdot Y +|Y|^2$ we get that 

$$\begin{aligned}
R_N(X,Y)& = \left(|X|^2+2X\cdot Y +|Y|^2\right)^{\frac{N}{2}}-|X|^N - N |X|^{N-2} X\cdot Y \\
&  \ge \frac{N}{2} |X|^{N-2}|Y|^2 +   \left( 2 X \cdot Y + |Y|^2\right)^\frac{N}{2}. 
\end{aligned}
$$

\medskip
\noindent
Since $\frac{N}{2}$ is even, we get 

\begin{equation}\label{RNBelow}
R_N(X,Y)\ge \frac{N}{2}|X|^{N-2} |Y|^2. 
\end{equation}

\medskip
\noindent
Similarly, if $\frac{N}{2}$ is odd (hence $N\ge 6$), we can apply inequality \eqref{Ineq2} (with $k, a$ and $b$ as above) to get:

$$\begin{aligned}
R_N(X,Y)&  \ge \frac{N}{2} |X|^{N-2}|Y|^2+   \frac{N}{2}|X|^2 \left( 2 X \cdot Y + |Y|^2\right)^{\frac{N}{2}-1} +   \left( 2 X \cdot Y + |Y|^2\right)^\frac{N}{2}  \\
& =  \frac{N}{2} |X|^{N-2}|Y|^2  +  \left( 2 X \cdot Y + |Y|^2\right)^{\frac{N}{2}-1} \left( \frac{N}{2}|X|^2  + 2 X \cdot Y + |Y|^2\right)
\\
& \ge  \frac{N}{2} |X|^{N-2}|Y|^2  +  \left( 2 X \cdot Y + |Y|^2\right)^{\frac{N}{2}-1} |X+Y|^2.
\end{aligned}
$$

\medskip
\noindent
Since $\frac{N}{2}-1$ is even, we get again \eqref{RNBelow}. 
Hence \eqref{RNBelow} holds whenever $\frac{N}{2}$ is integer and the conclusion follows by taking $X = \nabla v_N$ and $Y = \nabla u$.
  
\end{proof}

\begin{remark}
If $N \geq 3$ and $u \in W^{1,N}(\mathbb{R}^N)$, in general the integrability condition $|\nabla u|^2 |\nabla v_N|^{N-2} \in L^1(\mathbb{R}^N,dx)$ of Lemma \ref{W1nmun} is not satisfied. Indeed, let 

$$ \displaystyle u(x):= \sum_{k=2}^{\infty}  \frac{1}{k \sqrt{\ln k } } w(2(|x| - k)), $$

\medskip
\noindent
where

$$ \displaystyle
w(r):=
\begin{cases}
1+r & \text{if } r \in [-1,0[, \\
1-r & \text{if } r \in [0,1], \\
0   & \text{otherwise}.
\end{cases}
$$

\medskip
\noindent
A direct computation shows that  the function $u$  satisfies 
 
\begin{equation}\label{example}
u \in W^{1,N}(\mathbb{R}^N) \quad \text{ and } \quad |\nabla u|^2 |\nabla v_N|^{N-2} \not\in L^1(\mathbb{R}^N,dx). 
\end{equation} 

\medskip
\noindent
In particular, if $N$ is even, Lemma \ref{HBelow} yields

$$
\int_{\R^N} H(u,\mu_N) \,dx = +\infty.  
$$

\medskip
\noindent
To prove \eqref{example}  we observe that the functions $x\rightarrow w(2(|x|-k))$ have disjoint supports as $k$ varies in $\N$ with $k\ge 2$. Then 

$$ \displaystyle |u(x)|^N= \sum_{k=2}^{\infty}  \frac{1}{k^N (\ln k)^{\frac{N}{2}} }  (w(2(|x| - k)))^N.$$

\medskip
\noindent
Integrating in polar coordinates and using that $0\le w\le 1$, we find 

$$
\begin{aligned}
 \int_{\mathbb{R}^N} |u|^N \, dx 
& \leq \omega_{N-1} \sum_{k=2}^{\infty}  \frac{1}{k^N (\ln k)^{\frac{N}{2}} } \int_{k - \frac{1}{2}}^{k + \frac{1}{2}} \rho^{N-1} \, d \rho \\
& = \frac{\omega_{N-1}}{N} \sum_{k=2}^{\infty}  \frac{1}{k^N (\ln k)^{\frac{N}{2}} } \left[ \left(k + \frac{1}{2}   \right)^N - \left(k - \frac{1}{2}   \right)^N  \right] \\
& \leq c_1(N) \sum_{k=2}^{\infty} \frac{1}{k (\ln k)^{\frac{N}{2}} } < + \infty
\end{aligned}
$$

\medskip
\noindent
since $\frac{N}{2}>1$. Similarly, we have 

$$ \displaystyle |\nabla u(x)|=2 \sum_{k=2}^{\infty}  \frac{1}{k \sqrt{\ln k } }  | w'(2(|x| - k))|,$$

\medskip
\noindent
and 

$$ \displaystyle \int_{\mathbb{R}^N} |\nabla u(x)|^N \, dx =  c_2(N) \sum_{k=2}^{\infty}  \frac{1}{k^N (\ln k)^{\frac{N}{2}} } \int_{k - \frac{1}{2}}^{k + \frac{1}{2}} \rho^{N-1} \, d \rho < + \infty.$$

\medskip
\noindent
Now, we observe that $|\nabla v_N(x)| \sim \frac{N^2}{N-1} \frac{1}{|x|}$ if $|x| >> 1$, hence, taking $k_0$ sufficiently large, we obtain

\begin{align*} \displaystyle 
                        \int_{\mathbb{R}^N} |\nabla u |^2 |\nabla v_N |^{N-2} \, dx 
& \geq  c_3(N) + c_4(N) \sum_{k=k_0}^{\infty} \frac{1}{k^2 \ln k } \int_{k - \frac{1}{2}}^{k + \frac{1}{2}} \frac{1}{\rho^{N-2}} \rho^{N-1} \, d \rho \\
& \geq                  c_5(N) \left( 1 + \sum_{k=k_0}^{\infty} \frac{1}{k \ln k } \right) =                     + \infty.
\end{align*}

\end{remark}

\medskip
\noindent
To extend \eqref{EuclideanOnofriNdimensional} to $W_{\mu_N}(\mathbb{R}^N)$, we show that smooth compactly supported functions of $\mathbb{R}^N$ are dense in $W_{\mu_N}(\mathbb{R}^N)$ with respect to the norm   \eqref{norm}.

\medskip
\noindent
Let us define  $$\mathcal{X}:= \overline{C_0^{\infty} (\mathbb{R}^N)}^{\lVert \cdot \rVert_{\mu_N}}.$$ Clearly, $ C_0^{\infty} (\mathbb{R}^N) \subset W_{\mu_N}(\mathbb{R}^N)$, so that $\mathcal{X} \subseteq W_{\mu_N}(\mathbb{R}^N).$ Our aim is to prove that the two spaces coincide. 
We will prove this in two steps. Firstly, by using a similar idea of Tintarev and Fieseler in \cite[Remark 2.2]{TF}, we show that bounded functions in $W_{\mu_N}(\mathbb{R}^N)$  belong to $\mathcal{X}.$ 

\begin{lemma}\label{funzionilimitate}

Let $u \in W_{\mu_N}(\mathbb{R}^N) \cap L^{\infty}(\mathbb{R}^N)$. Then $u \in \mathcal{X}.$
\end{lemma} 

\begin{proof}
For any $k \in \mathbb{N},$ $k \geq 2$, let us denote by

 $$  k^*:= \frac{1}{\left( 1 - \frac{1}{\sqrt[k]{k}} \right)^k}.$$ 
 
\medskip
\noindent
We define

\begin{equation*}
\eta_k(x):=
\begin{cases}
 \displaystyle \frac{1}{\sqrt[k]{k}} & \text{ if } |x| \leq 1, \bigskip \\
\displaystyle  \frac{1}{\sqrt[k]{|x|}} + \frac{1}{\sqrt[k]{k}} - 1 & \text{ if }  1 < |x| \leq k^*, \bigskip \\
               0 & \text{ if }  k^* < |x| .
\end{cases}
\end{equation*}

\medskip 
\noindent
Let $u \in W_{\mu_N}(\mathbb{R}^N) \cap L^{\infty}(\mathbb{R}^N)$. We define $u_k:=u \, \eta_k.$
Observe that

\begin{equation*}
\nabla u_k(x):=
\begin{cases}
 \displaystyle \frac{1}{\sqrt[k]{k}} \nabla u(x) & \text{ if } |x| < 1, \bigskip \\
\displaystyle   \left( \frac{1}{\sqrt[k]{k}} + \frac{1}{\sqrt[k]{|x|}} - 1 \right)\nabla u(x) -\frac{1}{k} \frac{u(x)}{|x|^{\frac{1}{k}+1}} \frac{x}{|x|} & \text{ if }  1 < |x| < k^*, \bigskip \\
               0 & \text{ if }  k^* < |x| .
\end{cases}
\end{equation*}

\noindent
First of all, since $u \in L^{\infty}(\mathbb{R}^N)$ and $\eta_k(x) \to 1 $ for any $x \in \mathbb{R}^N$ as $k \to + \infty,$  we can apply   dominated convergence theorem to say that 

$$ \displaystyle \int_{\mathbb{R}^N} |u - u_k(x)| \, d \mu_N  \to 0 \quad \text{ as } k \to + \infty.$$

\medskip
\noindent
Now, we prove 

$$ \displaystyle \lVert \nabla ( u - u_k )  \rVert_N \to 0 \quad \text{ as } k \to + \infty.$$

\medskip
\noindent
We will denote by $o(1)$ any quantity that goes to $0$ as $k \to + \infty.$

\begin{align*}
\displaystyle 
\int_{\mathbb{R}^N} | \nabla u - \nabla u_k|^N \, dx 
& = \int_{|x| < 1} | \nabla u - \nabla u_k|^N \, dx +  \int_{1 < |x| < k^*} | \nabla u - \nabla u_k|^N \, dx + \int_{ |x| > k^*} | \nabla u |^N \, dx.
\end{align*}

\medskip
\noindent
We have

$$ \displaystyle \int_{|x| < 1} | \nabla u - \nabla u_k|^N \, dx= \left( 1 - \frac{1}{\sqrt[k]{k}}  \right)^N \int_{|x| < 1} | \nabla u|^N \, dx=o(1).$$

\medskip
\noindent
Since $\lVert \nabla u \rVert_N < + \infty$ and $k^* \to + \infty$ as $k \to + \infty,$ we easily obtain 

$$ \displaystyle \int_{|x| > k*} | \nabla u - \nabla u_k|^N \, dx=  \int_{|x| > k^*} | \nabla u|^N \, dx=o(1).$$

\medskip
\noindent
Since $u \in L^{\infty}(\mathbb{R}^N),$ there exists a constant $M_u$ such that $|u(x)| \leq M_u$ for almost every $ x \in \mathbb{R}^N$.

\begin{align*}
\displaystyle
    \int_{1 < |x| < k^*} | \nabla u - \nabla u_k|^N \, dx  
& =    \int_{1 < |x| < k^*} \left| \left(  2 - \frac{1}{\sqrt[k]{k}} - \frac{1}{\sqrt[k]{|x|}}  \right) \nabla u(x)  + \frac{1}{k} \frac{u(x)}{|x|^{\frac{1}{k}+1}} \frac{x}{|x|}    \right|^N  \, dx \\
& \leq C \int_{1 < |x| < k^* }  \left(  2 - \frac{1}{\sqrt[k]{k}} - \frac{1}{\sqrt[k]{|x|}}  \right)^N  \left| \nabla u(x) \right|^N \, dx \\
& + C \frac{M_u^N}{k^N} \int_{1 < |x| < k^*} \frac{1}{|x|^{\frac{N}{k}+N}} \, dx.
\end{align*}

\noindent
Now, we observe that 

$$
\displaystyle
 \frac{M_u^N}{k^N} \int_{ 1 < |x| < k^*} \frac{1}{|x|^{\frac{N}{k}+N}} \, dx < \frac{c_1}{k^N} \int_1^{+ \infty} \frac{1}{\rho^{\frac{N}{k}+1}} \, d \rho = \frac{c_2}{k^{N-1}}=o(1).
$$

\medskip
\noindent
Moreover, we can apply dominated convergence theorem to say that

$$ \displaystyle \int_{1 < |x| < k^* }  \left(  2 - \frac{1}{\sqrt[k]{k}} - \frac{1}{\sqrt[k]{|x|}}  \right)^N  \left| \nabla u(x) \right|^N \, dx \to 0 \quad \text{ as } k \to + \infty.$$

\medskip
\noindent
Hence

$$ \displaystyle \lVert \nabla ( u - u_k )  \rVert_N \to 0 \quad \text{ as } k \to + \infty.$$

\medskip
\noindent
It remains to prove, in the case $N \geq 3$, that

$$ \displaystyle \left( \int_{\mathbb{R}^N}  |\nabla u-  \nabla u_k |^2 |\nabla v_N|^{N-2} \, dx \right)^{\frac{1}{2}} \to 0 \quad \text{ as } k \to + \infty.$$

\noindent
We observe that $\displaystyle | \nabla v_N(x)| \sim \frac{1}{|x|}$ as $|x| \to + \infty,$ by which there exists a constant $ r_N > 1$ such that $ \displaystyle  |\nabla v_N (x)| \leq \frac{2}{|x|} $ if $|x| \geq r_N.$ Since $k \to + \infty$ implies $k^* \to + \infty,$ we can assume from the beginning that $k$ is such that $k^* > r_N$.
Moreover, we recall that $|\nabla v_N|$ is bounded from above.

\begin{align*}
\displaystyle
&       \int_{\mathbb{R}^N}  |\nabla u - \nabla u_k|^2 |\nabla v_N|^{N-2} \, dx \\
& =     \int_{ |x| < 1 }  |\nabla u - \nabla u_k|^2 |\nabla v_N|^{N-2} \, dx +  \int_{ 1 < |x| < r_N }  |\nabla u - \nabla u_k|^2 |\nabla v_N|^{N-2} \, dx \\
& +    \int_{ r_N < |x|  < k^* }  |\nabla u - \nabla u_k|^2 |\nabla v_N|^{N-2} \, dx + \int_{ k^* < |x|  }  |\nabla u|^2 |\nabla v_N|^{N-2} \, dx.
\end{align*}

\noindent
Reasoning in a similar way to the previous case, we obtain 
$$ \displaystyle \int_{ |x| < 1 }  |\nabla u - \nabla u_k|^2 |\nabla v_N|^{N-2} \, dx=o(1),$$
$$ \displaystyle   \int_{1 < |x| < r_N} | \nabla u - \nabla u_k|^2 |\nabla v_N|^{N-2} \, dx  =o(1),$$
and
$$ \displaystyle \int_{ k^* < |x|  }  |\nabla u |^2 |\nabla v_N|^{N-2} \, dx=o(1).$$
Now we observe that
\begin{align*}
\displaystyle
&       \int_{r_N < |x| < k^*} | \nabla u - \nabla u_k|^2 |\nabla v_N|^{N-2} \, dx  \\
& =    \int_{r_N < |x| < k^*} \left| \left(  2 - \frac{1}{\sqrt[k]{k}} - \frac{1}{\sqrt[k]{|x|}}  \right) \nabla u(x)  +\frac{1}{k} \frac{u(x)}{|x|^{\frac{1}{k}+1}} \frac{x}{|x|}    \right|^2 |\nabla v_N|^{N-2} \, dx \\
& \leq C \int_{r_N < |x| < k^*}  \left(  2 - \frac{1}{\sqrt[k]{k}} - \frac{1}{\sqrt[k]{|x|}}  \right)^2  \left| \nabla u(x) \right|^2 |\nabla v_N|^{N-2}  \, dx \\
& + C \frac{1}{k^2} \int_{r_N < |x|} \frac{1}{|x|^{\frac{2}{k}+N}} \, dx . \\
\end{align*}

\noindent
Clearly, 
$$ \displaystyle \frac{1}{k^2} \int_{r_N < |x|} \frac{1}{|x|^{\frac{2}{k}+N}} \, dx=o(1).$$

\noindent
Moreover, taking into account that $ \displaystyle |\nabla u|^2 |\nabla v_N|^{N-2} \in L^1(\mathbb{R}^N,dx),$  we can apply dominated convergence theorem to say that 
$$ \displaystyle  \int_{r_N < |x| < k^*}  \left(  2 - \frac{1}{\sqrt[k]{k}} - \frac{1}{\sqrt[k]{|x|}}  \right)^2  \left| \nabla u(x) \right|^2 |\nabla v_N|^{N-2}  \, dx=o(1).$$

\noindent
Finally, we have 
$$ \displaystyle \left( \int_{\mathbb{R}^N}  |\nabla u-  \nabla u_k |^2 |\nabla v_N|^{N-2} \, dx \right)^{\frac{1}{2}} \to 0 \quad \text{ as } k \to + \infty.$$

\noindent
We have proved that $\{u_k \}_k $ is a sequence of compactly supported functions such that $\lVert u -u_k \rVert_{\mu_N} \to 0$ as $k \to +\infty$.\\
We want to take a sequence $\{\tilde{u}_k \}_k \subset C_0^{\infty}(\mathbb{R}^N)$ such that $\lVert u -\tilde{u}_k \rVert_{\mu_N} \to 0$ as $k \to +\infty$.\\
Let us denote by $B_{k^*}$ the open ball $B_{k^*}(O)$.\\
Since ${u_{k}}_{|_{B_{k^*}}} \in W_0^{1,N}(B_{k^*}),$ there exists $\tilde{u}_k \in C_0^{\infty}(B_{k^*})$ such that 

$$ \displaystyle \left\lVert {u_{k}}_{|_{B_{k^*}}}- \tilde{u}_k  \right\rVert_{L^N(B_{k^*})} + \left\lVert \nabla \left({u_{k}}_{|_{B_{k^*}}}- \tilde{u}_k \right) \right\rVert_{L^N(B_{k^*})} < \frac{1}{k}.$$

\medskip
\noindent
We can extend to $0$ in $\mathbb{R}^N \setminus B_{k^*}$  and consider $\tilde{u}_k \in C_0^{\infty}(\mathbb{R}^N)$. Hence

$$ \displaystyle \lVert u_k - \tilde{u}_k  \rVert_N + \lVert \nabla \left(u_k - \tilde{u}_k \right) \rVert_N  < \frac{1}{k}.$$

\noindent
Since
\begin{align*}
\displaystyle
    \gamma_N:=\int_{\mathbb{R}^N} \frac{1}{\left( 1+ |x|^{\frac{N}{N-1}} \right)^{\frac{N^2}{N-1}}}  \, dx < + \infty  ,
\end{align*}

\noindent
we can apply H\"older inequality and say that

\begin{align*}
 \displaystyle 
       \int_{\mathbb{R}^N} |u_k(x) - \tilde{u}_k(x)| \, d \mu_N  < \gamma_N^{\frac{N-1}{N}} \lVert u_{k}- \tilde{u}_k  \rVert_N < \frac{\gamma_N^{\frac{N-1}{N}}}{k}=o(1).
\end{align*}

\noindent
Now, let $k$ be such that $k^* > r_N.$\\
Since

\begin{align*}
\displaystyle 
    \int_{0 < |x| \leq k^*} |\nabla v_N|^N \, dx \leq C \left( \ln k^* + 1   \right),
\end{align*}

\medskip
\noindent
we have
\begin{align*}
\displaystyle  \int_{\mathbb{R}^N}  |\nabla(u_k - \tilde{u}_k)|^2 |\nabla v_N|^{N-2} \, dx  <    C \frac{ \left( \ln k^* + 1    \right)^{\frac{N-2}{N}} }{k^2}.
\end{align*}

\medskip
\noindent
Now, observe that

$$ \displaystyle \lim_{k \to + \infty} \ln \, k^* = \lim_{k \to + \infty}  - k \ln \left( 1 - \frac{1}{\sqrt[k]{k}} \right) = + \infty,$$

\medskip
\noindent
and

\begin{align*}
\displaystyle
\ln \left( 1 - \frac{1}{\sqrt[k]{k}} \right) 
= \ln \left( 1 - e^{-\frac{\ln \, k}{k}} \right) 
= \ln \left(\frac{ \ln \, k}{k} + o \left( \frac{\ln \, k}{k} \right)     \right) 
= - ( 1 + o(1)) \ln \, k .
\end{align*}

\medskip
\noindent
Therefore, we have

\begin{align*}
\displaystyle
 \lim_{k \to + \infty} \frac{(\ln k^*+1)^{\frac{N-2}{N}}}{k^2} 
= \lim_{k \to + \infty} \frac{(\ln k^*)^{\frac{N-2}{N}}}{k^2}  
=  \lim_{k \to + \infty} \frac{ \left[  \ln \, k \right]^{1-\frac{2}{N}}}{k^{1+ \frac{2}{N}}} = 0.
\end{align*}

\medskip
\noindent
Hence

$$ \displaystyle \int_{\mathbb{R}^N}  |\nabla(u_k - \tilde{u}_k)|^2 |\nabla v_N|^{N-2} \, dx \to 0.$$

\medskip
\noindent
Finally, we have

\begin{align*}
\displaystyle
  \lVert u - \tilde{u}_k \rVert_{\mu_N} \leq  \lVert u -  u_k \rVert_{\mu_N} +  \lVert u_k - \tilde{u}_k \rVert_{\mu_N} \to 0 \quad \text{as } k \to + \infty.  
\end{align*}

\end{proof}

\noindent
Thanks to previous lemma, it remains to prove that any function $u \in W_{\mu_N}(\mathbb{R}^N)$ can be approximated by a sequence of bounded functions in $ W_{\mu_N}(\mathbb{R}^N)$.

\begin{lemma}\label{convergenzatroncate}\

\medskip
\noindent
Let $u \in W_{\mu_N}(\mathbb{R}^N)$. There exists a sequence $\{u_k\}_{k \in \mathbb{N}} \subset W_{\mu_N}(\mathbb{R}^N) \cap L^{\infty}(\mathbb{R}^N) $ such that $\lVert u-u_{k} \rVert_{\mu_N} \to 0$ as $k \to + \infty$. 
\end{lemma}

\begin{proof}

Let $ u \in W_{\mu_N}(\mathbb{R}^N).$ Taking into account \cite[Lemma A.4, Chapter 2]{KS}, we have $|u| \in  W_{\mu_N}(\mathbb{R}^N)$ with

\begin{equation*}
\nabla |u|=
\begin{cases}
 \nabla u  & \text{ a.e. in } \{ x \in \mathbb{R}^N  \; | \; u(x)>0\},\\
 0         & \text{ a.e. in } \{ x \in \mathbb{R}^N  \; | \; u(x)=0\},\\
 -\nabla u & \text{ a.e. in } \{ x \in \mathbb{R}^N  \; | \; u(x)<0\}.
\end{cases}
\end{equation*}

\medskip
\noindent
For any $\lambda > 0$, we can define the \textit{truncated function of $u$ of parameter $\lambda$} as $\displaystyle u_{\lambda}:=\max\{-\lambda,\min\{u,\lambda\}\},$ that is

\begin{equation*}
u_{\lambda}(x)=
\begin{cases}
   \lambda         & \text{ a.e. in } \{ x \in \mathbb{R}^N  \; | \; u(x)\geq  \lambda \},\\
    u              & \text{ a.e. in } \{ x \in \mathbb{R}^N  \; | \; - \lambda < u(x) < \lambda \},\\
 - \lambda         & \text{ a.e. in } \{ x \in \mathbb{R}^N  \; | \; u(x) \leq - \lambda \}.
\end{cases}
\end{equation*}

\medskip
\noindent
Let $k \in \mathbb{N},$ $k \geq 1,$ and let $\{u_k\}_{k \in \mathbb{N}}$ be a sequence of truncated function of $u$ of parameter $k$.\\ Clearly,  $\{u_k\}_{k \in \mathbb{N}} \subset W_{\mu_N}(\mathbb{R}^N) \cap L^{\infty}(\mathbb{R}^N) $.\\
Let us denote by $\displaystyle E_k:=\{ x \in \mathbb{R}^N \; | \; |u(x)| \geq k \}.$\\
Hence 

$$ \displaystyle \int_{\mathbb{R}^N} |u - u_k| \, d \mu_N = \int_{E_k} |u - u_k | \, d \mu_N \leq \int_{E_k} |u| \, d \mu_N. $$ 

\medskip
\noindent
Since $u \in L^1(\mathbb{R}^N, d \mu_N),$ by Chebyshev's inequality we know that 

$$|E_k|_{\mu_N} \to 0 \qquad \text{as } k \to + \infty.$$

\medskip
\noindent
Moreover, we can apply absolute continuity of the integral and say that 

$$ \displaystyle \int_{E_k} |u| \, d \mu_N \to 0 \qquad \text{ as } k \to + \infty.$$

\medskip
\noindent
Therefore, we deduce $u_k \to u$ in $L^1(\mathbb{R}^N, d\mu_N)$.

\medskip
\noindent
Let now $\varepsilon > 0$ and let us consider the set

$$ \displaystyle \left\{ x \in \mathbb{R}^N \; : \; \left| \frac{ \partial u_k}{ \partial x_i}(x) - \frac{ \partial u}{ \partial x_i}(x)   \right| \geq \varepsilon  \right\}. $$ 

\medskip
\noindent
Since

$$ \displaystyle \left\{ x \in \mathbb{R}^N \; : \; \left| \frac{ \partial u_k}{ \partial x_i}(x) - \frac{ \partial u}{ \partial x_i}(x)   \right| \geq \varepsilon  \right\} \subset E_k, $$ 

\medskip
\noindent
we deduce that 

$$ \displaystyle \frac{ \partial u_k}{ \partial x_i} \to \frac{ \partial u}{ \partial x_i} \quad \text{ in measure with respect to $\mu_N$.}$$

\medskip
\noindent
Hence there exists a subsequence $\{u_{k_j}\}_{j \in \mathbb{N}}$ such that $ \nabla u_{k_j}(x) \to \nabla u(x) $ almost everywhere in $\mathbb{R}^N$ with respect to the measure $\mu_N$. Since the Lebesgue measure $dx$ is absolutely continuous with respect to the measure $d\mu_N$, we deduce that
 $ \nabla u_{k_j} \to \nabla u $ almost everywhere in $\mathbb{R}^N$ also with respect to the Lebesgue measure.

\medskip
\noindent
Finally, since 

$$ \displaystyle \left| \nabla u(x) - \nabla u_{k_j}(x) \right|^N \leq \left| \nabla u(x) \right|^N \quad \text{a.e. } x \in \mathbb{R}^N,$$

\medskip
\noindent
and

$$ \displaystyle  |\nabla u(x)-  \nabla u_{k_j}(x) |^2 |\nabla v_N(x)|^{N-2} \leq |\nabla u(x) |^2 |\nabla v_N(x)|^{N-2} \quad \text{a.e. } x \in \mathbb{R}^N,$$

\medskip
\noindent
we can apply dominated convergence theorem to say that

$$ \lVert \nabla u - \nabla u_{k_j} \rVert_N \to 0 \quad \text{ as } j \to + \infty,$$

\medskip
\noindent
and

$$ \displaystyle \left( \int_{\mathbb{R}^N}  |\nabla u-  \nabla u_{k_j} |^2 |\nabla v_N|^{N-2} \, dx \right)^{\frac{1}{2}} \to 0 \quad \text{ as } j \to + \infty.$$

\end{proof}

\medskip
\noindent
Taking into account Lemma \ref{funzionilimitate} and Lemma \ref{convergenzatroncate}, we have proved Theorem \ref{densita}, that is 

$$ \displaystyle W_{\mu_N}(\mathbb{R}^N)=\mathcal{X}=\overline{C_0^{\infty} (\mathbb{R}^N)}^{\lVert \cdot \rVert_{\mu_N}}.$$

\medskip
\noindent
Now, we can extend the $N$-dimensional Euclidean Onofri inequality to the weighted Sobolev space $W_{\mu_N}(\mathbb{R}^N)$.

\medskip
\noindent
{\mbox {\it Proof of Theorem~\ref{EuclideanOnofriNdim}.~}}  Let $u \in W_{\mu_N}(\mathbb{R}^N)$. Since $W_{\mu_N}(\mathbb{R}^N)= \overline{C_0^{\infty} (\mathbb{R}^N)}^{\lVert \cdot \rVert_{\mu_N}},$ there exists a sequence $\displaystyle \{u_m\}_{m \in \mathbb{N}} \subset C_0^{\infty}(\mathbb{R}^N)$ such that $\displaystyle \lVert u-u_m \rVert_{\mu_N} \to 0.$ By \eqref{EuclideanOnofriNdimensional}, for any  $ m \in \mathbb{N}$ we have

$$ \displaystyle \ln \left( \int_{\mathbb{R}^N} e^{u_m} \, d\mu_N \right) \leq  \frac{1}{\widetilde{\omega_N}} \int_{\mathbb{R}^N} H_N(u_m,\mu_N) \, dx + \int_{\mathbb{R}^N} u_m \, d\mu_N. $$

\medskip
\noindent
The condition $\displaystyle \lVert u-u_m \rVert_{\mu_N} \to 0$ implies that there exists a subsequence, that we recall $\{u_m\}_{m \in \mathbb{N}}$, and a nonnegative function 
$w \in L^N(\mathbb{R}^N,dx)$ with $w^2 |\nabla v_N|^{N-2} \in L^1(\mathbb{R}^N,dx)$, such that:\\

\begin{itemize}
\item[•] $u_m(x) \to u(x)$ a.e. in $\mathbb{R}^N$ as $m \to + \infty$; \medskip
\item[•] $\nabla u_m(x) \to \nabla u(x)$ a.e. in $\mathbb{R}^N$ as $m \to + \infty$; \medskip
\item[•] for any $m$, $|\nabla u_m(x) | \leq w(x)$ a.e. in $\mathbb{R}^N$. 
\end{itemize}

\medskip
\noindent
We easily obtain that

$$ \displaystyle \int_{\mathbb{R}^N} u_m \, d\mu_N  \to \int_{\mathbb{R}^N} u \, d\mu_N \quad \text{as } m \to + \infty.$$

\medskip
\noindent
It is immediate to observe that

$$ \displaystyle H_N(u_m(x),\mu_N(x)) \to H_N(u(x),\mu_N(x)) \quad \text{ a.e. in } \mathbb{R}^N  \quad \text{ as   } \; m \to + \infty.$$

\medskip
\noindent
Taking into account Lemma \ref{stima}, we have

\begin{align*}
\displaystyle
  0 \leq H_N(u_m(x),\mu_N(x)) 
& \leq \frac{c_N}{2} \left( |\nabla u_m(x) |^N + |\nabla u_m(x) |^2|\nabla v_N(x) |^{N-2}   \right) \\
& \leq \frac{c_N}{2} \left( w(x)^N + w(x)^2|\nabla v_N(x) |^{N-2}   \right), 
\end{align*}

\medskip
\noindent
almost everywhere in $\mathbb{R}^N$ and for any $m$.\\
Since $w^N + w^2|\nabla v_N |^{N-2} \in L^1(\mathbb{R}^N,dx)$, we can apply dominated convergence theorem to say that

$$ \displaystyle \int_{\mathbb{R}^N} H_N(u_m,\mu_N) \, dx  \to \int_{\mathbb{R}^N} H_N(u,\mu_N) \, dx .$$

\medskip
\noindent
Finally, by Fatou's Lemma we have 

\begin{equation*}
\displaystyle
\int_{\mathbb{R}^N} e^u \, d\mu_N 
\leq \liminf_{m \to + \infty} \left( \int_{\mathbb{R}^N} e^{u_m} \, d\mu_N \right) 
\leq  e^{\frac{1}{\widetilde{\omega_N}} \int_{\mathbb{R}^N} H_N(u,\mu_N) \, dx + \int_{\mathbb{R}^N} u \, d\mu_N},
\end{equation*}

\medskip
\noindent
and passing to logarithm we can conclude that 

$$ \displaystyle \ln \left( \int_{\mathbb{R}^N} e^u \, d\mu_N \right) \leq  \frac{1}{\widetilde{\omega_N}} \int_{\mathbb{R}^N} H_N(u,\mu_N) \, dx + \int_{\mathbb{R}^N} u \, d\mu_N. $$

\qed

\medskip
\noindent
Notice that the functional  $I:W_{\mu_N}(\R^N) \to \mathbb{R}$, associated to \eqref{EuclideanOnofriNdim} and defined as

\begin{equation*}
\displaystyle I(u):=   \frac{1}{\widetilde{\omega_N}} \int_{\mathbb{R}^N} H_N(u,\mu_N) \, dx + \int_{\mathbb{R}^N} u \, d\mu_N- \ln \left(\int_{\mathbb{R}^N} e^u \, d\mu_N \right) 
\end{equation*} 

\medskip
\noindent
satisfies

$$ \displaystyle \inf_{W_{\mu_N}(\R^N)} I  = 0,$$

\medskip
\noindent
where the infimum is realized for any constant in $W_{\mu_N}(\R^N)$.

\bigskip

\section{Logarithmic Moser-Trudinger inequality on the balls of $\mathbb{R}^N$}\label{sezioneLogMTball}

\bigskip
\noindent
As previously noted, Moser \cite{M} investigated the borderline case of Sobolev embeddings for $W_0^{1,N}(\Omega)$, where $\Omega$ is a bounded domain in $\mathbb{R}^N$, refining the results of Trudinger \cite{Tru} (see also Yudovich \cite{Yudo} and Pohozaev \cite{poho}). More precisely, Moser established the following theorem:

\begin{theorem}[\textbf{Moser-Trudinger inequality on bounded domain of $\boldsymbol{\mathbb{R}^N}$}]\label{TheoremaMoserdomini}\

\medskip
\noindent
Let $\Omega$ be a bounded domain in $\R^N$ with $N \geq 2$. Let $\alpha_N:= N \, \omega_{N-1}^\frac{1}{N-1}$, where $\omega_{N-1}$ denotes the measure of the surface of the unit sphere of $\mathbb{R}^N$. Denoting by $\mathcal{H}:=\{u \in W_0^{1,N}(\Omega) \ | \  \lVert \nabla u \rVert_N \leq 1 \}$, there exists a constant $C_N > 1,$ such that

\begin{equation}\label{supsudominilimitati}
\sup_{u \in \mathcal{H}} \,  \frac{1}{|\Omega|} \int_{\Omega} e^{\alpha \left|u \right|^\frac{N}{N-1}}\,dx \le C_N \qquad \forall \, \alpha \leq  \alpha_N.
\end{equation}

\medskip
\noindent
The integral on the left is actually finite for any positive $\alpha$ and for any $u \in \mathcal{H}$, but if $\alpha > \alpha_N$ it can be made arbitrarily large by an appropriate choice of $u$, hence the supremum on $\mathcal{H}$ is equal to $+ \infty.$
\end{theorem}

As a straightforward consequence of previous result and Young's inequality, one can prove the following
weaker inequality:

\begin{proposition}[\textbf{Logarithmic M-T inequality on bounded domains of $\boldsymbol{\mathbb{R}^N}$}]\label{propositionlogMTinequalitydomains}\

\medskip
\noindent
Let $\Omega$ be a bounded domain of $\R^N$ with $N \geq 2$, and let $C_N$ be the constant appearing in \eqref{supsudominilimitati}.

\medskip
\noindent
For any $ u \in W_0^{1,N}(\Omega)$, we have

\begin{equation}\label{logMTinequalitydomains}
\displaystyle \ln \left( \frac{1}{|\Omega|} \int_{\Omega} e^u \, dx \right) \leq  \frac{1}{\widetilde{\omega_N}} \int_{\Omega} | \nabla u |^N \, dx + \ln C_N,
\end{equation}

\medskip
\noindent
where $ \displaystyle \widetilde{\omega_N}:= N^N \left( \frac{N}{N-1} \right)^{N-1} \omega_{N-1}.$

\end{proposition}

\begin{proof}

Let $u \in W^{1,N}_0(\Omega)$, $u \neq 0$, and let  $\displaystyle b^{\frac{N}{N-1}}:=\frac{N}{N-1}\frac{\alpha_N}{\lVert \nabla u \rVert^{\frac{N}{N-1}}_N}$.

\medskip
\noindent
By Young's inequality with conjugated exponents $N$ and $\frac{N}{N-1}$, we infer that 

$$ u \leq \frac{1}{b} \, \left(b|u|\right)  \leq  \frac{1}{Nb^N} + \frac{N-1}{N}b^{\frac{N}{N-1}}|u|^{\frac{N}{N-1}}=  \frac{\lVert \nabla u \rVert^N_N}{\widetilde{\omega_N}} + \alpha_N \left(\frac{\left| u \right|}{\lVert \nabla u \rVert_N} \right)^{\frac{N}{N-1}}. $$

\medskip
\noindent
Considering \eqref{supsudominilimitati}, we get

$$ \displaystyle \frac{1}{|\Omega|} \int_{\Omega} e^u \, dx \leq e^{\frac{\lVert \nabla u \rVert^N_N}{\widetilde{\omega_N}}} \frac{1}{|\Omega|} \int_{\Omega} e^{\alpha_N \left(\frac{ \left| u \right| }{\lVert \nabla u \rVert_N} \right)^{\frac{N}{N-1}}} \, dx \leq e^{\frac{\lVert \nabla u \rVert^N_N}{\widetilde{\omega_N}}} \, C_N . $$

\medskip
\noindent
Passing to logarithm, we obtain

$$ \displaystyle  \ln \left( \frac{1}{|\Omega|} \int_{\Omega} e^u \, dx \right) \leq  \frac{1}{\widetilde{\omega_N}} \int_{\Omega} | \nabla u |^N \, dx + \ln C_N.
$$

\medskip
\noindent
If $u \equiv 0$, inequality \eqref{logMTinequalitydomains} still holds.

\end{proof}

\medskip
\noindent
A refinement of inequality \eqref{logMTinequalitydomains} for the case where $\Omega$ is the unit ball $B_1$ of $\mathbb{R}^N$ follows from the work of Carleson and Chang \cite{CarlesonChang}. More precisely, they proved the following result \cite[Lemma 1]{CarlesonChang}:

\begin{lemma}\label{CarlesonChanglemma} Let $\delta >0$ and let

$$ \displaystyle \mathcal{K}_{\delta}:=  \left\lbrace  \, C^1-\text{ function $w$ defined in } 0 \leq t < \infty, \; w(0)=0, \; w'(t) \geq 0, \;  \text{ and } \int_0^{+ \infty} w'(t)^N \, dt \leq \delta    \right\rbrace.$$

\medskip
\noindent
For any $c>0$, we have

\begin{equation}\label{stimaCC}
\displaystyle 
\sup_{w \in \mathcal{K}_{\delta}} \int_0^{+ \infty} e^{c \,  w(t) - t} \, dt < e^{\left(\frac{N-1}{N} \right)^{N-1} \left( \frac{c^N \delta}{N} \right)} e^{1 + \frac{1}{2} + ... + \frac{1}{N-1}}.
\end{equation}

\medskip
\noindent
Moreover, when $c^N \delta \to + \infty$, the inequality in \eqref{stimaCC} tends asymptotically to an equality.
\end{lemma}

\begin{proposition}[\textbf{Logarithmic M-T inequality on the unit ball $B_1$ of $\boldsymbol{\mathbb{R}^N}$}]\label{propositionlogMTinequalityCC}\

\medskip
\noindent
For any $u \in  W^{1,N}_0(B_1)$, we have

\begin{equation}\label{logCC}
\displaystyle \ln \left( \frac{1}{V_N} \int_{B_1} e^{u}   dx  \right) < \frac{1}{\widetilde{\omega_N}} \int_{B_1} |\nabla u|^N  dx + \sum_{k=1}^{N-1} \frac{1}{k}.
\end{equation}

\end{proposition}

\begin{proof}
Let $u \in  W_0^{1,N}(B_1)$. We may assume that $u \geq 0$ since we can replace $u$ by $|u|$ without increasing the integral of the gradient. Also, by exploiting the density of the space $C^1(B_1)$ in $W_0^{1,N}(B_1)$,  we may assume that $u$ has compact support and is in $C^1(B_1)$. 

\medskip
\noindent
We recall that the symmetric decreasing rearrangement of $u$ is the unique right-continuous radially symmetric and decreasing function $u^* : B_1 \to [0,+ \infty)$ such that

$$ \displaystyle \left| \{ x \in B_1 \, : u^*(x) > t \}  \right| =  \left| \{ x \in B_1 \, : u (x) > t \}  \right| \qquad \text{ for every }  t \geq 0.$$

\medskip
\noindent
Among the properties of $u^*$, we recall (see \cite{LiebLoss} for details) that

$$ \displaystyle \int_{B_1} \left| \nabla u^* \right|^N \, dx \leq \int_{B_1} \left| \nabla u \right|^N \, dx $$

\medskip
\noindent
and

$$ \displaystyle \int_{B_1} e^{ u^*} dx = \int_{B_1} e^{u} dx.$$

\medskip
\noindent
So if the thesis holds with $u^*$, then it will hold also for $u$.\\
Therefore let us assume  $u \in C^1(B_1)$ be radially symmetric and decreasing.\\
Since $u(x)$ is a function that depends just on $\left| x \right|$, it is of the form

$$ \displaystyle u(x)= g(\left| x \right|) \qquad \text{ for some } g: [0, 1]  \to \mathbb{R} \;  \text{ of class } C^1.$$

\medskip
\noindent
Let us now introduce the change of variable  $\displaystyle t= - N \ln |x| ,$ by which $  e^{-t}= \left| x \right|^N. $

\medskip
\noindent
Now, we denote by $ \displaystyle \tilde{c}:= N^{\frac{N-1}{N}} \omega_{N-1}^{\frac{1}{N}} $ and we set $\displaystyle w(t):= \tilde{c} u(x),$ by which 

$$ \displaystyle u(x) = \frac{ w(t)}{\tilde{c} }= \frac{ w( - N \ln |x| )}{ \tilde{c} } \qquad \text{ for any } x \in B_1 \setminus \{O\}.$$

\medskip
\noindent
Observe that 

$$ w(+ \infty):= \lim_{ t \to + \infty} w(t)= \tilde{c} u(O) \in \mathbb{R}^+.$$

\medskip
\noindent
The function $w$ is monotone increasing in $t$. Moreover,

\begin{align*}
\displaystyle 
\frac{1}{V_N} \int_{B_1} e^{u(x)}   dx = \frac{N}{\omega_{N-1}} \int_{B_1} e^{ \frac{w(- N \ln \left| x \right| )}{\tilde{c}}} \, dx = N \int_0^1 e^{ \frac{ w(- N \ln \rho) }{\tilde{c}}} \rho^{N-1} \, d \rho = \int_0^{+ \infty} e^{ \frac{w(t)}{\tilde{c}} - t}  \, d t
\end{align*} 

\medskip
\noindent
and

\begin{align*}
\displaystyle \int_{B_1} \left| \nabla u(x)  \right|^N \, dx = N \int_0^1 w'(- N \ln \rho )^N \frac{1}{\rho^N} \rho^{N-1} \, d \rho  = \int_0^{+ \infty} w'(t)^N \, dt.
\end{align*}

\medskip
\noindent
By applying Lemma \ref{CarlesonChanglemma} with $\displaystyle c:=1/\tilde{c}$ and $ \displaystyle \delta:=\int_0^{+ \infty} w'(t)^N \, dt= \int_{B_1} \left| \nabla u(x) \right| ,$ we get

\begin{align*}
\displaystyle \frac{1}{V_N} \int_{B_1} e^{u(x)}   dx = \int_0^{+ \infty} e^{ c \, w(t) - t}  \, d t  < e^{\frac{1}{\widetilde{\omega_N}} \int_{B_1} |\nabla u|^N  dx } e^{1 + \frac{1}{2} + ... + \frac{1}{N-1}},
\end{align*}

\medskip
\noindent
and passing to logarithm we conclude that

\begin{equation*}
\displaystyle \ln \left( \frac{1}{V_N} \int_{B_1} e^{u}   dx  \right) < \frac{1}{\widetilde{\omega_N}} \int_{B_1} |\nabla u|^N  dx + \sum_{k=1}^{N-1} \frac{1}{k}.
\end{equation*}

\end{proof}

\medskip
\noindent
Notice that by Lemma \ref{CarlesonChanglemma} we also deduce that inequality \eqref{logCC} is sharp, that is the functional $J:W^{1,N}_{0}(B_1) \to \mathbb{R}$ defined as 

\begin{equation*}
\displaystyle J(u) := \frac{1}{\widetilde{\omega_N}} \int_{B_1} | \nabla u |^N \, dx  - \ln \left( \frac{1}{V_N} \int_{B_1} e^u \, dx \right),
\end{equation*}

\medskip
\noindent
satisfies 

$$
\inf_{W^{1,N}_{0}(B_1)} J = - \sum_{k=1}^{N-1} \frac{1}{k}.
$$

\medskip
\noindent
Furthermore, as both the $L^N$-norm of the gradient and the integral average of the exponential are invariant under translation and dilation, inequality \eqref{logCC} holds for any ball in $\mathbb{R}^N$. While translation invariance is immediate, the invariance under dilation is verified below:

\begin{lemma}\label{invariancedilation}
Let $r>0$ and let us denote by $B_r:=B_r(O) \subset \mathbb{R}^N$, where $N \geq 2$.\\
If $u \in W_0^{1,N}(B_1)$ then $u_r  \in W_0^{1,N}(B_r)$, where $ \displaystyle u_r(x):=u\left(\frac{x}{r}\right).$ Moreover, 

$$ \displaystyle \frac{1}{|B_r|} \int_{B_r} e^{u_r(x)} \, dx = \frac{1}{|B_1|} \int_{B_1} e^{u(y)} \, dy \qquad \text{ and } \qquad \int_{B_r} \left| \nabla u_r(x) \right|^N \, dx = \int_{B_1} \left| \nabla u(y) \right|^N \, dy .$$
\end{lemma}

\begin{proof}
Let us consider the following diffeomorphism:

\begin{align*}
\displaystyle B_r & \stackrel{{\normalfont\mbox{$\varphi$}}}{\longrightarrow} B_1 \\
             x  &   \longmapsto   y=\varphi(x):= \frac{x}{r}. 
\end{align*}

\medskip
\noindent
Note that

$$
J_{\varphi}(x)= \left( {\begin{array}{ccccc}
    \frac{1}{r} & 0 & \cdots &  & 0  \\
    0 & \frac{1}{r} & 0 & \cdots & 0 \\
                                    \\
    \vdots & & & \ddots & \vdots    \\
    0 & \cdots & & & \frac{1}{r}     \\
  \end{array} } \right)  
  \qquad |\det J_{\varphi}(x)| = \frac{1}{r^N}, \qquad dy=\frac{1}{r^N} dx.
$$

\medskip
\noindent
We have

$$ \displaystyle \int_{B_r} |u_r(x)|^N \, dx = \int_{B_r} \left|u\left(\frac{x}{r}\right)\right|^N \, dx = r^N \int_{B_1} |u(y)|^N \, dy < + \infty,$$ 

\medskip
\noindent
and

$$ \displaystyle \frac{1}{|B_r|} \int_{B_r} e^{u_r(x)} \, dx = \frac{1}{V_N r^N} \int_{B_r}  e^{u \left(\frac{x}{r} \right)} \, dx = \frac{1}{|B_1|} \int_{B_1} e^{u\left(y \right)} \, dy.$$ 

\medskip
\noindent
Observe now that $\displaystyle \nabla u_r(x) = \nabla (u \circ \varphi)(x)= \nabla u (\varphi(x) ) J_\varphi (x) = \frac{1}{r} \nabla u \left( \frac{x}{r} \right).$\\
Hence

$$ \displaystyle \int_{B_r} | \nabla u_r(x) |^N \, dx = \frac{1}{r^N} \int_{B_r} \left| \nabla u \left( \frac{x}{r} \right)  \right|^N \, dx = \int_{B_1} \left| \nabla u(y) \right|^N \, dy  $$

\end{proof}

\medskip
\noindent
A generalization of inequality \eqref{logCC} for $N=2$ was provided in \cite{CLMP} (see also \cite{CCL}). Specifically, if $\Omega \subseteq \mathbb{R}^2$ is a bounded domain with smooth boundary, then for any $u \in H^1_0(\Omega)$ we have 
 
\begin{equation}\label{OnofriDomain}
\ln \left(\frac{1}{|\Omega|}\int_{\Omega} e^u \,dx \right) \le  \frac{1}{16\pi}\int_{\Omega} |\nabla u|^2 dx  + 1 + {4\pi \sup_{\Omega}\gamma(x) + \ln \frac{\pi}{|\Omega|}},
\end{equation}

\medskip
\noindent
where $\gamma$ is Robin's function of $\Omega$, which is defined for  $x\in \Omega$ as $\gamma (x) = H_x(x)$, where $H_x$ is the unique solution of 

$$
\begin{cases}
\Delta H_x = 0  & \text{ in }\Omega,\\
H_x(y) = \frac{1}{2\pi }\ln|x-y|  & \text{ on }\partial \Omega. 
\end{cases}
$$

\medskip
\noindent
Differently from \eqref{logCC},  \eqref{OnofriDomain} is  not always sharp. Indeed, denoting by

$$
F(\Omega):= \inf_{u\in H^1_0(\Omega)} \left( \frac{1}{16\pi}\int_{\Omega} |\nabla u|^2 dx - \ln \left(\frac{1}{|\Omega|}\int_{\Omega} e^u \,dx \right)\right),
$$

\medskip
\noindent
there are domains $\Omega$ for which 

\begin{equation}\label{nonsharp}
 F(\Omega)>- 1- {4\pi \sup_{\Omega}\gamma(x) - \ln \frac{\pi}{|\Omega|}}.
\end{equation}

\medskip
\noindent
The domains for which \eqref{OnofriDomain} is sharp are known in the literature as domains of first kind, while the ones for which  \eqref{nonsharp} holds are known as domains of second kind. Inequality \eqref{OnofriDomain} plays a crucial role in the study of the mean field equation 

\begin{equation}\label{mean}
-\Delta u  = \frac{\rho e^u}{\int_{\Omega}e^u \,dx} \quad \text{ in }\Omega,
\end{equation}

\medskip
\noindent
which appears in the statistical mechanic description of vortex formations in 2d-models for turbulent flows in fluid dynamics (see e.g. \cite{JM,CLMP}). We refer to \cite{CL1,CL2,Malchiodi,Bart,BarMal} about existence results  for  \eqref{mean} and related problems.

\section{Equivalence between Euclidean Onofri Inequality and Logarithmic Moser-Trudinger Inequality}\label{SezEquiv}  

\medskip
\noindent
This Section is devoted to the proof of Theorem \ref{equivalence}. Let us recall that $J:W^{1,N}_{0}(B_1) \to \mathbb{R}$ and $I:W_{\mu_N}(\R^N) \to \mathbb{R}$ are the functionals defined in \eqref{Jfunctional} and \eqref{Ifunctional}, as

\begin{equation*}
\displaystyle J(u) := \frac{1}{\widetilde{\omega_N}} \int_{B_1} | \nabla u |^N \, dx  - \ln \left( \frac{1}{V_N} \int_{B_1} e^u \, dx \right),
\end{equation*}

\medskip
\noindent
and

\begin{equation*}
\displaystyle I(u):=   \frac{1}{\widetilde{\omega_N}} \int_{\mathbb{R}^N} H_N(u,\mu_N) \, dx + \int_{\mathbb{R}^N} u \, d\mu_N- \ln \left(\int_{\mathbb{R}^N} e^u \, d\mu_N \right). 
\end{equation*}

\medskip
\noindent
We will prove that 

$$
\inf_{W_{\mu_N}(\R^N)} I = \inf_{W^{1,N}_0(B_1)}  J +\sum_{k=1}^{N-1}\frac{1}{k}.
$$

\medskip
\noindent
{\mbox {\it Proof of Theorem~\ref{equivalence}.~}}
We begin with the proof of $ `` \geq "$.\\
Let $u \in W_{\mu_N}(\mathbb{R}^N)$. Since $W_{\mu_N}(\mathbb{R}^N)= \overline{C_0^{\infty} (\mathbb{R}^N)}^{\lVert \cdot \rVert_{\mu_N}}$, we can assume $u \in  C_0^{\infty}(\mathbb{R}^N)$ and the thesis will follow by density.\\
Let us consider the following cutoff function $\Psi \in C^{\infty}_0(\mathbb{R}^N)$ such that $0 \leq \Psi \leq 1$, $ \Psi \equiv 1$ if $ \displaystyle |x| \leq \frac{1}{2},$ and $\Psi \equiv 0 $ if $\displaystyle  |x| \geq 1$.\\
Let $r>0$ and let us define $\displaystyle \Psi_r(x):=\Psi\left(\frac{x}{r}\right)$, so that $ \Psi_r \equiv 1$ if $ \displaystyle |x| \leq \frac{r}{2},$ and $\Psi_r \equiv 0 $ if $\displaystyle  |x| \geq r.$ Obviously, $\displaystyle \Psi_r(x) \to 1$ as $r \to +\infty,$ for any $ x \in \mathbb{R}^N.$\\
For simplicity of notation, let's say $B_r:=B_r(O);$ moreover, we will denote by $o(1)$ any function that goes to zero as $r \to + \infty.$\\
Using the characteristic function $\chi$, we can write

$$ \displaystyle \int_{B_r} e^{u(x)\Psi_r(x) + v_N(x)}   \, dx =\int_{\mathbb{R}^N} e^{u(x)\Psi_r(x)+ v_N(x)}\chi_{B_r}(x)  \, dx.$$

\medskip
\noindent
Now we observe that $|u|$ is bounded by a constant $M_{u}$, so for any $ x \in \mathbb{R}^N$ and for any $r>0,$

$$ \displaystyle \left| e^{u(x)\Psi_r(x) + v_N(x)}\chi_{B_r}(x) \right| \leq e^{M_{u}+v_N(x)} =e^{M_{u}} \mu_N(x) \in L^1(\mathbb{R}^N).$$

\medskip
\noindent
Moreover, for any $x \in \mathbb{R}^N$,

$$ \displaystyle  e^{u(x)\Psi_r(x) + v_N(x)}\chi_{B_r}(x) \to e^{u(x) + v_N(x)} \quad  \text{ as } r \to + \infty. $$

\medskip
\noindent
Therefore we can apply dominated convergence theorem and say that, as $r \to + \infty$,

$$ \displaystyle  \int_{B_r} e^{u(x)\Psi_r(x) + v_N(x)}   \, dx  \to  \int_{\mathbb{R}^N} e^{u(x)+v_N(x)}\, dx  = \int_{\mathbb{R}^N} e^u \, d\mu_N. $$

\medskip
\noindent
We can write

\begin{equation*}
\displaystyle \ln \left( \int_{\mathbb{R}^N} e^u \, d\mu_N \right)= \lim_{r \to + \infty} \ln \left(  \int_{B_r} e^{u(x)\Psi_r(x) + v_N(x)}   \, dx \right),
\end{equation*}

\medskip
\noindent
therefore we can estimate

$$\displaystyle  \ln \left(  \int_{B_r} e^{u(x)\Psi_r(x) + v_N(x)}   \, dx \right).$$

\medskip
\noindent
In general the function $u\Psi_r + v_N \notin W^{1,N}_0(B_r)$.\\
If we define $w:= u\Psi_r + v_N - v_N(r)$, (where with $v_N(r)$ we make a slight abuse of notation to denote a constant function on $B_r$ whose value is $v_N(x)$ as $|x|=r$), we have $w \in  W^{1,N}_0(B_r)$.

\begin{align*}
\displaystyle 
    \ln \left( \int_{B_r} e^{u(x)\Psi_r(x) + v_N(x)}   \, dx   \right) 
& = \ln \left( \int_{B_r} e^{w(x)+v_N(r)}   \, dx   \right) \\
& = \ln \left( \mu_N(r) |B_r| \frac{1}{|B_r|}\int_{B_r} e^{w(x)}   \, dx   \right) \\
& = \ln \left( \frac{1}{V_N \left( 1+ r^{\frac{N}{N-1}} \right)^N} V_N r^N \right) + \ln \left( \frac{1}{|B_r|} \int_{B_r} e^{w(x)}   \, dx \right)  \\
& = N\ln \left( \frac{r}{  1+ r^{\frac{N}{N-1}}  }\right) + \ln \left( \frac{1}{|B_r|} \int_{B_r} e^{w(x)}   \, dx \right) \\
& = -\frac{N}{N-1} \ln r -N \ln \left( 1 + \frac{1}{r^{\frac{N}{N-1}}} \right) + \ln \left( \frac{1}{|B_r|} \int_{B_r} e^{w(x)}   \, dx \right)\\
& = -\frac{N}{N-1} \ln r + \ln \left( \frac{1}{|B_r|} \int_{B_r} e^{w(x)}   \, dx \right) + o(1).
\end{align*}

\medskip
\noindent
Since $w \in  W^{1,N}_0(B_r)$, we have 

$$ \displaystyle \inf_{W^{1,N}_0(B_1)}  J <  \frac{1}{\widetilde{\omega_N}} \int_{B_r} | \nabla w |^N \, dx - \ln \left( \frac{1}{|B_r|} \int_{B_r} e^{w(x)}   \, dx \right). $$

\medskip
\noindent
Therefore

\begin{align*}
\displaystyle 
\ln \left( \int_{B_r} e^{u\Psi_r} \, d \mu_N   \right) < \frac{1}{\widetilde{\omega_N}} \int_{B_r} | \nabla w |^N \, dx  -\frac{N}{N-1} \ln r - \inf_{W^{1,N}_0(B_1)}  J  + o(1). 
\end{align*}

\medskip
\noindent
We need to compute $ \displaystyle \frac{1}{\widetilde{\omega_N}} \int_{B_r} | \nabla w |^N \, dx$.\\
Since

$$ \displaystyle \nabla w = \nabla(u \Psi_r) + \nabla v_N$$

\medskip
\noindent
and

\begin{align*}
 \displaystyle
H_N(u \Psi_r, \mu_N) = | \nabla w|^N - |\nabla v_N |^N - N |\nabla v_N |^{N-2} \nabla  v_N  \cdot \nabla(u \Psi_r),
\end{align*}

\medskip
\noindent
we have

$$ \displaystyle | \nabla w|^N= H_N(u \Psi_r, \mu_N)  + |\nabla v_N |^N + N |\nabla v_N |^{N-2} \nabla  v_N  \cdot \nabla(u \Psi_r),$$

\medskip
\noindent
by which

\begin{align*}
\displaystyle
     \int_{B_r} | \nabla w |^N \, dx 
& =  \int_{B_r} H_N(u \Psi_r, \mu_N) \, dx \\
& +  \int_{B_r} |\nabla v_N |^N \, dx + \int_{B_r} |\nabla v_N |^{N-2} \nabla  v_N  \cdot \nabla(u \Psi_r) \, dx.
\end{align*}

\medskip
\noindent
By Green's formula and taking into account  \eqref{GaussGreen}, we get

\begin{align*}
\displaystyle
 \frac{N}{\widetilde{\omega_N}} \int_{B_r} |\nabla v_N |^{N-2} \nabla  v_N  \cdot \nabla(u \Psi_r) \, dx 
&= \frac{N}{\widetilde{\omega_N}} \left( \int_{ \partial B_r} u \Psi_r |\nabla v_N|^{N-2} \nabla v_N \cdot  \nu \, d\sigma -\int_{B_r} u \Psi_r \left( \Delta_N v_N \right) \, dx \right) \\
& = \frac{N}{\widetilde{\omega_N}} \int_{B_r} u \Psi_r (-\Delta_N v_N)  \, dx \\
& = \frac{N^{N+1}\left( \frac{N}{N-1} \right)^{N-1} V_N }{\widetilde{\omega_N}} \int_{B_r} u \Psi_r \mu_N \, dx \\
& = \int_{B_r} u \Psi_r  \mu_N \, dx .
\end{align*}

\medskip
\noindent
Now, we need to compute $\displaystyle \frac{1}{\widetilde{\omega_N}} \int_{B_r} |\nabla v_N |^N \, dx.$

\begin{align}\label{contovn}
\displaystyle
\frac{1}{\widetilde{\omega_N}} \int_{B_r} |\nabla v_N |^N \, dx 
& = \frac{N^{2 N}}{(N-1)^N \widetilde{\omega_N}}\int_{B_r} \frac{|x|^{\frac{N}{N-1}}}{\left(1+|x|^{\frac{N}{N-1}}\right)^N } \, dx  \nonumber \\
&= \frac{N^{2 N} \omega_{N-1}}{(N-1)^N \widetilde{\omega_N}} \int_0^r \frac{\rho^{\frac{N}{N-1}}}{\left(1+\rho^{\frac{N}{N-1}} \right)^N} \rho^{N-1} \, d \rho \nonumber\\
&= \frac{N}{N-1}  \int_0^r \frac{\rho^{N}}{\left(1+\rho^{\frac{N}{N-1}} \right)^N} \rho^{\frac{1}{N-1}} \, d \rho  \nonumber \\
& = \int_0^{r^{\frac{N}{N-1}}} \frac{t^{N-1}}{\left( 1+t \right)^N}  \, dt \\
& = \int_0^{r^{\frac{N}{N-1}}} \frac{1}{ 1+t}  \, dt - \sum_{k=0}^{N-2} \binom{N-1}{k} \int_0^{r^{\frac{N}{N-1}}} \frac{t^k}{\left( 1+t \right)^N}  \, dt \nonumber\\
& =  \frac{N}{N-1} \ln r - \sum_{k=0}^{N-2} \binom{N-1}{k} \int_0^{r^{\frac{N}{N-1}}} \frac{t^k}{\left( 1+t \right)^N}  \, dt + o(1). \nonumber
\end{align}

\medskip
\noindent
By induction, we prove that

\begin{equation}\label{ind}
\displaystyle  \binom{N-1}{k} \int_0^{r^{\frac{N}{N-1}}} \frac{t^k}{\left( 1+t \right)^N}  \, dt= \frac{1}{N-k-1} + o(1)
\end{equation}
for any  $k=0,...,N-2.$

\medskip
\noindent
If $k=0,$ we have

\begin{align*}
\displaystyle
\binom{N-1}{k} \int_0^{r^{\frac{N}{N-1}}} \frac{t^k}{\left( 1+t \right)^N}  \, dt
= \int_0^{r^{\frac{N}{N-1}}} \frac{1}{\left( 1+t \right)^N}  \, dt 
= \frac{1}{1-N} \left[\frac{1}{(1+t)^{N-1}} \right]_0^{r^{\frac{N}{N-1}}} 
= \frac{1}{N-1} + o(1).
\end{align*}

\medskip
\noindent
Now, let $N \geq 3,$ $k \geq 1,$  and let assume $\eqref{ind}$ to hold for $k-1.$

\begin{align*}
\displaystyle
&   \binom{N-1}{k} \int_0^{r^{\frac{N}{N-1}}} \frac{t^k}{\left( 1+t \right)^N}  \, dt \\
& = \binom{N-1}{k} \left[ \int_0^{r^{\frac{N}{N-1}}} t^k \left( 1+t \right)^{-N}  \, dt  \right]  \\
& = \binom{N-1}{k} \left[ \left[ t^k \frac{{\left( 1+t \right)^{1-N}}}{1-N} \right]_0^{ r^{\frac{N}{N-1}} }  - \frac{k}{1-N}\int_0^{r^{\frac{N}{N-1}}} t^{k-1}\left( 1+t \right)^{1-N}  \, dt      \right] \\
& = \binom{N-1}{k} \left[ o(1) + \frac{k}{N-1} \int_0^{r^{\frac{N}{N-1}}} \frac{t^{k-1}(1+t)}{\left( 1+t \right)^N}  \, dt       \right] \\
& = \binom{N-1}{k} \left[ o(1) + \frac{\frac{k}{N-1}}{\binom{N-1}{k-1}} \binom{N-1}{k-1} \int_0^{r^{\frac{N}{N-1}}} \frac{t^{k-1}}{\left( 1+t \right)^N}  \, dt   + \frac{k}{N-1}\int_0^{r^{\frac{N}{N-1}}} \frac{t^k}{\left( 1+t \right)^N}  \, dt      \right] \\
& = \binom{N-1}{k} \left[ o(1) + \frac{\frac{k}{N-1}}{\binom{N-1}{k-1}} \left( \frac{1}{N-k} + o(1)  \right)   + \frac{k}{N-1}\int_0^{r^{\frac{N}{N-1}}} \frac{t^k}{\left( 1+t \right)^N}  \, dt      \right] \\
& = \frac{\binom{N-1}{k}  }{ \binom{N-1}{k-1}  } \frac{k}{(N-1)(N-k)} + \frac{k}{N-1} \binom{N-1}{k} \int_0^{r^{\frac{N}{N-1}}} \frac{t^k}{\left( 1+t \right)^N}  \, dt + o(1) \\
& = \frac{1}{N-1}  +  \frac{k}{N-1} \binom{N-1}{k} \int_0^{r^{\frac{N}{N-1}}} \frac{t^k}{\left( 1+t \right)^N}  \, dt + o(1).
\end{align*}

\medskip
\noindent
Therefore

$$ \displaystyle  \frac{N-k-1}{N-1} \, \binom{N-1}{k}\int_0^{r^{\frac{N}{N-1}}} \frac{t^k}{\left( 1+t \right)^N}  \, dt = \frac{1}{N-1} + o(1), $$

\medskip
\noindent
by which

$$ \displaystyle  \binom{N-1}{k} \int_0^{r^{\frac{N}{N-1}}} \frac{t^k}{\left( 1+t \right)^N}  \, dt= \frac{1}{N-k-1} + o(1). $$

\medskip
\noindent
By $\eqref{ind}$ we deduce that

\begin{equation}\label{sum}
\displaystyle \sum_{k=0}^{N-2} \binom{N-1}{k} \int_0^{r^{\frac{N}{N-1}}} \frac{t^k}{\left( 1+t \right)^N}  \, dt= \sum_{k=0}^{N-2} \frac{1}{N-k-1} + o(1)= \sum_{k=1}^{N-1} \frac{1}{k} + o(1) .
\end{equation}

\medskip
\noindent
Finally,

\begin{equation}\label{contoGradvnallaN}
\displaystyle \frac{1}{\widetilde{\omega_N}} \int_{B_r} |\nabla v_N |^N \, dx =\frac{N}{N-1} \ln r - \sum_{k=1}^{N-1} \frac{1}{k} + o(1). 
\end{equation}

\medskip
\noindent
Summing up, we have shown that

\begin{align*}
\displaystyle 
 \inf_{W^{1,N}_0(B_1)}  J + \sum_{k=1}^{N-1} \frac{1}{k} <  I(u \Psi_r) + o(1).
\end{align*}

\medskip
\noindent
Now, we pass to the limit as $r \to + \infty.$ 

\medskip
\noindent
We observe that $|u|$ is bounded by a constant $M_{u}$, so for any $ x \in \mathbb{R}^N$ and for any $r>0,$

$$ \displaystyle \left|  u(x) \Psi_r(x)  \mu_N(x) \right| \leq M_{u} \mu_N(x) \in L^1(\mathbb{R}^N).$$

\medskip
\noindent
Moreover, for any $x \in \mathbb{R}^N$,

\medskip
\noindent
$$ \displaystyle   u(x) \Psi_r(x)  \mu_N(x) \to u(x)  \mu_N(x) \quad  \text{ as } r \to + \infty. $$

\medskip
\noindent
We can apply dominated convergence theorem to say that

$$ \displaystyle   \int_{\mathbb{R}^N} u \Psi_r  \, d\mu_N  \to  \int_{\mathbb{R}^N} u  d\mu_N. $$

\medskip
\noindent
Now we observe that, for any $x \in \mathbb{R}^N$,

$$ \displaystyle H_N(u(x) \Psi_r(x), \mu_N(x))  \to H_N(u(x), \mu_N(x)) \qquad \text{ as } r \to + \infty.$$

\medskip
\noindent
Moreover, considering Lemma \ref{stima} and $u \in  C_0^{\infty}(\mathbb{R}^N)$, there exists a function $v \in  L^1(\mathbb{R}^N, dx)$ such that, for any $ x \in \mathbb{R}^N$ and for any $r>1$, we have 

\begin{align*}
\displaystyle
        0 \leq  H_N(u(x) \Psi_r(x), \mu_N(x)) \leq v(x).
\end{align*}

\medskip
\noindent
In fact

\begin{align*}
\displaystyle
        0 
& \leq  H_N(u(x) \Psi_r(x), \mu_N(x)) \\
& \leq  \frac{c_N}{2} \left( |\nabla \left( u(x) \Psi_r(x) \right) |^N + |\nabla \left( u(x)\Psi_r(x) \right)|^2 |\nabla v_N(x) |^{N-2}   \right) \\
& \leq  C \left( \left( |u(x)| +  |\nabla u(x)| \right)^N + \left( |u(x)| +  |\nabla u(x)| \right)^2 |\nabla v_N(x) |^{N-2}   \right).
\end{align*}

\medskip
\noindent
Since  $u \in  C_0^{\infty}(\mathbb{R}^N)$, we have that 

$$ \displaystyle v(x):=C \left( \left( |u(x)| +  |\nabla u(x)| \right)^N + \left( |u(x)| +  |\nabla u(x)| \right)^2 |\nabla v_N(x) |^{N-2}   \right) \in L^1(\mathbb{R}^N, dx).$$

\medskip
\noindent
Therefore we can apply dominated convergence theorem and say that

$$ \displaystyle \frac{1}{\widetilde{\omega_N}} \int_{\mathbb{R}^N} H_N(u \Psi_r, \mu_N) \, dx  \to \frac{1}{\widetilde{\omega_N}} \int_{\mathbb{R}^N} H_N(u, \mu_N) \, dx .$$

\medskip
\noindent
Finally,

\begin{align*}
\displaystyle 
 \inf_{W^{1,N}_0(B_1)}  J + \sum_{k=1}^{N-1} \frac{1}{k} \leq  I(u),
\end{align*}

\medskip
\noindent
by which

$$
 \inf_{W^{1,N}_0(B_1)}  J +\sum_{k=1}^{N-1}\frac{1}{k} \leq \inf_{C_0^{\infty}(\R^N)} I.
$$

\medskip
\noindent
By density, that is by the proof of Theorem \ref{EuclideanOnofriNdim}, we obtain

$$
 \inf_{W^{1,N}_0(B_1)}  J +\sum_{k=1}^{N-1}\frac{1}{k} \leq \inf_{W_{\mu_N}(\R^N)} I.
$$

\medskip
\noindent
Now, we prove $ `` \leq "$.\\
Let $u  \in W^{1,N}_0(B_1)$ and let $r>0$. We define the function $u_r: \mathbb{R}^N \to \mathbb{R}$ as 

$$ \displaystyle 
u_r(x)=
\begin{cases}
\displaystyle u \left( \frac{x}{r} \right) - v_N(x)+v_N(r) & \text{if $|x| \leq r$,} \bigskip \\
0 & \text{if $|x| > r$.}
\end{cases}
$$

\medskip
\noindent
By Lemma \ref{W1nmun}, we have $u_r \in W_{\mu_N}(\mathbb{R}^N)$, by which

$$ \displaystyle  \inf_{W_{\mu_N}(\R^N)} I \leq I(u_r),
$$

\medskip
\noindent
where

\begin{equation}\label{calcoloMoserTrudinger}
\displaystyle I(u_r) = \frac{1}{\widetilde{\omega_N}} \int_{\mathbb{R}^N} H_N(u_r,\mu_N) \, dx + \int_{\mathbb{R}^N} u_r \, d\mu_N - \ln \left( \int_{\mathbb{R}^N} e^{u_r} \, d\mu_N \right).
\end{equation}

\medskip
\noindent
We need to compute the terms in \eqref{calcoloMoserTrudinger}. Firstly, we observe that

$$
\displaystyle 
\ln \left( \int_{\mathbb{R}^N} e^{u_r} \, d\mu_N \right) = \ln \left( \mu_N(r) \int_{B_r} e^{u \left( \frac{x}{r} \right)} \, dx + \int_{\mathbb{R}^N \setminus B_r}  \mu_N(x) \, dx \right).
$$

\medskip
\noindent
In particular,

\begin{align*}
\displaystyle
   \mu_N(r) \int_{B_r} e^{u \left( \frac{x}{r} \right)} \, dx 
& = \frac{1}{V_N \left( 1+r^{\frac{N}{N-1}}  \right)^N} \int_{B_r} e^{u \left( \frac{x}{r} \right)} \, dx \\
& = \frac{r^N}{\left( 1+r^{\frac{N}{N-1}}  \right)^N} \frac{1}{V_N} \int_{B_r} e^{u \left( \frac{x}{r} \right)} \frac{1}{r^N} \, dx \\
& = \left( \frac{r}{1+r^{\frac{N}{N-1}}}  \right)^N \frac{1}{V_N} \int_{B_1} e^{u(y)}  \, dy;
\end{align*}

\medskip
\noindent
and

\begin{align*}
\displaystyle
    \int_{\mathbb{R}^N \setminus B_r}  \mu_N(x) \, dx = \int_{\mathbb{R}^N \setminus B_r}  \frac{1}{V_N \left( 1+|x|^{\frac{N}{N-1}}  \right)^N}  \, dx  
=   \left[ \left( \frac{t}{1+t} \right)^{N-1} \right]_{r^{\frac{N}{N-1}}}^{+ \infty}  
= 1 - \frac{r^N}{\left( 1+r^{\frac{N}{N-1}}  \right)^{N-1}}.
\end{align*}

\medskip
\noindent
Therefore

\begin{align*}
\displaystyle
 \ln \left( \int_{\mathbb{R}^N} e^{u_r} \, d\mu_N \right) 
& = \ln \left( \int_{B_r} e^{u \left( \frac{x}{r} \right)} \mu_N(r) \, dx + \int_{\mathbb{R}^N \setminus B_r}  \mu_N(x) \, dx \right)\\
& = \ln \left( \left( \frac{r}{1+r^{\frac{N}{N-1}}}  \right)^N \frac{1}{V_N} \int_{B_1} e^{u}  \, dx + 1 - \frac{r^N}{\left( 1+r^{\frac{N}{N-1}}  \right)^{N-1}} \right) \\
&= \ln \left(   \left( \frac{r}{1+r^{\frac{N}{N-1}}}  \right)^N  \left[ \frac{1}{V_N} \int_{B_1} e^{u}  \, dx + \left( \frac{1+r^{\frac{N}{N-1}}}{r}  \right)^N - \left( 1 + r^{\frac{N}{N-1}}  \right)  \right]  \right) \\
& = N \ln \left( \frac{r}{1+r^{\frac{N}{N-1}}}  \right) + \ln \left( \frac{1}{V_N} \int_{B_1} e^{u}  \, dx + \frac{ (N-1)r^N + \sum_{k=0}^{N-2} \binom{N}{k} r^{\frac{Nk}{N-1}}}{r^N}    \right) \\
& = -\frac{N}{N-1} \ln r + \ln \left( \frac{1}{V_N} \int_{B_1} e^{u}  \, dx + N-1 + o(1) \right) + o(1).\\
\end{align*}

\medskip
\noindent
Now, we compute $ \displaystyle \frac{1}{\widetilde{\omega_N}} \int_{\mathbb{R}^N} H_N(u_r,\mu_N) \, dx.$

\medskip
\noindent
Since

$$ \displaystyle H_N(u, \mu_N) = | \nabla u + \nabla v_N |^N - |\nabla v_N |^N - N |\nabla  v_N |^{N-2} \nabla  v_N \cdot \nabla u $$

\medskip
\noindent
and

$$ \displaystyle
\nabla u_r(x)=
\begin{cases}
\displaystyle \frac{1}{r} \nabla_y u \left( \frac{x}{r} \right) - \nabla v_N(x) & \text{if $|x| \leq r$,} \bigskip \\
0 & \text{if $|x| > r$,}
\end{cases}
$$

\medskip
\noindent
we have

$$ \displaystyle
H_N(u_r, \mu_N)=
\begin{cases}
\displaystyle \frac{1}{r^N} \left|\nabla_y u \left( \frac{x}{r} \right) \right|^N + (N-1) |\nabla v_N|^N - \frac{N}{r} |\nabla  v_N |^{N-2} \nabla  v_N \cdot \nabla_y u \left( \frac{x}{r} \right) & \text{if $|x| \leq r$,} \bigskip \\
0 & \text{if $|x| > r$.}
\end{cases}
$$

\medskip
\noindent
Therefore

\begin{align*}
\displaystyle 
&   \frac{1}{\widetilde{\omega_N}} \int_{\mathbb{R}^N} H_N(u_r,\mu_N) \, dx \\
& = \frac{1}{\widetilde{\omega_N}} \int_{B_r} \left| \nabla_y u \left( \frac{x}{r} \right) \right|^N \frac{1}{r^N} \, dx + \frac{N-1}{\widetilde{\omega_N}} \int_{B_r} |\nabla v_N|^N \, dx - \frac{N}{ \widetilde{\omega_N} r } \int_{B_r} |\nabla  v_N |^{N-2} \nabla  v_N \cdot \nabla_y u \left( \frac{x}{r} \right) \, dx \\
& = \frac{1}{\widetilde{\omega_N}} \int_{B_1} |\nabla u |^N \, dy + \frac{N-1}{\widetilde{\omega_N}} \int_{B_r} |\nabla v_N|^N \, dx - \frac{N}{ \widetilde{\omega_N}} \int_{B_r} |\nabla  v_N |^{N-2} \nabla  v_N \cdot \nabla_x u \left( \frac{x}{r} \right) \, dx. \\
\end{align*}

\medskip
\noindent
Considering \eqref{contoGradvnallaN}, we have 

$$ \displaystyle \frac{N-1}{\widetilde{\omega_N}} \int_{B_r} |\nabla v_N|^N \, dx= N \ln r - (N-1)\sum_{k=1}^{N-1} \frac{1}{k}  + o(1). $$

\medskip
\noindent
Moreover, by Green's formula and taking into account \eqref{GaussGreen}, we obtain

\begin{align*}
\displaystyle
&   - \frac{N}{ \widetilde{\omega_N}} \int_{B_r} |\nabla  v_N |^{N-2} \nabla  v_N \cdot \nabla_x u \left( \frac{x}{r} \right) \, dx \\
& = -\frac{N}{\widetilde{\omega_N}} \left( \int_{ \partial B_r} u \left( \frac{x}{r} \right) |\nabla v_N|^{N-2} \nabla v_N \cdot  \nu \, d\sigma - \int_{B_r} u \left( \frac{x}{r} \right) \left( \Delta_N v_N \right) \, dx \right)\\
& = \frac{N}{\widetilde{\omega_N}} \int_{B_r} u \left( \frac{x}{r} \right) (\Delta_N v_N)  \, dx \\
& = -\frac{N^{N+1}\left( \frac{N}{N-1} \right)^{N-1} V_N }{\widetilde{\omega_N}} \int_{B_r} u \left( \frac{x}{r} \right) \mu_N \, dx \\
& = -\int_{B_r} u \left( \frac{x}{r} \right) \, d \mu_N  .\\
\end{align*}

\noindent
Therefore 

\begin{align*}
\displaystyle  
&   \frac{1}{\widetilde{\omega_N}} \int_{\mathbb{R}^N} H_N(u_r,\mu_N) \, dx \\
& = \frac{1}{\widetilde{\omega_N}} \int_{B_1} |\nabla u|^N dx -\int_{B_r} u \left( \frac{x}{r} \right)  \, d\mu_N  + N \ln r - (N-1)\sum_{k=1}^{N-1} \frac{1}{k}  + o(1).
\end{align*}

\noindent
It remains to compute the term $ \displaystyle \int_{\mathbb{R}^N} u_r \, d\mu_N $ in \eqref{calcoloMoserTrudinger} .

\begin{align*}
\displaystyle
  \int_{\mathbb{R}^N} u_r \, d\mu_N  = \int_{B_r} u \left( \frac{x}{r} \right) \, d\mu_N  - \int_{B_r} v_N(x) \, d\mu_N  + v_N(r) \int_{B_r} \,  d\mu_N . 
\end{align*}

\medskip
\noindent
In particular,

\begin{align*}
\displaystyle
    v_N(r) \int_{B_r}  d\mu_N  
&= \ln \left(\frac{1}{V_N \left( 1+ r^{\frac{N}{N-1}} \right)^N} \right)\left[ \left( \frac{t}{1+t} \right)^{N-1} \right]_0^{r^{\frac{N}{N-1}}} \\
& = - \frac{ r^N \ln V_N }{ \left( 1+ r^{\frac{N}{N-1}} \right)^{N-1}  } - N \frac{ r^N\ln \left( 1+ r^{\frac{N}{N-1}}  \right)}{ \left( 1+ r^{\frac{N}{N-1}}  \right)^{N-1}}  \\
& = - \ln V_N - \frac{N^2}{N-1} \ln r + o(1), \\
\end{align*}

\medskip
\noindent
and

\begin{align*}
\displaystyle
 - \int_{B_r}   v_N(x)\, d \mu_N 
& = - \int_{B_r} \frac{1}{V_N \left( 1 + |x|^{\frac{N}{N-1}}  \right)^N} \ln \left( \frac{1}{V_N \left( 1 + |x|^{\frac{N}{N-1}}  \right)^N} \right) \, dx\\
&= \frac{\omega_{N-1}}{V_N} \int_0^r \frac{\ln \left(V_N \left( 1 + \rho^{\frac{N}{N-1}} \right)^N \right)}{\left( 1 + \rho^{\frac{N}{N-1}}  \right)^N}  \rho^{N-1} \, d \rho \\
& = N \ln V_N \int_0^r \frac{\rho^{N-1}}{ \left( 1 + \rho^{\frac{N}{N-1}}  \right)^N} \, d \rho + N^2 \int_0^r \frac{ \rho^{N-1} \ln \left( 1 + \rho^{\frac{N}{N-1}} \right)}{\left( 1 + \rho^{\frac{N}{N-1}}  \right)^N} \, d \rho \\
& = \ln V_N + N^2 \int_0^r \frac{ \rho^{N-1} \ln \left( 1 + \rho^{\frac{N}{N-1}} \right)}{\left( 1 + \rho^{\frac{N}{N-1}}  \right)^N} \, d \rho + o(1).
\end{align*}

\medskip
\noindent
Considering \eqref{contovn} and \eqref{sum}, we have
 
\begin{align*}
\displaystyle
&   N^2 \int_0^r \frac{ \rho^{N-1} \ln \left( 1 + \rho^{\frac{N}{N-1}} \right)}{\left( 1 + \rho^{\frac{N}{N-1}}  \right)^N} \, d \rho \\
& = N(N-1) \int_0^{r^{\frac{N}{N-1}}} \frac{ t^{N-2} \ln \left( 1 + t \right)}{\left( 1 + t \right)^N} \, dt \\
& = N(N-1) \int_0^{r^{\frac{N}{N-1}}} \frac{1}{\left( 1+t \right)^2} \left( \frac{t}{1+t}  \right)^{N-2} \ln \left( 1+t \right) \, dt \\
& = N\left(   \left[ \left( \frac{t}{1+t}  \right)^{N-1} \ln \left( 1+t \right) \right]_0^{r^{\frac{N}{N-1}}} -  \int_0^{r^{\frac{N}{N-1}}}  \frac{t^{N-1}}{\left( 1+t \right)^N} \, dt   \right) \\
& = N \left( \frac{ r^N- \left( 1 + r^{\frac{N}{N-1}}  \right)^{N-1}    }{\left(1+r^{\frac{N}{N-1}}  \right)^{N-1}} \ln \left( 1+r^{\frac{N}{N-1}} \right)   +     \sum_{k=1}^{N-1} \frac{1}{k} + o(1)  \right) \\
& = N \sum_{k=1}^{N-1} \frac{1}{k} + o(1).
\end{align*}

\medskip
\noindent
Therefore

$$ \displaystyle  \int_{\mathbb{R}^N} u_r \, d\mu_N = \int_{B_r} u \left( \frac{x}{r} \right) \,  d\mu_N  - \frac{N^2}{N-1} \ln r + N \sum_{k=1}^{N-1} \frac{1}{k} + o(1). $$ 

\medskip
\noindent
Summing up, we have shown that 

\begin{align*}
\displaystyle 
 \inf_{W_{\mu_N}(\R^N)} I  
 \leq \frac{1}{\widetilde{\omega_N}} \int_{B_1} |\nabla u|^N  dx + \sum_{k=1}^{N-1} \frac{1}{k} - \ln \left( \frac{1}{V_N} \int_{B_1} e^{u}   dx + N-1 + o(1) \right) + o(1). 
\end{align*} 

\medskip
\noindent
Passing to limit as $r \to + \infty$, we obtain 

$$ \displaystyle \inf_{W_{\mu_N}(\R^N)} I  \leq \frac{1}{\widetilde{\omega_N}} \int_{B_1} |\nabla u|^N \, dx + \sum_{k=1}^{N-1} \frac{1}{k} - \ln \left( \frac{1}{V_N} \int_{B_1} e^{u}  \, dx + N-1 \right), $$

\medskip
\noindent
by which 

$$ \displaystyle \inf_{W_{\mu_N}(\R^N)} I  < J(u) + \sum_{k=1}^{N-1} \frac{1}{k}.$$

\medskip
\noindent
Finally,

$$ \displaystyle \inf_{W_{\mu_N}(\R^N)} I  \leq  \inf_{W^{1,N}_0(B_1)}  J + \sum_{k=1}^{N-1} \frac{1}{k}.$$

\qed

\section{Insight on $N=2$: connection with $\mathbb{S}^2$.}\label{SezN=2}

\bigskip
\noindent
Among the reasons that led Moser to improve Trudinger's result, it is important to mention  Nirenberg's problem on $\Sm^2$, which consists in describing the set of smooth functions that can be obtained as Gaussian curvature of a Riemannian metric conformally equivalent to the standard spherical metric. 

Let $\Sigma$ be a smooth surface, we recall that two Riemannian  metrics $g_0$ and $g$ on $\Sigma$ are said to be pointwise conformally equivalent if there exists  $u\in C^\infty(\Sigma)$ such that $g = e^u g_0$ on $\Sigma$. In this setting, the following identity holds:

$$
\Delta_{g_0} u + 2K_g e^{u} = 2 K_{g_0},
$$ 

where $\Delta_{g_0}$ denotes the Laplace-Beltrami operator on $\Sigma$, and $K_g$ and $K_{g_0}$ are the Gaussian curvatures of $\Sigma$ corresponding respectively to $g$ and $g_0$. When $\Sigma = \Sm^2$ and $g_0$ is the standard spherical metric, then a function $K \in C^\infty(\Sm^2)$ is the Gaussian curvature of a metric which is pointwise conformally equivalent to $g_0$ if and only if there exists $u\in C^\infty(\Sm^2)$ satisfying 

\begin{equation}\label{EqS2}   
\Delta_{g_0}u + 2K e^u=2 \quad \mbox{ on } \quad (\Sm^2,g_0).
\end{equation}

\medskip
\noindent
It is possible to give a variational formulation of \eqref{EqS2} by considering the functional 

$$ \displaystyle \Phi(u) = \frac{1}{2} \int_{\Sm^2} |\nabla_{g_0} u|^2_{g_0} \, dv_{g_0} + 2\int_{\Sm^2} u  \, dv_{g_0}  - 8 \pi \ln \left( \frac{1}{4\pi}\int_{\Sm^2} K e^u\, dv_{g_0} \right),$$ 

\medskip
\noindent
which is well defined in $$ \mathcal{A}=: \left\{ u \in W^{1,2}(\mathbb{S}^2,g_0) \quad | \quad \int_{\Sm^2} K e^u \, dv_{g_0} >0   \right\},$$

\medskip
\noindent
provided ${\displaystyle \max_{\Sm^2} K >0}$. A function $u \in W^{1,2}(\mathbb{S}^2,g_0)$ is a weak solution to \eqref{EqS2} if and only if it is a critical point of $\Phi$ and 

$$
\int_{\Sm^2} K e^{u}dv_{g_0}  = 4\pi.
$$

\medskip
\noindent
As an analog of Theorem \ref{TheoremaMoserdomini} for the Sobolev space $W^{1,2}(\mathbb{S}^2,g_0)$, Moser in \cite{M} proved the following  

\begin{theorem}[\textbf{Moser-Trudinger inequality on $(\mathbb{S}^2,g_0)$ }]\label{TheoremaMoserS2}\

\medskip
\noindent
Let $(\mathbb{S}^2,g_0)$ be the $two$-dimensional unit sphere endowed with standard metric $g_0$.

\medskip
\noindent
Let us denote with $ \displaystyle \bar{u}:= \frac{1}{4 \pi} \int_{\Sm^2} u \,d \nu_{g_0}$ and with $\bar{\mathcal{H}}:= \{u \in W^{1,2}(\mathbb{S}^2,g_0) \ | \ \bar{u}=0, \ \lVert \nabla_{g_0} u \rVert_2 \leq 1 \}$.\\
Then there exists a constant $S > 1$ such that

\begin{equation}\label{supsuS2}
\sup_{u \in \bar{\mathcal{H}} } \,  \frac{1}{4\pi}\int_{\Sm^2} e^{4\pi u^2 }\,d \nu_{g_0} \le S.
\end{equation}
\end{theorem}

\medskip
\noindent
As a straightforward consequence of previous result and Young's inequality, one can prove the following
weaker inequality:

\begin{proposition}[\textbf{Logarithmic M-T inequality on $(\mathbb{S}^2,g_0)$}]\label{propositionlogMTinequalityonS2}\

\medskip
\noindent
Let $S$ be the constant appearing in \eqref{supsuS2}. For any $\displaystyle u \in W^{1,2}(\mathbb{S}^2,g_0),$ we have

\begin{equation}\label{logMTinequalityS2}
\displaystyle  \ln \left( \frac{1}{4 \pi} \int_{\mathbb{S}^2} e^u \, dv_{g_0} \right) \leq  \frac{1}{16 \pi} \int_{\mathbb{S}^2} | \nabla_{g_0} u |_{g_0}^2 \, dv_{g_0} + \frac{1}{4 \pi} \int_{\mathbb{S}^2}u \, dv_{g_0} + \ln S.
\end{equation}

\end{proposition}

\begin{proof}

Let $u \in W^{1,2}(\mathbb{S}^2,g_0)$, $\lVert \nabla_{g_0} u \rVert_{2} \neq 0$, and let  $\displaystyle b^2=\frac{8 \pi}{\lVert \nabla_{g_0} u \rVert^2_2}$.

\medskip
\noindent
By Young's inequality, we infer that

$$ u \leq |u| \leq \frac{1}{2b^2} + \frac{1}{2}b^2u^2  = \frac{1}{16 \pi} \lVert \nabla_{g_0} u \rVert^2_2 +  4 \pi \frac{u^2}{\lVert \nabla_{g_0} u \rVert^2_2}, $$

\medskip
\noindent
by which

$$ \displaystyle \frac{1}{4 \pi} \int_{\mathbb{S}^2} e^u \, d \nu_{g_0} \leq e^{\frac{1}{16 \pi} \lVert \nabla_{g_0} u \rVert^2_2} \frac{1}{4 \pi} \int_{\mathbb{S}^2} e^{ 4 \pi \left( \frac{u}{\lVert \nabla_{g_0} u \rVert_2} \right)^2 } \, d \nu_{g_0}. $$

\medskip
\noindent
Applying previous inequality  to $u - \bar{u}$ and considering also \eqref{supsuS2}, we obtain

$$ \displaystyle \frac{1}{4 \pi} \int_{\mathbb{S}^2} e^{u- \bar{u}} \, d \nu_{g_0} \leq e^{\frac{1}{16 \pi} \lVert \nabla_{g_0} u \rVert^2_2} \frac{1}{4 \pi} \int_{\mathbb{S}^2} e^{ 4 \pi \left( \frac{u-\bar{u}}{\lVert \nabla_{g_0} u \rVert_2} \right)^2} \, d \nu_{g_0} \leq e^{\frac{1}{16 \pi} \lVert \nabla_{g_0} u \rVert^2_2} \, S . $$

\medskip
\noindent
Passing to logarithm, we obtain

$$ \displaystyle \ln \left( \frac{1}{4 \pi} \int_{\mathbb{S}^2} e^u \, d\nu_{g_0} \right) \leq  \frac{1}{16 \pi} \int_{\mathbb{S}^2} | \nabla_{g_0} u |_{g_0}^2 \, d\nu_{g_0} + \frac{1}{4 \pi} \int_{\mathbb{S}^2}u \, d\nu_{g_0} + \ln S.
$$

\medskip
\noindent
If $u \equiv c$ for any $c \in \mathbb{R},$ inequality \eqref{logMTinequalityS2} still holds.

\end{proof}

\medskip
\noindent
In particular,  \eqref{logMTinequalityS2} shows that $\Phi$ is bounded from below in $\mathcal{A}$. While this is not sufficient to apply standard minimization techniques, the existence of critical points for $\Phi$ can be obtained via improved versions of \eqref{logMTinequalityS2} under suitable assumptions on $K$. For example, in \cite{M2} Moser proves that if $u\in W^{1,2}(\mathbb{S}^2,g_0)$ is even, then \eqref{logMTinequalityS2} can be replaced by the stronger inequality:  

\begin{equation}\label{ConsMEven}
\ln \left(\frac{1}{4\pi}\int_{\Sm^2} e^u \, d\nu_{g_0} \right) \leq  \frac{1}{32 \pi} \int_{\Sm^2} | \nabla_{g_0} u |_{g_0}^2 \, d\nu_{g_0} + \frac{1}{4 \pi} \int_{\Sm^2}u \, d\nu_{g_0} + C ,
\end{equation}

\medskip
\noindent
for some $C\in \R$. In particular, if $K\in C^\infty(\Sm^2)$ is even, the restriction of $\Phi$ to the space of even functions in $W^{1,2}(\mathbb{S}^2,g_0)$ with zero average is coercive, and one can find a critical point of $\Phi$ by means of Palais symmetric criticality principle. This leads to one  of first general existence results for the Nirenberg problem: if $K\in C^\infty(\Sm^2)$ is an even function and $\displaystyle{\max_{\Sm^2} K>0}$, then there exists a solution of \eqref{EqS2}. We refer to {\cite{Berger,KW,CY1,CY2,StruweCurv} and references therein for further results concerning Nirenberg's problem on $\Sm^2$.}

\medskip
\noindent
Another interesting aspect of \eqref{logMTinequalityS2} is that while the coefficient  $\frac{1}{16\pi}$ is sharp, the constant $\ln S$ is not. Indeed, Onofri in \cite{O} proved that

$$
\inf_{W^{1,2}(\mathbb{S}^2,g_0)} \Phi  = 0 
$$

\medskip
\noindent
and he classified all minimum points of $\Phi$. Namely, we have the so called \textit{logarithmic Onofri-Moser-Trudinger inequality} or \textit{Onofri inequality} 

\begin{equation}\label{bestOnofri}
\ln \left(\frac{1}{4\pi}\int_{\Sm^2} e^u \, dv_{g_0} \right) \leq  \frac{1}{16 \pi} \int_{\Sm^2} | \nabla_{g_0} u |_{g_0}^2 \, dv_{g_0} + \frac{1}{4 \pi} \int_{\Sm^2}u \, dv_{g_0}, 
\end{equation}

\medskip
\noindent
for any $\, u \in W^{1,2}(\mathbb{S}^2,g_0)$, and equality  holds if and only if $u  = \ln |J_\psi| + C$, where $\psi:\Sm^2 \rightarrow \Sm^2$ is a conformal diffeomorphism of $(\Sm^2,g_0)$, $J_\psi$ denotes the Jacobian of $\psi$, and $C\in \R$. Inequality \eqref{bestOnofri} was proved in 1982 by Onofri in \cite{O}  using conformal invariance and an earlier result by Aubin  \cite{Aubin}. Indeed, inequality \eqref{bestOnofri} can be obtained via several methods, including various limiting procedures, mass transportation methods in the radial case, and rigidity methods associated with a related nonlinear flow (we refer to \cite{DEJ,DET} and references therein for complete proof and details). Onofri's inequality \eqref{bestOnofri} plays a role in spectral analysis for the Laplace-Beltrami operator. Indeed, thanks  to Polyakov's formula (see \cite{Poly1,Poly2,OPS,OPS2}), we have 

$$
\frac{\det (-\Delta_g)}{\det(-\Delta_{g_0})} =e^{-\frac{1}{24}\Phi(u)},
$$

\medskip
\noindent
for any metric $g = e^u g_0$ with $u \in C^\infty(\Sm^2)$. In particular,  Onofri's inequality \eqref{bestOnofri} implies that, in the conformal class of $g_0$, the determinant of the Laplace-Beltrami operator is maximized by  $g_0$. We mention also important applications; for instance, in chemotaxis, it appears in the study of the Keller-Segel model (see \cite{CALVCORR,Dolbeault-Perthame, GAJZAK}).

\medskip
\noindent
Now, let us recall that the two-dimensional Euclidean Onofri inequality states that

\begin{equation}\label{EuclideanOnofri2D}
\displaystyle  \ln \left(\int_{\mathbb{R}^2} e^u \, d \mu_2 \right) \leq  \frac{1}{16 \pi} \int_{\mathbb{R}^2} | \nabla u |^2 \, dx + \int_{\mathbb{R}^2} u \, d \mu_2
\end{equation}

\medskip
\noindent
for any $ u \in  W_{\mu_2}(\mathbb{R}^2):= \{ u \in L^1(\mathbb{R}^2,d\mu_2) \, : \, |\nabla u| \in L^2(\mathbb{R}^2,dx) \}$, where

$$ \displaystyle \mu_2(x):= \frac{1}{ \pi \left( 1+ |x|^2 \right)^2} \qquad \text{and} \quad d \mu_2:=\mu_2(x) \, dx. $$

\medskip
\noindent
In particular, we have the following:

\begin{proposition}\label{equivS2OE2}
Inequalities  \eqref{bestOnofri} and \eqref{EuclideanOnofri2D} are equivalent.
\end{proposition}

\begin{proof} Let $\pi_{\mathcal{N}} : \mathbb{S}^2 \setminus \lbrace {\mathcal{N}} \rbrace \to \mathbb{R}^2$ be the stereographic projection from the north pole ${\mathcal{N}}$, and let us denote by $\pi_{\mathcal{N}}^{-1}:\mathbb{R}^2 \to \mathbb{S}^2 \setminus \lbrace {\mathcal{N}} \rbrace$ its inverse.

\medskip
\noindent
The standard metric $g_0$ on $\mathbb{S}^2$ is represented by the matrix
 
$$
g_0(x)=\begin{pmatrix}
\lambda(x) & 0\\
0 & \lambda(x)
\end{pmatrix}, \qquad \text{where} \quad \lambda(x)= \frac{4}{(1+|x|^2)^2}.$$

\medskip
\noindent
If $v: \mathbb{S}^2 \to \mathbb{R}$ is a real valued function on $\mathbb{S}^2$, denoting with $\tilde{v}:= v \circ \pi_{\mathcal{N}}^{-1}: \mathbb{R}^2 \to \mathbb{R}$, by definition we have

$$ \displaystyle \int_{\mathbb{S}^2} v  \ d \nu_{g_0} := \int_{\mathbb{R}^2} \tilde{v}(x) \sqrt{\det g_0(x)}  \, dx = \int_{\mathbb{R}^2} \tilde{v}(x) \lambda(x) \, dx.$$

\medskip
\noindent
Therefore we have

$$ \frac{1}{4 \pi} \int_{\mathbb{S}^2} e^u \, d\nu_{g_0} = \frac{1}{4 \pi} \int_{\mathbb{R}^2} e^{\tilde{u}(x)} \lambda(x)  \, dx
=\int_{\mathbb{R}^2} e^{\tilde{u}(x)} \, d \mu_2 $$

\medskip
\noindent
and

$$\frac{1}{4 \pi} \int_{\mathbb{S}^2} u \, d\nu_{g_0} =  \int_{\mathbb{R}^2} \tilde{u}(x)  \lambda(x) \, dx  =\int_{\mathbb{R}^2} \tilde{u} \, d \mu_2.$$

\medskip
\noindent
Now we recall that if $p \in \mathbb{S}^2 \setminus{ \{ \mathcal{N} \} }$, the gradient of $v$ in $p$ with respect to the metric $g_0$ is given by 

$$ \displaystyle \nabla_{g_0} v (p) = \frac{1}{\lambda(x)} \nabla \tilde{v}(x) \bigg|_{x=\pi_{\mathcal{N}}(p)}. $$

\medskip
\noindent
In particular, we get

\begin{align*}
\displaystyle | \nabla_{g_0} u(p) |_{g_0}^2 = \nabla_{g_0} u(p) 
\begin{pmatrix}
\lambda(x) & 0\\
0 & \lambda(x)
\end{pmatrix}  
(\nabla_{g_0} u(p) )^T
= \frac{1}{\lambda(x)} \left| \nabla \tilde{v}(x) \right|^2 \bigg|_{x=\pi_{\mathcal{N}}(p)},
\end{align*}

\medskip
\noindent
whence

$$ \displaystyle  \int_{\mathbb{S}^2} | \nabla_{g_0} v |_{g_0}^2 \, d\nu_{g_0} = \int_{\mathbb{R}^2} | \nabla \tilde{v} |^2 \, dx. $$

\medskip
\noindent
Therefore

$$ \displaystyle \frac{1}{16 \pi} \int_{\mathbb{S}^2} | \nabla_{g_0} u |_{g_0}^2 \, d\nu_{g_0} = \frac{1}{16 \pi} \int_{\mathbb{R}^2} | \nabla \tilde{u} |^2 \, dx. $$

\medskip
\noindent
It remains to show that the space 

$$W^{1,2}(\mathbb{S}^2,g_0):=\left\lbrace u : \mathbb{S}^2 \to \mathbb{R} \, : \,\int_{\mathbb{S}^2} |u|^2 \, d\nu_{g_0} < + \infty \quad \text{ and } \quad  \int_{\mathbb{S}^2} | \nabla_{g_0} u |_{g_0}^2 \, d\nu_{g_0} < + \infty \right\rbrace$$

\medskip
\noindent
coincides with the space

$$ \displaystyle W:= \left\lbrace u : \mathbb{S}^2 \to \mathbb{R} \, : \,\int_{\mathbb{S}^2} |u| \, d\nu_{g_0} < + \infty \quad \text{ and }  \quad \int_{\mathbb{S}^2} | \nabla_{g_0} u |_{g_0}^2 \, d\nu_{g_0} < + \infty \right\rbrace.$$

\medskip
\noindent
It is clear that $W^{1,2}(\mathbb{S}^2,g_0) \subseteq W$. Conversely, if $u \in W$ then $u \in W^{1,1}(\mathbb{S}^2,g_0)$, and so by Sobolev embeddings on compact manifolds we get $u \in L^2(\mathbb{S}^2,g_0)$, i.e. $u \in W^{1,2}(\mathbb{S}^2,g_0)$.

\end{proof}

\bigskip

\section*{Appendix D}
\addcontentsline{toc}{section}{Appendix D}
\label{app:D}

\bigskip

\begin{lemma}\label{computationNlaplacian}

Let $ N \in \mathbb{N}$,  $N \geq 2$. Let $V_N:= \omega_{N-1}/N,$ where $\omega_{N-1}:=|\mathbb{S}^{N-1}|,$ and let us define 

$$ \displaystyle \mu_N(x):= \frac{1}{V_N\left( 1+ |x|^{\frac{N}{N-1}} \right)^N} \qquad \text{and} \qquad v_N(x):= \ln \mu_N(x).$$

\medskip
\noindent
Then

$$ \displaystyle   \Delta_N v_N(x) =- N^N \left( \frac{N}{N-1} \right)^{N-1} V_N \mu_N(x).$$

\end{lemma}

\begin{proof}

\medskip
\noindent
We recall that the $N$-Laplacian operator is defined as

$$ \displaystyle \Delta_N \, u := \text{div} \left( \left| \nabla u   \right|^{N-2} \nabla u     \right).$$ 

\medskip
\noindent
Observe that

\begin{equation*}
	\nabla v_{N}(x)= -  \frac{N^2}{N-1}  \frac{ \left| x   \right|^{\frac{1}{N-1} }}{ 1+|x|^{\frac{N}{N-1}}}   \frac{x}{\left| x \right| },
\qquad \text{and} \qquad 
	\left| \nabla v_{N}(x)  \, \right| =\frac{N^2}{N-1}  \frac{ \left| x   \right|^{\frac{1}{N-1} }}{ 1+|x|^{\frac{N}{N-1}}}.
\end{equation*}

\medskip
\noindent
Hence

\begin{align*}
\displaystyle \left|   \nabla v_{N}(x) \right|^{N-2}   \nabla v_{N}(x) = -  \left(    \frac{N^2}{N-1} \right)^{N-1} \frac{ x}{ \left(  1+|x|^{\frac{N}{N-1}}  \right)^{N-1}   } .   
\end{align*}

\medskip
\noindent
Now, observe that

\begin{align*}
\displaystyle
\frac{\partial}{ \partial x_i} \left(  \frac{ x_i}{ \left(  1+|x|^{\frac{N}{N-1}}  \right)^{N-1}   }         \right) 
& = \frac{ \displaystyle \left(  1+|x|^{\frac{N}{N-1}}  \right)^{N-1}  - N \left(  1+|x|^{\frac{N}{N-1}}  \right)^{N-2} |x|^{\frac{1}{N-1}} \,  \frac{x_i^2}{|x|} }{\left(  1+|x|^{\frac{N}{N-1}}  \right)^{2(N-1)} } \\
& = \frac{1}{\left(  1+|x|^{\frac{N}{N-1}}  \right)^{N-1}  } - N \frac{ \displaystyle |x|^{\frac{1}{N-1}} \,  \frac{x_i^2}{|x|} }{\left(  1+|x|^{\frac{N}{N-1}}  \right)^{N} } .
\end{align*}

\medskip
\noindent
Therefore, we have 

\begin{align*}
\displaystyle 
\Delta_N v_N(x)
& =  -  \left(    \frac{N^2}{N-1} \right)^{N-1} \sum_{i=1}^{N}  \frac{\partial}{ \partial x_i} \left(  \frac{ x_i}{ \left(  1+|x|^{\frac{N}{N-1}}  \right)^{N-1}   }         \right) \\
& =  -  \left(    \frac{N^2}{N-1} \right)^{N-1}  N \left[  \frac{1}{\left(  1+|x|^{\frac{N}{N-1}}  \right)^{N-1}  } -    \frac{ \displaystyle |x|^{\frac{N}{N-1}}  }{\left(  1+|x|^{\frac{N}{N-1}}  \right)^{N}} \right] \\
& = -N^N \left(    \frac{N}{N-1} \right)^{N-1}   V_N \mu_N(x).
\end{align*}

\end{proof}

\newpage

\thispagestyle{empty}

\newgeometry{paperwidth=10cm, paperheight=25cm, textheight=237mm, textwidth=98mm, top=40mm, left=28mm, right=28mm}

\section*{\center{\LARGE{Acknowledgments}}}

\medskip
\begin{center}
\emph{``Tu ne quaesieris, scire nefas...}''
\end{center}
\begin{center}
\emph{``Do not ask, it is not permitted to know...}''
\end{center}
\begin{center}
\emph{``Non chiedere, non è lecito sapere...}''
\end{center}

\medskip
\noindent
With these famous words begins the celebrated `Carpe Diem' ode by the Latin poet Quintus Horatius Flaccus. I chose these words because I believe that experiences of such human intensity should remain private. For this reason, an Acknowledgments section is absent from both my Bachelor’s and Master’s theses. Nevertheless, I feel it is right for a sign of that intensity to be visible, even if it only shows one percent of what lies beneath.

First and foremost, I must thank my supervisor, Professor Silvia Cingolani, for helping me grow. Under her guidance, my PhD in Mathematics was far more than just a PhD in Mathematics. She introduced me to a new, diverse, and complex branch of Mathematics, supporting me every step of the way. She guided me not as a vessel for internal knowledge, but as one that must know how to reflect, through its external form, the essence of what it contains. Even during my most challenging moments, her wide-ranging experience acted as a guiding light, providing a steady path for me to follow. I must thank again my supervisor, Professor Silvia Cingolani, because She taught me not just the substance of what to think, but the wisdom of how and when to do so.

Secondly, I would like to express my sincere gratitude to Professor David Arcoya Álvarez for his warm hospitality during my two visits to the University of Granada. He provided me with the opportunity not only to broaden my research scope, but also to present the results achieved. He encouraged me not only to understand the deep mathematical and historical roots of today’s research, but also to improve how I communicate and present my results.
Together with my supervisor, we initiated a collaboration that has culminated in a paper for which we are confidently awaiting the peer review results. 
The insightful discussions regarding the understanding of the topics of the paper have been of fundamental importance in my development as a researcher.

I would like to express my gratitude to Professor Gabriele Mancini. Our collaboration, which began together with my supervisor immediately after I completed a doctoral exam held by him, has been fundamental not only from a scientific perspective but also on a personal level.

Special thanks go to Professor Giusi Vannella. The collaboration we started together with my supervisor, provided invaluable scientific results and human support.

I would also like to thank my internal PhD commission, composed of Professor Giusi Vaira and Professor Felice Iavernaro, who followed my progress throughout these years, always providing wise advice and appreciation for the results achieved year after year.

I would also like to extend my deepest thanks to the referees of this thesis, Professor Alberto Farina and Professor José Carmona Tapia, for the time and effort they dedicated to evaluating my doctoral work.

I would like to thank both the former coordinator, Professor Francesca Mazzia, for her guidance, availability, and especially her appreciation throughout the three years of my PhD, and the current coordinator, Professor Giovanna Castellano, for her guidance and support during the final steps of this thesis.

I must also express my indebtedness to Professor Roberto Bellotti, the PI of the PNRR grant that funded this research. Without him, the words you are now reading would simply not exist.

\newpage

\thispagestyle{empty}

The list of people I should thank is—if I may be allowed to use an expression that is both mathematically and philosophically incorrect— `countably infinite'. I cannot but continue this list by thanking Professor Antonio Azzollini, my supervisor for both my Bachelor’s and Master’s degrees, who continued to take an interest in my progress and granted me the opportunity to hold a seminar at the University of Basilicata in Potenza. I must also thank Professor Benedetta Pellacci, who gave me the opportunity to hold a seminar at the Department of Mathematics of the University of Campania `Luigi Vanvitelli'. I shall not forget the taste of a pre-seminar pizza for lunch, but I certainly will never forget the beauty of their mathematical library. Special thanks are also due to Professor Berardino Sciunzi for the opportunity to present a seminar at a school held at the University of Cosenza. I would like to thank also Professor Jaroslaw Mederski, Professor Bartosz Bieganowski and Dr. Jacopo Schino for granting me the opportunity to present my work on several occasions, both in Warsaw and in Będlewo.

I would also like to extend my gratitude to everyone in the academic community I met during my time abroad in Spain. I owe special thanks to Rubén Fiñana Aránega for our conversations, which spanned from mathematics and linguistics to diverse subjects of a different nature.

I am also indebted to Dr. Marco Gallo for his support in various situations. I especially appreciate him giving a seminar here in Bari, where he even featured me as a character in one of his slides!

I would like to thank all my PhD colleagues at the Department of Mathematics. In particular, my heartfelt gratitude goes to Antonio Lagioia and Alessandro Cannone. Antonio welcomed me with his genuine kindness when I was still a stranger to everyone—more specifically, when I had not yet even figured out where the Department entrance was. Alessandro, on the other hand, has always been a source of inspiration to me. I have always admired his incredible balance; indeed, he even managed to convince me to drink (water) during lunch, besides teaching me many other important things on a human level. I must also thank Grazia Gargano (though I am not permitted here to use the name I have always called her), who had a significant impact on me at an institutional level. To Antonio, Alessandro, and Grazia, I have successfully taught the true meaning of Martial’s expression: \emph{`Parcere personis, dicere de vitiis}'.

I would also like to thank Giuseppe Rago; between one smile and the next, we enjoyed discussions ranging from the Greek god Hermès to D’Annunzio’s aestheticism.

I am profoundly grateful to Professor Sandra Lucente for perceiving and describing me, since our first meeting, as a man of the eighteenth and nineteenth centuries.

If I were to name each of you individually, I would have to write a second thesis just for that purpose. Therefore, I `limit' myself to saying that I truly want to thank the entire Department of Mathematics of the University of Bari Aldo Moro. I hope I managed to bring a smile to your face during our $\pi$-Day celebrations.

Finally, I wish to thank Professor Antonio Masiello. The words He told me at the beginning of my PhD— \emph{`Nella Matematica ci vuole fantasia}' translated in \emph{`Mathematics requires imagination}' — opened a fascinating door in my mind.\\

\noindent
As for my friends and my family, I am sorry \emph{``Tu ne quaesieris, scire nefas...}''. However, I wish to leave only one statement here. Even today, especially when I think of my sister Lucia, I have serious doubts about whether my cradle was swapped with someone else's at birth.

\bigskip
\bigskip

\begin{center}
 \textasciitilde \emph{ Non omnis moriar} \textasciitilde 
\end{center}

\newpage

\cleardoublepage
\phantomsection
\addcontentsline{toc}{chapter}{Bibliography}
\bibliographystyle{amsplain}

\begin{thebibliography}{999}


\bibitem{Adams} D. R. Adams, \textit{A sharp inequality of J. Moser for higher order derivatives}, Ann. of Math. \textbf{128} (1988), no. 2, 385--398.



\bibitem{AFTPAC} A. Aftalion and F. Pacella, \textit{Morse index and uniqueness for positive solutions of radial $p$-Laplace equations}, Trans. Amer. Math. Soc. \textbf{356} (2004), no. 11, 4255--4272.

\bibitem{ABG} M. Agueh, S. Boroushaki, and N. Ghoussoub, \textit{A dual Moser–Onofri inequality and its extensions to higher dimensional spheres}, Ann. Fac. Sci. Toulouse Math. (6) \textbf{26} (2017), no. 2, 217--233.

\bibitem{ALP} S. Ahmad, A. C. Lazer, and J. L. Paul, \textit{Elementary critical point theory and perturbations of elliptic boundary value problems at resonance}, Indiana Univ. Math. J. \textbf{25} (1976), no. 10, 933--944.

\bibitem{AH} W. Allegretto and H. Y. Xi, \textit{A Picone's identity for the p-Laplacian and applications}, Nonlinear Anal. \textbf{32} (1998), no. 7, 819--830.

\bibitem{ad} S. Almi and M. Degiovanni, \textit{On degree theory for quasilinear elliptic equations with natural growth conditions}, Recent Trends in Nonlinear Partial Differential Equations II: Stationary Problems, Contemp. Math. \textbf{595} (2013), 1--20.

\bibitem{amann_zehnder} H. Amann and E. Zehnder, \textit{Nontrivial solutions for a class of nonresonance problems and applications to nonlinear differential equations}, Ann. Sc. Norm. Super. Pisa Cl. Sci. (4) \textbf{7} (1980), no. 4, 539--603.

\bibitem{AMBROSETTIMANCINI} A. Ambrosetti and G. Mancini, \textit{Existence and multiplicity results for nonlinear elliptic problems with linear part at resonance: the case of the simple eigenvalue}, J. Differential Equations \textbf{28} (1978), no. 2, 220--245.

\bibitem{AL} A. Anane, \textit{Simplicité et isolation de la premiere valeur propre du p-laplacien avec poids}, C. R. Acad. Sci. Paris Sér. I Math. \textbf{305} (1987), no. 16, 725--728.

\bibitem{ANGOS} A. Anane and J.-P. Gossez, \textit{Strongly nonlinear elliptic problems near resonance: a variational approach}, Comm. Partial Differential Equations \textbf{15} (1990), no. 8, 1141--1159.

\bibitem{ANANETSOULI} A. Anane and N. Tsouli, \textit{On the second eigenvalue of the $p$-Laplacian}, in: Nonlinear partial differential equations (Fès, 1994), Pitman Res. Notes Math. Ser., 343, Longman, Harlow, 1996, 1--9.

\bibitem{ACCFM} C. A. Antonini, A. Cianchi, G. Ciraolo, A. Farina, and V. Maz'ya, \textit{Global second-order estimates in anisotropic elliptic problems}, Proc. Lond. Math. Soc. (3) \textbf{130} (2025), no. 1, e70034.

\bibitem{ACF} C. A. Antonini, G. Ciraolo, and A. Farina, \textit{Interior regularity results for inhomogeneous anisotropic quasilinear equations}, Math. Ann. \textbf{387} (2023), no. 3, 1745--1776.

\bibitem{ABC} D. Arcoya, N. Borgia and S. Cingolani, \textit{On the de Thélin eigenvalue problem and Landesman-Lazer conditions for quasilinear systems}, J. Differential Equations. {\bf 474} (2026), pp. 114513.
 \href{https://doi.org/10.1016/j.jde.2026.114513}{https://doi.org/10.1016/j.jde.2026.114513}


\bibitem{AG} D. Arcoya and J. L. Gámez, \textit{Bifurcation theory and related problems: anti-maximum principle and resonance}, Comm. Partial Differential Equations \textbf{26} (2001), no. 9-10, 1879--1911.

\bibitem{AO} D. Arcoya and L. Orsina, \textit{Landesman-Lazer conditions and quasilinear elliptic equations}, Nonlinear Anal. \textbf{28} (1997), no. 10, 1623--1632.

\bibitem{ASSCING} L. Asselle, S. Cingolani, and M. Starostka, \textit{Morse homology for a class of elliptic partial differential equations}, Commun. Contemp. Math. (2025). \href{https://doi.org/10.1142/S0219199726500082}{ https://doi.org/10.1142/S0219199726500082}


\bibitem{AM} G. Astarita and G. Marucci, \textit{Principles of Non-Newtonian Fluid Hydromechanics}, McGraw Hill, New York, 1974.

\bibitem{Aubin} T. Aubin, \textit{Meilleures constantes dans le théorème d'inclusion de Sobolev et un théorème de Fredholm non linéaire pour la transformation conforme de la courbure scalaire}, J. Funct. Anal. \textbf{32} (1979), no. 2, 148--174.

\bibitem{Bart} D. Bartolucci, \textit{Non-degeneracy, Mean Field Equations and the Onsager Theory of 2D Turbulence}, Arch. Ration. Mech. Anal. \textbf{230} (2018), 397--426.

\bibitem{BarMal} D. Bartolucci and A. Malchiodi, \textit{Mean field equations and domains of first kind}, Rev. Mat. Iberoam. \textbf{38} (2022), no. 4, 1067--1086.

\bibitem{bartsch_dancer2009} T. Bartsch and N. Dancer, \textit{Poincaré-Hopf type formulas on convex sets of Banach spaces}, Topol. Methods Nonlinear Anal. \textbf{34} (2009), 213--229.




\bibitem{benci1991} V. Benci, \textit{A new approach to the Morse-Conley theory and some applications}, Ann. Mat. Pura Appl. (4) \textbf{158} (1991), 231--305.

\bibitem{Berger} M. S. Berger, \textit{Riemannian structures of prescribed Gaussian curvature for compact 2-manifolds}, J. Differential Geom. \textbf{5} (1971), 325--332.


\bibitem{BDO} L. Boccardo, A. Dall'Aglio, and L. Orsina, \textit{Existence and regularity results for some elliptic equations with degenerate coercivity}, Atti Sem. Mat. Fis. Univ. Modena \textbf{46} (1998), 51--82.

\bibitem{boccardodefiguerido} L. Boccardo and D. G. de Figueiredo, \textit{Some remarks on a system of quasilinear elliptic equations}, NoDEA Nonlinear Differential Equations Appl. \textbf{9} (2002), no. 3, 309--323.

\bibitem{BDK} L. Boccardo, P. Drábek, and M. Kučera, \textit{Landesman--Lazer conditions for strongly nonlinear boundary value problems}, Comment. Math. Univ. Carolin. \textbf{30} (1989), no. 3, 411--427.

\bibitem{BDMS} V. Bögelein, F. Duzaar, P. Marcellini, and C. Scheven, \textit{Boundary regularity for elliptic systems with p, q-growth}, J. Math. Pures Appl. (9) \textbf{159} (2022), 250--293.



\bibitem{BON} A. Bonnet, \textit{A deformation lemma on a $C^1$ manifold}, Manuscripta Math. \textbf{81} (1993), no. 1, 339--359.

\bibitem{BCM} N. Borgia, S. Cingolani, and G. Mancini, \textit{New Analytical and Geometrical Aspects on Trudinger-Moser Type Inequality in 2D}, in: Singularities, Asymptotics, and Limiting Models, Springer INdAM Ser., vol. 64, Springer, Singapore, 2025, 235--254. \href{https://doi.org/10.1007/978-981-96-3584-9\_9}{https://doi.org/10.1007/978-981-96-3584-9\_9}

\bibitem{BCM2} N. Borgia, S. Cingolani, and G. Mancini, \textit{On the equivalence between an Onofri-type inequality by Del Pino–Dolbeault and the sharp logarithmic Moser–Trudinger inequality}, Calc. Var. Partial Differential Equations \textbf{64} (2025), no. 3, Art. 89, 22 pp.

\bibitem{BCV} N. Borgia, S. Cingolani, and G. Vannella, \textit{Nontrivial solutions for resonance quasilinear elliptic systems}, Adv. Nonlinear Anal. \textbf{13} (2024), no. 1, Art. 20230155.

\bibitem{BCV2} N. Borgia, S. Cingolani, and G. Vannella, \textit{A Poincaré–Hopf formula for functional associated to quasilinear elliptic system}, Nonlinear Anal. Real World Appl. \textbf{87} (2026), Art. 104443, 12 pp.

\bibitem{BCV3} N. Borgia, S. Cingolani and G. Vannella, \textit{Uniform $L^{\infty}$-boundedness for solutions of anisotropic quasilinear systems}, Adv. Nonlinear Stud. {\bf 26} (2026),  no. 1, pp. 127-143.\\
 \href{https://doi.org/10.1515/ans-2023-0203}{https://doi.org/10.1515/ans-2023-0203}

\bibitem{BM} Y. Bozhkov and E. Mitidieri, \textit{Existence of multiple solutions for quasilinear systems via fibering method}, J. Differential Equations \textbf{190} (2003), no. 1, 239--267.


\bibitem{BREZIS} Brezis, H. (2011). Functional Analysis, Sobolev Spaces and Partial Differential Equations. Springer New York.

\bibitem{browder1983} F. E. Browder, \textit{Fixed point theory and nonlinear problems}, Bull. Amer. Math. Soc. (N.S.) \textbf{9} (1983), no. 1, 1--39.



\bibitem{CLMP} E. Caglioti, P.-L. Lions, C. Marchioro, and M. Pulvirenti, \textit{A special class of stationary flows for two-dimensional Euler equations: a statistical mechanics description. Part II}, Comm. Math. Phys. \textbf{174} (1995), no. 2, 229--260.


\bibitem{CALVCORR} V. Calvez and L. Corrias, \textit{The parabolic-parabolic Keller-Segel model in $\mathbb{R}^2$}, Commun. Math. Sci. \textbf{6} (2008), no. 2, 417--447.

\bibitem{candela} A. M. Candela, E. Medeiros, G. Palmieri, and K. Perera, \textit{Weak solutions of quasilinear elliptic systems via the cohomological index}, Topol. Methods Nonlinear Anal. \textbf{36} (2010), no. 1, 1--18.

\bibitem{CSS} A. Candela, A. Salvatore, and C. Sportelli, \textit{Existence and multiplicity results for a class of coupled quasilinear elliptic systems of gradient type}, Adv. Nonlinear Stud. \textbf{21} (2021), no. 2, 461--488.



\bibitem{CarlesonChang} L. Carleson and A. Chang, \textit{On the existence of an extremal function for an inequality of J. Moser}, Bull. Sci. Math. (2) \textbf{110} (1986), no. 2, 113--127.

\bibitem{CCMV} J. Carmona, S. Cingolani, P. Martinez-Aparicio, and G. Vannella, \textit{Regularity and Morse index of the solutions to critical quasilinear elliptic systems}, Comm. Partial Differential Equations \textbf{38} (2013), no. 10, 1675--1711.


\bibitem{chang1981} K.-C. Chang, \textit{Variational methods for nondifferentiable functionals and their applications to partial differential equations}, J. Math. Anal. Appl. \textbf{80} (1981), no. 1, 102--129.

\bibitem{chang1981-cpam} K.-C. Chang, \textit{Solutions of asymptotically linear operator equations via Morse theory}, Comm. Pure Appl. Math. \textbf{34} (1981), no. 5, 693--712.

\bibitem{chang} K.-C. Chang, \textit{Infinite-dimensional Morse theory and multiple solution problems}, Progr. Nonlinear Differential Equations Appl., vol. 6, Birkhäuser Boston, Inc., Boston, MA, 1993.

\bibitem{CCL} S.-Y. A. Chang, C.-C. Chen, and C.-S. Lin, \textit{Extremal functions for a mean field equation in two dimension}, in: Lectures on Partial Differential Equations in Honor of Louis Nirenberg's 75th Birthday, Ser. Anal., vol. 2, Int. Press, Somerville, MA, 2003, 33--93.

\bibitem{chang-wang} K.-C. Chang and Z.-Q. Wang, \textit{Multiple non semi-trivial solutions for elliptic systems}, Adv. Nonlinear Stud. \textbf{12} (2012), no. 2, 363--381.

\bibitem{CY1} S.-Y. A. Chang and P. C. Yang, \textit{Prescribing Gaussian curvature on $S^2$}, Acta Math. \textbf{159} (1987), no. 1-2, 215--259.

\bibitem{CY2} S.-Y. A. Chang and P. C. Yang, \textit{Conformal deformation of metrics on $S^2$}, J. Differential Geom. \textbf{27} (2008), no. 2, 259--296.



\bibitem{CL1} C.-C. Chen and C.-S. Lin, \textit{Sharp Estimates for Solutions of Multi-Bubbles in Compact Riemann Surfaces}, Comm. Pure Appl. Math. \textbf{55} (2002), no. 6, 728--771.

\bibitem{CL2} C.-C. Chen and C.-S. Lin, \textit{Topological Degree for a Mean Field Equation on Riemann Surfaces}, Comm. Pure Appl. Math. \textbf{56} (2003), no. 12, 1667--1727.

\bibitem{CY} L. Chen and Y. Yang, \textit{Talenti comparison results and rigidity results for anisotropic p-Laplacian operator with Robin boundary conditions}, Adv. Nonlinear Stud. (2025). doi:10.1515/ans-2023-0194.

\bibitem{CM} A. Cianchi and V. G. Maz'ya, \textit{Optimal second-order regularity for the p-Laplace system}, J. Math. Pures Appl. (9) \textbf{132} (2019), 41--78.


\bibitem{CD} S. Cingolani and M. Degiovanni, \textit{Nontrivial solutions for p-Laplace equations with right hand side having p-linear growth at infinity}, Comm. Partial Differential Equations \textbf{30} (2005), no. 7-9, 1191--1203.

\bibitem{CD1} S. Cingolani and M. Degiovanni, \textit{On the Poincaré-Hopf Theorem for functionals defined on Banach spaces}, Adv. Nonlinear Stud. \textbf{9} (2009), no. 4, 679--699.

\bibitem{CDS} S. Cingolani, M. Degiovanni, and B. Sciunzi, \textit{Weighted Sobolev spaces and Morse estimates for quasilinear elliptic equations}, J. Funct. Anal. \textbf{286} (2024), no. 8, Art. 110346, 45 pp.

\bibitem{CDV} S. Cingolani, M. Degiovanni, and G. Vannella, \textit{Amann-Zehnder type results for p-Laplace problems}, Ann. Mat. Pura Appl. (4) \textbf{197} (2018), no. 2, 605--640.





\bibitem{CV} S. Cingolani and G. Vannella, \textit{Critical groups computations on a class of Sobolev Banach spaces via Morse index}, Ann. Inst. H. Poincaré C Anal. Non Linéaire \textbf{20} (2003), no. 2, 271--292.

\bibitem{CV2} S. Cingolani and G. Vannella, \textit{Morse index and critical groups for p-Laplace equations with critical exponents}, Mediterr. J. Math. \textbf{3} (2006), no. 3-4, 495--512.

\bibitem{CV3} S. Cingolani and G. Vannella, \textit{Marino-Prodi perturbation type results and Morse indices of minimax critical points for a class of functionals in Banach spaces}, Ann. Mat. Pura Appl. (4) \textbf{186} (2007), no. 1, 157--185.

\bibitem{CV4} S. Cingolani and G. Vannella, \textit{Multiple positive solutions for a critical quasilinear equation via Morse theory}, Ann. Inst. H. Poincaré C Anal. Non Linéaire \textbf{26} (2009), no. 2, 397--413.

\bibitem{CLARKSON} J. A. Clarkson, \textit{Uniformly convex spaces}, Trans. Amer. Math. Soc. \textbf{40} (1936), no. 3, 396--414.

\bibitem{CDM2} P. Clément, D. G. de Figueiredo, and E. Mitidieri, \textit{A priori estimates for positive solutions of semilinear elliptic systems via Hardy-Sobolev inequalities}, in: Differential Equations, Asymptotic Analysis, and Mathematical Physics (Potsdam, 1996), Pitman Res. Notes Math. Ser., 366, Longman, Harlow, 1997, 73--91.

\bibitem{CDM1} P. Clément, D. G. de Figueiredo, and E. Mitidieri, \textit{Positive solutions of semilinear elliptic systems}, in: Djairo G. de Figueiredo - Selected Papers (ed. D. Cassani et al.), Springer, Cham, 2014, 369--386.

\bibitem{C} J. N. Corvellec, \textit{Quantitative deformation theorems and critical point theory}, Pacific J. Math. \textbf{187} (1999), no. 2, 263--279.

\bibitem{CH} J. N. Corvellec and A. Hantoute, \textit{Homotopical stability of isolated critical points of continuous functionals}, Set-Valued Anal. \textbf{10} (2002), no. 2-3, 143--164.

\bibitem{CFV} M. Cozzi, A. Farina, and E. Valdinoci, \textit{Monotonicity formulae and classification results for singular, degenerate, anisotropic PDEs}, Adv. Math. \textbf{293} (2016), 343--381.



\bibitem{CUESTA} M. Cuesta, \textit{Eigenvalue problems for the p-Laplacian with indefinite weights}, Electron. J. Differential Equations \textbf{2001} (2001), no. 33, 1--9.


\bibitem{DeY} D. G. de Figueiredo and J. Yang, \textit{A priori bounds for positive solutions of a non-variational elliptic system}, in: Djairo G. de Figueiredo - Selected Papers (ed. D. Cassani et al.), Springer, Cham, 2014, 483--499.

\bibitem{DeFP} C. De Filippis and M. Piccinini, \textit{Borderline global regularity for nonuniformly elliptic systems}, Int. Math. Res. Not. IMRN \textbf{2023} (2023), no. 20, 17324--17376.

\bibitem{morais-souto} D. C. de Morais Filho and M. A. S. Souto, \textit{Systems of p-Laplacian equations involving homogeneous nonlinearities with critical Sobolev exponent degrees}, Comm. Partial Differential Equations \textbf{24} (1999), no. 7-8, 1537--1553.

\bibitem{denapolimarani} P. De Napoli and M. C. Mariani, \textit{Quasilinear elliptic systems of resonant type and nonlinear eigenvalue problems}, Abstr. Appl. Anal. \textbf{7} (2002), no. 3, 155--167.

\bibitem{DT} F. de Thélin, \textit{First eigenvalue of a nonlinear elliptic system}, C. R. Acad. Sci. Paris Sér. I Math. \textbf{311} (1990), no. 10, 603--606.

\bibitem{DD1} M. del Pino and J. Dolbeault, \textit{Best constants for Gagliardo–Nirenberg inequalities and applications to nonlinear diffusions}, J. Math. Pures Appl. (9) \textbf{81} (2002), no. 9, 847--875.

\bibitem{DD2} M. del Pino and J. Dolbeault, \textit{Nonlinear diffusions and optimal constants in Sobolev type inequalities: asymptotic behaviour of equations involving the p-Laplacian}, C. R. Math. Acad. Sci. Paris \textbf{334} (2002), no. 5, 365--370.

\bibitem{DD3} M. del Pino and J. Dolbeault, \textit{The Euclidean Onofri inequality in higher dimensions}, Int. Math. Res. Not. IMRN \textbf{2013} (2013), no. 15, 3600--3611.

\bibitem{degiovanni2009} M. Degiovanni, \textit{Critical groups of finite type for functionals defined on Banach spaces}, Commun. Appl. Anal. \textbf{13} (2009), no. 3, 395--410.

\bibitem{D} M. Degiovanni, \textit{On topological Morse theory}, J. Fixed Point Theory Appl. \textbf{10} (2011), no. 2, 197--218.

\bibitem{DL} M. Degiovanni and S. Lancelotti, \textit{Linking over cones and nontrivial solutions for p-Laplace equations with p-superlinear nonlinearity}, Ann. Inst. H. Poincaré C Anal. Non Linéaire \textbf{24} (2007), no. 6, 907--919.

\bibitem{DM} M. Degiovanni and M. Marzocchi, \textit{A critical point theory for nonsmooth functionals}, Ann. Mat. Pura Appl. (4) \textbf{167} (1994), 73--100.

\bibitem{DEM} K. Deimling, \textit{Nonlinear Functional Analysis}, Springer-Verlag, Berlin-Heidelberg, 1985.


\bibitem{diazthelin} J. I. Díaz and F. de Thélin, \textit{On a nonlinear parabolic problem arising in some models related to turbulent flows}, SIAM J. Math. Anal. \textbf{25} (1994), no. 4, 1085--1111.

\bibitem{DIAZSAA} J. I. Díaz and J. E. Saá, \textit{Existence et unicité de solutions positives pour certaines équations elliptiques quasilinéaires}, C. R. Acad. Sci. Paris Sér. I Math. \textbf{305} (1987), no. 12, 521--524.

\bibitem{dibenedetto1983} E. DiBenedetto, \textit{$C^{1+\alpha}$ local regularity of weak solutions of degenerate elliptic equations}, Nonlinear Anal. \textbf{7} (1983), no. 8, 827--850.

\bibitem{dinca_jebelean_mawhin2001} G. Dinca, P. Jebelean, and J. Mawhin, \textit{Variational and topological methods for Dirichlet problems with p-Laplacian}, Port. Math. (N.S.) \textbf{58} (2001), no. 3, 339--378.

\bibitem{dingxiao} L. Ding and S. W. Xiao, \textit{Multiple positive solutions for a critical quasilinear elliptic system}, Nonlinear Anal. \textbf{72} (2010), no. 5, 2592--2607.

\bibitem{DEJ} J. Dolbeault, M. J. Esteban, and G. Jankowiak, \textit{The Moser-Trudinger-Onofri inequality}, Chin. Ann. Math. Ser. B \textbf{36} (2015), no. 5, 777--802.

\bibitem{DET} J. Dolbeault, M. J. Esteban, and G. Tarantello, \textit{The role of Onofri type inequalities in the symmetry properties of extremals for Caffarelli-Kohn-Nirenberg inequalities, in two space dimensions}, Ann. Sc. Norm. Super. Pisa Cl. Sci. (5) \textbf{7} (2008), no. 2, 313--341.

\bibitem{Dolbeault-Perthame} J. Dolbeault and B. Perthame, \textit{Optimal critical mass in the two dimensional Keller-Segel model in $\mathbb{R}^2$}, C. R. Math. Acad. Sci. Paris \textbf{339} (2004), no. 9, 611--616.

\bibitem{DLL} M. Dong, N. Lam, and G. Lu, \textit{Sharp weighted Trudinger–Moser and Caffarelli–Kohn-Nirenberg inequalities and their extremal functions}, Nonlinear Anal. \textbf{173} (2018), 75--98.


\bibitem{DrabekRobinson} P. Drábek and S. B. Robinson, \textit{Resonance problems for the p-Laplacian}, J. Funct. Anal. \textbf{169} (1999), no. 1, 189--200.

\bibitem{drabek} P. Drábek, N. M. Stavrakakis, and N. B. Zographopoulos, \textit{Multiple nonsemitrivial solutions for quasilinear elliptic systems}, Differential Integral Equations \textbf{16} (2003), no. 12, 1519--1531.


\bibitem{perera} S. El Manouni and K. Perera, \textit{Existence and nonexistence results for a class of quasilinear elliptic systems}, Bound. Value Probl. \textbf{2007} (2007), Art. ID 85621, 5 pp.

\bibitem{ESE} S. Esedoglu and S. J. Osher, \textit{Decomposition of images by the anisotropic Rudin-Osher-Fatemi model}, Comm. Pure Appl. Math. \textbf{57} (2004), no. 12, 1609--1626.

\bibitem{esposito} P. Esposito, \textit{A classification result for the quasi-linear Liouville equation}, Ann. Inst. H. Poincaré C Anal. Non Linéaire \textbf{35} (2018), no. 3, 781--801.

\bibitem{fadell_rabinowitz1977} E. R. Fadell and P. H. Rabinowitz, \textit{Bifurcations for odd potential operators and an alternative topological index}, J. Funct. Anal. \textbf{26} (1977), no. 1, 48--67.

\bibitem{fadell_rabinowitz1978} E. R. Fadell and P. H. Rabinowitz, \textit{Generalized cohomological index theories for Lie group actions with an application to bifurcation questions for Hamiltonian systems}, Invent. Math. \textbf{45} (1978), no. 2, 139--174.


\bibitem{Fontana} L. Fontana, \textit{Sharp borderline Sobolev inequalities on compact Riemannian manifolds}, Comment. Math. Helv. \textbf{68} (1993), no. 3, 415--454.

\bibitem{furtadodepaiva} M. Furtado and F. O. V. de Paiva, \textit{Multiplicity of solutions for resonant elliptic systems}, J. Math. Anal. Appl. \textbf{319} (2006), no. 2, 435--449.

\bibitem{GAJZAK} H. Gajewski and K. Zacharias, \textit{Global behaviour of a reaction-diffusion system modelling chemotaxis}, Math. Nachr. \textbf{195} (1998), 77--114.


\bibitem{GIGA} Y. Giga, \textit{Surface Evolution Equations. A Level Set Approach}, Monographs in Mathematics, vol. 99, Birkhäuser Verlag, Basel, 2006.

\bibitem{glow} R. Glowinski and J. Rappaz, \textit{Approximation of a nonlinear elliptic problem arising in a non-Newtonian fluid flow model in glaciology}, ESAIM Math. Model. Numer. Anal. \textbf{37} (2003), no. 1, 175--186.

\bibitem{gromoll_meyer1969-t} D. Gromoll and W. Meyer, \textit{On differentiable functions with isolated critical points}, Topology \textbf{8} (1969), 361--369.

\bibitem{GV} M. Guedda and L. Véron, \textit{Quasilinear elliptic equations involving critical Sobolev exponents}, Nonlinear Anal. \textbf{13} (1989), no. 8, 879--902.

\bibitem{GPZ} Z. Guo, K. Perera, and W. Zou, \textit{On critical p-Laplacian systems}, Adv. Nonlinear Stud. \textbf{17} (2017), no. 4, 641--659.

\bibitem{GURTIN} M. E. Gurtin, \textit{Thermomechanics of Evolving Phase Boundaries in the Plane}, Oxford Mathematical Monographs, Oxford University Press, New York, 1993.

\bibitem{HESS} P. Hess, \textit{On a theorem by Landesman and Lazer}, Indiana Univ. Math. J. \textbf{23} (1973/74), no. 9, 827--829.

\bibitem{HO} K. Ho, Y. H. Kim, P. Winkert, and C. Zhang, \textit{The boundedness and Hölder continuity of weak solutions to elliptic equations involving variable exponents and critical growth}, J. Differential Equations \textbf{313} (2022), 503--532.




\bibitem{IM} S. Iula and G. Mancini, \textit{Extremal functions for singular Moser–Trudinger embeddings}, Nonlinear Anal. \textbf{156} (2017), 215--248.

\bibitem{JM} G. Joyce and D. Montgomery, \textit{Negative temperature states for the two-dimensional guiding-center plasma}, J. Plasma Phys. \textbf{10} (1973), 107--121.

\bibitem{K} G. Katriel, \textit{Mountain pass theorems and global homeomorphism theorems}, Ann. Inst. H. Poincaré C Anal. Non Linéaire \textbf{11} (1994), no. 2, 189--209.

\bibitem{KW} J. L. Kazdan and F. W. Warner, \textit{Curvature functions for compact 2-manifolds}, Ann. of Math. (2) \textbf{99} (1974), no. 1, 14--47.

\bibitem{kim_wang1989} S. K. Kim and T. X. Wang, \textit{The generalized Morse lemma and the Euler characteristic on Banach manifolds}, Topology Appl. \textbf{32} (1989), no. 1, 13--23.

\bibitem{KS} D. Kinderlehrer and G. Stampacchia, \textit{An Introduction to Variational Inequalities and Their Applications}, Classics in Applied Mathematics, vol. 31, Society for Industrial and Applied Mathematics (SIAM), Philadelphia, PA, 2000.

\bibitem{LWO} S. S. Lai, G. T. Will, and S. Otani, \textit{Model for inelastic biaxial bending of concrete members}, J. Structural Engineering \textbf{110} (1984), no. 11, 2563--2584.


\bibitem{6} S. Lancelotti, \textit{Morse index estimates for continuous functionals associated with quasilinear elliptic equations}, Adv. Differential Equations \textbf{7} (2002), no. 1, 99--128.

\bibitem{LANDLAZ} E. M. Landesman and A. C. Lazer, \textit{Nonlinear perturbations of linear elliptic boundary value problems at resonance}, J. Math. Mech. \textbf{19} (1970), no. 7, 609--623.

\bibitem{lazer_solimini1988} A. C. Lazer and S. Solimini, \textit{Nontrivial solutions of operator equations and Morse indices of critical points of min-max type}, Nonlinear Anal. \textbf{12} (1988), no. 8, 761--775.




\bibitem{li_li_liu2005} C. Li, S. Li, and J. Liu, \textit{Splitting theorem, Poincaré-Hopf theorem and jumping nonlinear problems}, J. Funct. Anal. \textbf{221} (2005), no. 2, 439--455.




\bibitem{LiebLoss} E. H. Lieb and M. Loss, \textit{Analysis}, Graduate Studies in Mathematics, vol. 14, American Mathematical Society, Providence, RI, second ed., 2001.

\bibitem{lieberman} G. M. Lieberman, \textit{Boundary regularity for solutions of degenerate elliptic equations}, Nonlinear Anal. \textbf{12} (1988), no. 11, 1203--1219.



\bibitem{Malchiodi} A. Malchiodi, \textit{Topological methods for an elliptic equation with exponential nonlinearities}, Discrete Contin. Dyn. Syst. \textbf{21} (2008), no. 1, 277--294.


\bibitem{manamaw} R. Manásevich and J. Mawhin, \textit{Periodic solutions for nonlinear systems with p-Laplacian-like operators}, J. Differential Equations \textbf{145} (1998), no. 2, 367--393.





\bibitem{marcellini1} P. Marcellini, \textit{On the definition and the lower semicontinuity of certain quasiconvex integrals}, Ann. Inst. H. Poincaré C Anal. Non Linéaire \textbf{3} (1986), no. 5, 391--409.

\bibitem{MarPer} A. J. Marinho and K. Perera, \textit{Local and nonlocal critical growth anisotropic quasilinear elliptic systems}, Calc. Var. Partial Differential Equations \textbf{64} (2025), no. 4, Art. 102, 26 pp.

\bibitem{MW} G. Marino and P. Winkert, \textit{Global a priori bounds for weak solutions of quasilinear elliptic systems with nonlinear boundary condition}, J. Math. Anal. Appl. \textbf{482} (2020), no. 2, Art. 123555, 19 pp.




\bibitem{MasielloPisani} A. Masiello and L. Pisani, \textit{Asymptotically linear elliptic problems at resonance}, Ann. Mat. Pura Appl. (4) \textbf{171} (1996), 1--13.

\bibitem{mawhin_willem1989} J. Mawhin and M. Willem, \textit{Critical Point Theory and Hamiltonian Systems}, Applied Mathematical Sciences, vol. 74, Springer-Verlag, New York, 1989.

\bibitem{MINGIONE} G. Mingione, \textit{Bounds for the singular set of solutions to non linear elliptic systems}, Calc. Var. Partial Differential Equations \textbf{18} (2003), no. 4, 373--400.

\bibitem{MMSV} L. Montoro, L. Muglia, B. Sciunzi, and D. Vuono, \textit{Regularity and symmetry results for the vectorial p-Laplacian}, Nonlinear Anal. \textbf{251} (2025), Art. 113700, 20 pp.

\bibitem{M} J. Moser, \textit{A sharp form of an inequality by N. Trudinger}, Indiana Univ. Math. J. \textbf{20} (1970/71), no. 11, 1077--1092.

\bibitem{M2} J. Moser, \textit{On a nonlinear problem in differential geometry}, in: Dynamical Systems (M. M. Peixoto, ed.), Academic Press, New York, 1973, 273--280.

\bibitem{O} E. Onofri, \textit{On the positivity of the effective action in a theory of random surfaces}, Comm. Math. Phys. \textbf{86} (1982), no. 3, 321--326.

\bibitem{OPS} B. Osgood, R. Phillips, and P. Sarnak, \textit{Extremals of determinants of Laplacians}, J. Funct. Anal. \textbf{80} (1988), no. 1, 148--211.

\bibitem{OPS2} B. Osgood, R. Phillips, and P. Sarnak, \textit{Compact isospectral sets of surfaces}, J. Funct. Anal. \textbf{80} (1988), no. 1, 212--234.

\bibitem{Tang} Z. Ou and C. Tang, \textit{Existence and multiplicity of nontrivial solutions for quasilinear elliptic systems}, J. Math. Anal. Appl. \textbf{383} (2011), no. 2, 423--438.

\bibitem{PR} M. C. Pelissier and L. Reynaud, \textit{Étude d'un modèle mathématique d'écoulement de glacier}, C. R. Acad. Sci. Paris Sér. A \textbf{279} (1974), 531--534.

\bibitem{PERERA} K. Perera, \textit{Nontrivial critical groups in p-Laplacian problems via the Yang index}, Topol. Methods Nonlinear Anal. \textbf{21} (2003), no. 2, 301--309.

\bibitem{PAO} K. Perera, R. P. Agarwal, and D. O'Regan, \textit{Morse Theoretic Aspects of p-Laplacian Type Operators}, Mathematical Surveys and Monographs, vol. 161, American Mathematical Society, Providence, RI, 2010.

\bibitem{poho} S. I. Pohozaev, \textit{On the eigenfunctions of the equation $\Delta u + \lambda f(u)=0$}, Dokl. Akad. Nauk SSSR \textbf{165} (1965), 36--39; English transl., Soviet Math. Dokl. \textbf{6} (1965), 1408--1411.

\bibitem{Poly1} A. M. Polyakov, \textit{Quantum geometry of bosonic strings}, Phys. Lett. B \textbf{103} (1981), no. 3, 207--210.

\bibitem{Poly2} A. M. Polyakov, \textit{Quantum geometry of fermionic strings}, Phys. Lett. B \textbf{103} (1981), no. 3, 211--213.


\bibitem{QS} P. Quittner and P. Souplet, \textit{A priori estimates and existence for elliptic systems via bootstrap in weighted Lebesgue spaces}, Arch. Ration. Mech. Anal. \textbf{174} (2004), no. 1, 49--81.

\bibitem{RAB} P. H. Rabinowitz, \textit{Some minimax theorems and applications to nonlinear partial differential equations}, Nonlinear Anal. \textbf{2} (1978), no. 2, 161--177.

\bibitem{RZ} V. D. Rădulescu and B. Zhang, \textit{Morse theory and local linking for a nonlinear degenerate problem arising in the theory of electrorheological fluids}, Nonlinear Anal. Real World Appl. \textbf{17} (2014), 311--321.


\bibitem{rybakowski1987} K. P. Rybakowski, \textit{The Homotopy Index and Partial Differential Equations}, Universitext, Springer-Verlag, Berlin, 1987.


\bibitem{spa} E. H. Spanier, \textit{Algebraic Topology}, McGraw-Hill Book Co., New York-Toronto, Ont.-London, 1966.

\bibitem{S} G. Stampacchia, \textit{Le problème de Dirichlet pour les équations elliptiques du second ordre à coefficients discontinus}, Ann. Inst. Fourier (Grenoble) \textbf{15} (1965), no. 1, 189--257.

\bibitem{SZ} N. M. Stavrakakis and N. B. Zographopoulos, \textit{Bifurcation results for quasilinear elliptic systems}, Adv. Differential Equations \textbf{8} (2003), no. 3, 315--336.


\bibitem{STRUWE} M. Struwe, \textit{Variational Methods}, third ed., Ergebnisse der Mathematik und ihrer Grenzgebiete, vol. 34, Springer-Verlag, Berlin, 2000.

\bibitem{StruweCurv} M. Struwe, \textit{A flow approach to Nirenberg’s problem}, Duke Math. J. \textbf{128} (2005), no. 1, 19--64.




\bibitem{SZU} A. Szulkin, \textit{Ljusternik-Schnirelmann theory on $C^1$-manifolds}, Ann. Inst. H. Poincaré Anal. Non Linéaire \textbf{5} (1988), no. 2, 119--139.

\bibitem{TAYLOR1} J. E. Taylor, \textit{Crystalline variational problems}, Bull. Amer. Math. Soc. \textbf{84} (1978), no. 4, 568--588.

\bibitem{TAYLOR2} J. E. Taylor, J. W. Cahn, and C. A. Handwerker, \textit{Geometric models of crystal growth}, Acta Metall. Mater. \textbf{40} (1992), no. 7, 1443--1474.



\bibitem{TF} K. Tintarev and K. H. Fieseler, \textit{Concentration Compactness: Functional-Analytic Grounds and Applications}, Imperial College Press, London, 2007.

\bibitem{tolksdorf1983} P. Tolksdorf, \textit{On the Dirichlet problem for quasilinear equations in domains with conical boundary points}, Comm. Partial Differential Equations \textbf{8} (1983), no. 7, 773--817.

\bibitem{tolksdorf1984} P. Tolksdorf, \textit{Regularity for a more general class of quasilinear elliptic equations}, J. Differential Equations \textbf{51} (1984), no. 1, 126--150.


\bibitem{Tru} N. S. Trudinger, \textit{On imbeddings into Orlicz spaces and some applications}, J. Math. Mech. \textbf{17} (1967), no. 5, 473--483.

\bibitem{U} K. Uhlenbeck, \textit{Morse theory on Banach manifolds}, J. Funct. Anal. \textbf{10} (1972), no. 4, 430--445.

\bibitem{Vann1} G. Vannella, \textit{Uniform $L^\infty$-estimates for quasilinear elliptic systems}, J. Math. \textbf{2023} (2023), Art. ID 2891340, 11 pp.


\bibitem{VAZ} J. L. Vázquez, \textit{A strong maximum principle for some quasilinear elliptic equations}, Appl. Math. Optim. \textbf{12} (1984), no. 3, 191--202.

\bibitem{velin-dethelin} J. Velin and F. de Thélin, \textit{Existence and nonexistence of nontrivial solutions for some nonlinear elliptic systems}, Rev. Mat. Univ. Complut. Madrid \textbf{6} (1993), no. 1, 153--194.



\bibitem{WULFF} G. Wulff, \textit{Zur Frage der Geschwindigkeit des Wachstums und der Auflösung der Kristallflächen}, Z. Krist. \textbf{34} (1901), 449--530.



\bibitem{Yudo} V. I. Yudovich, \textit{Some estimates connected with integral operators and with solutions of elliptic equations}, Dokl. Akad. Nauk SSSR \textbf{138} (1961), no. 4, 805--808; English transl., Soviet Math. Dokl. \textbf{2} (1961), 746--749.

\end{thebibliography}

\fancyhead[LE, RO]{\thepage}
\fancyhead[RE]{\nouppercase{\leftmark}}
\fancyhead[LO]{\nouppercase{\rightmark}}

\fancyhead[L]{ \ifthenelse{\isodd{\value{page}}}{\nouppercase{\leftmark}}{\thepage} }
\fancyhead[R]{ \ifthenelse{\isodd{\value{page}}}{\thepage}{\nouppercase{\leftmark}} }

\bigskip
\bigskip
\bigskip
\bigskip

\textbf{Mathematics Subject Classification:} 

\bigskip

\begin{tabular}{llllllll}
26D10
&
35B06
&
35B34
&
35B45
&
35B65
&
35J50
&
35J60
&
35J62
\\
35J70
&
35J92
&
46E30
&
46E35
&
47H11
&
47J30
&
58E05
&
58E35
\end{tabular}

\bigskip
\bigskip
\bigskip

\textbf{Key words:} 

\bigskip

\begin{tabular}{lll}
Asymptotically $(p,q)$ linear
&
Banach spaces
&
Critical groups
\\
Critical growth
&
Degree theory
&
Eigenvalue problem
\\
Fadell--Rabinowitz index
&
Generalized Stampacchia lemma
&
Morse index
\\
Poincaré--Hopf formula
&
Quasilinear elliptic systems
&
Resonant problems
\\
Trudinger--Moser inequality
&
Two dimension
&
Uniform $L^{\infty}$-boundedness 
\\
\end{tabular}

\bigskip
\bigskip
\bigskip

\textbf{Contacts:} \href{mailto: natalino.borgia@uniba.it}{natalino.borgia@uniba.it}, Dipartimento di Matematica, UNIBA (Bari).
 
\restoregeometry

\end{document}